\documentclass[11pt,leqno]{article}

\usepackage{amsmath,amsfonts,amscd,amssymb,theorem}

\long\def\comment#1\endcomment{}

\makeatletter
\begingroup
\gdef\th@dotted{\normalfont\itshape
  \def\@begintheorem##1##2{%
        \item[\hskip\labelsep \theorem@headerfont ##1\ ##2.]}%
\def\@opargbegintheorem##1##2##3{%
   \item[\hskip\labelsep \theorem@headerfont ##1\ ##2\ (##3).]}}
\endgroup
\makeatother

\theoremstyle{dotted}

\newtheorem{theorem}{Theorem}[section]
\newtheorem{lemma}[theorem]{Lemma}

\newtheorem{prop}[theorem]{Proposition}
\newtheorem{corr}[theorem]{Corollary}

\makeatletter
\begingroup
\gdef\th@upshape{\normalfont
  \def\@begintheorem##1##2{%
        \item[\hskip\labelsep \theorem@headerfont ##1\ ##2.]}%
\def\@opargbegintheorem##1##2##3{%
   \item[\hskip\labelsep \theorem@headerfont ##1\ ##2\ (##3).]}}
\endgroup
\makeatother

\theoremstyle{upshape}

\newtheorem{defn}[theorem]{Definition}
\newtheorem{remark}[theorem]{Remark}
\newtheorem{exa}[theorem]{Example}

\makeatletter
\renewcommand{\subsection}{\@startsection{subsection}{2}{0pt}{-3ex
plus -1ex minus -0.2ex}{-2mm plus -0pt minus
-2pt}{\normalfont\bfseries}} 
\renewcommand{\subsubsection}{\@startsection{subsubsection}{3}{0pt}{-3ex
plus -1ex minus -0.2ex}{-2mm plus -0pt minus
-2pt}{\normalfont\bfseries}} 
\makeatother

\makeatletter
\@addtoreset{equation}{section}
\makeatother

\newcommand{\cntrct}                
{\hspace{2pt}\raisebox{1pt}{\text{$\lrcorner$}}\hspace{2pt}}

\newcommand{\proof}[1][Proof.]{\smallskip\noindent{\em #1}}
\def\endproof{\hfill\ensuremath{\square}\par\medskip}

\renewcommand{\labelenumi}{{\normalfont(\roman{enumi})}}

\def\eqref#1{\thetag{\ref{#1}}}

\let\latexref=\ref
\def\ref#1{{\normalfont{\latexref{#1}}}}

\newcommand{\wt}{\widetilde}

\newcommand{\dg}{\dagger}

\newcommand{\idot}{{\:\raisebox{1pt}{\text{\circle*{1.5}}}}}
\newcommand{\hdot}{{\:\raisebox{3pt}{\text{\circle*{1.5}}}}}
\newcommand{\eps}{\varepsilon}
\renewcommand{\phi}{\varphi}

\newcommand{\vH}{\check{H}}
\newcommand{\vC}{\check{C}}

\newcommand{\bH}{\overline{H}}
\newcommand{\bC}{\overline{C}}

\newcommand{\Hom}{\operatorname{Hom}}
\newcommand{\Ext}{\operatorname{Ext}}
\newcommand{\RHom}{\operatorname{RHom}}
\newcommand{\Tor}{\operatorname{Tor}}

\newcommand{\Ker}{\operatorname{Ker}}

\newcommand{\id}{\operatorname{\sf id}}
\newcommand{\Id}{\operatorname{\sf Id}}
\newcommand{\Iso}{\operatorname{\sf Iso}}

\newcommand{\gr}{\operatorname{\sf gr}}
\newcommand{\tr}{\operatorname{\sf tr}}
\newcommand{\colim}{\operatorname{colim}}
\renewcommand{\lim}{\operatorname{lim}}

\newcommand{\D}{{\cal D}}
\newcommand{\C}{{\cal C}}

\newcommand{\K}{\mathcal{K}}

\newcommand{\hash}{\sharp}

\newcommand{\Sets}{\operatorname{Sets}}
\newcommand{\Aut}{{\operatorname{Aut}}}

\newcommand{\amod}{{\text{\rm -mod}}}
\newcommand{\bimod}{{\text{\rm -bimod}}}

\newcommand{\ppt}{{\sf pt}}

\newcommand{\lotimes}{\overset{\sf\scriptscriptstyle L}{\otimes}}

\newcommand{\Spec}{\operatorname{Spec}}

\newcommand{\cchar}{\operatorname{\sf char}}

\newcommand{\Z}{{\mathbb Z}}
\newcommand{\N}{{\mathbb N}}
\newcommand{\Q}{{\mathbb Q}}

\newcommand{\LR}{\Lambda R}

\newcommand{\Fun}{\operatorname{Fun}}

\newcommand{\E}{\mathcal{E}}

\newcommand{\Stab}{\operatorname{Stab}}

\newcommand{\Top}{\operatorname{Top}}
\newcommand{\Ho}{\operatorname{Ho}}

\newcommand{\hocolim}{\operatorname{\sf hocolim}}
\newcommand{\holim}{\operatorname{\sf holim}}
\newcommand{\ho}{\operatorname{\mathcal{H}{\it o}}}

\newcommand{\copr}{\sqcup}

\newcommand{\fFun}{\operatorname{{\mathcal F}{\mathit u}{\mathit n}}}
\newcommand{\sFun}{\operatorname{\sf Fun}}
\newcommand{\Sec}{\operatorname{Sec}}

\newcommand{\oFun}{\Fun^\otimes}
\newcommand{\bFun}{\Fun^\boxtimes}

\newcommand{\Tw}{\operatorname{\sf tw}}
\newcommand{\Ar}{\operatorname{\sf ar}}
\newcommand{\Sp}{\operatorname{Sp}}

\newcommand{\tw}{\operatorname{Tw}}

\newcommand{\Bi}{B_\infty}

\newcommand{\Cyl}{\operatorname{\sf C}}

\newcommand{\V}{{\sf V}}

\newcommand{\wCat}{\operatorname{Aug}}
\newcommand{\rCat}{\operatorname{Cat}}
\newcommand{\Refl}{\operatorname{Refl}}
\newcommand{\Alg}{\operatorname{Alg}}
\newcommand{\aAlg}{\operatorname{\mathcal{A}{\it lg}}}
\newcommand{\iCat}{\Cat^*}
\newcommand{\iMor}{\Mor^*}
\newcommand{\Mor}{\operatorname{\mathcal{M}{\it or}}}
\newcommand{\Cat}{\operatorname{\mathcal{C}{\it at}}}

\newcommand{\M}{\mathcal{M}}

\newcommand{\rAdj}{\operatorname{Adj}}
\newcommand{\aAdj}{\operatorname{\mathcal{A}{\it dj}}}

\newcommand{\adj}{\operatorname{\sf adj}}
\newcommand{\Adj}{\operatorname{\sf Adj}}
\newcommand{\eq}{\operatorname{\sf eq}}
\newcommand{\Eq}{\operatorname{\sf Eq}}
\newcommand{\nat}{\operatorname{\sf nat}}
\newcommand{\Nat}{\operatorname{\sf Nat}}

\newcommand{\Exp}{\operatorname{Exp}}
\newcommand{\Red}{\operatorname{Red}}

\newcommand{\Tr}{\operatorname{Tr}}

\newcommand{\ev}{\operatorname{\sf ev}}

\newcommand{\Av}{\operatorname{Av}}

\newcommand{\ogamma}{\overline{\gamma}}
\newcommand{\wgamma}{\wt{\gamma}}

\newcommand{\F}{\operatorname{\sf F}}
\newcommand{\Y}{\operatorname{\sf Y}}

\newcommand{\Tri}{\Tr_\infty}

\newcommand{\Dd}{\operatorname{\sf D}}

\newcommand{\Aa}{\operatorname{\sf A}}

\newcommand{\dm}{\diamond}

\newcommand{\Real}{\operatorname{\sf real}}
\newcommand{\Sing}{\operatorname{\sf sing}}

\newcommand{\Span}{\operatorname{\sf span}}

\renewcommand{\aa}{\mathfrak{a}}

\newcommand{\Ss}{\mathbb{S}}

\newcommand{\I}{\operatorname{\sf I}}

\newcommand{\HH}{\mathcal{H}\mathcal{H}}
\newcommand{\THH}{\mathcal{T}\mathcal{H}\mathcal{H}}
\newcommand{\WHH}{\mathcal{W}\mathcal{H}\mathcal{H}}

\newcommand{\Ff}{\mathbb{F}}

\newcommand{\X}{\mathcal{X}}

\newcommand{\Loc}{\operatorname{\sf Loc}}

\newcommand{\bI}{\overline{I}}

\newcommand{\nm}{{\displaystyle n\!\!-\!\!1}}
\newcommand{\np}{{\displaystyle n\!\!+\!\!1}}
\newcommand{\mm}{{\displaystyle m\!\!-\!\!1}}
\newcommand{\nl}{{\displaystyle n\!\!-\!\!l}}
\newcommand{\nlm}{{\displaystyle n\!\!-\!\!l\!\!-\!\!1}}

\newcommand{\CC}{CC}

\newcommand{\bK}{\K_\flat}

\newcommand{\bHP}{\overline{HP}}
\newcommand{\bCP}{\overline{CP}}

\title{Trace theories, B\"okstedt periodicity and Bott periodicity}

\author{D. Kaledin\thanks{Supported by the Russian Science
    Foundation, grant 21-11-00153.}}

\begin{document}

\maketitle

\tableofcontents

\section*{Introduction.}

The big miracle of Topological Hochschild Homology is that it is so
small: if $k$ is a perfect field of some positive characteristic
$p$, then $THH_\idot(k) = k[\sigma]$, the algebra of polynomials in
one variable $\sigma$ of homological degree $2$ known as the {\em
  B\"okstedt periodicity element}. If one only cares about
$THH_\idot(k)$ as a graded $k$-vector space, then there is a
spectral sequence converging to it whose $E_2$-page is easy to
compute. The sequence degenerates at $E_p$, so that if $p=2$, it
just degenerates (but one still has to do a separate argument if one
is interested in the multiplicative structure). For odd $p$, the
sequence does {\em not} degenerate: the $E_2$-page is rather large,
but the differential $d_{p-1}$ is non-trivial and cuts it down to
its proper size. Again, the multiplicative structure has to be studied
separately.

This miraculous behavior of $THH$ can be proved in several different
ways (we give a brief overview below in
Subsections~\ref{mcl.subs},~\ref{bok.subs}). However, none of the
proofs are easy. This paper arose as an attempt to answer a question
of L. Hesselholt: is there now, when Topological Hochschild Homology
is in its fourth decade, a simple and/or conceptually clear proof?

Let us state right away that our attempt mostly failed. The best we
can do with the spectral sequence argument is a modern repackaging
of the 1994 proof of \cite{slf}. We add strictly polynomial functors
that were not available at the time. The resulting argument is
reasonably short but needs one external input, namely, the fact that
the category of strictly polynomial functors of some fixed degree
has finite homological dimension. For multiplication, we can do
slightly better: we construct a ring map $\phi$ from $THH_\idot(k)$
to some ring such that $\phi(\sigma)$ is not nilpotent; this means
that $k[\sigma]$ is a subalgebra in $THH_\idot(k)$, and then it must
be the whole thing for dimension reasons. The map $\phi$ can be
obtained by a certain truncation of the cyclotomic structure map of
\cite{NS}, but it is actually simpler to construct it by general
nonsense.

To compensate for our failure at the main stated goal, we also prove
two comparison theorem. The first one gives a really simple and
purely algebraic expression of the Topological Hohchschild Homology
$THH_\idot(A,M)$ of a $k$-algebra $A$ with coefficients in an
$A$-bimodule $M$ in terms of the so-called Hochschild-Witt Homology
$WHH_\idot(A,M)$ introduced in \cite{witt}, \cite{hw} (see also an
overview in \cite{umn}). If $M=A$ is the diagonal bimodule, then we
also analyze the circle action on $THH_\idot(A) = THH_\idot(A,A)$,
and prove that the periodic Topological Cyclic Homology
$TP_\idot(A)$ of \cite{tp} coincides with the periodic version
$WHP_\idot(A)$ of $WHH_\idot(A)$. The second comparison theorem
concerns only the diagonal bimodule case: we prove that
$THH_\idot(A)$ coincides with zero term of the ``conjugate
filtration'' on the co-periodic cyclic homology $\bCP_\idot(A)$
introduced in \cite{coper}. In particular, if one inverts the
B\"okstedt generator $\sigma$, then $THH(A)$ becomes
$\bCP_\idot(A)$, and the identification sends $\sigma$ to the Bott
periodicity generator $u^{-1}$.

All of the above is contained in Section~\ref{thh.sec}, and at this
point, the reader might wonder what are we doing in the first 10
sections, and indeed why is the paper so long. The answer to this is
technological.

Our main technical tool for both comparison theorems is the notion
of a ``trace theory'' sketched in \cite{trace}. Roughly speaking,
the idea is to observe that $THH_\idot(A,M)$ for different $A$ and
$M$ are related by certain canonical isomorphisms, and one can use
it to show that $THH_\idot(k,M)$, for any $M$, completely determines
$THH_\idot(A,M)$ for any $A$ and $M$. In other words, if one works
with arbitrary coefficients, and takes account of the trace theory
structure, then one can forget about arbitrary $k$-algebras and only
consider $k$ itself. A toy version of such a reconstruction theorem
was proved in \cite{trace}. However, doing the whole theory properly
and in a reasonable generality requires quite a lot of work.

The second tool that we use heavily throughout the paper is the
technique of ``stabilization'', or ``additivization''. Namely, for
any functor $F(M)$ from some additive category such as that of
bimodules to the category of spaces, one can take its Goodwillie
derivative, thereby forcing the functor to land in connective
spectra and become additive with respect to $M$. We have found out
that this works really well for $THH_\idot(A,M)$ and related
theories: both the theories and relevant maps between them such as
our map $\phi$ can be obtained by stabilization from something
really simple and obvious (for $THH_\idot(A,M)$, this is the
``cyclic nerve'' of \cite{G.W}). Informally, it pays to only go
stable at the very end.  However, once again, developing the
stabilization machine that is strong enough and versatile enough
requires some space.

\medskip

The two stories, that of stabilization and that of trace theories,
are basicaly independent, so that the structure of the paper is
non-linear. A rough {\em leitfaden} would look something like this:
$$
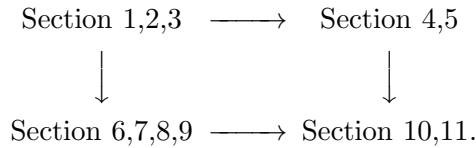
\begin{CD}
\text{Section 1,2,3} @>>> \text{Section 4,5}\\
@VVV @VVV\\
\text{Section 6,7,8,9} @>>> \text{Section 10,11}.
\end{CD}
$$
Before we explain in more detail what is being done where, let us
make a comment on methodology.

The natural context for trace theories is that of $2$-categories: a
trace theory is defined on a $2$-category, and a related notion of a
``trace functor'' makes sense for monoidal categories (that is,
$2$-categories with one object). The theory of $2$-categories is
notorious for horrible multidimensional diagrams, ``higher
associativity constraints'' and the like, with proofs being either
impossibly long, or incomplete. In our experience, the only way to
avoid this is to use systematically the usual category theory and
the machinery of Grothendieck fibrations. We thus replace a
$2$-category with its nerve, and treat it as a usual category
fibered over $\Delta$, modulo some conditions (but no extra
structure). Functors between $2$-categories are then functors over
$\Delta$, again modulo some conditions, and they are always
constructed by adjunction --- typically, as left or right Kan
extensions. All the ``higher constraints'' remain packaged via the
Grothendieck construction and do not appear explicitly; one never
has to check that two particular morphisms coincide.

For better or for worse, our approach to things homotopical is then
exactly the same: we rely on the usual category theory and
adjunction. We assume known that for any small category $I$,
localizing functors from $I$ to spaces gives a well-defined category
$\Ho(I)$, and for any functor $\gamma:I' \to I$, the pullback
functor $\gamma^*:\Ho(I) \to \Ho(I')$ has a left and right-adjoint
given by homotopy Kan extensions $\gamma_!,\gamma_*:\Ho(I') \to
\Ho(I)$. These facts can be established in any way one likes (for
example, by using Reedy model structures, as in \cite{DHKS}), and
then used as a black box. Everything else proceeds by using Kan
extensions, and the formalism is often identically the same as in
the non-homotopical setting.

This naive approach to homotopy theory has some obvious disadvantages ---
for example, we cannot works with $E_\infty$-ring spectra, nor with
equivariant and cyclotomic spectra --- but for our purposes in the
paper, it is sufficient.

\medskip

Let us now give a more detailed overview of the paper. Sections 1, 2
and 3 are preliminary. Section 1 is concerned with category theory
--- essentially, glorified bookkeeping --- so that nothing
whatsoever is new, and the main goal is to fix notation and
terminology (one slightly non-standard thing is the notion of a
``framing'' in Lemma~\ref{kan.le} that turns out to be very
convenient for computing Kan extensions). Section 2 deals with the
Grothendieck construction and Grothendieck fibrations. In fact,
since we have no need of model categories anywhere in the paper, we
take the liberty of restoring the original terminology of \cite{sga}
and speak about fibrations and cofibrations, with cofibrations
probably more important (for reasons explained in
Subsection~\ref{mon.subs}). Subsections~\ref{fib.subs} and
\ref{cart.subs} are fairly standard, Subsection~\ref{fun.subs}
contains a version of the right Kan extension construction for
cofibrations, and Subsection~\ref{ker.subs} contains technical
material that should probably be skipped until needed (``kernels''
of Definition~\ref{ker.def} are only used in
Subsection~\ref{ite.subs}, and ``reflections'' are only used to
construct the relative functor categories \eqref{ffun.r} and
\eqref{fun.ii} needed in Subsection~\ref{mono.subs}). Section 3 is
pure combinatorics: we fix notation for the category $\Delta$ of
finite non-empty ordinals and related small categories such as the
cyclic category $\Lambda$. Again, nothing is new (the treatment of
$\Lambda$ follows \cite[Subsection 1.2]{hw}).

Section 4 starts the homotopical part of the story. First, we
recall the standard Segal machine of \cite{seg}, and we show that
with very small modification, the same arguments provide a
stabilization functor for functors from some pointed category with
finite coproducts to spaces. We then get some additional results by
purely formal games with adjunction, and we treat multiplicative
structures, too (in a very naive way, see Remark~\ref{weak.rem}). In
Section 5, we compute some examples of stabilization; on one hand,
this illustrates the general machinery, and on the other hand,
provides the main technical results needed for Section 11. In
particular, we show how one can obtain the dual Steenrod algebra by
stabilization, and we also show that it has a natural filtered
counterpart obtained by stabilizing divided power functors. This
possibly extends all the way to a filtered version of the stable
homotopy category, see Remark~\ref{sthom.rem}, but we do not explore
this seriously at this point. The main technical tool is the notion
of ``truncated Tate cohomology'' of a finite group with respect to a
family of subgroups; this is the only thing that might possibly be
new (at least we do not know a ready reference).

In Section 6, we temporarily abandon homotopy theory and move to
$2$-categories. Monoidal structures and $2$-categories are treated
in Section 6 itself. Section 7 is devoted to adjunction and adjoint
pairs in general $2$-categories, and some results here might be new
-- in particular, this concerns a surprisingly simple description of
a free $2$-category generated by an $n$-tuple of composable adjoint
pairs given in Subsection~\ref{ite.subs}. We should mention though
that this description is only needed in
Subsection~\ref{fun.adj.subs}, so it can be skipped at first
reading. Moreover, the somewhat technical Proposition~\ref{tw.prop}
can be used as a black box. In Section 8, we construct and study one
specific $2$-category $\iMor(\C)$ of small categories enriched over
a fixed monoidal category $\C$. In Section 9, we introduce trace
theories, and we prove our main reconstruction theorem: a trace
theory on $\iMor(\C)$ is completely determined by the associated
trace functor on $\C$ (the actual statement is slightly stronger and
more precise, see Theorem~\ref{rec.thm}).

Sections 10 and 11 bring the two stories together. Section 10 is a
homotopical version of Section 9, and the gist of it is that
everything done for trace theories extends to homotopy trace
theories {\em verbatim}, simply by adding the adjective ``homotopy''
in appropriate places. We also combine this with Section 4, by
introducing the notion of a ``stable homotopy trace theory'' and
showing how to do stabilization in this context. With all this out
of the way, we reap the benefits: in Section 11, we apply the
machinery to THH and prove all the theorems mentioned above.

\section{Categorical preliminaries.}\label{cat.sec}

\subsection{Categories.}\label{cat.subs}

We distinguish between large and small categories but avoid more
advanced set theory. We denote by $\Sets$ the category of sets, and
we denote by $\rCat$ the category of small categories. For any
category $\C$, we write $c \in \C$ as a shorthand for ``$c$ is an
object in $\C$'', and for any $c,c' \in \C$, we write $\C(c,c')$ for
the set of maps from $c$ to $c'$ (we recall that by definition, this
is always a set). We denote by $\C^o$ the opposite category,
$\C^o(c',c) = \C(c,c')$, and for any morphism $f \in \C(c,c')$, we
denote by $f^o \in \C^o(c',c) = \C(c,c')$ the corresponding morphism
in the opposite category. We let $\C^>$ resp.\ $\C^<$ be the
category obtained by adding a new terminal resp.\ initial object $o$
to the category $\C$, and we note that we have $\C^{o<} \cong
\C^{>o}$, $\C^{o>} \cong \C^{<o}$. For any functor $\gamma:\C_0 \to
\C_1$, we denote by $\gamma^>:\C_0^> \to \C_1^>$, $\gamma^<:\C_0^<
\to \C_1^<$ its canonical extensions sending $o$ to $o$, and we let
$\gamma^o:\C_0^o \to \C_1^o$ be the opposite functor. For any two
categories $\C_0$, $\C_1$ equipped with functors $\gamma_0:\C_0 \to
\C$, $\gamma_1:\C_1 \to \C$ to a third category $\C$, we denote by
$\C_0 \ \raisebox{1pt}{$\mbox{}^{\gamma_0}$}\!\times_\C^{\gamma_1} \C_1$ the
category of triples $\langle c_0,c_1,\alpha \rangle$, $c_0 \in
\C_0$, $c_1 \in \C_1$, $\alpha:\gamma_0(c_0) \cong \gamma_1(c_1)$,
and we drop either or both of the superscripts $\gamma_0$,
$\gamma_1$ when they are clear from the context. If $\gamma_0$ is
clear but $\gamma_1$ is not, we will also denote $\C_0
\times^{\gamma_1}_\C \C_1$ by $\gamma_1^*\C_0$, and keep in mind
that it comes equipped with a projection $\gamma_1^*\C_0 \to \C_1$.

A subcategory $\C' \subset \C$ is {\em weakly full} if
$\C'(c,c')=\C(c,c')$ for any $c,c' \in \C' \subset \C$, and a weakly
full subcategory $\C' \subset \C$ is {\em full} if for any
isomorphism $c \cong c'$ in $\C$, $c \in \C'$ if and only if $c' \in
\C'$. We use set-theoretic notation for full subcategories: for a
full subcategory $\C' \subset \C$, the complement $\C' \setminus \C$
cosists of objects not in $\C'$, and for any two full subcategories
$\C',\C'' \subset \C$, we have the obvious full subcategories $\C'
\cap \C'',\C' \cup \C'' \subset \C$.

For any functor $\gamma:\C \to I$, the {\em left comma-fiber} $\C
/^\gamma i$ over an object $i$ in $I$ is the category of pairs
$\langle c,\alpha(c) \rangle$, $c \in \C$, $\alpha(c):\gamma(c) \to
i$ a map. The {\em fiber} $\C_i \subset \C/i$ is the full
subcategory spanned by $\langle c,\alpha(c) \rangle$ with invertible
$\alpha(c)$ (or equivalently, $\C_i \cong \C
{}^\gamma\times^{\eps(i)}_I \ppt$, where $\eps(i):\ppt \to I$ is
the embedding onto $i \in I$). Here we will also drop $\gamma$ from
notation when it is clear from the context. In particular, $I/i$
will stand for $I /^{\id} i$, that is, the category of objects $i'
\in i$ equpped with a map $i' \to i$. Dually, $i \setminus^\gamma \C
= (\C^o /^{\gamma^o} i)^o$ is the right comma-fiber of $\gamma$ over
$i$, it contains $\C_i$ as a full subcategory $\C_i \subset i
\setminus \C$, and $i \setminus I = i \setminus^{\id} I$ is the
category of objects $i' \in I$ equipped with a map $i \to i'$. For
any $i \in I$, we have the projection
\begin{equation}\label{comma.pr}
p_i:i \setminus^\gamma \C \to \C
\end{equation}
given by the forgetful functor, and for any map $f:i \to i'$,
composition with $f$ induced a functor
\begin{equation}\label{comma.tr}
f^*:i' \setminus^\gamma \C \to i \setminus^\gamma \C
\end{equation}
equipped with a canonical isomorphism $p_{i'} \circ f \cong p_i$.

A category $\C$ is a {\em groupoid} if all its morphisms are
invertible, {\em strongly discrete} if all its morphisms are
identity maps, and {\em discrete} if it is equivalent to a discrete
category. We say that a subcategory $\C' \subset \C$ in a category
$\C$ is {\em dense} if $\C$ and $\C'$ have the same objects (but
$\C'$ has less morphisms). Equivalently, a dense subcategory $\C_v
\subset \C$ is defined by a class $v$ of morphisms in $\C$ that is
closed under compositions and contains all the identity maps. For
any functor $\gamma:\C' \to \C$ and a dense subcategory $\C_v
\subset \C$ defined by a class of maps $v$, we denote by $\gamma^*v$
the class of maps $f$ in $\C'$ such that $\gamma(f) \in v$. For
every category $\C$, we have the strongly discrete dense subcategory
$\C_{\Id} \subset \C$ defined by the class $\Id$ of all identity
maps, and the dense {\em isomorphism groupoid} $\C_{\Iso} \subset
\C$ defined by the class $\Iso$ of all invertible maps. A functor
$\gamma:\C' \to \C$ from a groupoid $\C'$ by definition factors
uniquely through $\C_{\Iso}$, and if $\C'$ is strongly discrete, it
factors uniquely through $\C_{\Id}$. We denote by $\forall$ the
class of all maps in $\C$, so that $\C_\forall = \C$. As usual, a
functor $\gamma$ is {\em conservative} if $\gamma^*\Iso = \Iso$.

If $\C'$ is discrete, then any functor $\gamma:\C' \to \C$ is
isomorphic to a functor that factors through $\C_{\Id}$ but
unfortunately, this factorization is not unique (any two are
isomorphic as functors to $\C$ but not to $\C_{\Id}$). To aleviate
the problem, we say that a category $\C$ is {\em tight} if it has
exactly one object in each isomorphism class, we say that a {\em
  tighetening} of some $\C$ is a tight weakly full subcategory $\C'
\subset \C$ such that the embedding $\C' \to \C$ is an equivalence,
and we assume that tightenings exist (this requires slightly more of
the axiom of choice that one would like but that's life). If $\C$ is
tight, then any functor $\C' \to \C$ from a discrete $\C$ factors
uniquely through $\C_{\Id}$, and if $\C$ itself is discrete, it is
strongly discrete.

For any category $I$, its {\em category of arrows} $\Ar(I)$ has
arrows $f:i' \to i$ as objects, with morphisms from $f_0:i'_0 \to
i_0$ to $f_:i'_1 \to i_1$ given by commutative diagrams
\begin{equation}\label{ar.dia}
\begin{CD}
i'_0 @>{f_0}>> i_0\\
@V{g'}VV @VV{g}V\\
i'_1 @>{f_1}>> i_1.
\end{CD}
\end{equation}
The {\em twisted arrow category} $\Tw(I)$ has the same objects as
$\Ar(I)$, but morphisms are given by commutative diagrams
\begin{equation}\label{tw.dia}
\begin{CD}
i'_0 @>{f_0}>> i_0\\
@A{g'}AA @VV{g}V\\
i'_1 @>{f_1}>> i_1.
\end{CD}
\end{equation}
For any dense subcategory $I_v \subset I$ defined by a class of
morphisms $v$, we denote by $\Ar^v(I) \subset \Ar(I)$, $\Tw^v(I)
\subset \Tw(I)$ the full subcategories spanned by arrows in
$v$. Sending an arrow $i' \to i$ to its target $i$ gives functors
\begin{equation}\label{ar.tw.co}
t:\Ar(I) \to I, \qquad t:\Tw(I) \to I,
\end{equation}
and sending it to the source $i'$ gives functors
\begin{equation}\label{ar.tw}
s:\Ar(I) \to I, \qquad s:\Tw(I) \to I^o.
\end{equation}
We also have a functor
\begin{equation}\label{ar.eta}
\eta:I \to \Ar(I)
\end{equation}
sending an object $i \in I$ to $\id:i \to i$, and $\eta$ is fully
faithful, right-adjoint to $t$ and left-adjoint to $s$.

\begin{defn}\label{facto.def}
A {\em factorization system} on a category $I$ is a pair of classes
$\langle L,R \rangle$ defining dense subcategories $I_L,I_R \subset
I$ such that $L,R \supset \Iso$, any morphism $f:i' \to i$ in $I$
admits a factorization
\begin{equation}\label{l.r.facto}
\begin{CD}
i' @>{l}>> i'' @>{r}>> i
\end{CD}
\end{equation}
with $l \in L$, $r \in R$, and for any commutative square
$$
\begin{CD}
i_0 @>{f}>> i'_0\\
@V{l}VV @VV{r}V\\
i_1 @>{g}>> i'_1
\end{CD}
$$
in $I$ with $l \in L$, $r \in R$ there exists a unique morphism
$q:i_1 \to i_0'$ such that $f = q \circ l$ and $g= r \circ q$.
\end{defn}

\begin{exa}\label{facto.exa}
If $I = I_0 \times I_1$ is the product of two categories $I_0$,
$I_1$, with the projections $\pi_l:I \to I_l$, $l = 0,1$, then the
classes $\pi_0^*\Iso$ and $\pi_1^*\Iso$ form a factorization system
on $I$ (in either order).
\end{exa}

\begin{exa}\label{ar.facto.exa}
A factorization system $\langle L,R \rangle$ on a category $I$
defines a factorization system $\langle (s \times t)^*(L \times
L),(s \times t)^*(R \times R)\rangle$ on the arrow category $\Ar(I)$
that restricts to factorization systems on $\Ar^L(I),\Ar^R(I)
\subset \Ar(I)$. We have $\Ar^L(I)_{(s \times t)^*(L \times L)}
\cong \Ar(I_L)$ and $\Ar^R(I)_{(s \times t)^*(R \times R)} \cong
\Ar(I_R)$.
\end{exa}

For further details about factorization systems, we refer the reader
to \cite{bou}. In particular, Definition~\ref{facto.def} implies
that $L \cap R = \Iso$, either of the classes $L$, $R$ competely
determines the other one, and the factorization \eqref{l.r.facto} is
unique up to a unique isomorphism.

\begin{defn}\label{adm.def}
A full subcategory $I' \subset I$ is {\em left} resp.\ {\em
  right-admissible} if the embedding functor $\gamma:I' \to I$
admits a left resp.\ right-adjoint $\gamma_\dg:I \to I'$.
\end{defn}

\begin{exa}\label{adm.sub.exa}
If we have a left resp.\ right-admissible full subcategory $I'
\subset I$, and another subcategory $I'' \subset I$ that contains
$I'$ and the adjunction map $i \to \gamma(\gamma_\dg(i))$
resp.\ $\gamma(\gamma_\dg(i)) \to I$ for any $i \in I''$, then $I'
\subset I''$ is also left resp.\ right-admissible, with the same
$\gamma_\dg$.
\end{exa}

\begin{exa}\label{comma.exa}
Assume given a category $I$ equipped with a factorization system
$\langle L,R \rangle$, and for any $i \in I$, denote by $I /_R i
\subset I/i$ the full subcategory spanned by pairs $\langle
i',\alpha(i') \rangle$ with $\alpha(i') \in R$. Then $I /_R i
\subset I / i$ is left-admissible, with the adjoint functor sending
an arrow to the $R$-part of its decomposition \eqref{l.r.facto}, and
so is $I_L /_Ri \subset I_L / i$.
\end{exa}

\subsection{Functors.}

For any two categories $\C$, $\C'$ with $\C$ small, functors from
$\C$ to $\C'$ form a well-defined category that we denote by
$\Fun(\C,\C')$ (for example, if $\C=[1]$ is the single arrow
category with two objects $0$, $1$ and a single non-identity map $0
\to 1$, then $\Fun([1],\C') \cong \Ar(\C')$ is the arrow category of
the category $\C'$). If $\C$ is not small, one says that a functor
$\gamma:\C \to \C'$ is {\em continuous} if it commutes with filtered
colimits, and one says that a category $\C$ is {\em finitely
  presentable} if it has all filtered colimits and a set
$\overline{\C}$ of compact objects such that any $c \in \C$ is a
filtered colimit of objects in $\overline{\C}$. Then if $\C$ is
large but finitely presentable, continuous functors from $\C$ to
$\C'$ also form a well-defined category that we also denote by
$\Fun(\C,\C')$. For brevity, we will say that $\C$ is {\em bounded}
if it is either small or large and finitely presentable, and we will
say that a functor $\gamma:\C \to \C'$ is {\em bounded} if either
$\C$ is small or $\C$ is large and $\gamma$ is continuous, so that
in either case, $\Fun(\C,\C')$ denotes the category of bounded
functors. Whenever we consider a functor $\gamma:\C \to \C'$ with
bounded $\C$, it is assumed to be bounded unless indicated
otherwise. We have the {\em evaluation functor}
\begin{equation}\label{ev.eq}
\ev:\C \times \Fun(\C,\C') \to \C'
\end{equation}
with the obvious universal property: for any bounded category $\C''$
and bounded functor $\gamma:\C \times \C'' \to \C'$, there exists a
bounded functor $\gamma':\C'' \to \Fun(\C,\C')$ and an isomorphism
$\alpha:\gamma \cong \ev \circ (\Id \times \gamma')$, and the pair
$\langle \gamma',\alpha \rangle$ is unique up to a unique
isomorphism.

We say that a functor $\gamma:\C \to \C'$ {\em inverts a map $f$} if
$\gamma(f)$ is invertible, and we say that $\gamma$ is {\em locally
  constant} if it inverts all maps, and {\em constant} if up to an
isomorphism, it factors through the tautological projection $\C \to
\ppt$. For any class of maps $v$ in a bounded $\C$, we denote by
$\Fun^v(\C,\C') \subset \Fun(\C,\C')$ the full subcategory spanned
by functors that invert maps in $v$. In particular,
$\Fun^\forall(\C,\C')$ is the category of locally constant
functors. For any subclass $w \subset v$, we have $\Fun^v(\C,\C')
\subset \Fun^w(\C,\C')$, and we say that $w$ is {\em dense} in $v$
if the inclusion is an equivalence for any $\C'$.

For any bounded categories $I$, $I'$ and any category $\E$, a
bounded functor $\gamma:I \to I'$ induces a pullback functor
$\gamma^*:\Fun(I',\E) \to \Fun(I,\E)$, $F \mapsto F \circ
\gamma$. For any dense subcategory $I_v \in I$ defined by a class of
maps $v$, this restricts to a functor
\begin{equation}\label{pb.v}
\gamma^*:\Fun^v(I,\E) \to \Fun^{\gamma^*v}(I',\E).
\end{equation}
If $\gamma$ is essentially surjective, then \eqref{pb.v} is faithful
and conservative for any $\E$. It is often useful to strengthen this
in the following way.

\begin{defn}\label{loc.def}
A functor $\gamma:I' \to I$ between bounded categories $I$, $I'$ is
a {\em localization} if \eqref{pb.v} is an equivalence for any $v$
and $\E$. A category $I$ is {\em simply connected} if the
tautological projection $\tau:I \to \ppt$ is a localization.
\end{defn}

\begin{exa}\label{opp.exa}
If $\gamma:I_0 \to I_1$ is a localization of small categories, then
so is the opposite functor $\gamma^o:I_0^o \to I_1^o$ (just replace
$\E$ with $\E^o$ in \eqref{pb.v}). In particular, a small category
$I$ is simply connected if and only if so is $I^o$.
\end{exa}

\begin{exa}\label{adm.exa}
If a bounded full subcategory $I' \subset I$ in a bounded category
$I$ is left or right-admissible in the sense of
Definition~\ref{adm.def}, then the adjoint functor $\gamma_\dg:I'
\to I$ is a localization in the sense of
Definition~\ref{loc.def}. Moreover, any class of maps $v$ in $I'$
that contains the adjunction maps $i \to \gamma_\dg(\gamma(i))$
resp.\ $\gamma_\dg(\gamma(i)) \to i$ for any $i \in I'$ is dense in
$\gamma_\dg^*\Iso$.
\end{exa}

We of course assume known the notions of a limit and colimit, but to
fix notation, it is convenient to introduce the following.

\begin{defn}\label{f.aug.def}
An {\em augmentation} of a functor $E:I \to \E$ between some
categories $I$, $E$ is a functor $E_>:I^> \to \E$ equipped with an
isomorphism $E_>|_I \cong E$.  A {\em $e$-augmentation}, for some $e
\in \E$, is an augmentation $E_>$ equipped with an isomorphism
$E_>(o) \cong e$.  An augmentation $E_>$ is {\em exact} if for any
other augmentation $E_>'$, there exists a unique morphism of
augmentations $E_> \to E'_>$.
\end{defn}

Equivalently, an $e$-augmentation of a functor $E:I \to \E$ is given
by a morphism $E \to e_I$ to the constant functor $e_I:I \to \E$
with value $e$, or by a factorization of $E$ into a composition
\begin{equation}\label{aug.lft}
\begin{CD}
I @>>> \E / e @>{\phi}>> \E,
\end{CD}
\end{equation}
where $\phi:\E / e \to \E$ is the forgetful functor sending $e' \to
e$ to $e'$. By definition, $E$ admits an exact augmentation $E_>$ if
and only if the colimit $\colim_IE$ exists, and in this case, we
have $\colim_IE = E_>(o)$. A category $\E$ is {\em cocomplete} if
$\colim_IF$ exists for any bounded $I$ and $F \in \Fun(I,F)$.

\begin{exa}\label{pi0.exa}
The category $\Sets$ is cocomplete. For any bounded category $I$, we
denote $\pi_0(I) = \colim_I\ppt$, where $\ppt:I \to \Sets$ is the
constant functor that sends everything to a one-element
set. Informally, $\pi_0(I)$ is the set of connected components of
the category $I$, and $I$ is {\em connected} if
$\pi_0(I)\cong\ppt$.
\end{exa}

\subsection{Kan extensions.}

For any functor $\gamma:I \to I'$ between bounded categories, and
any category $\E$, the {\em left Kan extension} $\gamma_!E$ of a
functor $E \in \Fun(I,\E)$ is a functor $\gamma_!E \in \Fun(I',\E)$
equipped with a map $a:E \to \gamma^*\gamma_!E$ such that the
composition map
\begin{equation}\label{kan.dia}
\begin{CD}
\Hom(\gamma_!E,E') @>{\gamma^*}>>
\Hom(\gamma^*\gamma_!E,\gamma^*E) @>{- \circ a}>>
\Hom(E,\gamma^*E')
\end{CD}
\end{equation}
is an isomorphism for any $E' \in \Fun(I',\E)$. The left Kan
extension is unique and functorial in $E$, if it exists.

\begin{exa}\label{kan.loc.exa}
If $\gamma$ is a localization in the sense of
Definition~\ref{loc.def}, then $\gamma_!E$ exists for any $\E \in
\Fun^{\gamma^*v}(I,\E)$, for any class $v$, and $\gamma_!$ provides
the inverse equivalence to \eqref{pb.v}.
\end{exa}

If $\gamma_!E$ exists for any $E \in \Fun(I,\E)$, then
$\gamma_!:\Fun(I,\E) \to \Fun(I',\E)$ is the functor left-adjoint to
$\gamma^*$. A degenerate example is the tautological projection
$\tau:I \to \ppt$: in this case, $\Fun(\ppt,E) \cong \E$, providing
a pair $\langle e,a \rangle$ of an object $e \in E$ and a map $a:E
\to \tau^*e$ is equivalent to providing a $e$-augmentation $E_>:I^>
\to \E$ of the functor $E$, and $\langle e,a \rangle$ turns $e$ into
the a Kan extension if and only if $E_>$ is exact, so that we simply
have $\tau_!E = \colim_IE$. In general, $\gamma_!E$ is given by
\begin{equation}\label{kan.eq}
\gamma_!E(i) = \colim_{I /^\gamma i}E|_{I /^\gamma i}, \qquad E
\in \Fun(I,\E), \ i \in I',
\end{equation}
where $I /^\gamma i$ are the left comma-fibers, and it exists iff so
do the colimits in the right-hand side.

To check existence, it is often convenient to observe the
following. For any bounded category $I$ and left-admissible
subcategory $I' \subset I$, with the embedding functor $\gamma:I'
\to I$ and its left-adjoint $\gamma_\dg:I \to I'$, the pullback
$\gamma_\dg^*$ is right-adjoint to the pullback $\gamma^*$, so that
for any $E \in \Fun(I,\E)$, we have
\begin{equation}\label{adm.eq}
\colim_I E \cong \colim_{I'} E|_{I'}
\end{equation}
by adjunction, and the left-hand side exists iff so does the
right-hand side. In many practical cases, this helps to compute left
Kan extensions. Namely, for any functor $\gamma:I \to I'$, a
subcategory $I_0 \subset I'$ and an object $i \in I'$, say that an
$i$-augmentation of the induced functor $\gamma_0:I_0 \to I'$ is
{\em $\gamma$-admissible} if the embedding $I_0 \to I /^\gamma i$ of
\eqref{aug.lft} is left-admissible, and define a {\em framing} of
the functor $\gamma$ as a collection of subcategories $I(i) \subset
I$, one for each $i \in I'$, equipped with $\gamma$-admissible
$i$-augmentations $\gamma(i)_>$ of the functors $\gamma(i) =
\gamma|_{I(i)}$.

\begin{remark}\label{fr.rem}
If $\gamma:I \to I'$ is fully faithful, then for any $i \in I$, $I
/^\gamma \gamma(i) \cong I / i$ has a terminal object. Thus in this
case, to define a framing of $\gamma$, it suffices to consider
objects $i \in I'$ that are not in $\gamma(I) \subset I'$, and then
let $I(\gamma(i))$, $i \in I$ be the point category $\ppt$ embedded
onto $i \in I$.
\end{remark}

\begin{lemma}\label{kan.le}
Assume given bounded categories $I$, $I'$, a category $\E$, and a
bounded functor $\gamma:I \to I'$. Moreover, assume given a framing
of $\gamma$.
\begin{enumerate}
\item For any $E \in \Fun(I,\E)$, $\gamma_!E$ exists if and only if
  for any $i \in I$, $E|_{I(i)}:I(i) \to \E$ admits an exact augmentation.
\item For any $E \in \Fun(I',\E)$, $E$ itself with the identity map
  $\gamma^*E \to \gamma^*E$ is a left Kan extension
  $\gamma_!\gamma^*E$ if and only if for any $i \in I$, the
  augmentation $\gamma(i)_>^*E$ of the functor $\gamma(i)^*E$ is
  exact.
\end{enumerate}
\end{lemma}

\proof{} Clear. \endproof

\begin{exa}\label{cl.exa}
Say that a full embedding $\gamma:I' \to I$ is {\em right-closed} if
$I(i',i)$ is empty for any $i' \in I'$ and $i \in I \setminus I'$,
and dually, say that $\gamma$ is {\em left-closed} if $\gamma^o$ is
right-closed (equivalently, $\gamma$ is left resp.\ right-closed if
$I' \subset I$ is the fiber $I_0$ resp.\ $I_1$ of a functor $I \to
[1]$). Then for a right-closed $\gamma$, all the comma-fibers $I'
/^\gamma i$, $i \in I \setminus I'$ in \eqref{kan.eq} are
empty. Therefore by Remark~\ref{fr.rem}, the left Kand extension
$\gamma_!E$ exists for any target category $\E$ with an initial
object $o$, and coincides with the canonical extension $E^<:I \to
\E$ given by $E^<(i) = E(i)$ if $i \in I'$ and $o$ otherwise. A
common example of a right-closed embedding is the embedding $I
\subset I^<$ for any bounded $I$.
\end{exa}

Dually, the {\em right Kan extension} $\gamma_*E$ is a pair $\langle
\gamma_*E,a \rangle$, $\gamma_*E \in \Fun(I',\E)$,
$a:\gamma^*\gamma_*E \to E$ satisfying the universal property dual
to that of $\gamma_!E$, it is unique and functorial in $E$, and when
it exists for any $E$, it provides a functor $\gamma_*:\Fun(I',\E)
\to \Fun(I,\E)$ right-adjoint to $\gamma^*$. In fact, if $I$ and
$I'$ are small, so that there is no need to control boundedness, we
have $(\gamma_*)^o \cong \gamma_!E^o$, so that everything including
Lemma~\ref{kan.le} can be applied to the right Kan extensions simply
by passing to the opposite functors. For large bounded categories,
one can still compute right Kan extension $\gamma_*$ by using an
obvious dual version of \eqref{kan.eq} with limits and right
comma-fibers, possibly combined with a framing of the opposite
functor $\gamma^o$, but one has to check that the result is a
bounded functor. This is automatic in the situation of
Example~\ref{kan.loc.exa} (and in fact, in this case we have
$\gamma_!E \cong \gamma_*E$). This is also automatic in
Example~\ref{cl.exa}, with appropriate dualization: $\gamma:I' \to
I$ has to be leftt-closed, $o \in \E$ must be the terminal
object rather than the initial one, and we denote the canonical
extension by $E^>$.

Let us now mention that quite often, the isomorphism \eqref{adm.eq}
holds for full embeddings that are not left-admissible. Here is one
example.

\begin{defn}\label{cof.def}
A full subcategory $I' \subset I$ is {\em cofinal} if for any $i \in
I$, the right comma-fiber $i \setminus I'$ is simply connected in
the sense of Definition~\ref{loc.le}.
\end{defn}

\begin{remark}\label{cof.rem}
It is enough to check the condition of Definition~\ref{cof.def} for
$i \in I \setminus i$ --- indeed, $i \setminus I'$ for $i \in I'$
has an initial object $\id:i \to i$, thus it is tautologically
simply connected by Example~\ref{adm.exa}.
\end{remark}

A left-admissible embedding $\gamma:I' \to I$ is automatically
cofinal ($i \setminus I$ has an initial object given by the
adjunction map $i \to \gamma(\gamma_\dg(i))$). The converse is not
true. Nevertheless, for any cofinal full embedding $\gamma:I' \to
I$, the dual version of \eqref{kan.eq} immediately shows that
$\gamma_*E$ exists for any constant functor $E \in \Fun(I',\E)$, and
the adjunction map $E \to \gamma_*\gamma^*E$ is an isomorphism for
any constant functor $E \in \Fun(I,\E)$. Therefore for any $E \in
\Fun(I,\E)$, we still have the isomorphism \eqref{adm.eq}, and its
source exists iff so does its target. Even more generally,
$\gamma_*E$ exists for a locally constant $E$, and $\gamma_*$
induces an equivalence
\begin{equation}\label{forall.I}
\gamma_*:\Fun^\forall(I',\E) \cong \Fun^\forall(I,\E),
\end{equation}
so that $I$ is simply connected if and only if so is $I'$.

\begin{lemma}\label{cof.le}
Assume given full embeddings $I'' \subset I' \subset I$ of bounded
categories. Then $I'' \subset I$ is cofinal if and only if $I''
\subset I'$ and $I' \subset I$ are cofinal.
\end{lemma}

\proof{} If $I'' \subset I$ is cofinal, then $I'' \subset I'$ is
tautologically cofinal (the right comma-fiber $i \setminus I''$ for
$i \in I'$ does not change if we take $i$ as an object of $I$).  But
if $I'' \subset I'$ is cofinal, then the embedding $i \setminus I''
\subset i \setminus I'$ is cofinal for any $i \in I$ --- indeed, its
right comma-fibers are the same as for $I'' \subset I'$ --- and we
are done by \eqref{forall.I}.
\endproof

\begin{remark}
While a connected category $I$ need not be simply connected, the
pullback $\tau^*:\E \to \Fun(I,\E)$ with respect to the tautological
projection $\tau:I \to \ppt$ is still fully faithful for any target
category $\E$, so that $\colim_IE = \tau_!E$ and $\lim_IE = \tau_*E$
exist for any constant functor $E:I \to \E$, and we have $\tau_!E
\cong \tau_*E$, $E \cong \tau^*\tau_!E \cong \tau^*\tau_*E$. Thus if
one only wants \eqref{adm.eq}, one can get away with a weaker notion
of cofinality: it suffices to require that the right comma-fibers $i
\setminus I'$ are connected (as done in e.g.\ \cite[Chapter
  2.5]{KS}). However, we also need \eqref{forall.I}, so we use a
stronger definition. In fact, all the cofinal embeddings in the
paper will be also homotopy cofinal in the sense of
Subsection~\ref{ho.subs} below.
\end{remark}

\subsection{Homological algebra.}\label{hom.subs}

For any commutative ring $k$, we denote by $k\amod$ the category of
$k$-modules, and for any bounded category $I$, we simplify notation
by writing $\Fun(I,k)=\Fun(I,k\amod)$. This is an abelian category
with enough projective and injectives, and we denote its derived
category by $\D(I,k)$. We shorten $\D(\ppt,k)$ to $\D(k)$. We have a
natural functor
\begin{equation}\label{D.ho}
\D:\D(I,k) \to \Fun(I,\D(k)),
\end{equation}
and for any dense subcategory in $I$ defined by a class of morphisms
$v$, we denote by $\D^v(I,k) \subset \D(I,k)$ the full subcategory
spanned by objects $E$ with $\D(E) \in \Fun^v(I,\D(k))$. We recall
that explicitly, objects in $\D(I,k)$ are represented by chain
complexes $M_\idot$ in the abelian category $\Fun(I,k)$, and
$\D(I,k)$ has a natural $t$-structure whose term $\D^{\leq 0}(I,k)
\subset \D(I,k)$ is spanned by complexes concentrated in
non-negative homologocal degrees. For any bounded functor $\gamma:I'
\to I$, the pullback functor $\gamma^*$ descends to a functor
$\gamma^*:\D(I,k) \to \D(I',k)$ that has a left and a right-adjoint
functors $L^\hdot\gamma_!,R^\hdot\gamma_*:\D(I',k) \to \D(I,k)$
obtained by taking the derived functors of the left and right Kan
extensions. For any $i \in I$, with the embedding $\eps(i):\ppt \to
I$ onto $i$, the left Kan extension $\eps(i)_!$ is exact, and the
representable functor $k_i = \eps(i)_!k \in \Fun(I,k)$ is given by
\begin{equation}\label{k.i}
k_i(i') = k[I(i,i')], \qquad i' \in I,
\end{equation}
where for any set $S$, $k[S]$ is the free $k$-module generated by
$S$. Sending $i \in I$ to $k_i$ defines the fully faithful Yoneda
embedding
\begin{equation}\label{yo.k}
\Y:I^o \to \Fun(I,k).
\end{equation}
It may happen that the left Kan extension functor $\gamma_!$ itself
admits a left-ad\-joint; here is a useful example of such a
situation.

\begin{exa}\label{D.cl.exa}
Assume given a full embedding $\gamma:I' \to I$ right-closed in the
sense of Example~\ref{cl.exa}. Then $\gamma_!$ is given by extension
by $0$, thus $\gamma_!$ exact, and it in fact has a left-adjoint
$\gamma^!$ given by
\begin{equation}\label{ga.shriek}
\gamma^! = \Y_!\gamma^o_*\Y':\Fun(I,k) \to \Fun(I',k),
\end{equation}
where $\Y$, $\Y'$ are the Yoneda embeddings \eqref{yo.k} for the
categories $I$, $I'$, and the right Kan extension $\gamma^o_*$ is
again given by extension by $0$. The derived functor
$L^\hdot\gamma^!$ is then left-adjoint to $\gamma_! \cong
L^\hdot\gamma_!$.
\end{exa}

For any functor $M \in \Fun(I,k)$, the homology of the category $I$
with coefficients in $M$ is obtained by taking the total derived
functor $L^\hdot\colim_IM$ of the colimit functor $\colim_I$: we set
\begin{equation}\label{i.m}
C_\idot(I,M) = L^\hdot\colim_I M \in \D(k),
\end{equation}
and we denote by $H_\idot(I,M)$ the homology modules of the object
\eqref{i.m}. Both \eqref{adm.eq} and \eqref{kan.eq} have obvious
counterparts for homology and derived left Kan extensions, and so
does Lemma~\ref{kan.le}, so that one can compute derived Kan
extensions by choosing a framing. For a useful generalization of
\eqref{i.m}, say that a {\em $k$-valued bifunctor} on $I$ is a
functor $M \in \Fun(I^o \times I,k)$. Then for any map $f:i \to i'$
in $I$, we have a natural map
\begin{equation}\label{d.f}
d_f = M(f^o \times \id) \oplus (-M(\id \times f)):M(i',i) \to M(i,i)
\oplus M(i',i'),
\end{equation}
and we can define the {\em trace} $\Tr_I(M)$ of the bifunctor $M$ by
the exact sequence
\begin{equation}\label{tr.I}
\begin{CD}
\bigoplus_{f \in I(i,i')}M(i',i) @>{d}>> \bigoplus_{i \in
  I} @>>> \Tr_I(M) @>>> 0,
\end{CD}
\end{equation}
where the sum on the left is over all maps $f$ in $I$, and $d$ is
the sum of the maps \eqref{d.f} ($\Tr_IM$ is also known as the
``coend'' of the functor $M$, see e.g.\ \cite[IX.6]{Mcbook}). Then
$\Tr_I:\Fun(I^o \times I,k) \to k\amod$ is right-exact, with the
total derived functor $L^\hdot\Tr_I$, and we can define the {\em
  bifunctor homology object} $CH_\idot(I,M) \in \D(k)$ by
\begin{equation}\label{c.i.m}
CH_\idot(I,M) = L^\hdot\Tr_I(M).
\end{equation}
Bifunctor homology modules $HH_\idot(I,M) \in k\amod$ are then the
homology modules of the object \eqref{c.i.m}. In particular,
whenever we have two functors $N \in \Fun(I^o,k)$, $M \in
\Fun(I,k)$, we can define their box product $N \boxtimes_k M$ in
$\Fun(I^o \times I,k)$ by $N \boxtimes_k M(i \times i') = N(i)
\otimes_k M(i')$, and let
\begin{equation}\label{ot.I}
N \otimes_I M = \Tr_I(N \boxtimes_k M).
\end{equation}
Then deriving this in either of the two variables $N$, $M$ gives the
derived tensor product $N \lotimes_I M \in \D(k)$, and as soon as
either $N$ or $M$ is pointwise-flat --- that is, takes values in the
full subcategory $k\amod^{fl} \subset k\amod$ spanned by flat
$k$-modules --- we have a natural identification
\begin{equation}\label{lot.I}
N \lotimes_I M \cong CH_\idot(N \boxtimes_k M).
\end{equation}
Since any object in the derived category $\D(I,k)$ can be
represented by a complex of pointwise-flat $k$-modules, the product
\eqref{lot.I} extends to derived categories. The extended product
admits a covariant version of the adjunction isomorphism: for any
functor $\gamma:I' \to I$ from some bounded $I'$, we have a natural
identification
\begin{equation}\label{yo.cov}
E \otimes_I L^\hdot\gamma_!E' \cong \gamma^{o*}E \otimes_{I'} E'
\end{equation}
for any $E' \in \D(I',k)$, $E \in \D(I^o,k)$. In particular,
for every $M \in \Fun(I,k)$ and $i \in I$, we have
a natural identification
\begin{equation}\label{yo.M}
k_i \lotimes_I M \cong k_i \otimes_I M \cong M(i),
\end{equation}
where $k_i \in \Fun(I^o,k)$ is the representable functor of
\eqref{k.i} (that is automatically pointwise-flat). On the other
hand, the constant functor $k \in \Fun(I^o, k)$ with value $k$ is
also pointwise-flat, and we have
\begin{equation}\label{C.lot}
CH_\idot(I,\pi^*M) = k \lotimes_I M \cong C_\idot(I,M),
\end{equation}
where $\pi:I^o \times I \to I$ is the projection. In this way, the
functor homology \eqref{i.m} can be expressed as bifunctor homology
\eqref{c.i.m}. In the other direction, we have a canonical
quasiisomorphism
\begin{equation}\label{CH.Tw}
CH_\idot(I,M) \cong C_\idot(\Tw(I),(s \times t)^*M),
\end{equation}
where $\Tw(I)$ is the twisted arrow category, and $s \times t:\Tw(I)
\to I^o \times I$ is the product of the projections \eqref{ar.tw}
and \eqref{ar.tw.co}.

\begin{exa}
For a useful application of \eqref{C.lot}, assume given a bounded
category $I$ and a functor $I \to [1]$ with fibers $I_0,I_1 \subset
I$, as in Example~\ref{cl.exa}, and define the {\em homology of $I$
  with support in $I_0$} with coefficients in some $E \in \D(I,k)$
as
\begin{equation}\label{C.z}
C_\idot(I,I_1,E) = C_\idot(I,L^\hdot j_1^!E), \qquad
H_\idot(I,I_1,E) = H_\idot(I,L^\hdot j_1^!E),
\end{equation}
where $j_1^!$ is the functor \eqref{ga.shriek} of
Example~\ref{D.cl.exa} for the embedding $j_1:I_1 \to I$. Then we
have
\begin{equation}\label{C.z.L}
C_\idot(I,I_1,E) \cong j^o_{1*}k \lotimes_I E,
\end{equation}
where $k$ is the constant functor with value $k$, and $j^o_{1*}$ is
exact and given by extension by $0$. Moreover, if we denote by
$j_0:I_0 \to I$ the other embedding, then $j^o_{0!}$ is also exact
and given by extension by $0$, and we have a short exact sequence
$$
\begin{CD}
0 @>>> j^o_{0!}k @>>> k @>>> j^o_{1*} k @>>> 0
\end{CD}
$$
that gives rise to a distinguished triangle
\begin{equation}\label{aug.tria}
\begin{CD}
C_\idot(I_0,j_0^*E) @>>> C_\idot(I,E) @>>> C_\idot(I,I_1,E) @>>>
\end{CD}
\end{equation}
of homology complexes and the corresponding long exact sequence of
homology groups. For example, for any bounded category $I$, we have
the projection $I^> \to [1]$ sending $I \subset I^>$ to $0$ and the
new terminal object $o \in I^>$ to $1$; then for any augmented
functor $E \in \Fun(I^>,k)$, we have $C_\idot(I^>,E) \cong E(o)$,
the first map in \eqref{aug.tria} is the augmentation map, and it is
a quasiisomorphism if and only if $H_\idot(I^>,\{o\},E)=0$.
\end{exa}

We say that a category $I$ is {\em $k$-linear} if it is enriched
over $k\amod$ --- that is, all the sets $I(i,i')$ are equipped with
a $k$-module structure so that the composition maps are bilinear. We
say that a $k$-linear category $I$ is {\em flat} if so are all the
$k$-modules $I(i,i')$. For any $k$-linear bounded category $I$, we
denote by $\Fun_k(I,k) \subset \Fun(I,k)$ be the full subcategory
spanned by $k$-linear functors, and we let $\D(I)$ be its derived
category. More generally, a small {\em DG category} $I_\idot$ is a small
category $I$ enriched over the category $C_\idot(k\amod^{fl})$ of
chain complexes of flat $k$-modules --- that is, we are given a
complex $I_\idot(i,i')$ for any $i,i \in I'$ with
$I_0(i,i')=I(i,i')$ as sets, and the unital and associative
$k$-linear composition maps $I_\idot(i,i') \otimes_k I_\idot(i',i'')
\to I_\idot(i,i'')$ that extend the compositions in $I$. We assume
known that every DG category $I_\idot$ has a derived category
$\D(I_\idot)$ of DG modules with the standard properties that can be
found for example in \cite{kel}.

For any $k$-linear category $I$, we let $\Fun_k(I^o \times I,k)
\subset \Fun(I^o \times I,k)$ be the full subcategory spanned by
bifunctors that are $k$-linear in each of the two variables, and we
let $\Tr^k_I$ be the restriction of the functor \eqref{tr.I} to
$\Fun_k(I^o \times I,k)$. For any $M \in \Fun_k(I^o \times I,k)$,
the {\em $k$-linear bifunctor homology object} of $I$ with
coefficients in $M$ is defined by
$$
CH_\idot(I/k,M) = L^\hdot\Tr^k_I(M) \in \D(k),
$$
and the {\em $k$-linear bifunctor homology modules} $HH_\idot(I/k,M)
\in k\amod$ are its homology modules. The embedding $\Fun_k(I^o
\times I,k) \subset \Fun(I^o \times I,k)$ then induces a functorial
map
\begin{equation}\label{c.i.m.aug}
CH_\idot(I,M) \to CH_\idot(I/k,M)
\end{equation}
for any $M \in \Fun_k(I^o \times I,k)$. If $I$ is flat, we also have
the obvious $k$-linear version $- \otimes_{I/k} -$ of the product
\eqref{ot.I}, with its derived version $- \lotimes_{I/k} -$, and the
identification \eqref{lot.I} if $M$ or $N$ is pointwise-flat. In
particular, $I(-,i):I^o \to k\amod$ is then $k$-linear and
pointwise-flat, and we also have the $k$-linear version of the
Yoneda isomorphism \eqref{yo.M}.

\begin{exa}\label{tr.A.exa}
A $k$-linear category $I$ with one object is the same thing as an
associative unital $k$-algebra $A$, and it is flat iff so is
$A$. The category $\Fun_k(I^o \times I,k)$ is then the category
$A\bimod$ of $A$-bimodules --- that is, left modules over $A^o
\otimes_k A$ --- and $\Tr_I(M)$ coincides with $\Tr_A(M) = M/[A,M]$,
where $[A,M] \subset M$ is the $k$-submodule spanned by commutators
$am-ma$, $a \in A$, $m \in M$. The $k$-linear Yoneda isomorphism
\eqref{yo.M} reads as
\begin{equation}\label{yo.A}
L^\hdot\Tr_A(A^o \otimes_k M) \cong \Tr_A(A^o \otimes_k M) \cong M,
\end{equation}
for any flat $k$-algebra $A$ and left $A$-module $M$.
\end{exa}

\section{Grothendieck construction.}

\subsection{Fibrations and cofibrations.}\label{fib.subs}

We assume known the machinery of Grothendieck fibrations and
cofibrations of \cite{sga}. As a reminder, a functor $\gamma:\C \to
I$ is a {\em precofibration} if for any $i \in I$, the tautological
embedding $\C_i \subset \C/i$ admits a left-adjoint functor
\begin{equation}\label{beta.eq}
\zeta(i):\C/i \to \C_i.
\end{equation}
In this case, for any morphism $f:i' \to i$ in the category $I$, we
have an embedding $\C_{i'} \to \C/i$, $c \mapsto \langle c,f
\rangle$, and composing it with $\zeta(i)$ gives a functor
$f_!:\C_{i'} \to \C_i$ known as the {\em transition functor} of the
precofibration $\gamma$. A morphism $g:c' \to c$ in $\C$ is {\em
  cocartesian} over $I$ if $\zeta(\gamma(c))$ inverts its natural
lifting to a morphism in $\C/c$. For a composable pair of maps $f$,
$f'$, the adjunction provides a natural map
\begin{equation}\label{fib.compo}
(f \circ f')_! \to f_! \circ f'_!,
\end{equation}
and a precofibration is a {\em cofibration} if all these maps are
isomorphisms.

\begin{exa}\label{cyl.exa}
The {\em cylinder} $\Cyl(\gamma)$ of a functor $\gamma:\C_0 \to
\C_1$ is the category whose objects are those of $\C_0$ and $\C_1$,
and whose morphisms are given by
\begin{equation}\label{cyl.eq}
\Cyl(\gamma)(c,c') = \begin{cases}
\C_i(c,c'), &\quad c,c' \in \C_i, i = 0,1\\
\C_1(\gamma(c),c'), &\quad c \in \C_0,c' \in \C_1\\
\emptyset &\quad\text{otherwise}.
\end{cases}
\end{equation}
If $\tau:I \to \ppt$ is the tautological projection from some
category $I$ to the point, then $\Cyl(\tau) \cong I^>$, and if
$\id:\ppt \to \ppt$ is the identity functor from the point category
to itself, then $\Cyl(\id) \cong \ppt^> \cong [1]$. In general, we
have the projection $\chi:\Cyl(\gamma) \to [1]$ with fibers $\C_0$,
$\C_1$, and it is a cofibration with transition functor $\gamma$.
Conversely, for any cofibration $\C \to [1]$ with transition functor
$\gamma:\C_0 \to \C_1$, we have $\C \cong \Cyl(\gamma)$.
\end{exa}

\begin{exa}\label{gt.exa}
Assume given a cofibration $\pi:\C \to I$. Then its canonical
extension $\pi^>:\C^> \to I^>$ is also a cofibration; its fiber
$\C^>_o$ over the new terminal object $o \in I^>$ consists of the
terminal object $o \in \C^>$, and the transition functor $\C^>_i
\cong \C_i \to \C^>_o = \ppt$ corresponding to any $i \in I$ is the
tautological projection.

More generally, extending $\pi$ to a precofibration $\pi':\C' \to
I^>$ is equivalent to giving a category $\E = \C'_o$ and a functor
$\gamma = \zeta(o):\C \to \E$; for $\C^>$, one takes $\ppt$ as $\E$
and the tautological projection $\C \to \ppt$ as $\gamma$.
\end{exa}

A functor $\C \to I$ is a {\em fibration} if the opposite functor
$\C^o \to I^o$ is a cofibration. Thus a cofibration $\C \to I$
defines the opposite fibration $\C^o \to I^o$. It also defines the
transpose fibration $C^\perp \to I^o$ that has the same fibers
$\C^\perp_i =\C_i$, $i \in I$, the same transition functors $f^* =
f_!:\C_i \to \C_{i'}$ for any map $f:i \to i'$ in $I$, and the
isomorphisms \eqref{fib.compo} inverse to the same isomorphisms for
$\C$. Dually, a fibration $\C \to I$ defines the opposite
cofibration $\C^o \to I^o$ and the transpose cofibration $\C_\perp
\to I^o$. Moreover, a fibration $\C \to I$ can also be a
cofibration; this happens if all the transition functors $f^*$ have
left-adjoint functors $f_!$ (and these are then the cofibration
transition functors). In such a case, one says that $\C \to I$ is a
{\em bifibration}.

\begin{exa}\label{tw.exa}
  For any category $I$, the projections \eqref{ar.tw.co} are
  cofibrations, projection $\Ar(I) \to I$ of \eqref{ar.tw} is a
  fibration, and the projection $\Tw(I) \to I^o$ is its transpose
  cofibration.
\end{exa}

\begin{exa}\label{pb.exa}
For any functor $\gamma:I' \to I$, and any cofibration,
resp.\ fibration, resp.\ bifibration $\C \to I$, the pullback
$\gamma^*\C \to I'$ is a cofibration resp.\ fibration
resp.\ bifibration.
\end{exa}

\begin{exa}
A full embedding $I' \to I$ is a fibration resp.\ cofibration if and
only if it is left resp.\ right-closed in the sense of
Example~\ref{cl.exa}.
\end{exa}

For any category $I$ and categories $\C$, $\C'$ equipped with
functors $\pi:\C \to I$, $\pi':\C' \to I$, a {\em functor from $\C$
  to $\C'$ over $I$} is a pair of a functor $\gamma:\C \to \C'$ and
an isomorphism $\alpha(\gamma):\pi' \circ \gamma \cong \pi$. If
$\gamma$ has a left-adjoint $\gamma_\dg:\C \to \C'$, then we have a
map $\alpha_\dg(\gamma):\pi' \to \pi \circ \gamma_\dg$ adjoint to
$\alpha(\gamma)$, and we say that $\gamma_\dg$ is a {\em
  left-adjoint over $I$} if $\alpha_\dg(\gamma)$ is an
isomorphism. Dually, a right-adjoint $\gamma^\dg$ is a {\em
  right-adjoint over $I$} if $\alpha^\dg(\gamma):\pi \circ
\gamma^\dg \to \pi'$ is an isomorphism.

If $\C$ is bounded, and we are given another category $\phi:\E \to
I$ over $I$, then bounded functors from $\C$ to $\E$ over $I$ form a
well-defined category $\Fun_I(\C,\E)$. If $I$ itself is bounded, we
can let $\C = I$; then $\Sec(I,\E) = \Fun_I(I,\E)$ is the category
of sections $I \to \E$ of the projection $\phi$. If we have a functor
$\gamma:\C' \to \C$ over $I$ between bounded categories $\pi:\C \to
I$, $\pi':\C' \to I$, then just as in the absolute case, we have the
pullback functor $\gamma^*:\Fun_I(\C',\E) \to \Fun_I(\C,\E)$. For
any $E \in \Fun_I(\C,\E)$, the {\em left Kan extension $\gamma^I_!E$
  over $I$} is a pair of a functor $\gamma^I_!E \in \Fun_I(\C',\E)$
and a map $a:E \to \gamma^*\gamma^I_!E$ such that \eqref{kan.dia} is
an isomorphism for any $E' \in \Fun_I(\C',\E)$. Just as in the
absolute case, $\gamma_!^IE$ is unique and functorial, if it exists.

If $\phi:\E \to I$ is a cofibration, then $\gamma^I_!E$ can be
computed by a relative version of \eqref{kan.eq}. Namely, for any $c
\in \C'$, the functor $E$ defines a functor $E|_{\C /^\gamma c}:\C
/^\gamma c \to \E /^\phi \pi'(c)$ sending $\langle c',\alpha(c')
\rangle \in \C /^\gamma c$ to $E(c')$ equipped with the composition
map
$$
\begin{CD}
\phi(E(c')) @>{\alpha(E)(c')}>> \pi(c')
@>{\alpha(\gamma)(c')^{-1}}>> \pi'(\gamma(c'))
@>{\pi'(\alpha(c'))}>> \pi'(c),
\end{CD}
$$
and we then have
\begin{equation}\label{kan.I.eq}
\gamma^I_!E(c) = \colim_{\C /^\gamma c}\zeta(\pi'(c)) \circ
E|_{\C /^\gamma c}, \quad\quad E \in \Fun_I(\C,\E), \ c \in \C',
\end{equation}
where $\zeta(-)$ is the functor \eqref{beta.eq} for the cofibration
$\phi:\E \to I$. Moreover $\gamma^I_!E$ exists iff so do the
colimits in the right-hand side of \eqref{kan.I.eq}. We also have an
obvious relative counterpart of Lemma~\ref{kan.le}: for any framing
$\{\C(c) \subset \C,\gamma(c)_>:\C(c)^> \to \C'\}$, $c \in \C'$
of the functor $\gamma:\C \to \C'$, we have
\begin{equation}\label{fr.eq}
\colim_{C/^\gamma c}\zeta(\pi'(c)) \circ E|_{C /^\gamma
  c} \cong \colim_{\C(c)}\zeta(\pi'(c)) \circ E|_{\C(c)}
\end{equation}
for any $E \in \Fun(\C,\E)$, $c \in \C'$, where $\zeta(\pi'(c))$
acts on $E|_{\C(c)}$ through the factorization \eqref{aug.lft} of
the augmentation $\gamma(c)_>$, and for any $E' \in \Fun_I(\C',\E)$,
we have $E' \cong \gamma^I_!\gamma^*E'$ iff for any $c \in \C'$, the
augmented functor $\zeta(\pi'(c)) \circ \gamma(c)_>^*E$ is exact.

\begin{exa}\label{adj.I.exa}
If $\gamma:\C \to \C'$ admits a right-adjoint $\gamma^\dg:\C' \to
\C'$, a framing of $\gamma$ is obtained by taking $\C(c) = \ppt$
embedded onto $\gamma^\dg(c) \in \C$ and augmented by the adjunction
map $a:\gamma(\gamma^\dg(c)) \to c$. Then \eqref{kan.I.eq} and
\eqref{fr.eq} show that $\gamma_!^IE$ exists for any cofibration
$\phi:\E \to I$ and functor $E:\C \to \E$ over $I$, and it is given
by
\begin{equation}\label{adj.I.eq}
\gamma_!^IE(c) = \pi(a)_!E(\gamma^\dg(c)),
\end{equation}
where $\pi(a)_!$ is the transition functor for the cofibration $\C$
associated to the map $\pi(a):\pi(\gamma^\dg(c)) \cong
\pi'(\gamma(\gamma^\dg(c))) \to \pi'(c)$.
\end{exa}

A precofibration $\gamma:\C \to I$ is {\em bounded} resp.\ {\em
  discrete} if so are all its fibers $\C_i$, and {\em strongly
  discrete} if $\gamma^*\Id = \Id$. A discrete precofibration is of
course tautologically a cofibration. A functor $X:I \to \Sets$
defines a bounded strongly discrete cofibration $I[X] \to I$, where
$I[X]$ is the category of pairs $\langle i,x \rangle$, $i \in I$, $x
\in X(i)$, with maps $\langle i,x \rangle \to \langle i',x' \rangle$
given by maps $f:i \to i'$ such that $f(x) = x'$, and the projection
$I[X] \to I$ sending $\langle i,x \rangle$ to $i$. Every bounded
strongly discrete cofibration is of this form. More generally, for
any bounded confibration $\gamma:\C \to I$, we denote by
$\pi_0(\C/I) \to I$ the discrete cofibration corresponding to the
functor $\gamma_!\ppt:I \to \Sets$. Then we have a natural functor
$\gamma:\C \to \pi_0(\C/I)$, the functor is an equivalence if and
only if $C/I$ is discrete, and for any functor $\phi:\C \to \C'$
over $I$ to a bounded discrete cofibration $\C' \to I$ with the
corresponding equivalence $\gamma':\C' \to \pi_0(\C'/I)$, the
composition $\gamma' \circ \phi$ factors uniquely through $\gamma$.

\begin{exa}\label{tw.1.exa}
For any category $I$ with the twisted arrow category $\Tw(I)$,
sending $f:i' \to i$ to $i' \times i$ gives a discrete cofibration
$\Tw(I) \to I^o \times I$ corresponding to the $\Hom$-functor
$I(-,-):I^o \times I \to \Sets$; projecting further down to $I^o$,
we get the cofibration \eqref{ar.tw}.
\end{exa}

We will say that a cofibration $\gamma:\C \to I$ is {\em
  semidiscrete} if $\gamma^*\Iso = \Iso$, or equivalently, all its
fibers $\C_i$ are groupoids, or equivalently, all maps in $\C$ are
cocartesian over $I$. For any cofibration $\C \to I$, we have a
dense subcategory $\C_\natural \subset I$ defined by the class
$\natural$ of cocartesian maps, the induced projection $\C_\natural
\to I$ is a semidiscrete cofibration, and for any semidiscrete
cofibration $\C' \to I$, any functor $\C' \to \C$ cocartesian over
$I$ factors uniquely through $\C_\natural$. For any semidiscrete
cofibration $\gamma:\C' \to I$ and any functor $\pi:\C' \to \C$ such
that $\gamma \circ \pi:\C' \to I$ is a cofibration, $\pi$ itself is
a cofibration. A discrete cofibration is semidiscrete (although the
opposite is not true). Dually, a fibration $\C \to I$ is
semidiscrete if so is the opposite cofibration $\C^o \to I^o$, and
for any fibration $\C \to I$, we have the maximal semidiscrete
subfibration $\C_\ddag \subset \C$ defined by the class $\ddag$ of
all cartesian maps. Note that if a cofibration $\C \to I$ is
semidiscrete, then so is the transpose fibration $\C^\perp \to I^o$,
and we in fact have $\C^\perp \cong \C^o$.

\subsection{Cartesian functors.}\label{cart.subs}

If both $\pi:\C \to I$ and $\pi':\C' \to I$ are cofibrations, then a
functor $\gamma:\C \to \C'$ over $I$ is explicitly given by a
collection of functors $\gamma_i:\C_i \to \C_i'$ between their
fibers, and maps
\begin{equation}\label{fi.fu}
\gamma_f:f'_! \circ \gamma_i \to \gamma_{i'} \circ f_!,
\end{equation}
one for each morphism $f:i \to i'$ in $I$, subject to compatibility
conditions (where $f_!$, $f'_!$ are transition functors of the
cofibrations $\pi$, $\pi'$). The functor $\gamma$ is {\em
  cocartesian over $f$} iff \eqref{fi.fu} is an isomorphism, and it
is {\em cocartesian} if it is cocartesian over all maps in $I$. The
terminology for fibrations is dual (with ``cocartesian'' replaced by
``cartesian''). A functor is cocartesian if and only if it sends all
cocartesian maps to cocartesian maps. If $\C$ and $I$ are bounded,
then for any class of maps $v$ in $I$, we denote by
$\Fun_I^v(\C,\C') \subset \Fun_I(\C,\C')$ the full subcategory
spanned by functors cocartesian over all maps in $v$. We let
$\Sec^v(I,\C) = \Fun^v_I(I,\C)$ be the category of sections
cocartesian over maps in $v$, and dually, for any fibration $\C \to
I$, we let $\Sec_v(I,\C) \subset \Sec(I,\C)$ be the category of
section cartesian over maps in $v$. We note that speaking of
functors and sections cartesian resp.\ cocartesian over $v$ makes
sense even if $\C \to I$ is not a fibration resp.\ cofibration ---
it suffices to assume that it is such over $I_v \subset I$. For any
two cofibration $\C_0,\C_1 \to I$, a functor $\gamma:\C_0 \to \C_1$
cocartesian over $I$ induces a functor $\gamma^\perp:\C_0^\perp \to
\C_1^\perp$ between the transpose fibrations cartesian over $I$, and
vice versa (with $\gamma_i=\gamma^\perp_i$, $i \in I$, and the maps
\eqref{fi.fu} for $\gamma^\perp$ inverse to those for $\gamma$). For
any cofibration $\C \to I$, sending $E$ to $E^\perp$ provides an
equivalence
\begin{equation}\label{sec.perp}
\Sec^\forall(I,\C) \cong \Sec_\forall(I^o,\C^\perp).
\end{equation}
A convenient framing for a cocartesian functor $\gamma:\C \to \C'$
over $I$ is given by the left comma-fibers of the functors
$\gamma_i$: one takes $\C(c) = \C_{\pi'(c)} /^{\gamma_{\pi'(c)}} c$,
with the obvious augmentations. In particular, if $\C' = I$, $\pi' =
\id$, then $\gamma \cong \pi$ is automatically cocartesian, and has
a framing
\begin{equation}\label{fr.co}
\C(i) = \C_i
\end{equation}
given by its fibers. This immediately implies that if we have a
cartesian diagram of bounded categories
\begin{equation}\label{bc.sq}
\begin{CD}
\C' = \phi^*\C @>{\rho}>> \C\\
@V{\gamma'}VV @VV{\gamma}V\\
I' @>{\phi}>> I
\end{CD}
\end{equation}
and $\gamma$ is a cofibration, then for any functor $E \in
\Fun(\C,\E)$ to some category $\C$, $\gamma_!E$ exists iff so does
$\gamma'_!\rho^*E$, and the adjunction map
\begin{equation}\label{bc.eq}
\gamma'_!\rho^*E \to \phi^*\gamma_!E
\end{equation}
is an isomorphism. We call it the {\em base change isomorphism}
associated to a square \eqref{bc.sq}.

\begin{exa}\label{di.exa}
Assume given a bounded discrete cofibration $\gamma:I' \to I$, and a
functor $X':I' \to \Sets$. Then by \eqref{bc.eq}, $X = \gamma_!X'$
is given by
\begin{equation}\label{bc.di}
X(i) \cong \coprod_{i' \in I_i} X'(i'), \qquad i \in I,
\end{equation}
and this implies that $I'X' \cong IX$ (with the discrete cofibration
$IX \to I$ obtained by composing $\gamma$ with $I'X' \to I'$).
\end{exa}

\begin{exa}\label{av.exa}
As another application of \eqref{bc.eq}, assume given a small
category $I$ and a category $\E$, and note that we have an obvious
equivalence $\Fun^\forall(I,\E) \cong \Fun^\forall(I^o,\E)$. Then
the same equivalence can be realized by Kan extensions. Namely,
consider the twisted arrow category $\Tw(I^o)$, with the projections
$s:\Tw(I^o) \to I$, $t:\Tw(I^o) \to I^o$ of \eqref{ar.tw.co},
\eqref{ar.tw}, and for any $E \in \Fun(I,\E)$, let
\begin{equation}\label{av.eq}
\tw_I(E) = t_!s^*E
\end{equation}
whenever the Kan extension exist. Then $t$ is a cofibration with
fibers $\Tw(I^o)_i \cong i \setminus I$, $i \in I^o$, and for any $E
\in \Fun(I,\E)$, we have $s^*E|_{\Tw(I^o)_i} \cong p_i^*E$, where
$p_i$ is the projection \eqref{comma.pr}. Since $i \setminus I$ has
an initial object, this implies that as soon as $E$ is locally
constant, we have
$$
\tw_I(E)(i) \cong \colim_{i \setminus I}p_i^*E \cong E(i),
$$
so that $\tw_I(E)$ exists. Moreover, for any morphism $f:i \to i'$,
the morphism $\tw_I(E)(f)$ is inverse to $E(f)$, so that $\tw_I(E)$ is
locally constant. Thus $\tw_I$ gives the desired equivalence.
\end{exa}

\begin{lemma}\label{loc.le}
Assume given two bounded cofibrations $\C,\C' \to I$ and a
cocartesian functor $\gamma:\C \to \C'$ over $I$. Moreover, assume
that for any $i \in I$, the induced functor $\gamma_i:\C_i \to
\C'_i$ between the fibers \thetag{i} is full, or \thetag{ii} has a
right-adjoint $\gamma^\dg_i$, or \thetag{iii} is opposite to a
cofinal full embedding, or \thetag{iv} is a localization in the
sense of Definition~\ref{loc.def}. Then the same holds for $\gamma$.
\end{lemma}

\proof{} \thetag{i} is obvious. For \thetag{ii}, $\gamma^\dg$ has
components $\gamma^\dg_i$, with the maps \eqref{fi.fu} adjoint to
those for $\gamma$. For \thetag{iii}, the comma-fibers of $\gamma_i$
are left-admissible in the comma-fibers of $\gamma$, so the claim
follows from \eqref{forall.I}. Finally, for \thetag{iv}, it suffices
to check that the left Kan extension $\gamma_!E$ exists for any $E
\in \Fun^{\gamma^*\Iso}(\C,\E)$, the adjunction map $E \to
\gamma^*\gamma_!E$ is an isomorphism, and so id the adjunction map
$\gamma_!\gamma^*E' \to E'$ for any $E' \in \Fun(C,\E)$. By
\eqref{bc.eq}, all of this can be checked after restricting to the
fibers over all $i \in I$.
\endproof

\begin{lemma}\label{fib.le}
Let $\C \to I$ be a cofibration, and $\C' \subset \C$ a full
subcategory.
\begin{enumerate}
\item Assume that for any morphism $f:i \to i'$ in $I$, the transition
  functor $f_!:\C_i \to \C_{i'}$ sends $\C'_i \subset \C_i$ into
  $\C'_{i'} \subset \C_{i'}$. Then $\C' \to I$ is a cofibration, and
  the embedding $\C' \subset \C$ is cocartesian.
\item Assume that for any $i \in I$, the embedding $\C'_i \subset
  \C_i$ has a left-adjoint $\lambda_i:\C_i \to \C'_i$. Moreover, say
  that a map $g$ in $\C_i$, $i \in I$ is {\em $\lambda$-trivial} if
  $\lambda_i(g)$ is invertible, and assume that for any map $f$ in
  $I$, $f_!$ sends $\lambda$-trivial maps to $\lambda$-trivial
  maps. Then the embedding $\C'
  \subset \C$ admits a left-adjoint $\lambda:\C \to I$ over $I$.
\item Assume that the embedding $\C' \to \C$ admits a left-adjoint
  $\lambda:\C \to \C'$ over $I$. Then $\C' \to I$ is a cofibration,
  and $\lambda$ is cocartesian over $I$.
\end{enumerate}
\end{lemma}

\proof{} The first claim is obvious. For the second and the thrid,
the transition functors for $\C'$ are $f'_! = \lambda_{i'} \circ
f_!$ for any $f:i \to i'$ in $I$, and functor $\lambda$ in
\thetag{iii} is given by $\lambda_i$ with the maps \eqref{fi.fu}
given by adjunction. This shows that $\C' \to I$ is a
precofibration. Then either the fact that $\lambda$ is adjoint to
the embedding, or the final condition in \thetag{iii} is sufficient
to conclude the maps \eqref{fib.compo} and \eqref{fi.fu} are
isomorphisms.
\endproof

\begin{exa}\label{ar.exa}
Assume given a category equipped with a factorization system
$\langle L,R \rangle$ in the sense of Definition~\ref{facto.def},
and consider the subcategory $\Ar^R(I) \subset \Ar(I)$ and the
cofibration \eqref{ar.tw.co}. Then by Example~\ref{comma.exa}, the
fibers $\Ar^R(I)_i \subset \Ar(I)_i$ are left-admissible
subcategories, with the adjoint functors $\lambda_i:\Ar(I)_i \to
\Ar^R(I)_i$ sending an arrow $f:i' \to i$ to the component $r:i'' \to
i$ of its decomposition \eqref{l.r.facto}, and a map $g$ is
$\lambda$-trivial iff $s(g) \in L$. Since the transition functors
for the cofibration \eqref{ar.tw.co} commute with $s$, we are in the
situation of Lemma~\ref{fib.le}~\thetag{ii}, so that the induced
projection $t:\Ar^R(I) \to I$ is a cofibration. Dually, $s:\Ar^L(I)
\to I$ induced by \eqref{ar.tw} is a fibration.
\end{exa}

By abuse of terminology, for any cofibration $\pi:\C \to I$, we will
say that a functor $\gamma:\E \to \C$ from some category $\E$ is
{\em cocartesian over $I$} if $\gamma(f)$ is cocartesian for any map
$f$ in $\E$, and {\em vertical over $I$} if $\pi \circ \gamma$ is
constant. If $\E$ is bounded, then bounded vertical functors form a
cofibration $\Fun(\E,\C/I)$ over $I$ with fibers $\Fun(\E,\C/I)_i
\cong \Fun(\E,\C_i)$, $i \in I$. We will also say that a functor
$\gamma:\C \to \E$ is {\em cocartesian over $I$} if it inverts all
maps in $\C$ cocartesian over $I$ (or equivalently, if $\gamma
\times \id:\C \to \E \times I$ is cocartesian over $I$, or
equivalently, if the corresponding precofibration $\C' \to I^>$ of
Example~\ref{gt.exa} is a cofibration). For fibrations, the
terminology is again dual.

As a complement to Example~\ref{ar.exa}, note that for any
cofibration $\pi:\C \to I$, any map $f:c \to
c'$ in $\C$ has a factorization
\begin{equation}\label{c.v.facto}
\begin{CD}
c @>{c}>> c'' @>{v}>> c'
\end{CD}
\end{equation}
with cocartesian $c$ and $v \in \pi^*\Iso$, and \eqref{c.v.facto} is
unique up to a unique isomorphism, so that $\langle
\natural,\pi^*\Iso \rangle$ is a factorization system on $\C$. We
have $\Ar^\natural(\C) \cong \C \times^s_I \Ar(I)$, where $s$ is the
fibration \eqref{ar.tw}.

\begin{exa}\label{ar.1.exa}
In the situation of Example~\ref{ar.exa}, maps in $\Ar^R(I)$
cocartesian with respect to $t:\Ar^R(I) \to I$ correspond to
diagrams \eqref{ar.dia} with $g' \in L$. Then $s:\Ar^R(I) \to I$
restricts to a projection
\begin{equation}\label{s.R}
s:\Ar^R(I)_\natural \to I_L,
\end{equation}
and the dense embedding $\beta:I_L \to I$ factors as
\begin{equation}\label{be.t}
\begin{CD}
I_L @>{\eta}>> \Ar^R(I)_\natural @>{t}>> I,
\end{CD}
\end{equation}
where $\eta$ is induced by \eqref{ar.eta} and left-adjoint to the
projection \eqref{s.R}. In particular, we have $\beta_! \cong t_!
\circ \eta_! \cong t_! \circ s^*$, whenever these Kan extensions
exist.
\end{exa}

\subsection{Functor categories.}\label{fun.subs}

For any cofibration $\C \to I$, a functor $E:\C \to \E$ cocartesian
over $I$ defines a cocartesian functor $(E \times \id):\C \to \E
\times I$ over $I$, thus a transpose functor $(E \times
\id)^\perp:\C^\perp \to (\E \times I)^\perp \cong \E \times I^o$,
and then $(E \times \id)^\perp = E^\perp \times \id$ for a unique
functor $E^\perp:\C^\perp \to \E$ cartesian over $I^o$. If $I$ is
small and $\C$ is bounded, then $\C^\perp$ is bounded, and the
correspondence $E \mapsto E^\perp$ provides an equivalence
\begin{equation}\label{na.na}
\Fun^\natural(\C,\E) \cong \Fun^\ddag(\C^\perp,\E),
\end{equation}
where $\Fun^\natural(\C,\E) \subset \Fun(\C,\E)$
resp.\ $\Fun^\ddag(\C^\perp,\E) \subset \Fun(\C^\perp,\E)$ are the
full subcategories spanned by cocartesian resp.\ cartesian
functors. Alternatively, as in Example~\ref{av.exa}, \eqref{na.na}
can be realized by Kan extensions. Namely, consider the twisted
arrow category $\Tw(I)$ of Example~\ref{tw.exa}, and let $\Tw(\C/I)
= \Tw(I) \times_I \C$, $\Tw^\perp(\C/I) = \C^\perp \times_{I^o}
\Tw(I)$. For any $f:i' \to i$ considered as an object in
$\Tw(I)$, we have the transition functor $f_!:\C_{i'} \to \C_i$, and
taken together, these functors provide a functor $q:\Tw^\perp(\C/I)
\to \Tw(\C/I)$ over $\Tw(I)$. We then have the diagram
\begin{equation}\label{rtl.dia}
\begin{CD}
  \C @<{r}<< \Tw(\C/I) @<{q}<< \Tw^\perp(\C/I) @>{l}>> \C^\perp,
\end{CD}
\end{equation}
where $l$ and $r$ are the projection functors, and for any
cocartesian $E:\C \to \E$, we have $E^\perp \cong l_!q^*r^*E$.

For any cofibration $\C \to I$ over a category $I$ with an initial
object $o \in I$, the embedding $\C_o \to \C$ uniquely extends to a
functor
\begin{equation}\label{s.C}
\sigma:\C_o \times I \to \C
\end{equation}
cocartesian over $I$ (explicitly, $\sigma$ sends $c \times i$ to
$f(i)_!c$, where $f:o \to I$ is the unique map). This observation
has the following useful generalization. Let $\phi:I' \to I$ be a
bounded functor between bounded categories. Then for any cofibration
$\C \to I'$, we can define a cofibration $\phi_{**}\C \to I$ with
fibers
\begin{equation}\label{st.st}
\phi_{**}\C_i = \Sec(i \setminus^\phi I',p_i^*\C), \qquad i \in I,
\end{equation}
where $p_i$ is the projection \eqref{comma.tr}, and transition
functors $f_! = (f^*)^*$, where $f^*$ is the functor
\eqref{comma.tr}. For any class $v$ of maps in $I'$, we let
\begin{equation}\label{st.st.v}
\phi_{**}^v\C \subset \phi_{**}\C
\end{equation}
be the subcofibration spanned by sections cocartesian over maps in
$p_i^*v$, $i \in I$, and we let $\phi_*\C = \phi^{\forall}_{**}\C$
be the subcofibration spanned by cocartesian sections. Then since
for any $i \in I$, the right comma-fiber $i \setminus I$ has the
initial object, \eqref{s.C} induces a canonical equivalence $\id_*\C
\cong \C$ for any cofibration $\C \to I$, and then the pullback
functors $p_i^*$, $i \in I$ provide a functor
\begin{equation}\label{adj.2.1}
\C \to \phi_*\phi^*\C
\end{equation}
cocartesian over $I$. On the other hand, for any cofibration $\C \to
I'$, the evaluation functors \eqref{ev.eq} provide a functor
\begin{equation}\label{adj.2.2}
\ev:\phi^*\phi_{**}\C \to \C
\end{equation}
over $I'$ whose restriction to $\phi^*\phi_*\C \subset
\phi^*\phi_{**}\C$ is cocartesian over $I'$. By virtue of the
universal property of the functor \eqref{ev.eq}, for any bounded
cofibration $\C' \to I$ and any cofibration $\C \to I'$, the
functors \eqref{adj.2.1} and \eqref{adj.2.2} provide an equivalence
\begin{equation}\label{adj.2.eq}
\Fun_{I'}(\phi^*\C',\C) \cong \Fun_I^\natural(\C',\phi_{**}\C),
\end{equation}
where $\Fun_I^\natural(-,-) \subset \Fun_I(-,-)$ is the full
subcategory spanned by cocartesian functors. Moreover, for any class
of maps $v$ in $I'$, a functor $E:\phi^*\C' \to \C$ is cocartesian
over maps in $v$ if and only if the corresponding functor $E':\C'
\to \phi_{**}\C$ factors through $\phi^v_{**}\C \subset
\phi_{**}\C$, and in particular, \eqref{adj.2.eq} induces an
equivalence
\begin{equation}\label{adj.2.eq.1}
\Fun_{I'}^\forall(\phi^*\C',\C) \cong \Fun_I^\forall(\C',\phi_*\C),
\end{equation}
a sort of an adjunction property for $\phi^*$ and $\phi_*$. By
adjunction, for a composable pair $\phi:I' \to I$, $\phi':I'' \to I'$
of bounded functors, we have a natural equivalence $\phi_* \phi'_*\C
\cong (\phi \circ \phi')_*\C$ for any cofibration $\C \to I''$, and
if $\phi:I' \to I$ has a left-adjoint $\phi_\dg:I \to I'$, then we
have a natural equivalence $\phi_*\C \cong \phi_\dg^*\C$ for any
cofibration $\C \to I'$. In particular, if we have a
right-admissible full subcategory $I' \subset I$ in some bounded
$I$, with the embedding functor $\gamma:I' \to I$, and a cofibration
$\C \to I$, then \eqref{adj.2.eq.1} provides an equivalence
\begin{equation}\label{sec.adj}
\Sec^\forall(I,\C) \cong \Sec^\forall(I',\gamma^*\C).
\end{equation}
This implies that just as for usual Kan extensions, one can compute
the fibers of $\phi_*\C$ by choosing a framing $I'(i)$ for the
opposite functor $\phi^o$ in the sense of Lemma~\ref{kan.le}, and
replacing the fiber categories $i \setminus^\phi I'$ in
\eqref{st.st} with $I'(i)$. If $\phi$ is a fibration, one can use
the framing \eqref{fr.co}. Moreover, \eqref{sec.adj} is actually
induced by the pullback functor $\gamma^*$, with the inverse
equivalence provided by the pullback $\gamma_\dg^*$ with respect to
the right-adjoint functor $\gamma_\dg:I \to I'$, and then more
generally, for any class $v$ of maps in $I'$, the pullbacks
$\gamma^*$ and $\gamma^*_\dg$ provide an equivalence
\begin{equation}\label{sec.adj.v}
\Sec^{\gamma_\dg^*v}(I,\C) \cong \Sec^v(I',\gamma^*\C).
\end{equation}
In good situations, this allows to use framings to compute
$\phi^v_{**}$ for more general classes $v$.

\begin{exa}\label{id.exa}
Assume given a bounded category $I$ equipped with a factorization
system $\langle L,R \rangle$ in the sense of
Definition~\ref{facto.def}, and consider the projections
$s,t:\Ar^L(I) \to I$ of Example~\ref{ar.exa}. Then for any
cofibration $\C \to I$, the pullback functor $t^*$ induces a functor
\begin{equation}\label{id.st}
\Id^R_{**}\C \to s^\ddag_{**}t^*\C
\end{equation}
over $I$, where $\ddag = t^*(R)$ is the class of maps in $\Ar^L(I)$
cartesian over $I$, and \eqref{sec.adj.v} for the admissible
subcategories of Example~\ref{comma.exa} immediately shows that
\eqref{id.st} is an equivalence.
\end{exa}

\begin{exa}\label{ffun.exa}
For any bounded cofibration $\pi:\C \to I$ over a bounded category
$I$, and another cofibration $\C'/I$, let
\begin{equation}\label{ffun.ii}
\Fun^\forall(\C/I,\C'/I) = \pi^\natural_{**}\pi^*\C',
\end{equation}
where $\natural$ is the class of maps cocartesian over $I$. Then
$\Fun^\forall(\C/I,\C'/I) \to I$ is a cofibration with fibers
$$
\Fun^\forall(\C/I,\C'/I)_i \cong \Fun^{p_i^*\natural}_{i \setminus
  I}(p_i^*\C,p_i^*\C'), \qquad i \in I,
$$
where $p_i$ is the projection \eqref{comma.pr}, and
\eqref{adj.2.eq.1} provides an equivalence
\begin{equation}\label{ffun.eq}
\Fun_I^\forall(\C,\C') \cong
\Sec^\forall(I,\Fun^\forall(\C/I,\C'/I)).
\end{equation}
Informally, the category of categories cofibered over $I$ and
cocartesian functors between them is cartesian-closed, with
\eqref{ffun.ii} as the mapping category.
\end{exa}

Slightly more generally, for any bounded precofibration $\pi:\C \to
I$, and for any category $\E$, we can define a precofibration
$\Fun(\C/I,\E) \to I^o$ with fibers $\Fun(\C/I,\E)_i =
\Fun(\C_i,\E)$, transition functors $f^o_! = (f_!)^*$ induced by the
transition functors $f_!$ of the cofibration $\pi$, and the maps
\eqref{fib.compo} induced by the corresponding maps for $\pi$. If
$\pi$ is a cofibration, then $\Fun(\C/I,\E)$ is also a cofibration,
and we in fact have
\begin{equation}\label{ffun.st}
\Fun(\C/I,\E) \cong \pi^{\perp\ddag}_{**}\tau^*\E,
\end{equation}
where $\pi^\perp:\C^\perp \to I^o$ is transpose to $\pi$,
$\tau:\C^\perp \to \ppt$ is the tautological projection, and $\ddag$
is the class of maps cartesian over $I^o$, as in
Example~\ref{id.exa}. In this case, \eqref{adj.2.eq} identifies
functors from $\C$ to $\E$ and sections of the transpose fibration
$\Fun(\C/I,\E)^\perp \to I$, and this correspondence identifies
cocartesian functors with cartesian sections. By \eqref{na.na},
these in turn correspond to cocartesian sections of the cofibration
$\Fun(\C/I,\E) \to I^o$.

\subsection{Kernels and reflections.}\label{ker.subs}

In practical applications, it is useful to combine the pullback and
pushforward operations on cofibrations, as in
Example~\ref{id.exa}. Here is one specific construction of this type
that we will need.

\begin{defn}\label{ker.def}
For any bounded category $I$, an {\em $I$-kernel} is a category $\K$
equipped with a bounded cofibration $\K \to I^o \times I$. A {\em
  morphism} between $I$-kernels $\K$, $\K$ is a functor $\gamma:\K
\to \K'$ cocartesian over $I^o \times I$.
\end{defn}

For any bounded category $I$ with an $I$-kernel $\K$, the structural
cofibration $\K \to I^o \times I$ is the product $s \times t$ of
functors $s:\K \to I^o$, $\K:I' \to I$. Both are cofibrations, and
each of them is cocartesian with respect to the other one. We can
then consider the fibration $T = t^\perp:\K^\perp \to I^o$ transpose
to $t$, and being cocartesian, $s$ induces a cartesian functor
$S:\K^\perp \to I^o$ that is again a cofibration. Then for any
cofibration $\C \to I^o$, we let
\begin{equation}\label{ker.eq}
\K \otimes_I \C = T_{**}^\natural S^*\C,
\end{equation}
where $\natural$ is the class of maps in $\K^\perp$ cocartesian with
respect to $S$. This is a cofibration over $I^o$ with fibers
\begin{equation}\label{ker.exp}
(\K \otimes_I \C)_i \cong \Fun^\natural_{I^o}(\K_i,\C), \qquad i \in I,
\end{equation}
where $\K_i$ stands for the fiber of the cofibration $t$ cofibered
over $I^o$ by the cofibration $s$. For any morphism $\gamma:\K \to
\K'$ of $I$-kernels, the evaluation functor \eqref{adj.2.2} for
$\phi = \gamma^\perp$ then induces a functor
\begin{equation}\label{ker.mor}
\K' \otimes_I \C \to \K \otimes_I \C,
\end{equation}
cartesian over $I^o$. In terms of \eqref{ker.exp}, its component
over some $i \in I$ is the pullback $\gamma_i^*$ with respect to the
functor $\gamma_i:\K_i \to \K'_i$.

\begin{exa}\label{tw.yo.exa}
For any bounded category $I$, the twisted arrow category $\Tw(I)$
with its discrete cofibration $s \times t:\Tw(I) \to I^o \times I$
of Example~\ref{tw.1.exa} is an $I$-kernel. The category
$\Tw(I)^\perp$ is opposite to the arrow category $\Ar(I)^o$, with
the projections $S = s^o$, $T = t^o$ opposite to $s$ of
\eqref{ar.tw} and $t$ of \eqref{ar.tw.co}. Then \eqref{ar.eta}
provides their common section $\eta:I \to \Ar(I)$, and for any
cofibration $\C \to I^o$, we have a natural equivalence
$$
\Tw(I) \otimes_I \C \cong T_*S^*\C \cong T_*\eta_*\C \cong \C,
$$
where since $s:\Tw(I)_i \to I^o$ is discrete for any $i \in I$, we
can replace $T^\natural_{**}$ in \eqref{ker.eq} with $T_*$.
\end{exa}

\begin{exa}\label{sfun.exa}
For any two cofibrations $\C,\C' \to I$ with $\C$ bounded, let
$\pi:\C^\perp \to I^o$ be the transpose fibration, and let
\begin{equation}\label{sFun}
\sFun_I(\C,\C') = (\pi \times \id)^\ddag_{**}(\pi \times
\id)^*t^*\C',
\end{equation}
where $t:I^o \times I \to I$ is the projection, and $\ddag$ is the
class of maps in $\C' \times I$ cartesian over $I^o \times I$. Then
\eqref{sFun} is an $I$-kernel, and for any class of maps $v$ in $I$,
we have a natural identification
\begin{equation}\label{sfun.eq}
\Fun_{I^o \times I}^{t^*v}(\Tw(I),\sFun_I(\C,\C')) \cong
\Fun_I^v(\C,\C').
\end{equation}
In particular, morphisms of $I$-kernels $\Tw(I) \to \sFun_I(\C,\C')$
correspond to cocartesian functors $\C \to \C'$.
\end{exa}

In favourable circumstances, one can use Example~\ref{sfun.exa} to
establish a version of \eqref{ffun.eq} for functors that are only
cocartesian over a class of maps. To do this, we first observe that
the transpose fibration construction works in families. Namely, for
any cofibration $J \to I$, denote by $J^\lhd = J^{\perp o} \to I$
the cofibration with fibers $J_i^o$ and transition functors opposite
to those of $J \to I$, and note that a cocartesian functor
$\gamma:J_0 \to J_1$ between two cofibrations $J_0,J_1/I$ gives rise
to a cocartesian functor $\gamma^\lhd:J_0^\lhd \to J_1^\lhd$. Then
for any pair of cofibrations $J \to I$, $\gamma:\C \to J$, we can
define the {\em reflection} $(\C/J)^\lhd$ by
\begin{equation}\label{lhd.eq}
(C/J)^\lhd = (\C^\perp)_\perp,
\end{equation}
where the transpose fibration $\C^\perp$ is taken over $J$, and then
the transpose cofibration is taken over $I^o$. We have $(J/J)^\lhd =
J^\lhd$, and in general, $(C/J)^\lhd$ is equipped with a cofibration
$(C/J)^\lhd \to I$ and a functor $\gamma^\lhd:(C/J)^\lhd \to J^\lhd$
cocartesian over $I$. Over each $i \in I$, the fiber
$\gamma^\lhd_i:(\C/J)^\lhd_i \to J^\lhd_i \cong J_i^o$ is the
fibration transpose to $\gamma_i$. If we let $v$, $c$ be classes of
maps in $J$ and $J^\lhd$ vertical resp.\ cocartesian over $I$, then
$\gamma^\lhd$ is a cofibration over $J^\lhd_c$ and a fibration over
$J^\lhd_v$, and we have an equivalence
\begin{equation}\label{ref.perp}
\Sec^v(J,\C) \cong \Sec_v(J^\lhd,(\C/J)^\lhd),
\end{equation}
a relative version of \eqref{sec.perp}. For any commutative square
\begin{equation}\label{c.v.sq}
\begin{CD}
j_{00} @>{c^0}>> j_{01}\\
@V{v_0}VV @VV{v_1}V\\
j_{10} @>{c^1}>> j_{11}
\end{CD}
\end{equation}
in $J^\lhd$ with $c^0,c^1 \in c$, $v_0,v_1 \in v$, we have an
isomorphism $c^0_! \circ v_0^* \cong v_1^* \circ c^1_!$.

\begin{lemma}\label{refl.cof.le}
Assume given a cofibration $J \to I$ and a cofibration $\gamma:\C
\to J$ such that for any map $v:j \to j'$ in $J^\lhd$ vertical over
$I$, the transition functor $v^*:(\C/I)^\lhd_{j'} \to (\C/I^\lhd)_j$
admits a left-adjoint $v_!$, and for any square \eqref{c.v.sq}, the
map $v_{1!} \circ c^0_! \to c^1_!  \circ v_{0!}$ adjoint to the
isomorphism $c^0_! \circ v_0^* \cong v_1^* \circ c^1_!$ is itself
an isomorphism. Then $\gamma^\lhd$ is a cofibration.
\end{lemma}

\proof{} The first condition insures that $\gamma^\lhd$ is a
precofibration, with transition functor $f_! \cong v_! \circ c_!$
for any map $f = v \circ c$ in $J$, where $c$ is cocartesian and $v$
is vertical over $I$. The second condition then insures that the
maps \eqref{fib.compo} are isomorphisms, so that the precofibration
$\gamma^\lhd$ is a cofibration.
\endproof

\begin{lemma}\label{refl.restr.le}
Assume that $\pi:J^\lhd \to I$ has a fully faithful right-adjoint
$\eta:I \to J^\lhd$, and let $\C' = \eta^*(\C/J)^\lhd$. Then the
restriction functor
\begin{equation}\label{refl.restr.eq}
\eta^*:\Sec_v(J^\lhd,(\C/J)^\lhd) \to \Sec(I,\C')
\end{equation}
is an equivalence.
\end{lemma}

\proof{} Consider the category $\Ar^v(J)$, with the projections
$s,t:\Ar^v(J) \to J$ and the evaluation functor $\ev:[1] \times
\Ar^v(J) \to J$. Then $\ev^*\C \to [1]$ is a cofibration with fibers
$(\ev^*\C)_0 \cong s^*\C$, $(\ev^*\C)_1 \cong t^*\C$ and a
transition functor $q:s^*\C \to t^*\C$ cocartesian over $\Ar^v(J)$.
Since $t$ is a cofibration by Example~\ref{ar.exa}, $\Ar^v(J)$ is
cofibered over $I$, and we can consider the reflection $\C^\flat =
(\ev^*\C/[1] \times \Ar^v(J))^\lhd$. We have $([1] \times
\Ar^v(J))^\lhd \cong [1]^o \times \Ar^v(J^\lhd)$, with the
projections $s,t:\Ar^v(J^\lhd) \to J^\lhd$, and $\C^\flat$ is
fibered over $[1]^o = [1]$ with fibers $\C^\flat_0 \cong
s^*(\C/J)^\lhd$, $\C^\flat_1 \cong t^*(\C/J)^\lhd$ and the
transition functor $q^\lhd$. Since $\C^\flat \to [1]$ is a
fibration, we also have the functor $p:\C^\flat \to \C^\flat_0$
right-adjoint to the full embedding $\C^\flat_0 \subset \C^\flat$
given by $\id$ on $\C^\flat_0$ and $q^\lhd$ on $\C^\flat_1$.
We now observe that $t:\Ar^v(J^\lhd) \to J^\lhd$ has a fully
faithful right-adjoint $\eps:J^\lhd \to \Ar^v(J^\lhd)$ sending $j
\in J^\lhd$ to the adjunction map $j \to \eta(\pi(j))$, we have $s
\circ \eps \cong \id$ and $t \circ \eps \cong \eta \circ \pi$, so
that $\eps^*\C^\flat_0 \cong (\C/J)^\lhd$, $\eps^*\C^\flat_1 \cong
\pi^*\C'$, and $p$ restricts to a functor $\eps^*(p):\eps^*\C^\flat
\to (\C/J)^\lhd$. It then induces a functor
$$
\eps^*(p)^*:\Sec(J^\lhd,(\C/J)^\lhd) \to \Sec([1] \times
J^\lhd,\eps^*\C^\flat)
$$
sending a section $\sigma$ to a triple $\langle
\sigma_0,\sigma_1,\alpha \rangle$ of sections $\sigma_0 \in
\Sec(J^\lhd,(C/J)^\lhd)$, $\sigma_1 \in \Sec(J^\lhd,\pi^*\C')$ and a
map $\alpha:\sigma_0 \to \eps^*(q^\lhd) \circ \sigma_1$. We have
$\sigma_0 \cong \sigma$ and $\sigma_1 \cong
\pi^*\eta^*\sigma$. Moreover, $\eps^*(q^\lhd) \circ \sigma_1$ is
cartesian along $v$, and $\alpha$ is an isomorphism if and only if
so is $\sigma_0 = \sigma$. Therefore $\eps^*(q^\lhd) \circ \pi^*$ is
right-adjoint to $\eta^*$, and provides an inverse equivalence to
\eqref{refl.restr.eq}.
\endproof

Now assume given a category $I$ equipped with a factorization system
$\langle L,R \rangle$, consider the cofibration $t:\Ar^R(I) \to I$
of Example~\ref{ar.exa}, and let $\Tw^R(I) = \Ar^R(I)^\lhd$. The
cocartesian functor $\Ar(I) \to \Ar^R(I)$ left-adjoint to the
embedding $\Ar^(R)(I) \subset \Ar(I)$ induces a functor
\begin{equation}\label{pi.tw}
\pi^R:\Tw(I) \cong \Ar(I)^\lhd \to \Tw^R(I),
\end{equation}
and for any two cofibrations $\C,\C' \to I$ as in
Example~\ref{sfun.exa}, we can consider the cofibration
\begin{equation}\label{ffun.L}
\fFun^L(\C/I,\C'/I) = \pi^R_*(s^o \times t)^*\sFun_I(\C,\C') \to
\Tw^R(I).
\end{equation}
Then \eqref{adj.2.eq.1} and \eqref{sfun.eq} provide an equivalence
\begin{equation}\label{ffun.sec}
\Fun^L_I(\C,\C') \cong \Sec^{t^*L}(\Tw^R(I),\fFun^L(\C/I,\C'/I)).
\end{equation}
Explicitly, we have a factorization system $\langle t^*L \cap c,t^*R
\rangle$ on $\Ar^R(I)$, and \eqref{pi.tw} is a fibration over
$t^*R$, so that for any object in $\Ar^R(I)$ represented by an arrow
$r:i' \to i$, the embedding
\begin{equation}\label{tw.fr}
\Tw(i \setminus_L I) \cong r \setminus_{t^*L \cap c} \Tw(I) \subset
r \setminus \Tw(I) 
\end{equation}
is right-admissible by Example~\ref{comma.exa}. This gives a framing for
the functor opposite to \eqref{pi.tw}, and if we use this framing to
compute $\pi^L_*$ and apply \eqref{sfun.eq} for the comma-fiber $i
\setminus_L I$, we obtain an identification
\begin{equation}\label{ffun.r}
\fFun^L(\C/I,\C'/I)_r \cong \Fun^\forall_{i \setminus_L
  I}((f^*)^*p_{i'}^*\C,p_i^*\C'),
\end{equation}
where $f^*:i \setminus_L I \to i' \setminus_L I$ stands for the
transition functor of the fibration $s:\Ar^L(I) \to I$. In
particular, if $L = \Iso$, then $\fFun^L(\C/I,\C'/I)$ is the product
$\sFun_I(\C,\C') \times_{I^o \times I} \Tw(I))$, while at the other
extreme, if $L = \forall$, then \eqref{ffun.L} coincides with
\eqref{ffun.ii}.

We can then consider the reflection
$(\fFun^L(\C/,\C'/I)/\Tw^R(I))^\lhd$ of the cofibration
\eqref{ffun.L}. By definition, it is cofibered over $I$ and comes
equipped with a functor
\begin{equation}\label{ffun.p}
\gamma:(\fFun^L(\C/I,\C'/I)/\Tw^R(I))^\lhd \to \Ar^R(I) \cong
(\Tw^R(I))^\lhd
\end{equation}
cocartesian over $I$. Moreover, the projection $t:\Ar^R(I) \to I$
has a fully faithful right-adjoint $\eta:I \to \Ar^R(I)$ induced by
\eqref{ar.eta} that is cocartesian over $I_L$. Thus if we let
\begin{equation}\label{fun.ii}
\Fun^L(\C/I,\C'/I) = \eta^*(\fFun^L(\C/I,\C'/I)/\Tw^R(I))^\lhd \to I,
\end{equation}
then this is a cofibration over maps in $L$, and \eqref{ref.perp}
and Lemma~\ref{refl.restr.le} provide an equivalence
\begin{equation}\label{sec.ffun}
\Sec^{t^*L}(\Tw^R(I),\fFun^L(\C/I,\C'/I)) \cong
\Sec^L(I,\Fun^L(\C/I,\C',I)).
\end{equation}
Combining \eqref{sec.ffun} with \eqref{ffun.sec} then gives an
equivalence
\begin{equation}\label{fun.L.eq}
\Fun^L_I(\C,\C') \cong \Sec^L(I,\Fun^L(\C/I,\C'/I)).
\end{equation}
If $L=\forall$, this is the equivalence \eqref{ffun.eq} of
Example~\ref{ffun.exa}.

The ``favourable circumstances'' that we have mentioned are those
when the functor \eqref{ffun.p} satisfies the assumptions of
Lemma~\ref{refl.cof.le}. Then it is a cofibration over the whole
$I$, and so is its restriction \eqref{fun.ii}.

\section{Combinatorics.}\label{smcat.sec}

\subsection{Sets.}\label{set.sss}

A {\em pre-order} on a set is a binary relation $\leq$ that is
transitive and reflexive but not necessarily
antisymmetric. Equivalently, a pre-ordered set is a small category
that has at most one morphism between any two objects. We note that
for any pre-ordered set $J$, the categories $J^o$, $J^<$, $J^>$ are
also pre-ordered sets. The two opposite examples of a pre-order are
the maximal one (a morphism between any two objects) or a partial
order in the usual sense. For every $n \geq 0$, we denote by $[n]$
the partially ordered set $\{0,\dots,n\}$, with the standard order
(so that we have $[n]^< \cong [n]^> \cong [\np]$). In particular,
$[0]$ is the point category $\ppt$, and $[1]$ is the single arrow
category of Subsection~\ref{cat.subs}. We denote by $\V = \{0,1\}^<$
the partially ordered set with three elements $0,1,o$, and order
relations $0,1 \geq o$. Any set $S$ with the maximal pre-order is
denoted $e(S)$. More generally, for any preordered set $J$, define a
{\em $J$-augmented set} as a set $S$ equipped with a map $\pi:S \to
J$; then for any such $\langle S,\pi \rangle$, we let $e(S/J)$ be
the set $S$ with the pre-order induced by $\pi$ -- that is, $s \leq
s'$ iff $\pi(s) \leq \pi(s')$ (in particular, $e(S/\ppt) \cong
e(S)$, and $e(J/J) \cong J$). We say that a $J$-augmented set $S$ is
{\em proper} if $\pi:S \to J$ is surjective and we note that for a
proper $J$-augmented set $S$, $e(S/J) \to J$ is an equivalence of
categories.

We denote by $\Gamma$ the category of finite sets, and we let
$\Gamma_+$ be the category of finite sets and partially defined maps
between them -- that is, a map from $S_0$ to $S_1$ in $\Gamma_+$ is
a diagram
\begin{equation}\label{dom}
\begin{CD}
S_0 @<{i}<< \wt{S} @>{f}>> S_1
\end{CD}
\end{equation}
in $\Gamma$ with injective $i$. Equivalently, $\Gamma_+$ is the
category of finite pointed sets, with the equivalence sending $S$ to
its union $S_+ = S \cup \{o\}$ with an added distinguished
element. For clarity, we will always denote a finite set $S$
considered as an object in $\Gamma_+$ by $S_+$. A map in $\Gamma_+$
is an {\em anchor map} resp.\ a {\em structural map} if $f$
resp.\ $i$ in the diagram \eqref{dom} is an isomorphism. The
category $\Gamma_+$ is pointed, with the initial and terminal object
$o=\emptyset_+$ consisting of a single distinguished element, and
coproducts in $\Gamma$ are also coproducts in $\Gamma_+$,
traditionally written as $(S \copr S')_+ = S_+ \vee S'_+$. For any
$S,S' \in \Gamma$, the embeddings $i:S \to S \copr S'$, $i':S' \to S
\copr S'$ define anchor maps
\begin{equation}\label{S.spl}
a:S_+ \vee S'_+ \to S_+, \qquad a':S_+ \vee S'_+ \to S'_+
\end{equation}
given by the diagrams \eqref{dom} with $f=\id$ and $i = i_0$
resp.\ $i=i_1$. Cartesian product of finite sets is functorial with
respect to the diagrams \eqref{dom}, thus defines a product functor
\begin{equation}\label{sm.eq}
m:\Gamma_+ \times \Gamma_+ \to \Gamma_+.
\end{equation}
In terms of pointed sets, $m(S_+ \times S'_+)$ is the usual smash
product $S_+ \wedge S'_+ = (S_+ \times S'_+)/((S_+ \times o) \cup (o
\times S'_+))$.

\subsection{Ordinals.}\label{ord.subs}

As usual, we denote by $\Delta$ the category formed by non-empty
finite ordinals $[n]$, $n \geq 0$, and order-preserving maps between
them. Note that for any $[n]$, we have a unique isomorphism $[n]^o
\cong [n]$, and sending $[n]$ to $[n]^o$ gives an involution
$\iota:\Delta \to \Delta$. Adding the empty ordinal $\emptyset$ to
$\Delta$ gives the category $\Delta^<$, with $o = \emptyset$ (for
consistency, we will also denote it $[-1]$). We denote by
$V:\Delta^< \to \Sets$ the forgetful functor.  For any map $g:[n]
\to [m]$ in $\Delta^<$ with non-empty $[m]$, we can form the
cartesian square of partially ordered set
\begin{equation}\label{m.mb}
\begin{CD}
[n]_g @>{e_g}>> [n]\\
@VVV @VV{g}V\\
V([m]) @>>> [m],
\end{CD}
\end{equation}
where $V([m])$ is equipped with the discrete order, and then $[n]_g$
is the disjoint union of ordinals $[n_v]$ indexed by $v \in V([m]) =
\{0,\dots,m\}$ (non-empty if $g$ is surjective, possibly empty
otherwise). This provides a canonical equivalence
\begin{equation}\label{de.m}
\Delta^</[m] \cong \Delta^{<(m+1)}
\end{equation}
between the fiber category $\Delta^</[m]$ and the product of $m+1$
copies of $\Delta^<$.

For any $n \geq l \geq 0$, we denote by $s,t:[l] \to [n]$ the
(unique) embeddings onto the initial resp.\ terminal segment of the
target ordinal, and we note that we have a cocartesian square
\begin{equation}\label{seg.sq}
\begin{CD}
[0] @>{t}>> [l]\\
@V{s}VV @VV{s}V\\
[\nl] @>{t}>> [n].
\end{CD}
\end{equation}
The maps $s$ and $t$ for different $n \geq l \geq 0$ define dense
subcategories $\Delta_s,\Delta_t \subset \Delta$ (as abstract
categories, both are equivalent to the partially ordered set $\N$
of non-negative integers). We say that a map $f:[n] \to [m]$ is {\em
  special} resp.\ {\em antispecial} if $f(0)=0$ resp.\ $f(n)=m$,
and we denote by $\Delta_+,\Delta_- \subset \Delta$ the dense
subcategories defined by the classes of special resp.\ antispecial
maps. We note that we have $\Delta_s \subset \Delta_+$ and $\Delta_t
\subset \Delta_-$. A map is {\em bispecial} if it is both special
and antispecial, and $\Delta_{\pm} = \Delta_+ \cap \Delta_- \subset
\Delta$ is the corresponding dense subcategory. A map $f$ is special
iff $\iota(f)$ is antispecial, so that the involution $\iota$
establishes a equivalences $\iota:\Delta_+ \cong \Delta_-$,
$\iota:\Delta_{\pm} \cong \Delta_{\pm}$.

For another interpretation of special maps, consider the embedding
$\Delta \to \rCat$ sending an ordinal $[n]$ to the opposite ordinal
$[n]^o$ considered as a small category, and let
\begin{equation}\label{nu.dot}
\nu_\idot:\Delta_\idot \to \Delta
\end{equation}
be the corresponding cofibration. Explicitly, $\Delta_\idot$ is the
category of pairs $\langle [n],l \rangle$, $[n] \in \Delta$, $l \in
[n]$, with morphisms $\langle [n],l \rangle \to \langle [n'],l' \rangle$
given by maps $f:[n] \to [n']$ such that $f(l) \leq l'$. The
projection $\nu_\idot$ sends $\langle [n],l \rangle$ to $[n]$, it
has a right-adjoint $\nu_\dg:\Delta \to \Delta_\idot$ sending $[n]$
to $\langle [n],0 \rangle$, and this in turn has a right-adjoint
$\nu_\perp:\Delta_\idot \to \Delta$ sending $\langle [n],l \rangle$
to $[\nl] = \{l,\dots,n\} \subset [n]$. The dense subcategory in
$\Delta_\idot$ spanned by maps cocartesian over $\Delta$ is
equivalent to the right comma-fiber $[0] \setminus \Delta$, and a
map $f$ in $\Delta$ is special if and only if $\nu_\dg(f)$ is
cocartesian over $\Delta$. For antispecial maps, there is a similar
description using the cofibration $\iota^*\Delta_\idot$.

The embedding functor $\rho:\Delta_+ \to \Delta$ admits a
left-adjoint $\lambda:\Delta \to \Delta_+$ sending $[n]$ to
$[n]^<$. The composition $\kappa = \rho \circ \lambda$ is then an
endofunctor of the category $\Delta$, and we have the adjunction map
\begin{equation}\label{ka.adj}
a:\id \to \kappa
\end{equation}
equal to the embedding $t:[n] \to [\np] = \kappa([n])$ on any $[n]
\in \Delta$. The functor $\lambda$, hence also $\kappa$, extends to
$\Delta^< \supset \Delta$ by setting $\kappa(o)=\lambda(o) =
\emptyset^< = [0]$.

The embedding functor $\rho_\iota = \iota \circ \rho \circ
\iota:\Delta_- \to \Delta$ also admits a left-adjoint $\lambda_\iota
= \iota \circ \lambda \circ \iota:\Delta \to \Delta_-$, with the
composition $\kappa_\iota = \lambda_\iota \circ \rho_\iota \cong
\iota \circ \kappa \circ \iota$ and the adjunction map $a_\iota:\id
\to \kappa_\iota$, and the embedding
$\rho_\flat:\Delta_{\pm} \to \Delta$ admits a left-adjoint
$\lambda_\flat:\Delta \to \Delta_\pm$ with the composition
$\kappa_\flat = \lambda_\flat \circ \rho_\flat$ and the adjunction
map $\id \to \kappa_\flat$. We have $\kappa_\flat \cong \kappa \circ
\kappa_\iota \cong \kappa_\iota \circ \kappa$, and we have a
functorial cocartesian square
\begin{equation}\label{ka.sq}
\begin{CD}
\id @>{a}>> \kappa\\
@V{a_\iota}VV @VVV\\
\kappa_\iota @>>> \kappa_\flat.
\end{CD}
\end{equation}
Both $\kappa_\iota$ and $\kappa_\flat$ extend to $\Delta^<$, with
$\kappa_\iota(o)=\emptyset^>=[0]$ and
$\kappa_\flat(o)=\emptyset^{<>}=[1]$. The {\em concatenation} $[m]
\circ [n] \cong [m+n]$ of two ordinals $[m],[n] \in \Delta$ is their
disjoint union $[m] \copr [n]$ ordered left-to-right. Concatenation
is functorial with respect to either variable, and we have
$\kappa([n]) \cong [0] \circ [n]$, $\kappa_\iota([n]) \cong [n]
\circ [0]$ and $\kappa_\flat([n]) \cong [0] \circ [n] \circ [0]$.

A map $f:[n] \to [m]$ in $\Delta$ considered as a functor between
small categories has a right-adjoint $f^\dg:[m] \to [n]$ if and only
if it is special, and in this case $f^\dg$ is antispecial; for
left-adjoints $f_\dg$, the situation is dual. Thus sending $f$ to
its right resp. left-adjoint provides equivalences
\begin{equation}\label{p-m.eq}
\Delta_+ \cong \Delta_-^o, \qquad \Delta_- \cong \Delta_+^o.
\end{equation}
Moreover, for any map $f$ in $\Delta^<$, not only $\lambda(f)$ but
also the right-adjoint $\lambda(f)^\dg$ is special. Thus
$\lambda(f)^\dg$ is actually bispecial, and composing $\lambda$ with
the equivalence \eqref{p-m.eq} gives a functor
\begin{equation}\label{p-m-g.eq}
\Delta^< \to \Delta_{\pm}^o
\end{equation}
that also happens to be an equivalence. Restricting to $\Delta
\subset \Delta^<$, we obtain a canonical functor
\begin{equation}\label{theta.eq}
\theta:\Delta \to \Delta^o.
\end{equation}
We say that a map $f:[n] \to [m]$ in $\Delta$ is a {\em left-anchor
  map} resp.\ {\em right-anchor map} if $f=s$ resp.\ $f=t$, and we
say that $f$ is an {\em anchor map} iff $f:[n] \to [m]$ identifies
$[n]$ with some segment $\{l,l+1,\dots,l+n\} \subset [m]$ of the
ordinal $[m]$, or equivalently, if $f = s \circ t$ for some
right-anchor $t:[n] \to [n+l]$ and left-anchor $s:[n+l] \to [m]$.

By abuse of terminology, we say that a map $f$ in $\Delta^o$ is
special, anti-special, bispecial, anchor, left-anchor or
right-anchor if so is the opposite map $f^o$ in $\Delta$.

The classes of bispecial and anchor maps form a factorization system
$\langle \pm,a \rangle$ in $\Delta$ in the sense of
Definition~\ref{facto.def}, and we also have factorization systems
$\langle +,t \rangle$, $\langle -,s \rangle$. In particular, we have
the fibrations and cofibrations of Example~\ref{ar.exa}. Since
$\Delta_s \cong \Delta_t \cong \N$, we have $\Ar^s(\Delta) \cong
\Ar^t(\Delta) \cong \Delta_\idot$, and the cofibrations
$\Ar^s(\Delta),\Ar^t(\Delta) \to \Delta$ are both identified with
\eqref{nu.dot}. The fibrations
$\Ar(\Delta),\Ar^+(\Delta),\Ar^-(\Delta),\Ar^{\pm} \to \Delta$ have
a useful universl property with respect to cocartesian square
\eqref{seg.sq}. For any such square, the target of any map $f:[n]
\to [n']$ fits into a similar square with $l' = f(l)$; this is
functorial in $f$ and defines a cartesian diagram
\begin{equation}\label{pm.sq}
\begin{CD}
\Ar(\Delta)_{[n]} @>>> \Ar^+(\Delta)_{[n-l]}\\
@VVV @VVV\\
\Ar^-(\Delta)_{[l]} @>>> \Ar^\pm(\Delta)_{[0]},
\end{CD}
\end{equation}
where of course $\Ar^\pm(\Delta)_{[0]} \cong \ppt$, so that we
actually have a decomposition $\Ar(\Delta)_{[n]} \cong
\Ar^+(\Delta)_{[n-l]} \times \Ar^-(\Delta)_{[l]}$. If $f \in
\Ar(\Delta)_{[n]}$ is special resp.\ antispecial, then so are both
of the components of its decomposition, so that in particular,
\eqref{pm.sq} induces a decompositon
\begin{equation}\label{pm.eq}
\Ar^\pm(\Delta)_{[n]} \cong \Ar^\pm(\Delta)_{[l]} \times
\Ar^\pm(\Delta)_{[n-l]},
\end{equation}
where projections onto both factors are given by the transition
functors of the fibration $s:\Ar^\pm(\Delta) \to \Delta$.

Yet another useful factorization system $\langle p,i \rangle$ in
$\Delta$ is formed by the classes $p$ resp.\ $i$ of surjective
resp.\ injective maps. We have an embedding $\Ar^p(\Delta) \subset
\Ar^\pm(\Delta)$ cartesian with respect to the fibration
$s:\Ar^p(\Delta) \to \Ar^p(\Delta)$, and \eqref{pm.eq} induces a
corresponding decomposition for the fibration $\Ar^p(\Delta)$. Any
surjective map $p:[n] \to [m]$ considered as a functor between small
categories is both a fibration and a cofibration, and for any
injective map $i:[l] \to [m]$, there exists a cartesian square
\begin{equation}\label{p.i.sq}
\begin{CD}
[q] @>{i'}>> [n]\\
@V{p'}VV @VV{p}V\\
[l] @>{i}>> [m]
\end{CD}
\end{equation}
in $\Delta$ with some $[q]$, injective $i'$ and surjective $p'$
(this is obvious from \eqref{de.m}).

\subsection{Simplices.}

We denote by $\Delta^o\Sets$ the category of simplicial sets (that
is, functors from $\Delta^o$ to $\Sets$). For any simplicial set
$X$, its {\em category of simplices} $\Delta X$ is the corresponding
discrete fibration $\Delta X \to \Delta$. Explicitly, its objects
are pairs $\langle [n],x \rangle$, $[n] \in \Delta$, $x \in X([n])$,
and morphisms from $\langle [n],x \rangle$ to $\langle [n'],x'
\rangle$ are given by maps $f:[n] \to [n']$ such that
$X(f^o)(x')=x$. We let $\Delta^oX = (\Delta X)^o = \Delta^o[X]$ be
the opposite category. The {\em augmented category of simplices} is
given by $\Delta^< X = (\Delta X)^<$, with the opposite category
$\Delta^{o>}X = (\Delta^< X)^o = (\Delta X)^{o>}$.

Analogously, we denote by $\Sets_+$ the category of pointed sets,
and we let $\Delta^o\Sets_+$ be the category of simplicial pointed
sets. One distinguished simplicial pointed set is the {\em standard
  simplicial circle} $\Sigma$ obtained by taking the standard
simplicial interval $\Delta_1$, $\Delta_1([n]) = \Delta([n],[1])$,
and gluing together the two ends. The simplicial set $\Sigma$ is
finite, thus factors through a functor
\begin{equation}\label{Si}
\Sigma:\Delta^o \to \Gamma_+.
\end{equation}
A map $f$ in $\Delta$ is an anchor map iff $\Sigma(f^o)$ is an
anchor map in $\Gamma_+$.

A {\em contraction} of an augmented simplicial object
$c_>:\Delta^{o>} \to \C$ in some category $\C$ is a functor
$c_+:\Delta_+^o \to \C$ equipped with an isomorphism
$\lambda^{o*}c_+ \cong c^>$. Any simplicial set $X:\Delta^o \to
\Sets$ has a tautological augmentation $X_>$ sending $o \in \Delta^{o>}$
to the one-point set $\ppt$, and we say that $X$ is {\em
  contractible} if $X_>$ admits a contraction $X_+:\Delta_+^o \to
\Sets$. By the Grothendieck construction, such a contraction defines
a discrete fibration $\Delta_+X_+ \to \Delta_+$ with the opposite
discrete cofibration $\Delta_+^o X_+ = (\Delta_+ X_+)^o \to
\Delta_+^o$, and $\lambda$ extends to an embedding $\Delta^< X \to
\Delta_+X_+$.

\begin{defn}\label{aug.def}
For any contractible simplicial set $X \in \Delta^o\Sets$, an
augmentation $c^>:\Delta^{o>} X \to \C$ of a functor $c:\Delta^o X
\to \C$ to some category $\C$ is {\em contractible} if it further
extends to $\Delta^o_+X_+ \supset \Delta^{o>}X$ for some contraction
$X_+$ of the set $X$.
\end{defn}

\begin{lemma}\label{aug.le}
For any contractible simplicial set $X$ and category $\C$, a locally
constant functor $c:\Delta^o X \to \C$ is constant, and a
contractible augmentation $c_>$ of a functor $c:\Delta^o X \to \C$
is exact.
\end{lemma}

\proof{} The embedding $\lambda_X:\Delta X \to \Delta^< X \to
\Delta_+ X_+$ admits a left-adjoint $\rho_X:\Delta_+ X_+ \to \Delta
X$ given by the composition
$$
\begin{CD}
\Delta_+ X_+ @>{a^*}>> \rho^*\lambda^*\Delta_+ X_+ \cong
\rho^*\Delta X @>>> \Delta X,
\end{CD}
$$
where $a:\lambda \circ \rho \to \id$ is the adjunction
map. Therefore $\Delta^o X \subset \Delta^o_+X_+$ is a
left-admissible subcategory. On the other hand, by definition, the
set $X_+(o)$ consists of a single point, say $x$, and then $\langle
o,x \rangle \in \Delta_+X_+$ is the terminal object, thus defines a
left-admissible embedding $\ppt \to \Delta^o_+ X_+$. Then
\eqref{adm.eq} immediately shows that a contractible augmentation is
exact. Moreover, any locally constant $c_+:\Delta^o_+X_+ \to \C$ is
constant, and any locally constant $c:\Delta^oX \to \C$ inverts the
adjunction map $\id \to \rho_X^o \circ \lambda_X^o$, so that $c
\cong \lambda_X^{o*}\rho_X^{o*}c$, and $c$ is constant as well.
\endproof

As usual, the {\em nerve} $N(I)$ of a small category $I$ is the
simplicial set such that $N(I)([n])$, $[n] \in \Delta$ is the set of
functors $i_\idot:[n] \to I$. If $I$ is connected, its nerve $N(I)$
is connected, and if $I$ has an initial object, $N(I)$ is
contractible. Indeed, since $[0] \in \Delta_+$ is the initial
object, we have a natural projection $\rho^{o*}N(I) \to N(I)([0])$,
or in other words, a decomposition
\begin{equation}\label{ka.de.0}
\rho^{o*}N(I) = \coprod_{i \in I}N(I)_i, \qquad N(I)_i:\Delta_+^o
\to \Sets.
\end{equation}
Moreover, for any $i \in I$, we have $\lambda^{o*}N(I)_i \cong N(i
\setminus I)$. But if $i \in I$ is an initial object, the map $N(i
\setminus I) \to N(I)$ induced by the forgetful functor $i \setminus
I \to I$ is an isomorphism, so that $N(I)_i$ is a contraction of
$N(I)$.

The nerve is functorial in $I$ and provides a fully faithful
embedding $\rCat \to \Delta^o\Sets$ from the category $\rCat$ of small
categories to the category of simplicial sets. Its essential image
consists of simplicial sets $X:\Delta^o \to \Sets$ that send
cocartesian squares \eqref{seg.sq} to cartesian squares in $\Sets$
(this is called the {\em Segal condition}). The nerve embedding
commutes with limits but not with colimits. Colimits in $\rCat$ exist
but are hard to compute except in some special cases. Here is one
obvious example.

\begin{exa}\label{MV.exa}
If we have full subcategories $I_0,I_1 \subset I$ that are
left-closed in the sense of Example~\ref{cl.exa}, then $I_0 \cup
I_1$ is left-closed in $I$, $I_{01} = I_0 \cap I_1$ is left-closed
in $I_0$, $I_1$ and $I$, and we have $I_0 \cup I_1 \cong I_0
\copr_{I_{01}} I_1$. Giving a functor $I_0 \cup I_1 \to \E$ to some
category $\E$ is equivalent to giving functors $\gamma_0:I_0 \to
\E$, $\gamma_1:I_1 \to \E$, and an isomorphism between their
restrictions to $I_{01}$.
\end{exa}

Another example is path categories of quivers. Namely, a {\em
  quiver} $Q$ is a collection of two sets $V(Q)$, $E(Q)$ of its {\em
  vertices} and {\em edges}, and two maps $s,t:V(Q) \to E(Q)$
sending an edge to its source resp.\ target. Equivalently, a quiver
is a functor $Q:\Dd^o \to \Sets$, where $\Dd$ is the category with
two objects $0$, $1$ and two non-trivial maps $s,t:0 \to 1$, and we
take $V(Q) = Q(0)$, $E(Q) = Q(1)$. We have the embedding $\delta:\Dd
\to \Delta$ sending $1$ to $[1]$, $0$ to $[0]$ and $s$, $t$ to the
maps $s$, $t$ in $\Delta$, and the opposite functor $\delta^o$
factors as
\begin{equation}\label{P.a}
\begin{CD}
\Dd^o @>{\alpha}>> \Delta^o_a @>{\beta}>> \Delta^o,
\end{CD}
\end{equation}
where $\Delta_a \subset \Delta$ is the dense subcategory spanned by
anchor maps. Let $\Dd^o\Sets$ be the category of quivers, and let
$\Delta_a^o\Sets$ be the category of functors from $\Delta^o_a$
sets. Then the composition $P=\beta_!\alpha_*:\Dd^o\Sets \to
\Delta^o\Sets$ of the right Kan extension $\alpha_*$ and the left
Kan extension $\beta_!$ factors through $\rCat \subset \Delta^o\Sets$
and is left-adjoint to the restriction functor $\delta^{o*}:\rCat
\subset \Delta^o\Sets \to \Dd^o\Sets$. For any $Q \in \Dd^o\Sets$,
$P(Q)$ is called the {\em path category} of the quiver $Q$. By
adjunction, the functor $P:\Dd^o\Sets \to \rCat$ commutes with
colimits, and colimits of quivers can be computed pointwise.

\begin{remark}\label{qui.rem}
The presense of two different Kan extensions in the construction of
the path category $P(Q)$ is easy to understand if one observes that
for any $n \geq 0$, $[n]$ is the path category of a ``string
quiver'' $[n]_\delta$ with $V([n]_\delta) = V([n])$ and edges $e \in
E([n]_\delta) = \{1,\dots,n\}$ corresponding to intervals $\{l-1,l\}
\subset [n]$. Then maps between such quivers are exactly the anchor
maps, the squares \eqref{seg.sq} are induced by cocartesian squares
of quivers, so that the Segal condition makes sense for objects $X
\in \Delta^o_a\Sets$, and $\alpha_*$ in \eqref{P.a} identifies
$\Dd^o\Sets$ with the category of functors $X \in \Delta^o_a\Sets$
satifying this condition. Explicitly, it is given by
$$
\alpha_*Q([n]) \cong \Hom_{\Dd^o\Sets}([n]_\delta,Q), \qquad [n] \in
\Delta.
$$
The Kan extension $\beta_!$ can be then computed by \eqref{bc.di}
and the decomposition \eqref{be.t} for the anchor/bispecial
factorization system on $\Delta^o$, where one observes that since
all isomorphisms in $\Delta^o$ are identity maps, the cofibration
$\Ar^\pm(\Delta^o)_\natural \to \Delta^o$ is discrete. It is easy to
see from \eqref{pm.eq} that $\beta_!$ preserves the Segal
condition, and this creates the adjunction between $P$ and
$\delta^{o*}$.
\end{remark}

\subsection{Cylinders.}\label{cyl.subs}

Among other things, Lemma~\ref{aug.le} allows to prove some useful
results above nerves of small categories. To simplify notation, for
any small $I$, denote $\Delta I = \Delta N(I)$, and similarly for
$\Delta^o$, $\Delta^<$ and $\Delta^{o>}$. Note that $\Delta I^o
\cong \iota^*(\Delta I)$, where $\iota:\Delta \to \Delta$ sends
$[n]$ to $[n]^o$.

\begin{exa}\label{delta.n.exa}
For any $[n] \in \Delta$ considered as a small category, we have
$\Delta [n] \cong \Delta / [n]$, $\Delta^< [n] \cong \Delta^< / [n]$
and $\Delta^o [n] \cong [n] \setminus \Delta^o$, $\Delta^{o>}[n]
\cong [n] \setminus \Delta^{o>}$.
\end{exa}

As we saw, if $I$ has an initial object $i \in I$, then the nerve
$N(I)$ is contractible, so that by Lemma~\ref{aug.le}, $\Delta^o I$
is simply connected in the sense of Definition~\ref{loc.def}. The
first corollary of this is the following standard fact.

\begin{lemma}\label{shuf.le}
The diagonal embedding $\Delta^o \to \Delta^o \times \Delta^o$ is
cofinal in the sense of Definition~\ref{cof.def}.
\end{lemma}

\proof{} For each $[n] \times [m] \in \Delta^o$, the right comma-fiber
$([n] \times [m]) \setminus \Delta^o$ is $\Delta^o([n] \times [m])$
by Example~\ref{delta.n.exa}, and since $[n] \times [m]$ has an
initial object, it is simply connected by Lemma~\ref{aug.le}.
\endproof

\begin{remark}
A standard implication of Lemma~\ref{shuf.le} is that if a
cocomplete category $\E$ has finite products, then
$\colim_{\Delta^o}$ commutes with these products.
\end{remark}

\begin{lemma}\label{mu.le}
Assume given a full embedding $\gamma:I' \to I$ of small categories
that is either left or right-closed in the sense of
Example~\ref{cl.exa}. Then the induced functor
$\Delta^{o>}(\gamma):\Delta^{o>}I' \to \Delta^{o>}I$ admits a
left-adjoint functor
\begin{equation}\label{mu.adj.eq}
\mu:\Delta^{o>}I \to \Delta^{o>}I',
\end{equation}
and this functor is a cofibration.
\end{lemma}

\proof{} By Example~\ref{delta.n.exa} and \eqref{de.m}, we have
$\Delta^{o>}[1] \cong \Delta^{o>} \times \Delta^{o>}$, and then by
Example~\ref{facto.exa}, the classes $\pi_0^*\Iso$ and $\pi_1^*\Iso$
form a factorization system on $\Delta^{o>}[1]$, in either
order. Let us denote them by $0$ and $1$ and call them {\em
  $0$-special} and {\em $1$-special} maps (so that $l$-special maps
are those that are isomorphisms over $l \in [1]$). If $\gamma$ is
left-closed, then $I' = I_0$ for a functor $\chi:I \to [1]$, and
$\Delta^{o>}(\chi):\Delta^{o>}I \to \Delta^{o>}[1]$ is a discrete
cofibration. Thus if we say that a map $f$ in $\Delta^{o>}I$ is
$l$-special, $l=0,1$ whenever so is $\Delta^{o>}(\chi)(f)$, then
$0$-special and $1$-special maps also form a factorization system
$\langle 1,0 \rangle$ on $\Delta^{o>}I$. Consider the cofibration
$t:\Ar^0(\Delta^{o>}I) \to \Delta^{o>}$ of Example~\ref{ar.exa}, and
observe that $\Delta^{o>}(\gamma)^*\Ar^0(\Delta^{o>}I) \cong
\Delta^{o>}I$. Then $\mu = \Delta^{o>}(\gamma)^*t$ is our
cofibration \eqref{mu.adj.eq}, and the adjunction is obvious. For
a right-closed $\gamma$, interchange $0$ and $1$.
\endproof

Now take some $[n] \in \Delta$, and assume given a small category
$I$ and a functor $\pi:I \to [n]$, with the induced functor $\Delta
I \to \Delta [n] \cong \Delta/[n]$. Consider the full subcategories
$\Delta_p/[n] \subset \Delta_+/[n] \subset \Delta/[n]$ spanned by
surjective resp.\ special arrows $f:[m] \to [n]$, and let
\begin{equation}\label{del.p.pm.I}
\Delta_p(I/[n]) = \Delta I \times_{\Delta[n]} \Delta_p/[n], \qquad
\Delta_+(I/[n]) = \Delta I \times_{\Delta[n]} \Delta_+/[n].
\end{equation}
Say that $\pi:I \to [n]$ is an {\em iterated cylinder} if it is a
cofibration with non-empty fibers.

\begin{lemma}\label{cyl.le}
For any iterated cylinder $\pi:I \to [n]$, $\Delta_p(I/[n]) \subset
\Delta_+(I/[n])$ is cofinal. If $I$ has an initial object $i$, then
$\Delta_+(I/[n]) \subset \Delta I$ is also cofinal.
\end{lemma}

\proof{} The second claim follows from Lemma~\ref{aug.le} and
Lemma~\ref{cof.le}: since $I \cong i \setminus I$, we have a
left-admissible subcategory $\Delta_+N(I)_i \subset \Delta I$, and
since the fiber $I_0 \subset I$ is non-empty, we have $i \in I_0$,
so that $\Delta_+N(I)_i \subset \Delta_+(I/[n])$.

For the first claim, note that by definition, objects in $\Delta I$
are pairs $\langle [n],i_\idot \rangle$ of an object $[n] \in
\Delta$ and a functor $i_\idot:[n] \to I$, $l \mapsto i_l$. For any
$l \in [n]$, let $\Delta_lI \subset \Delta I$ be the full
subcategory consisting of pairs $\langle [m],i_\idot \rangle$ such
that the image of the map $\pi \circ i_\idot:[m] \to [n]$ contains
$s([l]) = \{0,\dots,l\} \subset [n]$. Then $\Delta_p(I/[n]) =
\Delta_0I$, $\Delta_+(I/[n]) = \Delta_nI$, and by induction and
Lemma~\ref{cof.le}, it suffices to prove that $\Delta_{l+1}I \subset
\Delta_lI$ is cofinal for any $l=0,\dots,n-1$. Choose such an $l$,
consider the embeddings $s:[l] \to [n]$, $t:[\nlm] \to [n]$, and
denote $I_{\leq l} = s^*I$, $I_{> l} = t^*I$. Then the full
embedding $I_{\leq l} \to I$ is left-closed, and the functor
opposite to the cofibration \eqref{mu.adj.eq} sends $\Delta_lI$ into
$\Delta I_{\leq l} \subset \Delta^<I_{\leq l}$. Therefore we have
the induced fibration $\Delta_lI \to \Delta I_{\leq l}$,
$\Delta_{l+1}I \subset \Delta_lI$ is a subfibration, and for any
object $\langle [m],i_\idot \rangle \in \Delta I_{\leq l}$, their
fibers are
\begin{equation}\label{mu.fib.eq}
(\Delta_{l+1}I)_{\langle [m],i_\idot \rangle} \cong \Delta_+ (i_m
  \setminus I_{>l}) \subset (\Delta_lI)_{\langle [m],i_\idot
    \rangle} \cong \Delta (i_m \setminus I_{>l}).
\end{equation}
But since $I \to [n]$ is a cofibration, the iterated cylinder $(i_m
\setminus I_{>l})/[\nlm]$ has an initial object. Therefore all the
embeddings \eqref{mu.fib.eq} are cofinal by the second claim, and we
are done by Lemma~\ref{loc.le}~\thetag{iii}.
\endproof

As an application of Lemma~\ref{cyl.le}, consider the arrow
categories $\Ar^p(\Delta)$, $\Ar^\pm(\Delta)$, and treat the
embedding $\delta:\Ar^p(\Delta)^o \to \Ar^\pm(\Delta)^o$ as a
functor over $\Delta^o$ with respect to the projections induced by
$t$ of \eqref{ar.tw.co}. Then for any object $c \in
\Ar^\pm(\Delta)^o$, we have an embedding
\begin{equation}\label{c.d.m}
\Ar^p(\Delta)^o /^\delta_{t^{o*}\Iso} c \subset
\Ar^p(\Delta)^o /^\delta c.
\end{equation}
Were $\delta$ cocartesian over $\Delta^o$, the embedding
\eqref{c.d.m} would have been left-admissible. While $\delta$ is
certainly not cocartesian --- its source and target are not even
cofibrations --- we still have the following result.

\begin{lemma}\label{cdm.le}
The embedding \eqref{c.d.m} is cofinal for any $c \in
\Ar^\pm(\Delta)^o$.
\end{lemma}

\proof{} By Example~\ref{ar.facto.exa}, the bispecial/anchor
factorization system $\langle \pm,a \rangle$ on $\Delta$ induces a
factorization system on $\Ar^\pm(\Delta)$, so that any map $c \to
c'$ in $\Ar^\pm(\Delta)$ uniquely factors as
$$
\begin{CD}
c @>{l}>> c'' @>{r}>> c'
\end{CD}
$$
with $l \in (s \times t)^*\pm$ and $r \in (s \times t)^*a$. One
immediately checks that if $c'$ lies in $\Ar^p(\Delta) \subset
\Ar^\pm(\Delta)$, then so does $c''$. Therefore by
Example~\ref{adm.sub.exa} and Example~\ref{comma.exa}, if we denote
by $\delta':\Ar^p(\Delta_\pm)^o \to \Ar(\Delta_\pm)$ the restriction
of the embedding $\delta$ to $\Ar^p(\Delta_\pm)^o \subset
\Ar^p(\Delta)^o$, then the embeddings
$$
\Ar^p(\Delta_\pm)^o /^{\delta'} c \subset \Ar^p(\Delta)^o /^\delta
c, \qquad \Ar^p(\Delta_\pm)^o /^{\delta'}_{t^{o*}\Iso} c \subset
\Ar^p(\Delta)^o /^\delta_{t^{*}\Iso} c
$$
are left-admissible. Then by Lemma~\ref{cof.le}, it suffices to
prove the claim with \eqref{c.d.m} replaced by the embedding
\begin{equation}\label{c.d.m.1}
\Ar^p(\Delta_\pm)^o /^{\delta'}_{t^{o*}\Iso} c \subset
\Ar^p(\Delta_\pm)^o /^{\delta'} c.
\end{equation}
To do this, take some object $c' \in \Ar^p(\Delta_\pm)^o /^{\delta'} c$,
and consider the right comma-fiber
$$
I(c,c') = c' \setminus (\Ar^p(\Delta_\pm)^o /^{\delta'}_{t^*\Iso} c)
$$
of the embedding \eqref{c.d.m.1}. Explicitly, $c$ is a bispecial map
$f:[n] \to [m]$, $c'$ is a pair of a surjective map $f':[n'] \to
[m']$ and a map $g:f \to f'$ in $\Ar^{\pm}(\Delta)$ with bispecial
components $g_n:[n] \to [n']$, $g_m:[m] \to [m']$, and $I(c,c')$ is
then opposite to the category of commutative diagrams
\begin{equation}\label{I.cc.dia}
\begin{CD}
[n] @>>> [n''] @>>> [n']\\
@V{f}VV @VV{p}V @VV{f'}V\\
[m] @= [m] @>{g_m}>> [m']
\end{CD}
\end{equation}
in $\Delta_{\pm}$ with surjective $p$ whose outer rectangle
represents the map $g$. We have to prove that $I(c,c')$ is simply
connected. Since bispecial maps admit the decomposition
\eqref{pm.eq} with respect to cartesian squares \eqref{seg.sq}, the
whole diagram \eqref{I.cc.dia} also functorially decomposes into
diagrams for $[l]$ and $[\nl]$, so that by induction, it suffices to
consider the case $[n]=1$. Then the left vertical arrow in
\eqref{I.cc.dia} carries no information and can be forgotten, and we
end up with an equivalence
$$
I(c,c') \cong \Delta_p^o(g_m^*[n']/[m]).
$$
But $[n']$ has an initial object, and since $f':[n'] \to [m']$ is
surjective, it is a cofibration and an iterated cylinder. Then
$g_m^*f':g_m^*[n'] \to [m]$ is also an iterated cylinder with an
initial object, and we are done by Lemma~\ref{cyl.le} and
Lemma~\ref{aug.le}.
\endproof

\subsection{Cycles.}\label{cycle.subs}

For any integer $n \geq 1$, we denote by $[n]_\Lambda$ the category
whose objects are residues $a \in \Z/n\Z$, with maps from $a$ to
$a'$ given by integers $l \geq 0$ such that $a' = a + l \mod n$
(equivalently, $[n]_\Lambda$ is the path category of the wheel
quiver $[n]_\lambda$ with $n$ vertices, and $l$ is the length of the
path). For any $a \in [n]_\Lambda$, we have a map $\tau_a:a \to a$
given by $n$, and for any functor $f:[n]_\Lambda \to [m]_\Lambda$,
there exists a unique integer $\deg f \geq 0$ such that $f(\tau_a) =
\tau_{f(a)}^{\deg f}$ for all $a \in [n]_\Lambda$. We say that $f$
is {\em non-degenerate} if $\deg f \geq 1$, and {\em horizontal} if
$\deg f = 1$. The {\em cyclic category $\Lambda$} of A. Connes is
the category with objects $[n]$, $n \geq 1$, and morphisms from
$[n]$ to $[m]$ given by horizontal functors $[n]_\Lambda \to
[m]_\Lambda$. For any $[n] \in \Lambda$, we denote by
$V([n]_\lambda)$ the set of objects or the category $[n]_\Lambda$
(or equivalently, the set of vertices of the corresponding wheel
quiver). Every horizontal functor $f:[n]_\Lambda \to [m]_\Lambda$
admits a left-adjoint $f_\dg:[m]_\Lambda \to [n]_\Lambda$, and
sending $[n]$ to $[n]$ and $f$ to $f^o_\dg$ provides an equivalence
of categories
\begin{equation}\label{la.du}
\theta:\Lambda \cong \Lambda^o.
\end{equation}
For any $[n]$ in $\Lambda$ and a horizontal functor $f:[n]_\Lambda
\to [1]_\Lambda$, we have a cartesian square of small categories
\begin{equation}\label{la.de}
\begin{CD}
[\nm] @>{\eps(f)}>> [n]_\Lambda\\
@VVV @VV{f}V\\
[0] @>{\eps}>> [1]_\Lambda,
\end{CD}
\end{equation}
where $[0],[\nm] \in \Delta$ are considered as small categories in
the usual way, and $\eps$ in the bottom line is the embedding onto
the unique object in $[1]_\Lambda$ (in terms of quivers,
\eqref{la.de} corresponds to removing one edge from a wheel quiver
to obtain a string quiver). Sending the rightmost colum in
\eqref{la.de} to the leftmost column defines a functor $\Lambda/[1]
\to \Delta$ that happens to be an equivalence; composing the inverse
equivalence with the forgetful functor $\Lambda/[1] \to \Lambda$
gives an embedding $j:\Delta \to \Lambda$. More generally, for any
$[n] \in \Lambda$, we can choose a map $f:[n] \to [1]$, and this
gives an identification
\begin{equation}\label{la.de.n}
\Lambda/[n] \cong \Delta/[\nm].
\end{equation}
Passing to the opposite categories and identifying $\Lambda \cong
\Lambda^o$ by \eqref{la.du}, we obtain a dual embedding
$j^o:\Delta^o \to \Lambda$ and an equivalence $\Lambda \setminus [n]
\cong \Delta^o \setminus [n]$ dual to \eqref{la.de.n} for any $[n]
\in \Lambda$. One can also see the equivalence $[1] \setminus
\Lambda \cong \Delta^o$ directly in terms of the equivalence
\eqref{p-m-g.eq}: if we let $\omega_1:[1] \to [1]_\Lambda$ be the
functor that sends $0,1 \in [1]$ to the unique object $0 \in
[1]_\Lambda$, and the map $0 \to 1$ to the map $\tau_0:0 \to 0$,
then for any injective bispecial map $b:[1] \to [n]$ in $\Delta$, we
have the cocartesian square of small categories
\begin{equation}\label{de.la}
\begin{CD}
[1] @>{\omega_1}>> [1]_\Lambda\\
@V{b}VV @VV{b'}V\\
[n] @>{\omega_n}>> [n]_\Lambda,
\end{CD}
\end{equation}
and the functor $b'$ is horizontal (to see that \eqref{de.la} is
cocartesian, note that it is induced by a cocartesian square that
glues the two end vertices of a string quiver to obtain a wheel). The
correspodence $b \mapsto b'$ is functorial and identifies $[1]
\setminus \Lambda$ with the full subcategory $\Delta^o \subset
\Delta_{\pm}$ spanned by injective bispecial maps. The
identifications given by \eqref{la.de} and \eqref{de.la} intertwine
the duality \eqref{la.du} and the embedding \eqref{theta.eq}.

\begin{lemma}\label{av.le}
For any cocomplete target category $\E$, the functor $\tw_\Lambda$
of \eqref{av.eq} takes values in $\Fun^\forall(\Lambda,\E) \cong
\Fun^\forall(\Lambda^o,\E) \subset \Fun(\Lambda^o,\E)$.
\end{lemma}

\proof{} For any $E \in \Fun(\Lambda,\E)$, $E'=\tw_\Lambda(E)$ exists
since $\E$ is cocomplete, and we have to check that for any map
$f:[m] \to [n]$ in $\Lambda$, $E'$ inverts $f^o$. It
suffices to choose a map $g:[1] \to [m]$ and prove that $E'$ inverts
$g^o$ and $(f \circ g)^o$, so we may assume right away that
$[m]=1$. Then by \eqref{bc.eq}, it suffices to check that the
transition functor $f^o_!:\Tw(\Lambda^o)_{[n]} \to
\Tw(\Lambda^o)_{[1]}$ of the cofibration $t:\Tw(\Lambda^o) \to
\Lambda^o$ is cofinal. But by \eqref{la.de.n} and \eqref{la.du},
this transition functor can be identified with the projection
$(\Delta/[\nm])^o \cong \Delta^o[\nm] \to \Delta^o$, and then its
right comma-fibers are of the form $\Delta^o([\nm] \times [l])$,
$[l] \in \Delta^o$. As in the proof of Lemma~\ref{shuf.le}, these
are simply connected by Lemma~\ref{aug.le}.
\endproof

The {\em cyclotomic category} $\LR \subset \Lambda$ has the same
objects $[n]$, $n \geq 1$, and all non-degenerate functors
$f:[n]_\Lambda \to [m]_\Lambda$ as morphisms. A morphism $f$ in
$\LR$ is {\em vertical} if it is a discrete bifibration. Any
morphism $f$ in $\LR$ factors as $f = v \circ h$ with horizontal $h$
and vertical $v$, and $\langle h,v \rangle$ is a factorization
system on $\LR$. For any integer $l \geq 1$, one defines the
category $\Lambda_l$ as the category of vertical arrows $v:[nl] \to
[n]$ in $\LR$ of degree $l$, with morphisms given by commutative
square
$$
\begin{CD}
[nl] @>{v}>> [n]\\
@V{h'}VV @VV{h}V\\
[n'l] @>{v'}>> [n']
\end{CD}
$$
in $\LR$ with horizontal $h$, $h'$. Such a square is automatically
cartesian. Sending an arrow to its source resp.\ target gives two
functors
\begin{equation}\label{i.pi}
i_l,\pi_l:\Lambda_l \to \Lambda.
\end{equation}
The functor $i_l$ is known as the {\em edgewise subdivision functor}
and goes back to \cite{edge}. The functor $\pi_l$ has no special
name; it is a bifibration whose fibers are connected groupoids
$\ppt_l$ with one object $o$, $\Aut(o)=\Z/l\Z$.

\section{Stabilization and additivization.}\label{stab.sec}

\subsection{Homotopy categories.}\label{ho.subs}

We denote by $\Top_+$ the category of pointed compactly generated
topological spaces, and we denote by $\Ho$ its homotopy category
obtained by localizing $\Top_+$ with respect to the class of weak
equivalences. For any bounded category $I$, we denote by $\Ho(I)$
the localization of the category $\Fun(I,\Top_+)$ with respect to
the class of pointwise weak equivalences. For any bounded functor
$\gamma:I \to I'$ between bounded categories, the pullback
$\gamma^*$ preserves weak equivalences, thus descends to a functor
$\gamma^*:\Ho(I') \to \Ho(I)$. This functor $\gamma^*$ has a left
and a right-adjoint $\gamma_!,\gamma_*:\Ho(I) \to \Ho(I')$ known as
the {\em homotopy Kan extensions}. If $I'$ is the point category
$\ppt$, then $\gamma_! = \hocolim_I$ is the homotopy colimit over
$I$, and $\gamma_* = \holim_I$ is the homotopy limit.

For any bounded cofibration $\C \to I$, the transition functors
$(f_!)^*$ of the transpose fibration $\Fun(\C/I,\Top_+)^\perp \to I$
are pullback functors, so that $\Fun(\C/I,\Top_+)^\perp$ descends to
a fibration $\ho(\C/I) \to I$ with fibers $\ho(\C/I)_i = \Ho(\C_i)$,
$i \in I$. This fibration is also a cofibration, with transition
functors given by homotopy Kan extensions $(f_!)_!$. Any cocartesian
functor $\gamma:\C' \to \C$ between two cofibrations over $I$
induces a cartesian pullback functor $\gamma^*:\ho(\C/I) \to
\ho(\C'/I)$ that has a left-adjoint $\gamma_!:\Ho(\C'/I) \to
\Ho(\C/I)$; for each $i \in I$, we have $(\gamma^*)_i \cong
\gamma_i^*$ and $(\gamma_!)_i \cong \gamma_{i!}$. For any bounded
cofibration $\C \to I$, we have a comparison functor
\begin{equation}\label{ho.Ho}
\ho:\Ho(\C) \to \Sec(I,\Ho(\C/I))
\end{equation}
that reduces to $\Ho(\C) \to \Fun(I,\Ho)$ if $\C = I$. This functor
is conservative and commutes with pullbacks (but it is certainly not
an equivalence unless $I$ is a point).

To compute homotopy Kan extensions, one actually needs to know very
little about them (apart from their sheer existence known from
e.g.\ \cite{DHKS} that can be used as a black box). Namely, just as
in the non-homotopical case, an adjunction between functors
$\gamma_0:I' \to I$, $\gamma_1:I \to I'$ induces an adjunction
between the pullback functors $\gamma_0^*$ and $\gamma_1^*$, and
this yields \eqref{adm.eq} (with colimits replaced by homotopy
colimits). Moreover, for a discrete cofibration $\pi:I' \to I$, the
left Kan extension $\pi_!:\Fun(I',\Top_+) \to \Fun(I,\Top_+)$
respects weak equivalences, thus descends to the left homotopy Kan
extension, and dually for discrete fibrations and $\pi_*$. Together
with \eqref{adm.eq}, this describes Kan extensions with respect to
embeddings $\ppt \to I$, and the description immediately yields
homotopical versions of \eqref{kan.eq}, Lemma~\ref{kan.le} and the
base change isomorphism \eqref{bc.eq} (where again, one replaces
$\colim$ with $\hocolim$ in all the statements, and uses $\gamma_!$
to denote the homotopy Kan extension).

For any bounded $I$ and a class of maps $v$ in $I$, we denote by
$\Ho^v(I) \subset \Ho(I)$ the full subcategory spanned by $X \in
\Ho(I)$ such that $\ho(X)$ lies in $\Fun^v(I,\Ho)$. We say that $X
\in \Ho(I)$ is {\em locally constant} if so is $\ho(X)$, so that
$\Ho^\forall(I) \subset \Ho(I)$ is the full subcategory spanned by
locally constant $X$. For any functor $\gamma:I' \to I$ and any
class $v$ of maps in $I$, \eqref{pb.v} induces a functor
\begin{equation}\label{ho.pb.v}
\gamma_v^*:\Ho^v(I) \to \Ho^{\gamma^*v}(I'),
\end{equation}
and as in Definition~\ref{loc.def}, we say that $\gamma$ is a {\em
  homotopy localization} if \eqref{ho.pb.v} is an equivalence for
any $v$. We note that the localizations of Example~\ref{adm.exa} are
also homotopy localizations. Somewhat weaker, we say that $\gamma$
is a {\em weak equivalence} if \eqref{ho.pb.v} is an equivalence for
the class $v=\forall$ of all maps. Such an equivalence is
automatically compatible with limits and colimits: for any $X \in
\Ho^\forall(I')$, the natural maps
\begin{equation}\label{ho.co.eq}
\hocolim_{I'}\gamma^*X \to \hocolim_IX
\end{equation}
and
\begin{equation}\label{ho.li.eq}
\holim_X \to \holim_{I'}\gamma^*X
\end{equation}
are both isomorphisms in $\Ho$. A homotopy localization is trivially
a weak equivalence, and so are its one-sided inverses; in
particular, it applies to left and right-admissible full
embeddings.

In general, a functor $\gamma:I' \to I$ between small categories
$I$, $I'$ is a weak equivalence if and only if so is the induced map
$|I'| \to |I|$ of the geometric realizations of their
nerves. However, this is a serious theorem that we will not need. We
restrict ourselves to the following trivial but very useful
observation whose effectiveness was demonstrated to us by
L. Hesselholt.

\begin{lemma}\label{Qui.le}
Assume given functors $\gamma_0:I' \to I$, $\gamma_1:I \to I'$
between bounded categories $I'$ and $I$ equipped with a map between
$\Id$ and $\gamma_0 \circ \gamma_1$ and a map between $\Id$ and
$\gamma_1 \circ \gamma_0$ (in either direction). Then
$\gamma_{0\forall}^*$ and $\gamma_{1\forall}^*$ are mutually inverse
equivalences.
\end{lemma}

\proof{} A one-sided inverse to a weak equivalence is trivially a
weak equivalence. For any $I$, consider the product $I \times [1]$,
with the embeddings $s,t:I \to I \times [1]$ onto $I \times 0$ and
$I \times 1$, and their common one-sided inverse $e:I \times [1] \to
I$ given by the projection onto the first factor. Then $e$ is a weak
equivalence by Example~\ref{adm.exa}, and therefoe so are $s$ and
$t$. Therefore for any two functors $\phi_0,\phi_1:I' \to I$, a map
$\phi_0 \to \phi_1$ induced an isomorphism $\phi_{0\forall}^* \cong
\phi_{1\forall}^*$. In the assumptions of the Lemma, this provides
isomorphisms $\Id \cong \gamma_{0\forall}^* \circ
\gamma_{1\forall}^*$ and $\Id \cong \gamma_{1\forall}^* \circ
\gamma_{0\forall}^*$.
\endproof

We will say that a bounded category $I$ is {\em homotopy
  contractible} if the tautological projection $\tau:I \to \ppt$ is
a weak equivalence, and we will say that a functor $\gamma:I' \to I$
between bounded categories is {\em homotopy cofinal} if its right
comma-fiber $i \setminus^\gamma I'$ is homotopy contractible for any
$i \in I$. For any contractible $I$, $\holim_I$ and $\hocolim_I$ are
isomorphic, and both are inverse equivalences to
$\tau^*_\forall$. As for \eqref{adm.eq}, by virtue of (the dual
homotopical version of) \eqref{kan.eq}, this implies that $\gamma_*
\circ \tau^* \cong (\tau \circ \gamma)^*$ for any cofinal functor
$\gamma:I' \to I$, and then by adjunction, this in turn implies that
\eqref{ho.co.eq} is an isomorphism for any $X \in \Ho(I)$.

A contractible $I$ is simply connected in the sense of
Definition~\ref{loc.le}, and a homotopy cofinal full embedding
$\gamma$ is cofinal in the sense of Definition~\ref{cof.le}. The
converse is not true; however, Lemma~\ref{cof.le},
Remark~\ref{cof.rem} and all the criteria of cofinality of
Section~\ref{cat.sec} still hold for homotopy cofinal embeddings,
and all the cofinal embeddings of Section~\ref{smcat.sec} are
homotopy cofinal. In particular, this is true for the embedding of
Lemma~\ref{shuf.le}, so that for any $X \in \Ho(\Delta^o \times
  \Delta^o)$, the natural map
\begin{equation}\label{suf}
\hocolim_{\Delta^o}\delta^*X \to \hocolim_{\Delta^o \times
  \Delta^o}X
\end{equation}
induced by the diagonal embedding $\delta:\Delta^o \to \Delta^o
\times \Delta^o$ is a homotopy equivalence. This implies that
$\hocolim_{\Delta^o}$ commutes with finite products.

\subsection{Basic stabilization.}\label{stab.subs}

We begin by recalling the standard Segal construction of \cite{seg}
(with slightly different notation). As in Subsection~\ref{set.sss},
let $\Gamma_+$ be the category of finite pointed sets. A {\em
  $\Gamma$-space} is an object $X \in \Ho(\Gamma_+)$.

\begin{defn}\label{seg.def}
A $\Gamma$-space $X$ is {\em special} if $X(o)$ is contractible, and
for any $S,S' \in \Gamma$, the map
\begin{equation}\label{a.X}
X(a) \times X(a'):X(S_+ \vee S'_+) \to X(S_+) \times X(S'_+)
\end{equation}
induced by the anchor maps \eqref{S.spl} is a weak homotopy
equivalence.
\end{defn}

It is well-known that any special $\Gamma$-space $X$ is a
commutative monoid object in $\Ho(\Gamma_+)$. Indeed, if we take
$S=S'$ and let $d:\Gamma_+ \to \Gamma_+$ be the functor sending
$S_+$ to $S_+ \wedge S_+$, then the anchor maps \eqref{S.spl} and
the codiagonal maps $m:S_+ \wedge S_+ \to S_+$ give rise to a
diagram
$$
\begin{CD}
X \times X @<{X(a) \times X(a')}<< d^*X @>{X(m)}>> X,
\end{CD}
$$
and the map on the left is an equivalence, so it can be inverted in
$\Ho(\Gamma_+)$. This defines a product map
\begin{equation}\label{G.prod}
X \times X \to X
\end{equation}
in $\Ho(\Gamma_+)$, and one checks that it is associative and
commutative. In particular, $\pi_0(X(\ppt_+))$ is a commutative
monoid, and $X$ is called {\em group-like} if this monoid is a
group.

\begin{defn}\label{stab.def}
For any bounded category $I$, an object $X \in \Ho(\Gamma_+ \times
I)$ is {\em stable} if for any $i \in I$, its restriction to
$\Gamma_+ \times \{i\}$ is a group-like special $\Gamma$-space.
\end{defn}

\begin{prop}\label{stab.prop}
For any bounded category $I$, denote by $\Ho^{st}(\Gamma_+,I)$ the
full subcategory in $\Ho(\Gamma_+ \times I)$ spanned by stable
objects. Then the embedding $\Ho^{st}(\Gamma_+, I) \subset
\Ho(\Gamma_+ \times I)$ admits a left-adjoint {\em stabilization
  functor}
$$
\Stab:\Ho(\Gamma_+ \times I) \to \Ho^{st}(\Gamma_+, I).
$$
\end{prop}

\proof{} Note that for any $S_+,S'_+ \in \Gamma_+$ and functor
$X:\Gamma_+ \to \Top_+$, we have a natural map
\begin{equation}\label{X.S}
\bigvee_{s \in S}X(i_{s+} \wedge \id):\bigvee_{s \in S} X(S'_+) =
S_+ \wedge X(S'_+) \to X(S_+ \wedge S'_+),
\end{equation}
where $i_s:\ppt \to S$ is the embedding onto $s \in S$, and the map
\eqref{X.S} is functorial in $S_+$, $S'_+$ and $X$. Consider the
product $\Delta^o \times \Gamma_+ \times I$, with the projection
$\pi:\Delta^o \times \Gamma_+ \times I \to \Gamma_+ \times I$ and
the functor $\beta:\Delta^o \times \Gamma_+ \times I \to \Gamma_+
\times I$ given by the composition
$$
\begin{CD}
\Delta^o \times \Gamma_+ \times I @>{\Sigma \times \id \times \id}>>
\Gamma_+ \times \Gamma_+ \times I @>{m \times \id}>> \Gamma_+ \times
I,
\end{CD}
$$
where $m$ is the product functor \eqref{sm.eq}, and $\Sigma:\Delta^o
\to \Gamma_+$ is the simplicial circle. For any $X \in \Ho(\Gamma_+
\times I)$, let $BX = \pi_!\beta^*X$. The maps \eqref{X.S} provide a
map $\Sigma \wedge \pi^*X \to \beta^*X$ that gives rise to a map
$$
\hocolim_{\Delta^o}(\Sigma \wedge \pi^*X) \cong
\hocolim_\Delta^o\Sigma \wedge X \cong \Sigma X \to BX,
$$
and by adjunction, we obtain a map $X \to \Omega B X$, where
$\Sigma$ and $\Omega$ stand for the suspension and loop space
functor applies pointwise. By induction, we obtain a map
$\Omega^nB^nX \to \Omega^{n+1}B^{n+1}X$ for any $n \geq 1$, and we
can set
\begin{equation}\label{stab.bo}
\Stab(X) = \colim_n\Omega^nB^nX.
\end{equation}
Then $\Stab$ obviously commutes with coproducts and Kan extensions
$\gamma_!$, $\gamma:I' \to I$, and it is the main result of
\cite{seg} that if $X$ is stable, $BX$ is also stable and $X \to
\Omega BX$ is a homotopy equivalence. Thus to show that $\Stab$ with
the natural map $\Id \to \Stab$ provides the required stabilization
functor, it suffice to show that $\Stab(X)$ is stable for any $X \in
\Ho(\Gamma_+ \times I)$.

To check this, we may assume that $I = \ppt$, so that $X$ is a
$\Gamma$-space. For any such $X$ and $S_+,S'_+ \in \Gamma_+$, denote
by $X(S_+,S'_+)$ the homotopy cofiber of the map \eqref{a.X}. Say
that $X$ is {\em $n$-connected}, $n \geq 0$, if $X(S_+)$ is
$n$-connected for any $S_+ \in \Gamma_+$, and {\em $n$-stable} if it
is group-like and $X(S_+,S'_+)$ is $n$-connected for any $S_+,S'_+
\in \Gamma_+$. Then if a $\Gamma$-space $X$ is $m$-stable, the
corresponding maps \eqref{a.X} induce isomorphisms on the homotopy
groups $\pi_i$ for $i < m$, and then for any $n \leq m$, the same
maps for $\Omega^nX$ are isomorphisms on $\pi_i$ for $i < m -
n$. Thus it suffices to prove that for any $X$ and $n \geq 0$,
$B^nX$ is $2n$-stable. By induction, it further suffices to prove
that $B X$ is $0$-connected and $1$-stable, and $(n+1)$-connected
and $(m+2)$-stable if $X$ is $n$-connected and $m$-stable. For
connectedness, recall that for any $Y:\Delta^o \to \Top_+$ such that
$Y([q])$ is $n$-connected for any $[q] \in \Delta^o$ and
contractible for $q=0$, $\hocolim_{\Delta^o}Y$ is $(n+1)$-connected,
and note that $Y:\Delta^o \to \Top_+$ sending $[q] \in \Delta^o$ to
$X(S_+ \wedge \Sigma([q]))$ satisfies this assumption for any $S_+
\in \Gamma_+$. For stability, note that
$$
BX(S_+,S'_+) = \hocolim_{\Delta^o}\delta^*Y,
$$
where $Y:\Delta^o \times \Delta^o \to \Top_+$ sends $[q]
\times[q']$ to $X(S_+ \wedge \Sigma([q]),S'_+ \wedge \Sigma([q']))$,
apply the isomorphism \eqref{suf}, and repeat the same argument
twice.
\endproof

\begin{remark}\label{hoco.stab.rem}
For any functor $\gamma:I' \to I$, the pullback $(\id \times
\gamma)^*$ sends stable objects to stable objects, so that by
adjunction, we have a natural map $\Stab \circ (\id \times \gamma)^*
\to (\id \times \gamma)^* \circ \Stab$, and the construction of
$\Stab$ given in Proposition~\ref{stab.prop} immediately shows that
this map is an isomorphism. If $I' = \Delta^o$ and $I = \ppt$, then
the adjoint map $(\id \times \gamma)_!  \circ \Stab \to \Stab \circ
(\id \times \gamma)_!$ is also an isomorphism. Indeed, in this case
$(\id \times \gamma)^*$ is fully faithful, so it suffices to check
that $(\id \times \gamma)_!$ sends stable objects to stable
objects. But since $\hocolim_{\Delta^o}$ commutes with finite
products, $(\id \times \gamma)_!$ sends special objects to special
objects, and it obviously also preserves the group-like condition.
\end{remark}

\begin{remark}\label{spectr.rem}
For any bounded category $I$, we have the forgetful functor
$e^*:\Ho^{st}(\Gamma_+, I) \to \Ho(I)$, where
\begin{equation}\label{ga.e}
e:I \times \Gamma_+
\end{equation}
is the embedding onto $\ppt_+ \times I$, and by virtue of
Proposition~\ref{stab.prop}, it has a left-adjoint {\em
  spectrification functor} $\Stab \circ e_!:\Ho(I) \to
\Ho^{st}(\Gamma_+, I)$.
\end{remark}

\begin{remark}\label{stab.0.rem}
Say that a {\em $\Gamma$-set} is a functor $X:\Gamma_+ \to \Sets$,
say that it is {\em special} if $X(o)=\ppt$ and the maps \eqref{a.X}
are isomorphisms, say that it is moreover {\em stable} if the
product \eqref{G.prod} turns it into a group, and for any category
$I$, say that $X:\Gamma_+ \times I \to \Sets$ is {\em stable} if so
is $X_i:\Gamma_+ \to \Sets$ for any $i \in I$. Then a stable
$\Gamma$-set is the same thing as an abelian group, so if $I$ is
bounded, we have a full embedding
\begin{equation}\label{stab.set}
\Fun(I,\Z) \cong \Fun^{st}(\Gamma_+ \times I,\Sets) \subset
\Fun(\Gamma_+ \times I,\Sets),
\end{equation}
where $\Fun^{st}$ stands for the full subcategory spanned by stable
objects. It is easy to see that the embedding \eqref{stab.set}
admits a left-adjoint stabilization functor $\Stab_0$. If we compose
it with the left Kan extension $e_!$, as in Remark~\ref{spectr.rem},
then $\Stab_0(e_! S) \cong \Z[S]$, the free abelian group generated
by the functor $S:I \to \Sets$. The functor $\pi_0:\Top_+ \to \Sets$
is left-adjoint to the embedding $\Sets \to \Top_+$ sending a set
$S$ to the discrete topological space $S_+$, the same holds in
families, and then by adjunction, we have $\pi_0(\Stab(X)) \cong
\Stab_0(\pi_0(X))$ for any $X \in \Ho(\Gamma^+ \times I)$ and
$\pi_0(\Stab(e_! X)) \cong \Z[\pi_0(X)]$ for any $X \in \Ho$.
\end{remark}

\subsection{Advanced stabilization.}\label{adv.stab.subs}

It turns out that purely by formal games with adjunction,
Proposition~\ref{stab.prop} yields several useful corollaries.  We
will need two of them. For the first one, say that a category $I$ is
{\em half-additive} if it is pointed and has finite coproducts. Then
for any bounded half-additive category $I$, we have a natural
functor
\begin{equation}\label{ga.m}
m:\Gamma_+ \times I \to I
\end{equation}
sending $S \times i$ to the coproduct of copies of $i$ numbered by
elements $s \in S$, and \eqref{ga.e} is a section of
\eqref{ga.m}. Say that $X \in \Ho(I)$ is {\em stable} if so is $m^*X
\in \Ho(\Gamma_+ \times I)$, and denote by $\Ho^{st}(I) \subset
\Ho(I)$ the full subcategory spanned by stable functors. Note that
since $I$ is pointed, we have natural maps $i \copr i' \to i$, $i
\copr ' \to i'$ for any $i,i' \in I$, and then for any $X \in
\Ho^{st}(I)$, the induced map $X(i \copr i') \to X(i) \times X(i')$
is a weak equivalence (if $i = i'$, this follows from
Definition~\ref{seg.def}, and if not, observe that $i \copr i'$ is a
retract of $(i \copr i') \copr (i \copr i')$). For any bounded
half-additive category $I$, the functor \eqref{ho.Ho} extends to a
functor
\begin{equation}\label{ho.st}
\ho^{st}:\Ho^{st}(I) \to \Fun(I,\Ho^{st}(\Gamma_+)), \quad E
\mapsto \ho(m^*E),
\end{equation}
where $\ho(m^*E)$ is taken with respect to the projection $\Gamma_+
\times I \to I$. A functor $\gamma:I' \to I$ between half-additive
categories $I$, $I'$ is {\em half-additive} if it preserves finite
coproducts, and if $I$, $I'$ are bounded, such a functor induces a
pullback functor $\gamma^*:\Ho^{st}(I) \to \Ho^{st}(I')$.

\begin{lemma}\label{stab.le}
For any bounded half-additive category $I$, the full embedding
$\Ho^{st}(I) \subset \Ho(I)$ admits a left-adjoint functor
$\Stab_I:\Ho(I) \to \Ho^{st}(I)$.
\end{lemma}

\proof{} Take a finite set $S$, and consider the product
$\Gamma_+^S$ of copies of the category $\Gamma_+$ indexed by
elements $s \in S$.  Since $\Gamma_+$ is pointed, the projection
$\tau_s:\Gamma_+^S \to \Gamma_+$ onto the component corresponding to
some $s \in S$ has a left and right-adjoint functor $i_s:\Gamma_+
\to \Gamma_+^S$, and since $\Gamma_+$ has finite coproducts, the
diagonal embedding $\delta_S:\Gamma_+ \to \Gamma_+^S$ has a
left-adjoint functor $\Sigma_S:\Gamma_+^S \to \Gamma_+$. Moreover,
we have tautological anchor maps $\Sigma_S \to \tau_s$, $s \in S$,
and and these induce a map
\begin{equation}\label{Sig}
\Sigma_S^*X \to \prod_{s \in S}\tau_s^*X
\end{equation}
for any $X \in \Ho(\Gamma_+)$ that is an isomorphism if $X$ is
stable. But $\delta_S^*$ is left-adjoint to $\Sigma_S^*$, and
$i_s^*$ is left-adjoint to $\tau_s^*$ for any $s \in S$. Therefore
by adjunction, \eqref{Sig} induces an isomorphism
\begin{equation}\label{Sig.st}
\Stab(\delta_S^*Y) \cong \prod_{s \in S}\Stab(i_s^*Y)
\end{equation}
for any $Y \in \Ho(\Gamma_+^S)$. In particular, if we consider the
endofunctor $\sigma_S = \Sigma_S \circ \delta_S$ of the category
$\Gamma_+$, then for any $Y \in \Ho(\Gamma_+)$, we have a natural
isomorphism
\begin{equation}\label{stab.S}
\Stab(\sigma_S^*Y) \cong \Stab(Y)^S.
\end{equation}
Now take a bounded half-additive category $I$, with the functors
\eqref{ga.m} and \eqref{ga.e}, and let $\Stab_I:\Ho(I) \to \Ho(I)$
be the composition $e^* \circ \Stab \circ m^*$, with a map $\Id
\cong e^* \circ m^* \to \Stab_I$ induced by the adjunction map $\Id
\to \Stab$. Then for any stable $X \in \Ho^{st}(I)$, the map $X \to
\Stab_I(X)$ is obviously an isomorphism, and for any $Y \in \Ho(I)$,
$\Stab_I(Y)$ is stable by \eqref{stab.S}.
\endproof

Slightly more generally, for any bounded category $I'$ and bounded
half-additive $I$, say that $X \in \Ho(I \times I')$ is stable if so
is its pullback $(m \times \id)^*X$ in $\Ho(\Gamma_+ \times I \times
I')$, and let $\Ho^{st}(I,I') \subset \Ho(I \times I'$ be the full
subcategory spanned by stable objects (if $I=\Gamma_+$, this reduces
to Definition~\ref{stab.def} and the notation of
Proposition~\ref{stab.prop}).

\begin{corr}\label{stab.corr}
For any bounded $I'$ and bounded half-additive $I$, the embedding
$\Ho^{st}(I,I') \subset \Ho(I \times I')$ admits a left-adjoint
stabilization functor $\Stab_I$, and for any bounded $I''$ and
functor $\gamma:I'' \to I'$, we have the base change isomorphism
$(\id \times \gamma)^* \circ \Stab_I \cong \Stab_I \circ (\id \times
\gamma)^*$. Moreover, for any bounded hald-additive $I''$ and
half-additive functor $\gamma:I'' \to I$, we have the base chaneg
isomorphism $(\gamma \times \id)^* \circ \Stab_I \cong \Stab_{I''}
\circ (\gamma \times \id)^*$.
\end{corr}

\proof{} As in Lemma~\ref{stab.le}, take $\Stab_I = (e \times \id)^*
\circ \Stab \circ (m \times \id)^*$. This commutes with base change
by Remark~\ref{hoco.stab.rem}, and then defines the desired adjoint
by the same argument as in Lemma~\ref{stab.le}.
\endproof

For the second corollary, again let $I$ be an arbitrary bounded
category, take some $n \geq 2$, and consider the product $\Gamma_+^n
\times I$. Then it can be decomposed as $\Gamma_+ \times
(\Gamma^{n-1}_+ \times I)$ in $n$ different ways. Say that $X \in
\Ho(\Gamma_+^n \times I)$ is {\em polystable} if it is stable with
respect to all $n$ decompositions, and let $\Ho^{st}(\Gamma_+^n, I)
\subset \Ho(\Gamma_+^n \times I)$ be the full subcategory spanned by
polystable objects. Then smash product in $\Gamma_+$ induces a
functor
\begin{equation}\label{ga.m.n}
\Gamma_+^n \times I \to \Gamma \times I
\end{equation}
over $I$, and for any splitting $\Gamma_+^n = \Gamma_+ \times
\Gamma_+^{n-1}$ and objects $S_+,S'_+ \in \Gamma_+$, $S_\idot \in
\Gamma_+^{n-1}$, we have $m_n((S_+ \vee S'_+) \times S_\idot) \cong
m_n(S_+ \times S_\idot) \vee m_n(S'_+ \times S_\idot)$. Therefore
$m_n^*$ sends stable objects to polystable objects and induces a
functor
\begin{equation}\label{m.n.pb}
m_n^*:\Ho^{st}(\Gamma_+, I) \to \Ho^{st}(\Gamma_+^n, I),
\end{equation}
We then have the following result.

\begin{lemma}\label{poly.st.le}
For any $n \geq 2$ and $I$, the functor \eqref{m.n.pb} is an
equivalence.
\end{lemma}

\proof{} The functor $m_n$ of \eqref{ga.m.n} is the composition of
functors $m_2$ for the categories $\Gamma_+^l \times I$, $0 \leq l <
n$, so by induction, it suffices to consider the case
$n=2$. Moreover, \eqref{m.n.pb} admits a left-adjoint $\Stab \circ
m_{n!}$, and we have to check that the adjunction maps $\Id \to
m_n^* \circ \Stab \circ m_{n!}$, $\Stab \circ m_{n!} \circ m_n^* \to
\Id$ are invertible. By Remark~\ref{hoco.stab.rem}, both commute
with base change with respect to the category $I$, so it suffices to
consider the case $I=\ppt$, when $m_2$ is the usual smash product
functor $m:\Gamma_+^2 \to \Gamma_+$.

Let $r:[1] \to \Gamma_+$ be the embedding sending $0$ to $\ppt_+$
and $1$ to $\emptyset_+$, with the unique anchor maps between them,
and let $t_l:\Gamma_+ \to \Gamma_+ \times [1]$ be the embedding onto
$\Gamma_+ \times l$, $l=0,1$.  Fix a splitting $\Gamma_+^2 =
\Gamma_+ \times \Gamma_+$, and consider the embedding $\id \times
r:\Gamma_+ \times [1] \to \Gamma_+ \times \Gamma_+ \cong
\Gamma_+^2$. Denote $q = m \circ (\id \times r):\Gamma_+ \times [1]
\to \Gamma_+$. Denote by $\Ho_0(\Gamma_+ \times [1]) \subset
\Ho(\Gamma_+ \times [1])$ and by $\Ho_0^{st}(\Gamma_+, [1]) \subset
\Ho^{st}(\Gamma_+ \times [1])$ the full subcategories spanned by
objects with contractible $t_1^*X$. Then $q$ and $\id \times r$
induce functors
\begin{equation}\label{q.r.mn}
\begin{CD}
\Ho^{st}(\Gamma_+) @>{q^*}>> \Ho^{st}_0(\Gamma_+, [1]) @<{(\id
  \times r)^*}<< \Ho^{st}(\Gamma_+^2),
\end{CD}
\end{equation}
and it suffices to prove that both are equivalences.

For $q^*$, for any $S_+ \in \Gamma_+$, consider the functor
$\eps(S_+):\V^o \to \Gamma_+ \times [1]$ sending $0$, $1$, $o$ to
$S_+ \times 0$, $S_+ \times 1$, $\emptyset_+ \times 1$, and note
that $q^o \circ \eps(S_+)^o:\V \to \Gamma^o_+$ sends $0$ to $S_+$
and the rest to the initial object $\ppt_+ \in \Gamma_+^o$, so that
it is naturally augmented by $S_+$. Moreover, this defines a framing
for $q^o$ in the sense of Lemma~\ref{kan.le}, and if we compute the
homotopy right Kan extension $q_*$ using this framing, we see that
$q_*X$ for any $X \in \Ho(\Gamma_+ \times [1])$ fits into a homotopy
cartesian square
\begin{equation}\label{q.sq}
\begin{CD}
q_* X @>>> t^*_0X\\
@VVV @VVV\\
X(\emptyset_+ \times 1) @>>> t_1^*X.
\end{CD}
\end{equation}
Then if $X$ is in $\Ho^{st}_0(\Gamma_+, [1])$, the bottom arrow in
\eqref{q.sq} is a homotopy equivalence, so that $q_*X \cong t_0^*X$
is stable, and $q_*$ then gives an equivalence inverse to $q^*$.

For $(\id \times r^*)$, let $\langle a,s \rangle$ be the
anchor/structural factorization system on the category
$\Gamma_+$. Then $r:[1] \to (\Gamma_+)_a$ is full, so that by
Example~\ref{comma.exa}, we have a framing for $r$ given by
categories $[1] /^r_s S_+$, $S_+ \in \Gamma_+$. Explicitly, we have
$[1] /^r_s S_+ \cong S^<$, the discrete category $S$ with the added
initial object $o$. This induces a framing for $\id \times r$, and
then for any $X \in \Ho_0(\Gamma_+ \times [1])$ and $S_+ \in
\Gamma_+$, we have
\begin{equation}\label{X.t.0}
(\id \times r)_!X |_{\Gamma_+ \times S_+} \cong \bigvee_{s \in S}
t_0^*X \cong t_0^*X \wedge S_+.
\end{equation}
But for any $\Gamma$-space $Z \in \Ho(\Gamma_+)$, the product $Z
\wedge S_+ = Z \vee \dots \vee Z$ is of the form $\delta^*Z^{\vee
  S}$ for the object
$$
Z^{\vee S} = \bigvee_{s \in S}\tau_s^*Z \in \Ho(\Gamma_+^S),
$$
and then as in the proof of Lemma~\ref{stab.le}, the isomorphism
\eqref{Sig.st} shows that we have
\begin{equation}\label{stab.Z}
\Stab(Z \wedge S_+) \cong \Stab(Z)^S.
\end{equation}
Therefore if $X$ is stable, $(\Stab \circ (\id \times r)_!X)
|_{\Gamma_+ \times S_+} \cong (t_0^*X)^S$ by \eqref{X.t.0}. Then
firstly, $\Stab((\id \times r)_!X)$ that {\em a priori} only lies in
$\Ho^{st}(\Gamma_+, \Gamma_+) \subset \Ho(\Gamma_+)$ is polystable,
secondly, the adjunction map $X \to (\id \times r)^* \circ \Stab
\circ (\id \times r)_!X$ is an isomorphism, and finally, the
adjunction map $\Stab((\id \times r)_!r^*Y) \to Y$ is an isomorphism
for any polystable $Y \in \Ho^{st}(\Gamma_+^2,\ppt) \subset
\Ho(\Gamma_+^2)$. Then $\Stab \circ (\id \times r)_!$ is the inverse
equivalence to $(\id \times r)^*$, and we are done.
\endproof

\begin{corr}
For a bounded half-additive category $I$, the pullback functor
$m^*:\Ho^{st}(I) \to \Ho^{st}(\Gamma_+,I)$ is fully faithful, with a
left-adjoint given by restricting $e^* \circ \Stab_I:\Ho(I \times
\Gamma_+) \to \Ho(I)$ to $\Ho^{st}(\Gamma_+,I) \subset \Ho(\Gamma_+
\times I)$.
\end{corr}

\proof{} Since stabilization commutes with finite product by
Remark~\ref{hoco.stab.rem}, the stabilization functor
$\Stab_I:\Ho(\Gamma_+ \times I) \to \Ho^{st}(I,\Gamma_+) \subset
\Ho(\Gamma_+ \times I)$ of Corollary~\ref{stab.corr} sends
$\Ho^{st}(\Gamma_+,I)$ into itself. Therefore if we denote by
$\Ho^{2st}(\Gamma_+,I) = \Ho^{st}(\Gamma_+,I) \cap
\Ho^{st}(I,\Gamma_+) \subset \Ho(\Gamma_+,I)$ the full subcategory
spanned by objects $X \in \Ho(\Gamma_+ \times I)$ that are stable in
each variable, the embedding $\Ho^{2st}(\Gamma_+,I) \to \Ho(\Gamma_+
\times I)$ admits a left-adjoint induced by $\Stab_I$, and it
suffices to prove that $m^*:\Ho^{st}(I) \to \Ho^{2st}(\Gamma_+,I)$
is an equivalence of categories. Since $e^* \circ m^* \cong \id$, it
further suffices to prove that $e^*:\Ho^{2st}(\Gamma_+,I) \to
\Ho^{st}(I)$ is an equivalence. However, it has a left-adjoint given
by $\Stab_I \circ \Stab \circ e_!$, and since both stabilizations
and $e_!$ satisfy base change with respect to $I$, it suffices to
consider the case $I=\Gamma_+$. This is Lemma~\ref{poly.st.le}.
\endproof

\subsection{Symmetric monoidal structures.}\label{mono.subs}

A convenient way to package all the data associated to a symmetric
monoidal category is the same Segal machine as in
Definition~\ref{seg.def}, with functors to spaces replaced by
Grothendieck cofibrations.

\begin{defn}\label{sym.1.def}
A {\em unital symmetric monoidal structure} on a category $\C$ is a
cofibration $\Bi\C \to \Gamma_+$ equipped with an equivalence
$\Bi\C_{\ppt_+} \cong \C$ such that $\Bi\C_{o}$ is the point
category $\ppt$, and for any $S,S' \in \Gamma$, the functor
$$
a_! \times a'_!:\Bi\C_{S_+ \vee S'_+} \to \Bi\C_{S_+} \times
\Bi\C_{S'_+}
$$
induced by the maps \eqref{S.spl} is an equivalence of categories.
\end{defn}

Informally, if we consider the set $\{0,1\}$ with two elements, then
the two anchor maps $\{0,1\}_+ \to \ppt_+$ identify the fiber
$\Bi\C_{\{0,1\}_+}$ with the product $\C^2$, and the unique
structural map $\{0,1\}_+ \to \ppt_+$ provides the product functor
$\mu:\C^2 \to \C$. All the associativity, commutativity and
unitality constraints are packaged into the maps \eqref{fib.compo}
for the cofibration $\Bi\C \to \Gamma_+$.

\begin{exa}\label{point.sym.exa}
For any bounded category $I$, and any category $\C$ equipped with a
unital symmetric mono\-i\-dal structure $\Bi\C$, the functor
category $\Fun(I,\C)$ carries a natural {\em pointwise} unital
symmetric monoidal structure given by $\Bi\Fun(I,\C) \cong
\Fun(I,\Bi\C/\Gamma_+)$.
\end{exa}

\begin{exa}\label{prod.mono.exa}
A category $\C$ that has finite cartesian products carries a unital
symmetric monoidal structure $\Bi\C$ given by these
products. Formally, take the embedding $e:\ppt \to \Gamma_+$ of
\eqref{ga.e}, observe that the fibration $\Bi'\C = (e^o_*\C^o)^o \to
\Gamma_+$ is also a cofibration, and take the preimage $\Bi\C =
\zeta(o)^{-1}(1) \subset \Bi'\C$ of the terminal object $1 \in \C =
(\Bi'\C)_{o}$ under the functor \eqref{beta.eq} for the terminal
object $o=\emptyset_+ \in \Gamma_+$. By
Lemma~\ref{fib.le}~\thetag{ii}, $\Bi\C \to \Gamma_+$ is a
cofibration, and the embedding $\Bi\C \subset \Bi'\C$ has a
left-adjoint functor $\lambda$ cocartesian over $\Gamma_+$.
\end{exa}

\begin{exa}\label{X.mono.exa}
A commutative monoid $X$ considered as a discrete category carries a
unital symmetric monoidal structures $\Bi X = \Gamma_+[X_\infty]$,
where $X_\infty:\Gamma_+ \to \Sets_+$ is the discrete special
$\Gamma$-space corresponding to $X$, with the unit element as the
distinguished point.
\end{exa}

\begin{defn}\label{sym.2.def}
A {\em lax monoidal structure} on a functor $F:\C' \to \C$ between
two categories $\C$, $\C'$ equipped with unital symmetric monoidal
structures $\Bi\C$ resp.\ $\Bi\C'$ is a functor $\Bi F:\Bi\C' \to
\Bi\C$ over $\Gamma_+$, cocartesian along anchor maps in $\Gamma_+$
and equipped with an isomorphism $\Bi F_{\ppt_+} \cong F$. A lax
monoidal structure $\Bi F$ is {\em monoidal} if it is cocartesian
over all maps.
\end{defn}

The essential structure for a lax monoidal functor $F$ is the map
$\mu \circ F^2 \to F \circ \mu'$ provided by \eqref{fi.fu} for the
functor $\Bi F$. The functor $F$ is monoidal iff this map is an
isomorphism.

\begin{exa}\label{inf.adj.exa}
Assume that a functor $\gamma:\C \to \C'$ between two symmetric
monoidal categories is equipped with a symmetric monoidal structure
$\Bi\gamma$, and at the same time admits a right-adjoint
$\gamma^\dg:\C' \to \C$. Then $\Bi\gamma$ admits a right-adjoint
$\Bi\gamma^\dg$ over $\Gamma_+$ that defines a lax symmetric
monoidal structure on $\gamma\dg$.
\end{exa}

\begin{exa}\label{oppa.exa}
For any category $\C$ equipped with a unital symmetric
mo\-no\-i\-dal structure $\Bi\C$, the opposite category $\C^o$
carries a unital symmetric monoidal structure $\Bi\C^o =
\Bi\C^{\perp o}$. For any two categories $\C$, $\C'$ with unital
symmetric monoidal structures $\Bi\C$, $\Bi\C'$, the product $\C
\times \C'$ carries a unital symmetric monoidal structure $\Bi (\C
\times \C') = \Bi\C \times_{\Gamma_+} \Bi\C'$. For any unital
symmetric monoidal category $\C$, the Yoneda pairing
\begin{equation}\label{yo.mono}
\Y:\C^o \times \C \to \Sets
\end{equation}
carries a natural lax monoidal structure $\Bi\Y$ with respect to the
cartesian product in $\Sets$.
\end{exa}

\begin{exa}\label{alg.exa}
For any object $A \in \C$ in a unital symmetric monoidal category
$\C$, giving a lax monoidal structure $\Bi i_A$ on the embedding
$i_A:\ppt \to \C$ is equivalent to turning $A \in \C$ into a unital
commutative associative algebra.
\end{exa}

\begin{exa}\label{ga.ga.exa}
By Example~\ref{prod.mono.exa} and Example~\ref{oppa.exa}, any
category $I$ that has finite coproducts is a symmetric unital
monoidal category with respect to coproducts. In particular, this
applies to $\Gamma$; the corresponding category $\Bi\Gamma$ is
$\Ar^s(\Gamma_+)$, where $s$ is the class of structural maps. For
the opposite category $\Gamma^o$, we have $\Bi\Gamma^o \cong
\Tw^s(\Gamma_+)$. The structural cofibration $\Bi\Gamma \cong
\Ar^s(\Gamma_+) \to \Gamma_+$ admits a right-adjoint $\eta:\Gamma_+
\to \Bi\Gamma$, and by Example~\ref{inf.adj.exa} and
Example~\ref{alg.exa}, this turns the one-point set $\ppt \in
\Gamma$ into an algebra object in $\C$. This algebra object is
universal: for any unital symmetric monoidal category $\C$ and
object $A \in \C$, a lax monoidal structure $\Bi i_A$ on the
embedding $i_A:\ppt \to \C$ factors uniquely as $\Bi i_A \cong \Bi
I_A \circ \eta$ for a unique monoidal functor $I_A:\Gamma \to
\C$. Explicitly, $\Bi I_A = \eta^{\Gamma_+}_!\Bi i_A$, where the
left Kan extension is given by \eqref{adj.I.eq}.
\end{exa}

\begin{exa}\label{ga.iso.exa}
The isomorphism groupoid $\overline{I}$ of a unital symmetric
mono\-i\-dal category $I$ carries a unital symmetric monoidal
structure $\Bi\overline{I} = (\Bi I)_\natural$. For $\Gamma$ with
$\Bi\Gamma$ of Example~\ref{ga.ga.exa}, the groupoid
$\overline{\Gamma}$ also has a universal property: for any object $c
\in \C$ in a unital symmetric monoidal category, there exists a
unique symmetric monoidal functor $\gamma:\overline{\Gamma} \to \C$
sending $\ppt$ to $c$. In particular, if we take the monoid $\N$ of
non-negative integers with respect to addition, and treat is a
discrete unital symmetric monoidal category as in
Example~\ref{X.mono.exa}, then $1 \in \N$ defines a monoidal
functor $c:\Gamma \to \N$ sending a set to its cardinality. Note
that $\Bi c$ has no sections.
\end{exa}

For any unital symmetric monoidal categories $I$, $\C$ with bounded
$I$, lax monoidal functors from $I$ to $\C$ form a category that we
denote $\Fun_\infty(I,\C)$. More generally, if we are given a
functor $\gamma:\C \to I$ equipped with a monoidal structure
$\Bi\gamma$, then sections of $\Bi\gamma$ that are lax monoidal
functors form a category that we denote $\Sec_\infty(I,\C)$. We
tautologically have $\Fun_\infty(I,\C) \cong \Sec_\infty(I,I \times
\C)$, and sections of the tautological projection $\C \to \ppt$ are
algebra objects in $\C$, as in Example~\ref{alg.exa}. It is useful
to observe that under favourable conditions, arbitrary lax monoidal
functors can be also described as algebra objects. Namely,
informally, for any small category $I$ and symmetric monoidal
category $\C$ with a product $- \otimes -$, the product in $\C$
induces a box product
\begin{equation}\label{box.eq}
- \boxtimes -:\Fun(I,\C)^2 \to \Fun(I^2,\C), \quad (F \boxtimes F')
(i \times i') = E(i) \otimes E(i')
\end{equation}
If $I$ itself is unital symmetric monoidal, we also have the product
functor $m:I^2 \to I$ and the pullback functor $m^*:\Fun(I,\C) \to
\Fun(I^2,\C)$. It might happen that $m^*$ has a left-adjoint $m_!$
that commutes with all products $- \otimes c$, $c \in \C$. In this
case, we can turn $\Fun(I,\C)$ into a symmetric monoidal tensor
category by considering the convolution product
\begin{equation}\label{conv.eq}
F \circ F' = m_!(F \boxtimes F').
\end{equation}
More formally, let $(a,s)$ be the anchor/structural factorization
system on $\Gamma_+$, and use the technology of
Subsection~\ref{ker.subs}. If denote by $\Bi^s I \to \Gamma$
the restriction the cofibration $\Bi I \to \Gamma_+$ to $\Gamma =
(\Gamma_+)_s \subset \Gamma_+$, then it immediately follows from
\eqref{ffun.r} that the cofibration
\begin{equation}\label{Bi.FF}
\Bi\Fun(\Bi^s I/\Gamma,\C) \cong \fFun^a(\Bi I/\Gamma_+,\Bi
\C/\Gamma_+)  \to \Gamma_+
\end{equation}
defines a symmetric monoidal structure on the category $\Fun(\Bi^s
I/\Gamma,\C)$ of \eqref{ffun.st}. Explicitly, ob\-jects in
$\Fun((\Bi I)_s/\Gamma,\C)$ are pairs $\langle S,F \rangle$, $S \in
\Gamma$, $F:I^S \to \C$, and symmetric monoidal structure is induced
by the box product \eqref{box.eq}: we have $\langle S,F \rangle
\circ \langle S',F' \rangle = \langle S \copr S',F \boxtimes F'
\rangle$. The projection $\Fun(\Bi^s I/\Gamma,\C) \to \Gamma^o$ is
monoidal with respect to the monoidal structure $\Bi\Gamma^o \cong
\Tw^s(\Gamma_+)$ of Example~\ref{ga.ga.exa}, and \eqref{ffun.sec}
reads as
\begin{equation}\label{mo.ff}
\Fun_\infty(I,\C) \cong \Sec_\infty^\natural(\Gamma^o,\Fun(\Bi^s
I/\Gamma,\C),
\end{equation}
where $\Sec_\infty^\natural(-,-)$ stands for the full subcategory
spanned by cocartesian sections. We also have the reflection of the
cofibration \eqref{Bi.FF} over $\Bi\Gamma^o$, and one immediately
checks by induction that if $m^*$ admits a left-adjoint $m_!$ such
that $m_!(a \otimes c) \cong m_!a \otimes c$, $a \in \Fun(I^2,\C)$,
$c \in \C$, then \eqref{ffun.p} satisfies the conditions of
Lemma~\ref{refl.cof.le}, so that \eqref{fun.ii} is a cofibration.
Then again, one checks by \eqref{ffun.r} that
\begin{equation}\label{Bi.F}
\Bi\Fun(I,\C) \cong \Fun^a((\Bi I/\Gamma_+,\Bi \C/\Gamma_+) \to
\Gamma_+
\end{equation}
defines a symmetric monoidal structure on $\Fun(I,\C)$. This is our
convolution monoidal structure, and \eqref{fun.L.eq} then provides an
equivalence
\begin{equation}\label{mo.f}
\Fun_\infty(I,\C) \cong \Sec_\infty(\ppt,\Fun(I,\C)),
\end{equation}
so that lax monoidal functors from $I$ to $\C$ are identified with
algebra objects in $\Fun(I,\C)$.

\begin{remark}\label{kan.mono}
A monoidal functor $\gamma:I' \to I$ between bounded symmetric
monoidal categories induces a cocartesian functor $\Bi^s\gamma:\Bi^s
I' \to \Bi^s I$, and the pullback functor $(\Bi^s\gamma)^*$ is
naturally monoidal with respect to the structures \eqref{Bi.FF}. If
the cofibrations \eqref{Bi.F} exist, then $\gamma^*:\Fun(I,\C) \to
\Fun(I',\C)$ is naturally lax monoidal (but only lax since
$\gamma^*$ commutes with $m^*$ but not necessarily with
$m_!$). However, if $\gamma^*$ admits a left-adjoint Kan extension
functor $\gamma_!$, then $\gamma_!$ does commute with $m_!$. In this
case, $\Bi\gamma_*$ admits a left-adjoint $\Bi\gamma_!$ over
$\Gamma_+$ that is cocartesian and defines a monoidal structure on
$\gamma_!$.
\end{remark}

\begin{exa}\label{symseq.exa}
Let $\C$ be a cocomplete unital symmetric monoidal category such
that $c \otimes -:\C \to \C$ preserves colimits for any $c \in \C$,
and let $\N$, $\overline{\Gamma}$ be the unital symmetric monoidal
categories of Example~\ref{ga.iso.exa}. Then $\Fun(\N,\C)$ is the
monoidal category of non-negatively graded objects in $\C$. Objects
in $\Fun(\overline{\Gamma},\C)$ are ``symmetric sequences'' of
\cite{HSS}. If $c:\overline{\Gamma} \to \N$ is the cardinality
functor, then the pullback $c^*$ is fully faithful, and the left Kan
extension $c_!:\Fun(\overline{\Gamma},\C) \to \Fun(\N,\C)$ is
monoidal by Remark~\ref{kan.mono}, so that $\Fun(\N,\C)$ is a
monoidal localization of $\Fun(\overline{\Gamma},\C)$. There are
other interesting localizations. For example, if $\C=\Z\amod$ is the
category of abelian groups, then we have the symmetric sequence
$\Z_\idot \in \Fun(\overline{\Gamma},\C)$ given by the sign
representations, the pointwise tensor product $\Z_\idot \otimes
\Z_\idot$ is $\Z$, and we have a full embedding
$c^*_\idot:\Fun(\N,\Z) \to \Fun(\overline{\Gamma},\Z)$ given by
$c^*_\idot(E) = c^*E \otimes \Z_\idot$. It has a left-adjoint
$c_!^\hdot$ sending $E_\idot \in \Fun(\overline{\Gamma},\Z)$ to
$c_!(E_\idot \otimes \Z_\idot)$, and by
Lemma~\ref{fib.le}~\thetag{ii}, this induces a new symmetric
monoidal product on $\Fun(\N,\Z)$ such that $c_!^\hdot$ is
monoidal. This new product is the old one but twisted by the usual
homological signs, so that commutative algebras in $\Fun(\N,\Z)$
with the new symmetric monoidal structure are graded-commutative
algebras over $\Z$.
\end{exa}

\subsection{Multiplicative $\Gamma$-spaces.}\label{mult.subs}

We now recall the lesser-known part of \cite{seg} that deals with
multiplications. We note that the smash product defines a unital
symmetric monoidal structure both on $\Gamma_+$ and on
$\Top_+$. Since smash product preserves weak equivalences, it also
induces a unital symmetric monoidal structure on $\Ho$.

\begin{defn}
A {\em multiplicative $\Gamma$-space} is a lax monoidal functor
$X$ from $\Gamma_+$ to $\Top_+$.
\end{defn}

By virtue of \eqref{mo.f}, multiplicative $\Gamma$-spaces correspond
to algebra objects in $\Fun(\Gamma_+,\Top_+)$ with respect to the
convolution product, but this not very useful since the convolution
product is not homotopy invariant. To construct a homotopy invariant
formalism, it is better to use \eqref{mo.ff}. For any bounded unital
symmetric monoidal category $I$, $\Fun(\Bi^s I/\Gamma,\Top_+) \to
\Gamma^o$ is a cofibration with fibers $\Fun(I^n,\Top_+)$, $n \geq
0$, and its transition functors respect weak equivalences. Therefore
it induces a cofibration $\Ho(\Bi^s I/\Gamma)$ with fibers
$\Ho(I^n)$, and since smash product in $\Top_+$ also preserves weak
equivalences, \eqref{Bi.FF} induces a symmetric monoidal structure
$\Bi\Ho(\Bi^s I/\Gamma)$ on $\Ho(\Bi^s I/\Gamma)$. We still have a
cofibration $\Bi\Ho(\Bi^s I/\Gamma) \to \Bi\Gamma^o$, and since
homotopy left Kan extensions commute with smash product $- \wedge X$
for any $X \in \Top_+$, the reflection $(\Bi\Ho(\Bi^s
I/\Gamma)/\Bi\Gamma^o)^\lhd$ again satisfies the conditions of
Lemma~\ref{refl.cof.le}. As in \eqref{fun.ii}, we can then restrict
it to the subcategory $\Gamma_+ \subset \Ar^s(\Gamma_+) =
\Bi\Gamma$, and obtain a cofibration
\begin{equation}\label{bi.ho}
\Bi\Ho(I) \to \Gamma_+
\end{equation}
that defines a convolution symmetric monoidal structure on the
homotopy category $\Ho(I)$. We still have the equivalence
\begin{equation}\label{ho.sec.mono}
\Sec^\natural_\infty(\Gamma^o,\Ho(\Bi^s I/\Gamma)) \cong
\Sec_\infty(\ppt,\Ho(I)),
\end{equation}
and we can treat an object in either of these two equivalent
categories as an enhanced version of a lax monoidal functor from $I$
to $\Ho$. By \eqref{mo.ff}, a geniune lax monoidal functor $I \to
\Top_+$ generates such an object by localization, and forgetting the
enhancement corresponds to considering the associated functor
\eqref{ho.Ho} that is lax symmetric monoidal with respect to the
monoidal structure \eqref{bi.ho}. We also note that as in
Remark~\ref{kan.mono}, the homotopy Kan extension $\gamma_!$ with
respect to a monoidal functor $\gamma:I' \to I$ carries a natural
monoidal structure.

Now take a bounded unital symmetric monoidal category $I$, and
consider the product $\Gamma_+ \times I$ with the product monoidal
structure $\Bi(\Gamma_+ \times I)$ of Example~\ref{oppa.exa}. Then
$\Ho(\Gamma_+ \times I)$ carries a unital symmetric monoidal
structure $\Bi\Ho(\Gamma_+ \times I)$ of \eqref{bi.ho}.

\begin{prop}\label{mult.prop}
Let $\Bi\Ho^{st}(\Gamma_+, I) \subset \Bi\Ho(\Gamma_+ \times I)$ be
the full subcategory spanned by $\Ho^{st}(\Gamma_+, I)^S \subset
\Ho(\Gamma_+ \times I)^S \cong \Bi\Ho(\Gamma_+ \times I)_{S_+}$, $S
\in \Gamma$. Then the induced projection $\Bi\Ho^{st}(\Gamma_+, I)
\to \Gamma_+$ is a cofibration that defines a unital symmetric
monoidal structure on $\Ho^{st}(\Gamma_+, I)$, and the stabilization
functor $\Stab$ of Proposition~\ref{stab.prop} extends to a monoidal
functor $\Bi\Stab:\Bi\Ho(\Gamma_+ \times I) \to
\Bi\Ho^{st}(\Gamma_+, I)$.
\end{prop}

\proof{} For any bounded category $I$, say that a map $f$ in the
category $\Ho(\Gamma_+ \times I)$ is {\em stably trivial} if
$\Stab(f)$ is invertible. Moreover, for any $S \in \Gamma$, say that
a map $g$ in $\Ho(\Gamma_+ \times I)^S$ is stably trivial if so is
each of its components $g_s$, $s \in S$. Then by
Lemma~\ref{fib.le}~\thetag{ii}, it suffices to prove that the
transition functors $f_!$ of the cofibration \eqref{bi.ho} send
stably trivial maps to stably trivial maps.

If $f$ is an anchor map, then $f_!$ commutes with stabilization by
definition, so by induction, it suffices to consider the case when
$f:\{0,1\}_+ \to \ppt_+$ is the structural map that induces the
product on $\Ho(\Gamma_+ \times I)$. Explicitly, if we denote by
$m:(\Gamma_+ \times I)^2 \to \Gamma_+ \times I$ the smash product on
$\Gamma_+ \times I$, then we have to check that for any stably
trivial map $g$ in $\Ho(\Gamma_+ \times I)$ and any $X \in
\Ho(\Gamma_+ \times I)$, $m_!(g \boxtimes X)$ is stably
trivial. Moreover, since $m^*$ sends stable objects to stable
objects, and the stabilization functor of
Proposition~\ref{stab.prop} commutes with base change with respect
to $I$, this amounts to checking that for any $Y \in \Ho$, $g \wedge
Y$ is stably trivial, and it suffices to check it for $I = \ppt$.

Indeed, if $Y=S_+$ is a finite pointed set, then the claim
immediately follows from \eqref{stab.Z}. But by adjunction, a
homotopy colimit of stably trivial maps is stably trivial. Therefore
the claim is clear for arbitrary pointed sets, and then for homotopy
colimits over $\Delta^o$ of arbitrary simplicial pointed sets. This
is the whole $\Ho$.
\endproof

\begin{remark}\label{weak.rem}
The category $\Ho^{st}(\Gamma_+)$ is of course just the category of
connective spectra, the connective part of the standard
$t$-structure on the stable homotopy category, and the monoidal
structure of Proposition~\ref{mult.prop} is the standard smash
product. A commutative associative unital algebra object in
$\Ho^{st}(\Gamma_+)$ is then what used to be called a ``commutative
associative unital connective ring spectrum'' at the time of
\cite{seg}, and Proposition~\ref{mult.prop} also shows that any
multiplicative $\Gamma$-space generates such a thing by
stabilization. The modern notion of an $E_\infty$-ring spectrum is of
course much stronger; however, the old naive version will be
sufficient for our purposes.
\end{remark}

\begin{remark}
Even if $X$ is a commutative ring spectrum in the weak sense of
Remark~\ref{weak.rem}, the homotopy groups $\pi_\idot(X)$ still form
a graded commutative $\Z$-algebra. To see this in the stabilization
formalism, one first observes that by the same argument as in
Proposition~\ref{mult.prop}, the stabilization functor $\Stab_0$ of
Remark~\ref{stab.0.rem} creates a unital symmetric monoidal
structure on $\Fun^{st}(\Gamma_+ \times I)$, and the equivalence
\eqref{stab.set} sends it to the convolution product on $\Fun(I,\Z)$
(in particular, if $I=\ppt$, we recover the standard tensor product
of abelian groups). Now take the circle $S^1 \in \Ho$, extend it to
a monoidal functor $\overline{\Gamma} \to \Ho$, $\ppt \mapsto S^1$
as in Example~\ref{ga.iso.exa}, and let
$\Sigma_\idot:\overline{\Gamma}^o \cong \overline{\Gamma} \to \Ho^o$
be the opposite functor. Then any algebra object $X$ in
$\Ho^{st}(\Gamma_+)$ defines a lax monoidal functor $\Ho(X):\Gamma_+
\to \Ho$, and we can consider the Yoneda pairing
$$
\pi_\idot(X) = \Y(\Sigma_\idot \times \Ho(X)):\overline{\Gamma}
\times \Gamma_+ \to \Sets
$$
of \eqref{yo.mono}. This is a lax monoidal functor, and moreover,
since $X$ is stable, $\pi_\idot(X)$ is stable as an object in
$\Fun^{st}(\Gamma_+ \times \overline{\Gamma})$. Thus by
\eqref{stab.set}, $\pi_\idot(X)$ is a symmetric sequence of abelian
groups that is an algebra with respect to the convolution product of
Example~\ref{symseq.exa}. One then checks that it actually lies in
the essential image of the full embedding $c_\idot^*$, thus gives a
graded commutative algebra over $\Z$.
\end{remark}

Assume now given a half-additive category $I$, and assume that $I$
is equipped with a unital symmetric monoidal structure. Then we say
that the structure is {\em distributive} if $0 \otimes i = 0$ for
any $i \in I$, where $0 \in I$ is the initial object, and the
natural map $(i_0 \otimes i_1) \copr (i_0' \otimes i_1) \to (i_0
\copr i_0') \otimes i_1$ is an isomorphism for any $i_0,i_0',i_1 \in
I$. In this case, the functors \eqref{ga.e} and \eqref{ga.m} are
monoidal, and we have the following corollary of
Proposition~\ref{mult.prop}.

\begin{corr}\label{mult.corr}
Assume given a bounded half-additive category $I$, and assume that
it is equipped with a distributive unital symmetric monoidal
structure. As in Proposition~\ref{mult.prop}, let $\Bi\Ho^{st}(I)
\subset \Bi\Ho(I)$ be the full subcategory spanned by $\Ho^{st}(I)^S
\subset \Ho(I)^S \cong \Bi\Ho(I)_{S_+}$, $S \in \Gamma$. Then the
induced projection $\Bi\Ho^{st}(I) \to \Gamma_+$ is a cofibration
that defines a unital symmetric monoidal structure on $\Ho^{st}(I)$,
and the stabilization functor $\Stab_I$ of Lemma~\ref{stab.le}
carries a monoidal structure $\Bi\Stab_I$.
\end{corr}

\proof{} By Lemma~\ref{fib.le}~\thetag{iii}, it suffices to
construct the left-adjoint functor $\Bi\Stab_I$. By
Lemma~\ref{stab.le}, $\Bi\Stab_I = \Bi t^* \circ \Bi \Stab \circ \Bi
m^*$ works.
\endproof

\subsection{Additivization.}\label{add.subs}

Let us now describe a linear version of the stabilization
formalism. Assume given a commutative ring $k$. For any bounded
category $I$, denote by $\Ho(I,k)$ the localization of the category
$\Fun(I,\Delta^o k\amod)$ with respect to weak equivalences. By
Dold-Kan, $\Delta^ok\amod$ is equivalent to the category $C_{\geq
  0}(k)$ of chain complexes of $k$-modules bounded from below by $0$
(in homological degrees), weak equivalences are quasiisomorphisms,
and $\Ho(I,k)$ is the connective part of the standard $t$-structure
on the derived category $\D(I,k)$ of Subsection~\ref{hom.subs},
while \eqref{ho.Ho} is induced by \eqref{D.ho}. A functor $\gamma:I'
\to I$ gives rise to a pullback functor $\gamma^*:\Ho(I,k) \to
\Ho(I,k)$, and we have the left and right homotopy Kan extensions
$\gamma_!$, $\gamma_*$ (where in homological terms, $\gamma_!$ is
the left-derived functor of the usual left Kan extension, and
$\gamma_*$ is the right-derived functor of the right Kan extension
composed with truncation with respect to the standard
$t$-structure). The isomorphisms \eqref{adm.eq}, \eqref{kan.eq},
\eqref{bc.eq} and Lemma~\ref{kan.le} also hold for the categories
$\Ho(-,k)$, and a weak equivalence $\gamma:I' \to I$ induces an
equivalence of categories $\gamma_\forall^*:\Ho^\forall(I,k) \cong
\Ho^\forall(I',k)$.

Since every $k$-module has a flat resolution, $\Ho(I,k)$ can
also be obtained by localizing the category
$\Fun(I,\Delta^ok\amod^{fl})$, where $k\amod^{fl} \subset k\amod$ is
the full subcategory spanned by flat $k$-modules. Then
$\Delta^ok\amod^{fl}$ is symmetric monoidal, the product respects
weak equivalences, and the same procedure as in
Subsection~\ref{mult.subs} provides a unital symmetric monoidal
structure on $\Ho(\Gamma_+,k)$. We have the forgetul functor
\begin{equation}\label{phi.eq}
\phi:k\amod^{fl} \to \Sets_+
\end{equation}
sending a $k$-module $M$ to its underlying set with $0 \in M$ as the
distinguished point, and it has a left-adjoint {\em reduced span}
functor $\Span_k:\Sets_+ \to k\amod^{fl}$ sending a pointed set
$S_+$ to the quotient $k[S_+]/k \cdot o$ of the free $k$-module
$k[S_+]$ generated by $S_+$ by the submodule $k \cdot o$ spanned by
the distinguished point $o \in S_+$. The functor $\Span_k$ is
monoidal with respect to the smash product on $\Sets_+$, so that
\eqref{phi.eq} is lax monoidal by adjunction, as in
Example~\ref{inf.adj.exa}. On the level of homotopy categories, for
any bounded $I$, \eqref{phi.eq} induces a functor $\phi:\Ho(I,k) \to
\Ho(I)$ sending a simplicial $k$-module to the geometric realization
of the underlying simplicial set, with its left-adjoint
$\Span_k:\Ho(I,k) \to \Ho(I)$, and these functors commute with
pullbacks $\gamma^*$ with respect to functors $\gamma:I' \to I$.

All the results about stabilization also have counterparts for the
$k$-linear homotopy categories $\Ho(-,k)$. Namely, we say that $X
\in \Ho(\Gamma_+ \times I,k)$ is stable if so is $\phi(X) \in
\Ho(\Gamma_+ \times I,k)$ and analogously, if $I$ is bounded and
half-additive, $X \in \Ho(I,k)$ is stable if so if $\phi(X)$. We
denote by
\begin{equation}\label{k.emb}
\Ho^{st}(\Gamma_+, I,k) \subset \Ho(\Gamma_+ \times I,k),
\qquad \Ho^{st}(I,k) \subset \Ho(I,k)
\end{equation}
the full subcategories spanned by stable objects, and the
stabilization functors of Proposition~\ref{stab.prop}
resp.\ Lemma~\ref{stab.le} also provide stabilization functors
$\Stab$, $\Stab_I$ left-adjoint to the embeddings \eqref{k.emb} (one
has to replace $\Top_+$ with $\Delta^ok\amod$, and use the truncated
homological shift instead of the loop functor $\Omega$). These
stablization functors are monoidal with respect to natural symmetric
monoidal structures. We also have natural isomorphisms
\begin{equation}\label{phi.stab}
\phi \circ \Stab \cong \Stab \circ \phi, \qquad \phi \circ \Stab_I
\cong \Stab_I \circ \phi,
\end{equation}
where $\phi$ are the forgetful functors induces by
\eqref{phi.eq}. In particular, the functor $\Stab_0$ of
Remark~\ref{stab.0.rem} is given by $\Stab_0(X) = \tau_{\leq
  0}(\Stab(\Z[X]))$, where $\tau_{\leq 0}$ is the truncation with
respect to the standard $t$-structure. A new feature in the
$k$-linear case is that finite products in $\Delta^ok\amod$ coincide
with coproducts, so that not only pullbacks but also left Kan
extensions preserve stability: for any bounded functor $\gamma:I'
\to I$ bewteen bounded categories, $(\id \times \gamma)^*$ and $(\id
\times \gamma)_!$ induce an adjoint pair of functors
\begin{equation}\label{k.st.pb}
\begin{aligned}
(\id \times \gamma)^*:&\Ho^{st}(\Gamma_+, I,k) \to
  \Ho^{st}(\Gamma_+, I',k),\\
(\id \times \gamma)_!:&\Ho^{st}(\Gamma_+, I',k) \to
  \Ho^{st}(\Gamma_+, I,k),
\end{aligned}
\end{equation}
and similarly for $\Ho^{st}(I,k)$, $\Ho^{st}(I',k)$ for a
half-additive $I$. Then $k$-linear stabilization commutes both with
pullbacks and left Kan extensions, so that \eqref{phi.stab} and the
isomorphism $\Stab \circ (\id \times \gamma)^* \cong (\id \times
\gamma)^* \circ \Stab$ of Remark~\ref{hoco.stab.rem} induce by
adjunction a natural isomorphism
\begin{equation}\label{stab.pb}
\Stab \circ (\id \times \gamma)_! \circ \phi \cong (\id \times \gamma_!)
\circ \phi \circ \Stab,
\end{equation}
and similarly for $\Stab_I$, $\Stab_{I'}$.

\begin{remark}\label{kawa.rem}
For any $\Gamma$-space $X \in \Ho(\Gamma_+)$ and any commutative
ring $k$, the isomorphism $\phi \circ \Stab \cong \Stab \circ \phi$
induces a map $\Stab(\Span_k(X)) \to \Span_k(\Stab(X))$, However,
this map is {\em not} an isomorphism. In effect, for any $X \in
\Ho$, $\Span_k(X) \in \Ho(k)$ is represented by the homology complex
of the space $X$ with coefficients in $k$, so that
$\Span_k(\Stab(X))(\ppt_+)$ is the homology of the infinite loop
space $\Stab(X)(\ppt_+)$. On the other hand,
$\Stab(\Span_k(X))(\ppt_+)$ is the homology of the corresponding
spectrum (this is more-or-less obvious from the construction, or see
a proof in \cite{ka.kawa}).
\end{remark}

Alternatively, one can construct the stabilization functors directly
by purely homological means; we only give a sketch, and refer the
reader to \cite[Section 3]{ka0} for details. We note that the
group-like requirement in the definition of stability is automatic
in the $k$-linear case, and again, the products and coproducts in
$\Ho(k)$ coincide. Thus in the $k$-linear context, ``stable'' just
means ``additive'', or rather, ``half-additive'' --- that is,
sending coproducts to coproducts. The category
$\Ho^{st}(\Gamma_+,k)$ is then equivalent to $\Ho(k)$, with the
equivalence $\Ho(k) \to \Ho^{st}(\Gamma_+,k) \subset
\Ho(\Gamma_+,k)$ sending $M \in \Ho(k)$ to $M \otimes_k T$, where $T
\in \Ho(\Gamma_+,k)$ is the reduced span functor, $S_+ \mapsto
\Span_k(S_+)$. Therefore the left-adjoint functor
$\Stab:\Ho(\Gamma_+,k) \to \Ho(k)$ is given by
\begin{equation}\label{stab.k}
\Stab(M) = M \lotimes_{\Gamma_+} T^o, \qquad T^o(S_+) =
\Hom_k(T(S_+),k),
\end{equation}
where $\lotimes_{\Gamma_+^o}$ is the derived tensor product
\eqref{lot.I} over the small category $\Gamma_+$. To compute it, one
needs to choose a projective resolution of $T^o$ in the abelian
category $\Fun(\Gamma_+^o,k)$. This is greatly simplified by first
constructing an equivalence
\begin{equation}\label{pdk.eq}
\Fun(\Gamma_+^o,k) \cong \Fun(\Gamma_-^o,k),
\end{equation}
where $\Gamma_-$ is the category of finite sets and surjective maps
between them, and doing it in such a way that $T^o$ is sent to $t
\in \Fun(\Gamma_-^o,k)$ with $t(S)=k$ for $S=\ppt$ and $0$
otherwise. Then the functors $s_n \in \Fun(\Gamma_-^o,k)$
represented by sets of cardinality $n$, $n \geq 0$ give a set of
projective generators of the category $\Fun(\Gamma_-^o,k)$, with
$s_1 = k$ being the constant functor with value $k$. We have an
exact sequence
$$
\begin{CD}
s_2 @>{\eta}>> s_1 @>>> t @>>> 0,
\end{CD}
$$
and it extends to a resolution $P_\idot$ of $t$ given by
\begin{equation}\label{Q.reso}
P_n = s_2^{\otimes n}, \qquad d = \sum_{1 \leq i \leq n}\id^{i-1}
\otimes \eta \otimes \id^{n-i} \text{ on } P_n,
\end{equation}
where one needs to check separately that $s_2^{\otimes n}$ is
projective for every $n$. As an additional bonus, the equivalence
\eqref{pdk.eq} sends the convolution monoidal structure on
$\Fun(\Gamma_+,k)$ to the pointwise monoidal structure on
$\Fun(\Gamma_-^o,k)$, so that the resolution \eqref{Q.reso} is
monoidal, and also describes the monoidal structure on $\Stab$.

\begin{remark}\label{stab.prod.rem}
Another advantage of the resolution \eqref{Q.reso} is that each term
$P_n$ is the sum of a finite number of representable objects; this
implies that additivization commutes with arbitrary products. More
generally, \eqref{stab.k} makes sense for any $M \in
\D(\Gamma_+,k)$, and $\Stab$ commutes with arbitrary products in
$\D^{\leq n}(\Gamma_+,k)$ for any fixed $n$.
\end{remark}

\begin{remark}
Although unknown to the author at the time, all the constructions of
\cite[Section 3]{ka0} were actually discovered much earlier by
T. Pirashvili. The general additivization construction is summarized
nicely in \cite{LP}, while \eqref{pdk.eq} is in \cite{P.dold} (and
should be called ``Pirashvili-Dold-Kan equivalence'').
\end{remark}

\section{Examples.}

\subsection{Basic examples.}\label{add.exa.subs}

To illustrate stabilization and additivization, we take a
commutative ring $k$ and consider functors $F$ from $k\amod^{fl}$ to
itself. The most basic example is $F = \Span_{k} \circ
\phi:k\amod^{fl} \to k\amod^{fl}$, the composition of the forgetful
functor \eqref{phi.eq} and its left-adjoint. By
Remark~\ref{kawa.rem}, $\Stab(F)$ sends a flat $k$-module $M$ to
$Q(M/k) = \Stab(F)(M)$ in $\Ho(k) \cong \Ho^{st}(k)$ that represents
its stable homology with coefficients in $k$ (that is, the homology
of the corresponding Eilenberg-Mac Lane spectrum). Since $F$ is
symmetric monoidal, $\Stab(F)$ is also symmetric monoidal, so that
for any commutative flat $k$-algebra $A$, $Q(A/k)$ is a commutative
ring object in $\Ho(k)$. If $k = \Ff_p$ is a prime field, then the
homology algebra $H_\idot(Q(k/k))$ is the dual Steenrod algebra. We
remind the reader that if $p$ is odd, this is the free
graded-commutative algebra
\begin{equation}\label{st.odd}
H_\idot(Q(k/k)) = k[\beta,\xi_i,\tau_i]
\end{equation}
generated by the Bokstein element $\beta$ of degree $1$, and
elements $\tau_i$, $\xi_i$, $i \geq 1$ of degrees $\deg \xi_i =
2p^i-2$, $\deg \tau_i = 2p^i-1$ (for a proof, see any textbook on
algebraic topology or Subsection~\ref{steen.subs} below). If $p=2$,
then the formula is the same but without $\tau_i$, $\deg \xi_i =
2^{i+1}-1$, and one has to remember that in $\cchar 2$,
``graded-commutative'' means ``commutative''.

If one computes $\Stab(F)$ by \eqref{stab.k} with the resolution
\eqref{Q.reso}, one can lift it to a functor $k\amod^{fl} \to
\Delta^o k\amod \cong C_{\geq 0}(k)$, $M \mapsto Q_\idot(M/k)$. If
$k=\Z$, this functor actually coincides on the nose with the cubical
construction $Q_\idot(M)$ of Eilenberg and Mac Lane. It is lax
monoidal (although not symmetric monoidal, since \eqref{Q.reso} only
has a non-commutative algebra structure), so that for any
$k$-algebra $A$, we have a DG algebra $Q_\idot(A/k)$.

Generalizing this further, for any additive flat $k$-linear category
$A$, one can define a $k$-linear DG category $Q_\idot(A/k)$ with the
same objects, and with morphism complexes given by
$Q_\idot(A/k)(a,a') = Q_\idot(A(a,a')/k)$. Moreover, one can extend
$Q_\idot(-/k)$ to a functor $\Delta^ok\amod^{fl} \to C_{\geq 0}(k)$
by applying it pointwise to a simplicial group, and then totalizing
the bicomplex obtained by the Dold-Kan equivalence from the
resulting object $Q_\idot(M/k) \in \Delta^oC_{\geq 0}(k)$. Then the
extended functor $Q_\idot$ can be applied pointwise to simplicial
objects in $\Fun(A,k)$ pointwise-flat over $k$, and this provides a
functor $Q:\Ho(A,k) \to \D(Q_\idot(A/k))$ to the derived category of
modules over the DG category $Q_\idot(A/k)$. By definition, it
factors through the stabilization functor $\Stab$, and in fact
induces an equivalence
\begin{equation}\label{Q.A.eq}
Q:\Ho^{st}(A,k) \cong \D_{\geq 0}(Q_\idot(A/k)),
\end{equation}
where $\D_{\geq 0}(Q_\idot(A/k)) \subset \D(Q_\idot(A/k))$ is
spanned by modules concentrated in non-negative homological
degrees. To recover the whole triangulated category
$\D(Q_\idot(A/k))$, take the closure of $\Ho^{st}(A,k) \subset
\Ho(A,k) = \D^{\leq 0}(A,k) \subset \D(A,k)$ with respect to
homological shifts.

\medskip

For other examples, consider polynomial functors $k\amod^{fl} \to
k\amod$. The simplest such is the tensor power functor $T_n(M) =
M^{\otimes n}$ for some fixed $n$, but this is not interesting:
\eqref{Sig.st} immediately implies that we have $\Stab(T_n)=0$ as
soon as $n \geq 2$. For an interesting example, fix a prime $p$, and
consider the {\em cyclic power functor} $C$ given by
\begin{equation}\label{C.M.eq}
C(M) = H^0(\Z/p\Z,M^{\otimes_k p}) = \tau^{\leq
  0}C^\hdot(\Z/p\Z,M^{\otimes_k p}),
\end{equation}
where the $\Z/p\Z$-action on $M^{\otimes p}$ is generated by the
longest cycle permutation, and $C^\hdot(\Z/p\Z,-) \in \D(k)$ stands
for the cohomology complex of the group $\Z/p\Z$. Denote $R =
\Stab(C)$.

\begin{lemma}\label{tate.le}
For any $M \in k\amod^{fl}$, we have a functorial isomorphism
\begin{equation}\label{R.M.eq}
R(M) \cong \tau^{\leq 0}\vC^\hdot(\Z/p\Z,M^{\otimes_k p}),
\end{equation}
where $\vC^\hdot(\Z/p\Z,-) \in \D(k)$ stands for the Tate cohomology
complex, and the composition of \eqref{R.M.eq} with the natural
map $C(M) \to R(M)$ is induced by the embedding
$C^\hdot(\Z/p\Z,M^{\otimes_k p}) \to \vC^\hdot(\Z/p\Z,M^{\otimes_k p})$.
\end{lemma}

\proof{} To compute the target of the hypothetical isomorphism
\eqref{R.M.eq}, choose a resolution $P_\idot$ of the trivial
$k[\Z/p\Z]$-module $k$ by free $k[\Z/p\Z]$-modules $P_i = V_i
\otimes_k k[\Z/p\Z]$, $V_i \in k\amod$, let $P'_\idot$ be the
shifted cone of the augmentation map $P_\idot \to k$, and let
$$
R_\idot(M) = H^0(\Z/p\Z,P'_\idot \otimes_k M^{\otimes_k p}).
$$
Then $R_\idot(M) \in C_{\geq 0}(k)$ represents the target of
\eqref{R.M.eq}, and the embedding $k \to P'_\idot$ induces a
functorial map
\begin{equation}\label{C.M}
C(M) \to R_\idot(M).
\end{equation}
For any $i \geq 0$, we have $R_{i+1}(M) \cong V_i \otimes_k
M^{\otimes p}$, and since $\Stab(T_p)=0$, we also have
$\Stab(R_{i+1})=0$, so that the map in $\Ho(k\amod,k)$ represented
by \eqref{C.M} is stably trivial. On the other hand, for any two flat
$k$-modules $M_0,M_1 \in k\amod^{fl}$, the quotient $(M_0 \oplus
M_1)^{\otimes_k p}/(M_0^{\otimes_k p} \oplus M_1^{\otimes_k p})$ is
a free $k[\Z/p\Z]$-module with trivial Tate cohomology, so
that the object in $\Ho(k\amod,k)$ given by $R_\idot(M)$ is
stable. Stabilizing \eqref{C.M}, we get \eqref{R.M.eq}.
\endproof

Up to now, all our example were $k$-linear. However, since
stabilization commutes with the forgetul functor \eqref{phi.eq},
even a non-linear map $F_0 \to F_1$ between functors $F_0,F_1:k\amod
\to k\amod$ induces a map between their stabilizations. For an
example of this, assume that $k$ is a perfect field of
characteristic $p$, and for any $k$-vector space $M \in k\amod$,
consider the map $M \to M^{\otimes p}$ sending $m \in M$ to
$m^{\otimes p}$. This map is $\Z/p\Z$-equivariant and functorial,
thus induces a map
\begin{equation}\label{sp.psi.eq}
\psi:\phi = \phi \circ \Id \to \phi \circ C,
\end{equation}
and this gives rise to a functorial map
\begin{equation}\label{st.psi.eq}
\Stab(\psi):\phi(M) \to \phi(R(M))
\end{equation}
in $\Ho^{st}(\Gamma_+)$ for any $k$-vector space $M \in
k\amod$. Lemma~\ref{tate.le} easily implies that $\Stab(\psi)$ is an
isomorphism on $\pi_0$, but $\phi(R(M))$ also has higher homotopy
groups. In fact, \eqref{R.M.eq} induces a functorial exact triangle
\begin{equation}\label{C.tria}
\begin{CD}
C(M) @>>> R(M) @>{a}>> C_\idot(\Z/p\Z,M^{\otimes p})[1] @>>>
\end{CD}
\end{equation}
in $\D_{\geq 0}(k)$, where $C_\idot(\Z/p\Z,-)$ stands for the group
homology complex, and composing the projection $a$ with
\eqref{st.psi.eq}, we obtain a map
\begin{equation}\label{pr.psi.eq}
\psi' = a \circ \psi:\phi(M) \to \phi(C_\idot(\Z/p\Z,M^{\otimes p})[1]).
\end{equation}
Both the source and the target of $\psi'$ are $k$-linear, the source
is discrete, and the target is $1$-connective, so were $\psi'$ to be
$k$-linear, it would vanish. It does not, even for the
one-dimensional space $M=k$. In principle, the whole map
\eqref{pr.psi.eq} can be computed explicitly in terms of the
Steenrod power operations, see \cite{NS}, but for our purposes, the
following it sufficient.

\begin{lemma}\label{boks.le}
The map $\overline{\psi}:\phi(k) \to \phi(k[1])$ obtained by
composing \eqref{pr.psi.eq} with the projection $C_\idot(\Z/p\Z,k)
\to H_0(\Z/p\Z,k) = k$ is the composition of the Frobenius
endomorphism $k \to k$ and the Bokstein map.
\end{lemma}

\proof{} More generally, for any $M \in k\amod$, let $C'(M) =
H_0(\Z/p\Z.M^{\otimes_k p})$, and note that the correspondence $m
\mapsto m^{\otimes p}$ used in \eqref{sp.psi.eq} also provides a
$k$-linear functorial map $\psi^*:M^{(1)} \to C'(M)$, where
$M^{(1)}$ is the Frobenius twist. If $M=k$, then $\psi^*:k \to k$ is
the Frobenius endomorphism. For any $M$, we can compose $\psi^*$
with the Bokstein map and obtain a functorial map
\begin{equation}\label{bo.eq}
\phi(M) = \phi(M^{(1)}) \to \phi(C'(M)[1])
\end{equation}
of connective spectra. On the other hand, composing
\eqref{pr.psi.eq} with the truncation map
$C_\idot(\Z/p\Z,M^{\otimes_k p}) \to C'(M)$ also gives a functorial
map $\overline{\psi}:\phi(M) \to \phi(C'(M)[1])$, and it suffices to
prove that it coincides with \eqref{bo.eq}.

To do this, note that, unlike the Bokstein map itself, the map
\eqref{bo.eq} admits a functorial cone given by the second {\em
  polynomial Witt vectors functor} $W_2$ of \cite{witt}. By
definition, this is a functor from $k$-vector spaces to modules over
the second Witt vectors ring $W_2(k)$ that fits into a functorial
short exact sequence
\begin{equation}\label{w2.seq}
\begin{CD}
0 @>>> C'(M) @>{V}>> W_2(M) @>{R}>> M^{(1)} @>>> 0.
\end{CD}
\end{equation}
If we apply the forgetful functor $\D_{\geq 0}(W_2(k)) \to
\Ho^{st}(\Gamma)$, this sequences induces an exact triangle of
spectra whose connecting differential is
\eqref{bo.eq}. Equivalently, \eqref{bo.eq} is induced by the
composition
\begin{equation}\label{Z.M.eq}
\begin{CD}
M^{(1)} @>{\overline{R}^{-1}}>> \overline{Z}(M) @>{\overline{a}}>>
C'(M)[1]
\end{CD}
\end{equation}
in $\D_{\geq 0}(W_2(k))$, where $\overline{Z}(M)$ is the cone of the
map $V$, the isomorphism $\overline{R}:\overline{Z}(M) \to M$ is
induced by $R$, and $\overline{a}$ is the natural
projection. Moreover, $W_2$ also fits into a functorial short exact
sequence
\begin{equation}\label{w2.2.seq}
\begin{CD}
0 @>>> M^{(1)} @>{C}>> W_2(M) @>{F}>> C(M) @>>> 0,
\end{CD}
\end{equation}
and we have a functorial {\em Teichm\"uller map} $T:\phi(M) \to
\phi(W_2(M))$ such that $\phi(R) \circ T = \id$ and $\phi(F) \circ
T:M \to C(M)$ is the map \eqref{sp.psi.eq}.

Now let $Z = \Stab(W_2)$ be the stabilization of the functor $W_2$,
and note that by \eqref{w2.2.seq} and the same argument as in
Lemma~\ref{tate.le}, it fits into a functorial exact triangle
\begin{equation}\label{w2.tria}
\begin{CD}
W_2(M) @>>> Z(M) @>{a}>> C_\idot(\Z/p\Z,M^{\otimes_k p})[1]
@>>>
\end{CD}
\end{equation}
in $\D_{\geq 0}(W_2(k))$ whose connecting differential is the
composition of the truncation map $C_\idot(\Z/p\Z,M^{\otimes_k })
\to C'(M)$ and the map $V$ of \eqref{w2.seq}. In particular, the
truncation $\tau_{\leq 1}R^2(M)$ is exactly $\overline{Z}(M)$ of
\eqref{Z.M.eq}, and $\tau_{\leq 1}(a)$ is the map $\overline{a}$.
Moreover, since $\phi(R) \circ T \cong \id$, the stablization
$\Stab(T)$ of the Teichm\"uller map induces an isomorphism $\phi(M)
\to \phi(\overline{Z}(M))$ inverse to $\phi(\overline{R})$. But on
the other hand, we have $\psi = \phi(F) \circ T$, so that
$\Stab(\psi) = \phi(\Stab(F)) \circ \Stab(T)$, and the triangles
\eqref{C.tria}, \eqref{w2.tria} induce a commutative diagram
$$
\begin{CD}
Z(M) @>>> \overline{Z}(M) @>{\overline{a}}>> C'(M)[1]\\
@V{\Stab(F)}VV @V{\tau_{\leq 1}\Stab(F)}VV @|\\
R(M) @>>> \tau_{\leq 1}R(M) @>{\tau_{\leq 1}a}>> C'(M)[1].
\end{CD}
$$
Combined with \eqref{Z.M.eq} and the definiton of $\overline{\psi}$,
this proves the claim.
\endproof

\subsection{Generalized Tate cohomology.}

It turns out that can one connect the two main examples of
addivization of Subsection~\ref{add.exa.subs} by computing
addivization of symmetric power functors and their divided power
counterparts. If $n=p$ is a prime, this can be done by exactly the
same argument as for the cyclic power functor \eqref{C.M.eq} of
Lemma~\ref{tate.le}. For other integers, we first need an
appropriate generalization of Tate cohomology.

Fix a commutative ring $k$. We will say that an {\em admissible
  family of subgroups} in a finite group $G$ is a collection $\X$ of
subgroups $H \subset G$ that is closed under conjugation and
intersections, and does not contain $G$ itself. We will also say
that a morphism $M \to N$ between $k[G]$-modules is {\em
  $\X$-surjective} if $M^H \to N^H$ is surjective for any $H \in
\X$, and a $k[G]$-module $P$ is {\em $\X$-projective} if $\Hom(P,-)$
sends $\X$-surjective maps to surjective maps. For example, for any
$H \in \X$, the induced $k[G]$-module $k[G/H] = k \otimes_{k[H]}
k[G]$ is $\X$-projective by adjunction, and for any $k[G]$-module
$M$ projective over $k$, the product $M \otimes_k k[G/H]$ is
$\X$-projective by the projection formula. A complex $M_\idot$ of
$k[G]$-modules is {\em $\X$-exact} if $M_\idot^H$ is acyclic for any
$H \in \X$. An {\em $\X$-resolution} of a $k[G]$-module $M$ is a
complex $P_\idot$ of $\X$-projective $k[G]$-modules, trivial in
negative homological degrees and equipped with a map $a:P_\idot \to
M$ whose cone $\wt{P}_\idot$ is $\X$-exact. By exactly the same
argument as in the usual case, a map between two $k[G]$-modules
equipped with $\X$-resolutions lifts to a map between resolutions,
and the liting is unique up to a chain homotopy. To construct an
$\X$-resolution of the trivial $k[G]$-module $k$, one can take an
$\X$-surjective cover $\eta:P_0 \to k$ --- for example, by setting
$P_0 = \bigoplus_{H \in \X}k[G/H]$ --- and then continue by taking
$P_i = P_0^{\otimes_k i+1}$ with the differential as in
\eqref{Q.reso}. Note that for this particular resolution, the cone
$\wt{P}_\idot$ is naturally a DG algebra.

\begin{defn}\label{ttate.def}
The {\em truncated Tate cohomology} $\bH_\idot(G,\X,M_\idot)$ of the
group $G$ with respect to the admissible family $\X$ and with
coefficients in a bounded complex $M_\idot$ of $k[G]$-modules is the
homology of the complex
\begin{equation}\label{bC.G}
\bC_\idot(G,\X,M_\idot) = (\wt{P}_\idot \otimes_k M_\idot)^G,
\end{equation}
where $\wt{P}_\idot$ is the cone of a $\X$-resolution $P_\idot \to
k$.
\end{defn}

By the lifting property of $\X$-resolutions, truncated Tate homology
is well-defined and independent of the choice of the $\X$-resolution
$P_\idot$, and the complex \eqref{bC.G} is well-defined as an object
in the derived category $\D(k)$. In particular, we can choose
$P_\idot$ in such a way that $\wt{P}_\idot$ is a DG algebra, and
this turns $\bH_\idot(G,\X,-)$ into a lax monoidal functor. It is
also somewhat cohomological, in that a short exact
sequence of complexes 
$$
\begin{CD}
0 @>>> M_\idot' @>>> M_\idot @>{a}>> M''_\idot @>>> 0
\end{CD}
$$
with termwise $\X$-surjective $a$ gives rise to a long exact
sequence of the truncated Tate cohomology groups. However, if $a$ is
simply surjective, this need not be true. In particular,
$\bH_\idot(G,\X,-)$ does not preserve quasiisomorphisms, and can be
non-trivial even if the complex $M_\idot$ is acyclic.

\begin{exa}
If $\X = \{e\}$ consists of the trivial subgroup $\{e\} \subset G$,
then $\X$-surjective is surjective, $\X$-exact is exact,
$\X$-projective is projective, and for any $k[G]$-module $M$, we
have
$$
\bH_\idot(G,\{e\},M) \cong \tau^{\leq 0}\vH^\hdot(G,M),
$$
where as in Lemma~\ref{tate.le}, the right-hand side is the
truncation of the usual Tate cohomology groups.
\end{exa}

\begin{exa}\label{norma.exa}
If we have a surjective map of groups $f:G \to W$ and an admissible
family $\X$ of subgroups in $W$, then the family $f^{-1}\X =
\{f^{-1}(H)|H \in \X\}$ is an admissible family of subgroups in
$G$. Any $k[W]$-module $V$ defines a $k[G]$-mo\-dule $f_*V$ by
restriction of scalars, and for any $\X$-resolution $P_\idot$ of
$k$, $f_*P_\idot$ is a $f^{-1}\X$-resolution, so that
$\bH_\idot(G,f^{-1}\X,f_*V) \cong \bH_\idot(W,\X,V)$.
\end{exa}

\begin{remark}
If one replaces $G$-invariants in \eqref{bC.G} with the full
cohomology complexes $C^\hdot(G,-)$, and takes the sum-total complex
of the resulting bicomplex, one arrives at the well-known notion of
{\em generalized Tate cohomology} $\vH^\hdot(G,\X,-)$ (see
e.g.\ \cite[Section 7]{proma}). Alternatively, we have
$$
\vH^\hdot(G,\X,M_\idot) \cong
\RHom^\hdot_{\D^b(k[G])/\D^b_\X(k[G])}(k,M_\idot),
$$
where $\D^b_\X(k[G]) \subset \D^b(k[G])$ is the Karoubi closure of
the full triangulated subcategory in $\D^b(k[G])$ generated by
$\X$-projective modules. We always have a map
$\bH_\idot(G,\X,M_\idot) \to \vH^\hdot(G,\X,M_\idot)$ but it need
not be injective even when $M_\idot$ is $k$ in degree $0$.  For
example, if one takes $k=\Ff_p$ and $G=\Z/p^2\Z$, with $\X$
consisting of $\Z/p\Z \subset \Z/p^2\Z$, then
$\bH_\idot(G,\X,k) \cong \bH_\idot(\Z/p\Z,k)$ by
Example~\ref{norma.exa}, while $\vH^\hdot(G,\X,k)=0$.
\end{remark}

To study truncated Tate cohomology, it is convenient to do the
following. Let $\Gamma_G$ be the category of finite $G$-sets ---
that is, finite sets $S$ equipped with an action of $G$ --- and let
$O_G \subset \Gamma_G$ be the full subcategory spanned by {\em
  $G$-orbits}, that is, $G$-sets $S$ such that the action is
transitive. For any subgroup $H \subset G$, the quotient $G/H$ is a
$G$-orbit that we denote by $[G/H] \in O_G \subset \Gamma_G$, and
all $G$-orbits are of this form. The category $\Gamma_G$ has finite
products and finite coproducts. Say that a functor $E \in
\Fun(\Gamma_G,k)$ is {\em additive} if for any $S,S' \in \Gamma_S$,
the natural map $E(S) \oplus E(S') \to E(S \copr S')$ is an
isomorphism, and let $\Fun_{add}(\Gamma_G,k) \subset
\Fun(\Gamma_G,k)$ be the full subcategory spanned by additive
functors. Then we have a natural equivalence
\begin{equation}\label{O.add}
\Fun(O_G,k) \cong \Fun_{add}(\Gamma_G,k)
\end{equation}
given by the pullback $\eps^*$ with respect to the embedding
$\eps:O_G \to \Gamma_G$, with the inverse equivalence given by the
left Kan extension $\eps_!$.

For any subgroup $H \subset G$, the left comma-fiber $\Gamma_G /
[G/H]$ is naturally identified with $\Gamma_H$, and the forgetful
functor $\psi^H:\Gamma_H \cong \Gamma_G / [G/H] \to \Gamma_G$ has a
right-adjoint $\mu_H:\Gamma_G \to \Gamma_H$ sending $S \in \Gamma_G$
to $S \times [G/H]$ with its natural projection to $[G/H]$. Then by
adjunction, we have
\begin{equation}\label{mu.psi}
\mu_H^* \cong \psi^H_!,
\end{equation}
so that $\psi^H_!$ is exact. The functor $\mu_H$ commutes with
finite coproducts, so that \eqref{mu.psi} sends
$\Fun_{add}(\Gamma_H,k) \subset \Fun(\Gamma_H,k)$ into
$\Fun_{add}(\Gamma_G,k) \subset \Fun(\Gamma_G,k)$. The functor
$\psi^H$ sends orbits to orbits, thus restricts to a functor
$\psi^H:O_H \to O_G$, and this is compatible with the equivalences
$\eps_!$ of \eqref{O.add}, so that the left Kan extension
$\psi^H_!:\Fun(O_H,k) \to \Fun(O_G,k)$ is also exact.

If $G$ is equipped with an admissible family of subgroups $\X$, then
we can consider the full subcategory $O_\X \subset O_G$ spanned by
orbits $[G/H]$ with $H$ in $\X$. Since by assumption, $\X$ does not
contain $G$ itself, we can extend the embedding $O_\X \to O_G$ to a
full embedding $O_\X^> \to O_G$ by sending the terminal object $o
\in O_\X^>$ to the one-point $G$-orbit $[G/G]$, and we let
$\Gamma_\X \subset \Gamma_G$ be the full subcategory spanned by
finite coproducts of orbits in $O^>_\X \subset O_G$. Then $\Gamma_\X
\subset \Gamma_G$ is closed under finite products. We still have the
notion of an additive functor in $\Fun(\Gamma_\X,k)$ and the
equivalence \eqref{O.add}, and for any $H$ in $\X$, we still have
the functor $\psi^H:O_H \to O_\X \subset O_\X^>$ such that the left
Kan extension $\psi^H_!:\Fun(O_H,k) \to \Fun(O_\X^>,k)$ is exact and
given by \eqref{mu.psi}.

Now, let $\sigma:\Gamma^o_G \to k[G]\amod$ be the tautological
functor sending a $G$-set $S$ to the space $k(S)$ of $k$-valued
functions on $S$. Restricting $\sigma$ to $\Gamma_\X^o \subset
\Gamma^o_G$ and taking its left Kan extension with respect to the
Yoneda embedding \eqref{yo.k} gives a right-exact functor
$\Y_!(\sigma):\Fun(\Gamma_\X,k) \to k[G]\amod$ that has a
right-adjoint $\Loc_\X:k[G]\amod \to \Fun(\Gamma_\X,k)$. Explicitly,
$\Loc_\X$ is given by
\begin{equation}\label{loc.G}
\Loc_\X(M)(S) = (M \otimes_k k[S])^G, \qquad S \in \Gamma_\X,M \in
k[G]\amod,
\end{equation}
so that in particular, it takes values in the subcategory
$\Fun_{add}(\Gamma_\X,k)$. We extend $\Loc_\X$ to bounded complexes
by applying it pointwise, so that for any bounded complex $M_\idot
\in C^b_\idot(k[G])$, we have a well-defined object
$\Loc_\X(M_\idot)$ in the bounded derived category $\D^b(O_\X^>,k)$.

\begin{lemma}\label{O.le}
For any bounded complex $M_\idot$ of $k[G]$-modules, we have
$$
\bH_\idot(G,\X,M_\idot) \cong H_\idot(O_\X^>,\{o\},\Loc_\X(M_\idot))
$$
where $H_\idot(-)$ is the homology with support of \eqref{C.z}.
\end{lemma}

\proof{} Let $\sigma^*:\Gamma_G \to k[G]\amod$ be the functor
sending $S \in \Gamma_G$ to the free $k$-module $k[S]$ (that is, to
the dual $k$-module to $k(S)$), with the left Kan extension $L =
\Y_!(\sigma^*):\Fun(O_\X^{>o},k) \to k[G]\amod$, and note that for
any $[G/H] \in O_\X^{<o}$ with the representable functor $k_{[G/H]}
\in \Fun(O_G^{>o},k)$, we have
\begin{equation}\label{Y.P}
L(k_{[G/H]}) \cong \sigma^*([G/H]) \cong k[G/H].
\end{equation}
Then \eqref{loc.G} provides a functorial identification
\begin{equation}\label{tate.sup}
L(k_{[G/H]}) \otimes_{O_G^>} \Loc(M_\idot) \cong (k[G/H] \otimes_k
M_\idot)^G
\end{equation}
To compute the homology with support, we need to choose a projective
resolution $P_\idot$ of the functor $j_{1*}^ok$ in \eqref{C.z.L},
and we can do it by taking $P_0 = k = k_{[G/G]}$, and then letting
all the $P_i$, $i \geq 1$ be the appropriate sums of representable
functors $k_{[G/H]}$, $[G/H] \in O^o_\X \subset O_\X^{>o}$. By
\eqref{Y.P}, we then have $L(P_0) \cong k$, and for any $i \geq 1$,
$L(P_i) \in k[G]\amod$ is $\X$-projective. Then by \eqref{tate.sup},
to finish the proof, it remains to show that $\Y(\sigma^*)(P_\idot)$
is $\X$-exact, and this amounts to checking that
$H^{\{o\}}(O_\X^>,\Loc_\X(k[G/H])) = 0$ for any $H \in \X$.
However, by \eqref{loc.G} and \eqref{mu.psi}, we have
$\Loc_\X(k[G/H]) \cong \psi^H_!\Loc_H(k[G/H])$, where $\Loc_H$ is
the functor \eqref{loc.G} for the category $\Gamma_H$, and by
\eqref{yo.cov}, it then suffices to observe that
$\psi^{Ho*}k_{\{o\}} = 0$.  \endproof

\subsection{Divided powers.}\label{div.subs}

Recall that we denote by $\Gamma_-$ the category of non-empty finite
sets and surjective maps. For any integer $n \geq 1$, let $[n] \in
\Gamma_-$ be the set of cardinality $n$, and let $\Sigma_n =
\Aut([n])$ be the corresponding permutation group. For any flat
$k$-module $V$, the $n$-th {\em symmetric power} $S^n(V)$ and {\em
  divided power} $D^n(V)$ are given by
\begin{equation}\label{S.D}
S^n(V) = T^n(V)_{\Sigma_n}, \qquad D^n(V) = T^n(V)^{\Sigma_n},
\end{equation}
where $T^n(V) = V^{\otimes_k n}$ is the tensor power functor. For
any $n,m \geq 1$, the isomorphism $T^m(T^n(V)) \cong T^{mn}(V)$
induces functorial maps
\begin{equation}\label{SD.mn}
S^m(S^n(V)) \to S^{mn}(V), \qquad D^{mn}(V) \to D^m(D^n(V)).
\end{equation}
For any map $f:[n] \to [m]$ in $\Gamma_-$, let $\Sigma_f \subset
\Sigma_n$ be the subgroup of automorphisms $a:[n] \to [n]$ such that
$f \circ a = f$. Explicitly, $f$ defines a partition $[n] = [n_1]
\copr \dots \copr [n_m]$ of $[n]$ into $[m]$ disjoint subsets $[n_i]
= f^{-1}(i)$, $i \in [m]$, and $\Sigma_f = \Sigma_{[n_1]} \times
\dots \times \Sigma_{[n_m]}$ consists of permutations that preserve
each $[n_i]$. For any $[n]$, let $\X_n$ be the family of subgroups
$\Sigma_f \subset \Sigma_n$ for all $f:[n] \to [m]$ in $\Gamma_n$
with $m \geq 2$. Then the family $\X_n$ is admissible. The augmented
orbit category $O_{\X_n}^>$ is naturally identified with the
category $\Gamma_n$ whose objects are maps $f:[n] \to [m]$ in
$\Gamma_-$, and whose maps from $f':[n] \to [m']$ to $f:[n] \to [m]$
are maps $g:[m'] \to [m]$ such that $g \circ f'$ factors through $f$
(but we do not specify a factorization, so that $\Gamma_n(f',f)
\cong \Ar(\Gamma_-)(f',f)/\Sigma_f$).

\begin{remark}
The standard notation for the divided power functor is $\Gamma^n$
and not $D^n$; we allow ourselves to change it to avoid confusion
with the notation for the category of finite sets.
\end{remark}

For any $n$, the divided power functor $D^n$ of \eqref{S.D} is a
functor from $k\amod^{fl}$ to $k\amod$, and it is lax monoidal by
adjunction, thus defines a lax monoidal functor
\begin{equation}\label{Q.n}
Q^n_\idot = \Stab(D^n):k\amod^{fl} \to C_{\geq 0}(k),
\end{equation}
where as in Subsection~\ref{add.exa.subs}, we can rigidify things by
computing $\Stab$ via \eqref{stab.k} with the resolution
\eqref{Q.reso}. For any $V \in k\amod^{fl}$, we denote by
$HQ^n_\idot(V)$ the homology of the complex
$Q_\idot^n(V)$. If we take $V=k$, then $Q^n_\idot(k)$ is a
DG algebra that defines a commutative ring object in $\Ho(k) \cong
\D_{\geq 0}(k)$, and $HQ_\idot(k)$ is a graded-commutative algebra
over $k$.

\begin{lemma}\label{D.le}
For any $n \geq 1$, we have a natural isomorphism
\begin{equation}\label{D.stab}
Q^n_\idot(V) \cong \bC_\idot(\Sigma_n,\X_n,T^n(V))
\end{equation}
of commutative ring objects in $\D_{\geq 0}(k)$ whose target is as
in \eqref{bC.G}.
\end{lemma}

\proof{} Choose an $\X_n$-resolution $P_\idot$ of the trivial
$k[\Sigma_n]$-module $k$ such that each $P_i$ is a sum of modules
$k[\Sigma_n/\Sigma_f]$, $\Sigma_f \in \X_n$, and then proceed
exactly as in Lemma~\ref{tate.le}. On one hand, for each $f:[n] \to
[m]$, we have
\begin{equation}\label{T.D}
(T^n(V) \otimes_k k[\Sigma_n/\Sigma_f])^{\Sigma_n} \cong
T^n(V)^{\Sigma_f} \cong \bigotimes_{1 \leq i \leq m}D^{n_i}(V),
\end{equation}
and since $m \geq 2$, this has trivial stabilization by
\eqref{Sig.st}. On the other hand, the quotient $T^n(V_0 \oplus
V_1)/(T^n(V_0) \oplus T^n(V_1))$ is $\X_n$-projective, so that the
right-hand side of \eqref{D.stab} is already stable.
\endproof

With a little bit of extra effort, one can show that \eqref{D.stab}
lifts to a quasiisomorphism between DG algebras, but we will not
need this: our main interest lies in the homology algebra
$HQ_\idot^n(k)$. To compute it, recall that Lemma~\ref{O.le}
provides an isomorphism
\begin{equation}\label{Si.Ho}
\bH_\idot(\Sigma_n,\X_n,T^n(V)) \cong
H_\idot(\Gamma_n,\{o\},L_n(V)),
\end{equation}
where $\{o\} \in \Gamma_n$ is the terminal object $[n] \to [1]$, and
we simplify notation by writing $L_n = \Loc_{\X_n} \circ T^n$. Note
that for any $p \geq 1$, the cartesian product of finite sets
induces a product functor
$$
\Gamma_n \times \Gamma_p \to \Gamma_{np},
$$
and in particular, taking the product with the terminal object $[p]
\to [1]$ gives a functor $\eps_p:\Gamma_n \to \Gamma_{pn}$. This is
a fully faithful embedding whose essential image consists of
partitions $[pn] = [d_1] \copr \dots \copr [d_l]$ such that every
$d_i$, $1 \leq i \leq l$ is divisible by $p$. It is also
right-closed in the sense of Example~\ref{cl.exa}. For any bounded
complex $M_\idot$ of flat $k$-modules, \eqref{loc.G}, \eqref{SD.mn}
and \eqref{T.D} provide a natural map
\begin{equation}\label{p.eq}
\eps_p^*L_{np}(M_\idot) \to L_n(D^p(M_\idot))
\end{equation}
that is an isomorphism if $M=k$ or $M=k[1]$.

Now let $I_\idot$ be the acyclic length-$2$ complex $\id:k \to k$
placed in homological degrees $0$ and $1$, so that fits into a short
exact sequence
\begin{equation}\label{I.seq}
\begin{CD}
0 @>>> k @>{b}>> I_\idot @>{a}>> k[1] @>>> 0,
\end{CD}
\end{equation}
and let $I_\idot^* = I_\idot[-1]$ be the same complex $\id:k \to k$
placed in degrees $0$ and $-1$. Moreover, assume from now on that
$k$ is annihilated by a prime $p$.

\begin{lemma}\label{D.I.le}
The complex $D^m(I_\idot)$ is acyclic unless $m = np$ is divisible
by $p$, and in the latter case, the map $\phi_n:D^{np}(I_\idot) \to
D^n(D^p(I_\idot))$ induced by \eqref{SD.mn} is a
quasiisomorphism. Moreover, if $p$ is odd, then this map is an
isomorphism, and $D^m(I^*_\idot)$ is acyclic for any $m$.
\end{lemma}

\proof{} By duality, we may replace $D^m$ with $S^m$, replace
$\phi_n$ with the dual map $\phi_n^*$, and swap $I_\idot$ and
$I^*_\idot$.  If $p$ is odd, then the symmetric algebra
$S^\hdot(I_\idot)$ resp.\ $S^\hdot(I^*_\idot)$ is the free
graded-commutative algebra $k[t,\xi]$ with $\deg t = 0$, $\deg \xi =
1$ resp. $-1$, and the differential $d\xi = t$ resp.\ $dt=\xi$. The
claim is then obvious. If $p=2$, then we still have
$S^\hdot(I_\idot^*) \cong k[t,\xi]$, $\deg \xi=-1$, $dt=\xi$, but
``graded-commutative'' now means commutative, so that $\phi_n^*$ is
not an isomorphism anymore. However, it is trivial to check that
$k[t,\xi]$ is acyclic in degrees other than $0$, and its homology in
degree $0$ is $k[t^2]$.
\endproof

\begin{remark}\label{D.rem}
If $m$ is not divisible by $p$, then Lemma~\ref{D.I.le} provides a
map $k \to D^m(I_\idot)_1 \subset T^m(I_\idot)_1$ that splits the
differential $D^m(I_\idot)_1 \to D^m(I_\idot)_0 \cong D^m(k) \cong
k$, so that $T^m(I_\idot)_1 \cong k \oplus M$ for some
$\Sigma_m$-module $M$ flat over $k$. Then $T^m(I_\idot) \cong
I_\idot \otimes_k S^\hdot(M[1])$ is $\Sigma_m$-equivariantly
contractible, so that $D^m(V \otimes I_\idot)$ is acyclic for any
flat $k$-module $V$.
\end{remark}

Since the embedding $\eps_p:\Gamma_n \to \Gamma_{pn}$ is
right-closed for any $n \geq 1$, the left Kan extension $\eps_{p!}$
is given by extension by $0$. Then by \eqref{loc.G} and \eqref{T.D},
Lemma~\ref{D.I.le} immediately implies that the adjunction map
$\eps_{p!}\eps_p^*L_{np}(I_\idot) \to L_{np}(I_\idot)$ is a
quasiisomorphism. Moreover, if $p$ is odd, then \eqref{p.eq} is an
isomorphism for $M_\idot=I_\idot$. The complex $\bI_\idot =
D^p(I_\idot)$ is isomorphic to $k \oplus k[1]$, and the isomorphism
inverse to \eqref{p.eq} induces by adjunction a quasiisomorphism
\begin{equation}\label{eps.L.p}
\eps_{p!}L_n(\bI_\idot) \to L_{np}(I_\idot)
\end{equation}
that gives rise to an isomorphism
\begin{equation}\label{eps.p}
HQ^n_\idot(\bI_\idot)) \to HQ^{np}_\idot(I_\idot)
\end{equation}
in $\D_{\geq 0}(k)$. If $p=2$, then $D^p(I_\idot)$ is
quasiisomorphic to $k$ placed in degree $0$, so we may let
$\bI_\idot=k$, construct a quasiisomorphism \eqref{eps.L.p} by
composing the adjunction map $\eps_{p!}L_n(k) \to L_{pn}(k)$ with
the map induced by the embedding $b:k \to I_\idot$ of \eqref{I.seq},
and still obtain an isomorphism \eqref{eps.p}. Moreover, if we turn
$I_\idot$ and $\bI_\idot$ into DG algebras with trivial
multiplication on $\bI_1$ and $I_1$, so that $b$ in \eqref{I.seq} is
a DG algebra map, then \eqref{eps.p} is an algebra isomorphism.

\begin{lemma}\label{p.corr}
For any commutative ring $k$ annhilated by an odd prime $p$, and any
$n \geq 1$, we have an isomorphism of graded-commutative algebras
\begin{equation}\label{p.corr.eq}
HQ_\idot^{np}(k) \cong HQ_\idot^n(k)[\xi,\tau],
\end{equation}
where the generators $\xi$, $\tau$ have degrees $\deg \xi =
2(pn-1)$, $\deg \tau = 2n-1$. If $p=2$, we have the same isomorphism
but without $\tau$ and with $\deg \xi = 2n-1$.
\end{lemma}

\proof{} For any $m \geq 1$, the sequence \eqref{I.seq} induces a
sequence
\begin{equation}\label{I.m.seq}
\begin{CD}
0 @>>> T^m(k) @>{b}>> T^m(I_\idot) @>{a}>> T^m(k[1]) @>>> 0
\end{CD}
\end{equation}
of complexes of $k[\Sigma_m]$-modules. This sequence is
termwise-split, and while it is not exact in the middle term, the
homology there is $\X_m$-projective in each degree. Therefore
\eqref{I.m.seq} gives rise to an exact triangle
\begin{equation}\label{I.m.tria}
\begin{CD}
Q_\idot^m(k) @>{\beta}>> Q_\idot^m(I_\idot) @>{\alpha}>>
Q_\idot^m(k[1]) @>{\delta}>>
\end{CD}
\end{equation}
in $D_{\geq 0}(k)$. All maps here are maps of
$Q_\idot^m(k)$-modules, and $\beta$ is also an algebra map. If
$p=2$, we have $T^m(k[1]) \cong T^m(k)[m]$, so that $Q_\idot^m(k[1])
\cong Q^m(k)[m]$. If $p$ is odd, we can shift the sequence
\eqref{I.seq} and repeat the argument to obtain an exact triangle
\begin{equation}\label{I.m.st.tria}
\begin{CD}
Q_\idot^m(k[-1]) @>>> Q_\idot^m(I^*_\idot) @>>>
Q_\idot^m(k) @>>>
\end{CD}
\end{equation}
in $D_{\geq 0}(k)$, where $L^m(I_\idot^*)$, hence also
$Q^m_\idot(I^*_\idot)$ vanishes by \eqref{loc.G}, \eqref{T.D} and
Lemma~\ref{D.I.le}. Therefore the connecting differential in the
triangle \eqref{I.m.st.tria} is an isomorphism, and since $T^m(k[1])
\cong T^m(k[-1])[2m]$, we have $Q^m_\idot(k[1]) \cong
Q^m_\idot(k[-1])[2m] \cong Q^m_\idot(k)[2m-1]$. If we take $m=pn$,
then we have the identification \eqref{eps.p} induced by the
quasiisomorphism \eqref{eps.L.p}, and in terms of this
identification, the map $\alpha$ in \eqref{I.m.tria} is induced by a
map $\eps_{p!}L^n(\bI_\idot) \to L^{pn}(k[1])$ adjoint to a map
$\alpha':L^n(\bI_\idot) \to \eps_p^*L_{pn}(k[1]) \cong
L_n(D^p(k[1]))$. If $p$ is odd, then $D^p(k[1])=0$, and if $p=2$ and
$D^p(k[1]) \cong k[p]$, we have $\alpha' = L_n(D^p(a \circ b))=0$
since $a \circ b = 0$. Thus in any case, $\alpha'$ vanishes, so does
$\alpha$, and the triangle \eqref{I.m.tria} splits after one
rotation. Altogether, we obtain a short exact sequence
$$
0 \longrightarrow HQ^{np}_{\idot-d}(k)
\overset{\delta}{\longrightarrow} HQ^{np}_\idot(k)
\overset{\beta}{\longrightarrow} HQ^n_\idot(\overline{I}_\idot)
\longrightarrow 0
$$
of graded $HQ^{np}_\idot(k)$-modules, where $d=2m-2=2(np-1)$ if $p$
is odd and $d=np-1=2n-1$ if $p=2$. Moreover, $\beta$ is an algebra
map. Since $\delta$ is a module map, it must be given by
multiplication by the element $\xi=\delta(1)$ of degree $d$, and if
$p=2$, $\bI_\idot=k$ and we are done. If $p$ is odd, then to finish
the proof, it remains to construct a commutative algebra isomorphism
$HQ^n_\idot(\overline{I}_\idot) \cong HQ^n_\idot(k)[\tau]$, where
$\tau$ is a generator of degree $2n-1$. To do this, note that
$\overline{I}_\idot$ also fits into an exact sequence of the form
\eqref{I.seq} that is moreover split, and then the triangle
\eqref{I.m.seq} for $m=n$ induced by this sequence is also split.
\endproof

\subsection{The dual Steenrod algebra.}\label{steen.subs}

Now assume that $k=\Ff_p$ is a prime field. Then for any $k$-vector
space $V$, we have a natural $k$-linear map $V \to S^p(V)$ sending
$v \in V$ to $v^{\otimes p}$, and for any $n \geq 1$, this can be
combined with \eqref{SD.mn} to give a functorial map
\begin{equation}\label{S.fr}
\begin{CD}
S^n(V) @>>> S^p(S^n(V)) @>>> S^{pn}(V).
\end{CD}
\end{equation}
If $V$ is finite-dimensional, we have the dual map
\begin{equation}\label{D.fr}
\begin{CD}
D^{pn}(V) @>>> D^p(D^n(V)) @>>> D^n(V),
\end{CD}
\end{equation}
and we extend it to all vector spaces by taking filtered
colimits. Being functorial, the map
\eqref{D.fr} induces maps
\begin{equation}\label{Q.fr}
Q_\idot^{pn}(V) \to Q_\idot^n(V), \qquad HQ^{pn}_\idot(V) \to
HQ^n_\idot(V)
\end{equation}
of the stabilizations \eqref{Q.n} of functors \eqref{S.D} and their
homology modules. Note that the map \eqref{D.fr} is compatible with
the lax monoidal structures, so that \eqref{Q.fr} is also a lax
monidal map.

\begin{lemma}
In terms of the isomorphism \eqref{p.corr.eq}, the map \eqref{Q.fr}
is obtained by sending the generators $\tau$ and $\xi$ to $0$.
\end{lemma}

\proof{} Let $I_\idot$ and $b$ be as in \eqref{I.seq}, and note that
for any $k$-vector space $V$, the map \eqref{D.fr} for $n=1$ factors
as
\begin{equation}\label{d.p.d}
\begin{CD}
D^p(V) @>{D^p(b)}>> D^p(V \otimes I_\idot) @>{d}>> V
\end{CD}
\end{equation}
for a certain functorial map $d$. Therefore we have a functorial
diagram
\begin{equation}\label{fr.dia}
\begin{CD}
L_{pn}(V) @>{L_{pn}(b)}>> L_{pn}(V \otimes I_\idot) @<{e}<<
\eps_{p!}\eps_p^*L_{pn}(V \otimes I_\idot)\\
@. @. @VVV\\
@. \eps_{p!}L_n(V) @<{L_n(d)}<<  \eps_{p!}L_n(D^p(V \otimes I_\idot))
\end{CD}
\end{equation}
of complexes in $\Fun(\Gamma_{pn},k)$, where the vertical map is
\eqref{p.eq}, and $e$ is the adjunction map. By Remark~\ref{D.rem},
$e$ is a quasiisomorphism, thus invertible in $\D_{\geq
  0}(\Gamma_{pn},k)$, so that \eqref{fr.dia} defines a map
$L_{pn}(V) \to \eps_{p!}(V)$ in $\D_{\geq 0}(\Gamma_{pn},k)$. After
evaluation at $\{o\} \in \Gamma_{pn}$, this map becomes
\eqref{D.fr}. Being functorial, the diagram \eqref{fr.dia} also
defines a diagram in $\Ho(k\amod \times \Gamma_{pn},k) \cong
\D_{\geq 0}(k\amod \times \Gamma_{pn},k)$, where $e$ is again
invertible. Then if we denote by $z:k\amod \to k\amod \times
\Gamma_{pn}$ the right-closed embedding onto $k\amod \times \{o\}$,
and apply the functor $L^\hdot z^!$ of \eqref{ga.shriek} together
with the identifications \eqref{Si.Ho} and \eqref{D.stab}, we obtain
a map $Q^{pn} \to Q^n$ in $\Ho(k\amod,k)$ that fits into a
commutative diagram
$$
\begin{CD}
D^{pn} @>>> D^n\\
@VVV @VVV\\
Q^{pn} @>>> Q^n,
\end{CD}
$$
where the top arrow is \eqref{D.fr}, and the vertical arrows are
stabilization maps. By the universal property of stabilization, the
bottom arrow then must coincide with \eqref{Q.fr}. To finish the
proof, it remains to evaluate at $V = k$, and compare \eqref{fr.dia}
with the construction of the isomorphism \eqref{p.corr.eq} given in
Corollary~\ref{p.corr}.
\endproof

Now, Lemma~\ref{p.corr} immediately implies by induction that
$Q_\idot^n(k)=0$ unless $n$ is a power of $p$, and then being a
$Q_\idot^n(k)$-module, $Q_\idot^n(V)$ must vanish for any $k$-vector
space $V$. Thus we might as well renumber the functors $Q_\idot^n$
by setting
\begin{equation}\label{Q.i}
Q^{(i)}_\idot = Q^{p^i}_\idot, \quad\quad HQ_\idot^{(i)} =
HQ_\idot^{p^i}, \qquad i \geq 0,
\end{equation}
and let $D^{(i)} = D^{p^i}$. Since $k$ is a prime field, the
commutative algebra $k(V)$ of all $k$-valued functions on a
$k$-vector space $V$ has the Frobenius endomorphism equal to the
identity, $f^p=f$, so that the tautological map $V^* \to k[V]$
extends to maps $S^n(V^*) \to k(V)$, $n \geq 1$ that are compatible
with the maps \eqref{S.fr}. Dually, we have maps $k[V] \to D^n(V)$
that give rise to a map
\begin{equation}\label{D.lim}
k[V] \to \lim_iD^{(i)}(V),
\end{equation}
where the limit is taken with respect to the
maps \eqref{D.fr}. Both its source and target are functorial in $V$,
so that it induces a map
\begin{equation}\label{Q.lim}
Q_\idot(k/k) \to R^\hdot\lim_iQ^{(i)}_\idot(k),
\end{equation}
where as in Subsection~\ref{add.exa.subs}, $Q_\idot(k/k)$ is the
stabilization of the linear span functor $V \mapsto k[V]$.

\begin{lemma}
The map \eqref{Q.lim} is a quasiisomorphism.
\end{lemma}

\proof{} To construct stabilizations in $\Ho(I,k)$ for some
half-additive category $I$, we only need to consider finite
coproducts in $I$. Therefore we can restrict our attention to the
subcategory $P(k) \subset k\amod$ spanned by finite-dimensional
$k$-vector spaces $V$. For any such $V$, denote by $\K_\idot(V) =
S^\hdot(V \otimes_k I_\idot)$ the total symmetric power of the
complex $V \otimes_k I_\idot$, where $I_\idot$ is as
\eqref{I.seq}. Then $\K_\idot(V)$ is a flat $S^\hdot(V)$-algebra
quasiisomorphic to $k$ (this is the standard Koszul
resolution). Moreover, any $k$-linear map $a:V \to S^\hdot(V)$
uniquely extends to an algebra map $\exp(a):S^\hdot(V) \to
S^\hdot(V)$, and we can consider the complex $\K_\idot(V,a) =
\K_\idot(V) \otimes_{S^\hdot(V)} S^\hdot(V)$, where $S^\hdot(V)$ is
a module over itself via the map $\exp(a)$. We have a map $\exp(a) =
\id \otimes \exp(a):\K_\idot(V) \to \K_\idot(V,a)$, and if $\exp(a)$
is flat, the complex $\K_\idot(V,a)$ only has homology in degree
$0$. In particular, this is the case if $a=\phi:V \to S^p(V) \subset
S^\hdot(V)$ is the map \eqref{S.fr}, and if $a=\phi-\id$ is its
difference with the identity map $\id:V \to V \subset
S^\hdot(V)$. Moreover, in the latter case, we can equip
$\K_\idot(V,\phi-\id)$ with a multiplicative increasing filtration
$F_\idot$ by assigning filtered degree $1$ to $V \subset S^\hdot(V)
= \K_0(V,\phi-\id)$ and $p$ to $V \otimes k \subset V \otimes
S^\hdot(V) = \K_1(V,\phi-\id)$. Then
$\gr_\idot^F\K_\idot(V,\phi-\id) \cong \K_\idot(V,\phi)$, so that
all the associated graded pieces also have homology only in degree
$0$, and then by the spectral sequence argument, the same holds for
the filtered pieces $F_n\K(V,\phi-\id)$, $n \geq 0$.

Since $V$ is finite-dimensional and $k$ is prime, all $k$-valued
functions on the dual vector space $V^*$ are polynomial, so that the
map $S^\hdot(V) \to k(V^*)$ is surjective. In effect, $k(V^*)$ is
the quotient of $S^\hdot(V)$ by the ideal generated by
$(\phi-\id)(V)$, and we have a quasiisomorphism
$\K_\idot(V,\phi-\id) \cong k(V^*)$. Dually, for any map $a:V^* \to
S^\hdot(V^*)$, let $\K^\hdot(V,a) = \K_\idot(V^*,a)^*$, with the
descreaing filtration $F^n\K^\hdot(V,a) = F_n\K_\idot(V^*,a)^*$;
then we have a functorial quasiisomorphism
\begin{equation}\label{K.lim}
k[V] \cong \K^\hdot(V,\phi-\id) \cong \lim_nF^n\K^\hdot(V,\phi-\id),
\end{equation}
where all the terms in the limit have homology concentrated in
degree $0$. The induced filtration on the group algebra $k[V]$ is
the filtration by the powers of the augmentation ideal, and since
$V$ is finite-dimensional, it goes to $0$ at some finite
step. Moreover, if we assign filtered degree $p^i$ to $D^{(i)}(V)$,
then \eqref{D.lim} becomes a filtered map, and by
Remark~\ref{stab.prod.rem}, additivization with respect to $V$
commutes with the limits in \eqref{D.lim} and in
\eqref{K.lim}. Therefore it suffices to prove that the associated
graded quotient of the map \eqref{D.lim} becomes a quasiisomorphism
after applying additivization.

To do this, it is convenient to lift \eqref{D.lim} to a map of
complexes. On one hand, we have the surjective map $\exp(\phi -
\id)^*:\K^\hdot(V,\phi-\id) \to \K^\hdot(V)$, it becomes filtered if
the rescale the filtration on its target by $p$, and $\K^\hdot(V)
\cong k$ has trivial additivization on all its graded pieces, so
that we may replace $\K^\hdot(V,\phi-\id)$ in \eqref{K.lim} with
$\bK^\hdot(V) = \Ker \exp(\phi-\id)^*$. On the other hand, by the
telescope construction, we have an exact sequence
\begin{equation}\label{D.exa}
0 \longrightarrow \lim_iD^{(i)}V \longrightarrow \prod_{i \geq
  0}D^{(i)}(V) \overset{\phi^* - \id}{\longrightarrow} \prod_{i \geq
  0}D^{(i)}(V) \longrightarrow 0,
\end{equation}
where $\phi^*:D^{(i+1)}(V) \to D^{(i)}(V)$ are the maps
\eqref{D.fr}, and \eqref{D.exa} becomes a filtered exact sequence if
we assign filtered degree $p^i$ resp.\ $p^{i+1}$ to $D^{(i)}(V)$ in
the middle resp.\ rightmost term. But for any $a$, we have a natural
projection
$$
\K^\hdot(V,a) \to \K^0(V,a) \cong \prod_nD^n(V) \to
\prod_iD^{(i)}(V),
$$
and it intertwines $\exp(\phi-\id)^*$ and the map $\phi^*-\id$ in
\eqref{D.exa}, thus induces a filtered map $\bK^\hdot(V) \to
\lim_iD^{(i)}(V)$. This is our map.

It now remains to observe that
$\gr_F^n\K^i(V,\phi-\id)=\gr_F^n\K^i(V)=0$ if $ip > n$, and
$\exp(\phi-\id)^*:\gr_F^n\K^i(V,\phi-\id) \to \gr_F^n\K^i(V)$ is an
isomorphism if $n=ip$, so that $\gr_F^n\bK^i(V)=0$ when $ip \geq
n$. If $0 < ip < n$, then $\gr_F^n\K^i(V,\phi-\id) \cong
\Lambda^i(V) \otimes D^{n-ip}(V)$ has trivial additivization by
\eqref{Sig.st}, and similarly for $\gr_F^n\K^i(V)$, and if $i=0$ but
$n$ is not a power of $p$, then both have trivial additivization by
Lemma~\ref{p.corr}. Finally, if $n=p^i$, then $\gr_F^n\bK^0(V)$ is
the $n$-th associated graded piece of \eqref{D.exa} on the nose.
\endproof

To see how \eqref{Q.lim} yields \eqref{st.odd}, note that by
Lemma~\ref{p.corr} and induction, we have
\begin{equation}\label{HQ.i}
HQ^{(i)}_\idot(k) \cong
k[\beta,\xi_1,\tau_1,\dots,\xi_{i-1},\tau_{i-1},\xi_i],
\end{equation}
where the degrees of the generators are the same as in
\eqref{st.odd}, there are no $\tau$ if $p=2$, and the transition
maps in \eqref{Q.lim} act on homology by sending the extra
generators to $0$. Therefore the inverse system of homology groups
stabilizes at a finite step in each degree, $R^1\lim$ vanishes, and
the homology $HQ_\idot(k/k)$ of the DG algebra $Q_\idot(k/k)$ is
exactly as in \eqref{st.odd}.

\begin{remark}\label{sthom.rem}
The dual Steenrod algebra is of course not just an algebra but a
Hopf algebra, with some comultiplication that is too non-linear to
admit a DG model. In terms of stabilization, one observes that by
adjunction, the linear span functor $V \mapsto k[V]$ is a comonad,
and then the endofunctor of the category $\Ho(k)$ given by its
stabilization $Q_\idot$ is also a comonad. Since the comonad is
non-linear, adding an enhancement to it requires some technology,
but whatever technology one uses, enhanced coalgebras over this
enhanced comonad are connective spectra, for more-or-less
tautological reasons (for example, if one uses ``stable model
pairs'', then this is \cite[Theorem 10.6]{kazh}). However, observe
that the whole projective system \eqref{D.lim} has a structure of a
comonad, with the structure maps \eqref{SD.mn}, and then so does its
stabilization \eqref{Q.lim}. In other words, the filtration on
$Q_\idot$ given by \eqref{Q.lim} is compatible with the comonad
structure, and then coalgebras over this filtered version of
$Q_\idot$ form a ``filtered'' version of the stable homotopy
category. It seems that this category has not been considered yet,
and it might be interesting. We hope to return to this elsewhere.
\end{remark}

\section{Monoidal structures and $2$-categories.}

\subsection{The Segal condition and $2$-categories.}

Let us now recall the description of symmetric monoidal structures
in terms of the category $\Gamma_+$ given in
Subsection~\ref{mono.subs}.  To encode non-symmetric monoidal
categori\-es, one can use the same Segal machine as in
Definition~\ref{sym.1.def} but with $\Delta^o$ instead of
$\Gamma_+$. It is convenient to start with a more general notion of
a $2$-category.

\begin{defn}\label{2cat.def}
A cofibration $\C \to \Delta^o$ {\em satisfies the Segal condition}
if for any $n \geq l \geq 0$, the functor
\begin{equation}\label{seg.eq}
\C_{[n]} \to \C_{[l]} \times_{\C_{[0]}} \C_{[n-l]}
\end{equation}
induced by \eqref{seg.sq} is an equivalence. A {\em $2$-category} is
a cofibration $\C \to \Delta^o$ with discrete $\C_{[0]}$ satisfying
the Segal condition. A {\em lax $2$-functor} between $2$-categories
$\C$, $\C'$ is a functor $\gamma:\C \to \C'$ over $\Delta^o$
cocartesian over anchor maps. A {\em $2$-functor} is a lax
$2$-functor that is cocartesian over all maps.
\end{defn}

\begin{remark}\label{seg.rem}
More generally, for any $n \geq l' \geq l \geq 0$, we have a
cocartesian square
\begin{equation}\label{seg.1.sq}
\begin{CD}
[l'-l] @>{t}>> [l']\\
@V{s}VV @VV{s}V\\
[\nl] @>{t}>> [n]
\end{CD}
\end{equation}
in $\Delta$, and if $\C \to \Delta^o$ satisfies the Segal condition,
then we also have
$$
\C_{[n]} \cong \C_{[l']} \times_{\C_{[l'-l]}} \C_{[n-l]}
$$
for any square \eqref{seg.1.sq} (just combine \eqref{seg.eq} for $n
\geq l' \geq 0$ and $l' \geq l \geq 0$).
\end{remark}

\begin{remark}\label{st.rem}
Assume given a cofibration $\pi:\C \to \Delta^o$, and say that a
functor $c:\V \to \C^o$ is {\em standard} if it is cartesian over
$\Delta$, and $\pi^o \circ c:\V \to \Delta^o$ is the top left part
of a square \eqref{seg.1.sq}. Then $\C$ satisfies the Segal condition
of Definition~\ref{2cat.def} if and only if for any standard $c:\V
\to \C^o$ there exists a colimit $\colim_{\V} c^o$, and the natural
map $\pi(\colim_{\V} c^o) \to \colim_{\V} (\pi^o \circ c) = [n]$ is
an isomorphism. Moreover, since $s:[l] \to [n]$ and $t:[\nl] \to
[n]$ are anchor maps, the opposite $\gamma^o$ to any lax $2$-functor
$\gamma:\C \to \C'$ to some $2$-category $\C'$ preserves the
colimits of standard functors $c:\V \to \C^o$. In addition to this,
say that $c:\V \to \C^o$ is {\em half-standard} if $c(o) \to c(0)$
is a cartesian lifting of an anchor map, and $c(o) \to c(1)$ is
vertical with respect to $\pi^o$. Then for such a $c$,
$\colim_{\V} c$ exists as soon as $\C$ satisfies the Segal
condition, and these colimits are also preserved by opposites
$\gamma^o$ to lax $2$-functors.
\end{remark}

\begin{exa}\label{E.exa}
For any category $\C$, let $\eps:\ppt \to \Delta^o$ be the embedding
onto $[0]$, and let $E\C = \eps_*\C$. Then the fibers of the
cofibration $E\C \to \Delta^o$ are $E\C_{[n]} \cong \C^{V([n])}$, so
that $E\C$ trivially satisfies the Segal condition, and if $\C$ is
discrete, then $E\C$ is a $2$-category. Moreover, for any
cofibration $\pi:\C \to \Delta^o$, we have a natural functor
\begin{equation}\label{nu.eq}
\nu:\C \to E\C_{[0]} = \eps_*\eps^*\C
\end{equation}
induced by \eqref{adj.2.1}. This functor $\nu$ is cocartesian over
$\Delta^o$, and if $\C$ is a $2$-category, $\nu$ is a $2$-functor.
\end{exa}

\begin{remark}\label{red.rem}
The second condition of Definition~\ref{2cat.def} -- namely, the
requirement that $\C_{[0]}$ is discrete -- can be always achieved by
the following trick. For any cofibration $\C \to \Delta^o$, consider
the discrete subcategory $\C_{[0],\Id} \subset \C_{[0]}$, and define
the {\em reduction} $\C^{red}$ by the cartesian square
$$
\begin{CD}
  \C^{red} @>>> \C\\
  @VVV @VV{\nu}V\\
  E\C_{[0],\Id} @>>> E\C_{[0]},
\end{CD}
$$
where $\nu$ is the functor \eqref{nu.eq}. Then if $\C$ satisfies the
Segal condition, so does $\C^{red}$, and $(\C^{red})_{[0]} =
\C_{[0],\Id}$ is discrete. However, this procedure has to be used
with caution, since $\C_{[0],\Id} \subset \C_{[0]}$, hence also
$\C^{red}$ depend on $\C$ on the nose, and not only on its
equivalence class. To aleviate the problem, it is better to first
replace $\C$ with its tightening.
\end{remark}

Definition~\ref{2cat.def} is a convenient packaging of the usual
notion of a (weak) $2$-category: objects are objects of $\C_{[0]}$,
and for any $c,c' \in \C_{[0]}$, the fiber $\C(c,c')$ of the
projection $s^o_! \times t^o_!:\C_{[1]} \to \C_{[0]}$ over $c \times
c'$ is the category of morphisms from $c$ to $c'$. For $m \geq 2$,
we have an equivalence
\begin{equation}\label{c.m.eq}
\prod_{i=1}^ma_{i!}^o:\C_{[m]} \cong \C_{[1]} \times_{\C_{[0]}} \dots
\times_{\C_{[0]}} \C_{[1]},
\end{equation}
with $m$ copies of $\C_{[1]}$ numbered by edges $i \in E([n]_\delta)
= \{1,\dots,m\}$ of the string quiver $[n]_\delta$, and $a_i:[1] \to
[m]$ standing for the embeddings of the edges. For every $c \in
\C_{[0]}$, we have the identity object $\id_c \in \C(c,c)$ induced
by the tautological projection $[1] \to [0]$, and the composition
functors $- o -:\C(c,c') \times \C(c',c'') \to \C(c,c'')$ are
induced by the map $m:[1] \to [2]$ sending $0$ to $0$ and $1$ to
$2$.

The {\em point $2$-category} $\ppt^2$ is $\Delta^o$ itself. As it
should, it has one object, $\ppt^2_{[0]} \cong \ppt$, and for any
object $c \in \C_{[0]}$ in a $2$-category $\C$, the embedding
$\eps(c):\ppt \to \C_{[0]}$ onto $c$ uniquely extends to a
$2$-functor $\eps(c):\ppt^2 \to \C$ that we also call the embedding
onto $c$. The tautological projection $\tau:\C \to \ppt^2$ is the
structural cofibration. If we have another $2$-category $\C'$, then
the {\em constant $2$-functor $\C' \to \C$ with value $c$} is the
composition of $\tau:\C' \to \ppt^2$ and $\eps(c):\ppt^2 \to
\C$. The {\em $2$-product} $\C \times^2 \C'$ is the product $\C
\times_{\Delta^o} \C'$. A $2$-category $\C$ is {\em discrete} if $\C
\cong \C_{[0]} \times \Delta$. For any $2$-category $\C$, applying
the involution $\iota:\Delta \to \Delta$, $[n] \mapsto [n]^o$ gives
a cofibration $\iota^*\C \to \Delta$, and this is also a
$2$-category; we call it the {\em opposite $2$-category} and denote
$\C^\iota$.

In keeping with our usage for ordinary categories, we say that a lax
$2$-functor $\gamma:\C' \to \C$ is {\em dense} if
$\gamma_{[0]}:\C'_{[0]} \to \C_{[0]}$ is an equivalence. For any lax
$2$-functor $\gamma$, we can define a $2$-category $\gamma^*\C$ by
the cartesian square
\begin{equation}\label{2.ga}
\begin{CD}
\gamma^*\C @>>> \C\\
@VVV @VV{\nu}V\\
E\C'_{[0]} @>>> E\C_{[0]},
\end{CD}
\end{equation}
where $\nu$ is as in Example~\ref{E.exa}. Then $\gamma$ factors as
\begin{equation}\label{d.f.facto}
\begin{CD}
\C' @>{\wgamma}>> \gamma^*\C @>{\ogamma}>> \C,
\end{CD}
\end{equation}
where $\wgamma$ is tautologically dense. We say that $\gamma$ is
{\em $2$-fully faithful} if $\wgamma$ is an equivalence; in
particular, $\ogamma$ in \eqref{d.f.facto} is
tautologically $2$-fully faithful.

\begin{exa}\label{sub.exa}
If $\gamma$ is the tautological embedding $\gamma:\C_{[0]} \times
\Delta^o \to \C$, then $\gamma^*\C \cong \C$. More generally, if we
have an embedding $S \subset \C_{[0]}$, with the corresponding
functor $\gamma:S \times \Delta^o \to \C$, then $\gamma^*\C$ can be
thought of a the full $2$-subcategory in $\C$ spanned by objects $c
\in S \subset \C_{[0]}$.
\end{exa}

It is useful to generalize Example~\ref{sub.exa} in the following
way. Assume given a functor $S:\C_{[0]} \to \Sets$ (that is, a set
$S_c$ for any $c \in \C_{[0]}$). Denote by $\pi:\C_{[0]}[S] \to
\C_{[0]}$ be the corresponding discrete cofibration, let
$\chi:\C_{[0]}[S] \times \Delta^o \to \C$ be the composition of the
$2$-functor $\pi \times \id:\C_{[0]}[S] \times \Delta \to \C_{[0]}
\times \Delta$ with the tautological embedding $\C_{[0]} \times
\Delta^o \to \C$, and let
\begin{equation}\label{C.S.eq}
\C[S] = \chi^*\C,
\end{equation}
where the right-hand side is as in \eqref{2.ga}.  Equivalently,
$\C[S]$ can by obtained by extending $S$ to a functor $S:\C \to
\Sets$ by the right Kan extension with respect to the embedding
$\C_{[0]} \to \C$, and taking the corresponding discrete cofibration
$\C[S] \to \C$. Explicitly, for any $[n] \in \Delta$ and $c \in
\C_{[n]}$, we have
$$
S(c) \cong \prod_{i=0}^nS(b_{i!}^oc),
$$
where $b_i:[0] \to [n]$ sends $0$ to $i$. One can also consider the
universal situation: the forgetful functor $\Sets_+ \to \Sets$ is a
discrete cofibration with fiber $S$ over any $S \in \Sets$, we can
consider the induced cofibration $E\Sets_+ \to E\Sets$, and then
$\C[S]$ fits into a cartesian square
\begin{equation}\label{C.S.sq}
\begin{CD}
\C[S] @>>> E\Sets_+\\
@V{\pi}VV @VVV\\
\C @>{E(S) \circ \nu}>> E\Sets,
\end{CD}
\end{equation}
where $\nu$ is the functor \eqref{nu.eq}. Note that up to an
isomorphism, the functor $S$ factors through $\Sets_{\Id}$, so that
we may replace $E\Sets_+ \to E\Sets$ in \eqref{C.S.sq} with the
induced cofibration over $E\Sets_{\Id} \subset E\Sets$ without
changing $\C[S]$ (since the cofibration $E\Sets_{\Id} \to \Delta^o$
is discrete, $E(S) \circ \nu$ in \eqref{C.S.sq} then also becomes a
cofibration). In either description, we obviously have $\C[\ppt]
\cong \C$, where $\ppt:\C_{[0]} \to \Sets$ sends everything to the
one-element set.

\subsection{Simplcial replacements and special maps.}

In any meaningful formalism, a usual category should define a
$2$-category. In the context of Definition~\ref{2cat.def}, this can
be achieved as follows. For any category $I$, denote
\begin{equation}\label{del.I}
\Delta^o\langle I \rangle = \Id^+_{**}(I \times \Delta^o),
\end{equation}
where $I \times \Delta^o \to \Delta^o$ is the trivial cofibration,
$\Id:\Delta^o \to \Delta^o$ is the identity functor, $+$ is the
class of special maps, and $\Id^+_{**}$ has the same meaning as in
Example~\ref{id.exa}. If we identify $\Ar^t(\Delta)$ and the
category $\Delta_\idot$ of \eqref{nu.dot}, then \eqref{id.st} and
\eqref{ffun.st} provide an identification $\Delta^o \langle I
\rangle \cong \Fun(\Delta_\idot/\Delta,I)$, so that for any $[n]
\in \Delta^o$, the fiber $\Delta^o\langle I \rangle_{[n]}$ of the
cofibration $\Delta(I) \to \Delta^o$ is the functor category
$\Fun([n],I)$. In particular, $\Delta^o\langle I \rangle$ obviously
satisfies the Segal condition.

\begin{defn}\label{DI.def}
The {\em simplicial replacement} $\Delta^o I$ of a category $I$ is
the reduction $\Delta(I)^{red}$ of the cofibration \eqref{del.I} in
the sense of Remark~\ref{red.rem}.
\end{defn}

Explicitly, objects in $\Delta^o I$ are pairs $\langle
[n],i_\idot\rangle$ of an object $[n] \in \Delta$ and a functor
$i_\idot:[n] \to I$, with maps from $\langle [n],i_\idot \rangle$ to
$\langle [n'],i'_\idot \rangle$ given by a map $f:[n'] \to [n]$ and
a map $i'_\idot \to f^*i_\idot$ that is pointwise an identity
map. The {\em augmented simplicial replacement} is the category
$\Delta^{o>} I = (\Delta^o I)^>$. If the category $I$ is small,
$\Delta^o I \cong \Delta^o N(I)$ and $\Delta^{o>}I =
\Delta^{o>}N(I)$ are the categories of simplicis of its nerve
$NI:\Delta^o \to \Sets$, as in Subsection~\ref{cyl.subs}. For any
$I$, $\Delta^o I$ is a $2$-category in the sense of
Definition~\ref{2cat.def}, with the structural cofibration $\Delta^o
I \to \Delta^o$ given by the forgetul functor $\langle [n],i_\idot
\rangle \mapsto [n]$, and the cofibration is discrete. Its extension
$\Delta^{o>}I \to \Delta^{o>} = \Delta^{<o}$ is also a discrete
cofibration by Example~\ref{gt.exa}. We have $\Delta^o \ppt \cong
\Delta \cong \ppt^2$, and $\Delta^oI^o \cong (\Delta^o I)^\iota$ is
the $2$-category opposite to $\Delta^oI$. A functor $\gamma:I' \to
I$ induces a $2$-functor $\Delta^o \gamma:\Delta^o I' \to \Delta^o
I$.

\begin{exa}\label{eS.exa}
If $S$ is a set, and $e(S)$ is $S$ with maximal preorder, as in
Subsection~\ref{set.sss}, then $\Delta^o e(S) \cong ES$, where $ES
\to \Delta^o$ is as in Example~\ref{E.exa}.
\end{exa}

The category of simplices $\Delta^o X$ of a simplicial set $X$
satisfies the Segal condition if and only so does $X$, and in this
case, $X \cong N(I)$ and $\Delta^o X \cong \Delta^o I$ for a unique
small category $I$. Thus any bounded $2$-category $\C$ such that the
cofibration $\C \to \Delta^o$ is discrete is actually the simplicial
replacement of a small $1$-category, $\C \cong \Delta^oI$. For any
bounded $2$-category $\C$ with the structural cofibration $\pi$, the
simplicial set $\pi_!\ppt:\Delta^o \to \Sets$ obviously satisfies
the Segal condition, so that $\pi_0(\C/\Delta^o) \cong
\Delta^o\tau(\C)$ for a unique small category $\tau(\C)$ that we
call the {\em truncation} of the $2$-category $\C$. Explicitly,
objects of $\tau(\C)$ are objects $c \in \C_{[0]}$, and morphisms
are connected components of the morphism categories in $\C$,
$\tau(\C)(c,c') = \pi_0(\C(c,c'))$.

\begin{remark}
The correspondence $I \mapsto \Delta^oI$ respects objects and
morphisms but loses the $2$-categorical structure. In particular,
for an equivalence $I' \to I$ between small categories, the induced
functor $\Delta^oI' \to \Delta^oI$ is not in general an
equivalence. In terms of \eqref{C.S.eq}, we have $\Delta^oI' \cong
\Delta^oI[S]$, where $S:I_{\Id} \to \Sets$ sends $i \in I$ to the
set of its preimages in $I'$.
\end{remark}

Another useful class of $2$-categories is {\em $(1,2)$-categories},
namely, $2$-cate\-gories $\C$ such that for any $c,c' \in \C_{[0]}$,
the category $\C(c,c')$ is a groupoid. Equivalently, one can require
that the structural cofibration $\C \to \Delta^o$ is semidiscrete.
For any $2$-category $\C$, the dense subcategory $\C_\natural
\subset \C$ spanned by cocartesian maps is then a $(1,2)$-category,
with the same object and isomorphism groupoids $\C(c,c')_{\Iso}
\subset \C(c,c')$ as categories of morphisms. If $\C_\natural$ is
bounded, we can consider its truncation $\tau(\C_\natural)$. Objects
in $\tau(\C_\natural)$ are still objects $c \in \C_{[0]}$, and maps
are isomorphism classes of $1$-morphisms in $\C$.

By abuse of terminology, we define a $2$-functor from a category $I$
to a $2$-category $\C$ as a $2$-functor $\gamma$ from $\Delta^o I$
to $\C$, and similarly for lax $2$-functors. In the other
directions, since $\Delta^o I$ is discrete, any lax $2$-functor from
$\C$ to $\Delta^oI$ is automatically a $2$-functor. It turns out
that these then correspond bijectively to usual functors $\C \to I$
of a special kind.

\begin{defn}\label{sp.def}
A morphism $f$ in a $2$-category $\C$ is {\em special} if it is a
cocartesian lifting of a special map in $\Delta^o$. A functor $E:\C
\to \E$ to some category $\E$ is {\em special} if it inverts all
special maps, and a cofibration $\C' \to \C$ is {\em special} if for
any special map $f$ in $\C$, the transition functor $f_!$ is an
equivalence.
\end{defn}

For any small category $I$ with simplicial replacement $\Delta^oI$,
the evaluation functor \eqref{adj.2.2} induces a special functor
\begin{equation}\label{xi.eq}
\xi:\Delta^o I \to I
\end{equation}
that sends $\langle [n],i_\idot \rangle$ to $i_\idot(0) \in I$. We
also have $\Delta^oI^o \cong \iota^*\Delta^oI$, so that
\eqref{xi.eq} for the category $I^o$ provides a functor
$\xi_\perp:\Delta^o I \to I^o$ sending $\langle [n],i_\idot \rangle$
to $i_\idot(n) \in I^o$. One can also combine $\xi$ and $\xi_\perp$
into a single functor
\begin{equation}\label{xi.fl}
\xi_\flat:\Delta^o I \to \Tw(I)
\end{equation}
sending $\langle [n],i_\idot \rangle$ to the arrow $i_\idot(0) \to
i_\idot(n)$. Then $\xi$ and $\xi_\perp$ are obtained by composing
\eqref{xi.fl} with the projections \eqref{ar.tw.co}
resp.\ \eqref{ar.tw}.

Now, for any $2$-category $\C$ with truncation $\tau(\C)$, we can
compose \eqref{xi.eq} for $\tau(\C)$ with the natural projection $\C
\to \Delta^o\tau(\C) = \pi_0(\C)$ to obtain a special functor
\begin{equation}\label{xi.C.eq}
\xi:\C \to \tau(\C).
\end{equation}
Then \eqref{xi.eq} and \eqref{xi.C.eq} enjoy the following universal
properties.

\begin{lemma}\label{sp.le}
Any special functor $E:\Delta^o I \to \E$ to some category $\E$
factors uniquely through the functor \eqref{xi.eq}, and any special
cofibration $\C' \to \Delta^o I$ is of the form $\C' \cong \xi^*\C$
for a unique cofibration $\C \to I$.
\end{lemma}

\proof{} The functor $\zeta([0])$ of \eqref{beta.eq} for the
cofibration $\Delta_+^oI = \rho^{o*}\Delta^oI \to \Delta_+$ gives a
projection $\Delta_+^oI \to I_{\Id}$, or in other words, a
decomposition
\begin{equation}\label{C.prod}
\Delta^o_+I \cong \coprod_{i \in I}(\Delta^o_+I)_i,
\end{equation}
a categorical version of \eqref{ka.de.0}, and then as in
Lemma~\ref{aug.le}, for any $i \in I$, we have an adjoint
pair of functors
\begin{equation}\label{rho.la.C}
\begin{aligned}
\lambda(i) &= \lambda:i \setminus^\xi \Delta^o I \cong
\lambda^*\rho^*(\Delta^o_+I)_i \to (\Delta^o_+I)_i,\\
\rho(i) &= a^*\rho^*:(\Delta^o_+ I)_i \to i \setminus^\xi \Delta^o I,
\end{aligned}
\end{equation}
where $a:\lambda \circ \rho \to \id$ is the adjunction
map. Therefore the subcategories $(\Delta^o_+I)_i^o \subset
(\Delta^o I)^o$ with the augmentations $\rho^o(i)=\rho(i)^o$ give a
framing for the functor $\xi^o:(\Delta^o I)^o \to I^o$ in the sense
of Lemma~\ref{kan.le}. Moreover, for any special functor
$E:\Delta^oI \to \E$, the opposite functor $E^o$ is locally constant
on ${(\Delta^o_+I)_i}$, thus constant since $(\Delta^o_+I)_i$ has an
initial object. Then by Lemma~\ref{kan.le}, $\xi^o_!E^o$ exists, and
the adjunction map $\xi^{o*}\xi^o_!E^o \to E^o$ is an isomorphism,
so that $E^o$ indeed factors through $\xi^o$ (and then $E$ factors
through $\xi$). For cofibrations, let $\C = \xi_*\C'$, and use the
adjunctions \eqref{rho.la.C} and equivalences \eqref{sec.adj} to
check that the functor $\xi^*\C \to \C'$ of \eqref{adj.2.1} is an
equivalence.
\endproof

\begin{corr}\label{sp.corr}
Any special functor $E:\C \to \E$ from a $2$-category $\C$ to a
category $\E$ factors uniquely through the functor \eqref{xi.C.eq}.
\end{corr}

\proof{} By \eqref{adj.2.eq}, $E$ factors through $\Delta(\E)$, and
then since $\C_{[0]}$ is discrete, it further factors as
$$
\begin{CD}
\C @>{E'}>> \Delta^o\E @>{\xi}>> \E,
\end{CD}
$$
where $E'$ is cocartesian over $\Delta^o$. Since $\Delta^o\E \to
\Delta^o$ is discrete, $E'$ further factors through $\pi_0(\C) =
\Delta^o\tau(\C)$, and we are done by Lemma~\ref{sp.le}.
\endproof

\begin{remark}
In terms of Definition~\ref{loc.def}, Corollary~\ref{sp.corr} can be
rephrased to say that \eqref{xi.C.eq} is a localization, and special
maps are dense in $\chi^*\Iso$.
\end{remark}

\begin{remark}
An obvious counterpart of Corollary~\ref{sp.corr} for special
cofibrations is completely wrong.
\end{remark}

\subsection{Cylinders and $2$-functors.}

For any two $2$-ca\-te\-go\-ries $\C$, $\C'$ with $\C'$ bounded,
bounded $2$-functors from $\C'$ to $\C$ form a full subcategory in
$\Fun_{\Delta^o}(\C',\C)$ that we denote by $\Fun^2(\C',\C) \subset
\Fun_{\Delta^o}(\C',\C)$. For any small category $I$ and
$2$-category $\C$, we simplify notation by writing $\Fun^2(I,\C) =
\Fun^2(\Delta^o I,\C)$. For any $[n] \in \Delta$, we tautologically
have $\Fun^2([n],\C) \cong \C_{[n]}$. More generally, if are given a
quiver $Q:\Dd^o \to \Sets$, a $2$-category $\C$ defines a
cofibration $\C(Q) \to \Dd^o$ with fibers $\C(Q)_i =
\C^{Q(i)}_{[i]}$, $i = 0,1$.

\begin{lemma}\label{qui.le}
For any quiver $Q$ and $2$-category $\C$, we have a natural
identification
\begin{equation}\label{P.fun.2}
\Fun^2(P(Q),\C) \cong \Sec^\natural(\Dd^o,\C(Q))
\end{equation}
between $2$-functors from the path category $P(Q)$ to $\C$ and
cocartesian sections of the cofibration $\C(Q) \to \Dd^o$.
\end{lemma}

In particular, if we have a cocartesian square $Q_\idot:[1]^2 \to
\Dd^o\Sets$ of quivers, so that $P(Q_\idot):[1]^2 \to \rCat$ is a
cocartesian square of small categories, then by Lemma~\ref{qui.le},
the corresponding square of categories $\Fun^2(P(Q_\idot),\C)$ is
cartesian. This of course includes the squares \eqref{seg.1.sq} but
there are other useful examples.

\proof{} Let $A(Q) = \alpha^o_*Q:\Delta_a^o \to \Sets$ be the right
Kan extension with respect to the functor \eqref{P.a}, and let
$\Dd^oQ \to \Dd^o$, $\Delta_a^oA(X) \to \Delta_a^o$ be the discrete
cofibrations corresponding to $Q$ and $A(Q)$. Then \eqref{P.a} lifts
to a diagram
\begin{equation}\label{P.a.Q}
\begin{CD}
\Dd^oQ @>{\alpha(Q)}>> \Delta_a^oA(Q) @>{\beta(Q)}>> \Delta^oP(Q)
\end{CD}
\end{equation}
of functors over $\Delta^o$, and \eqref{adj.2.eq.1} provides an
identification $\Sec^\natural(\Dd^o,\C(Q)) \cong
\Fun_{\Dd^o}^\natural(\Dd^oQ,\delta^{o*}\C)$. For any $E \in
\Fun_{\Dd^o}^\natural(\Dd^oQ,\delta^{o*}\C)$, the right Kan
extension $\alpha(Q)_*E$ with respect to \eqref{P.a.Q} can be
computed by (the dual version of) \eqref{kan.eq}, and since the
cofibration $\Delta_a^oA(X) \to \Delta_a^o$ is discrete, it
identifies the right comma-fibers of the functor $\alpha(Q)$ with
those of $\alpha$ of \eqref{P.a}. Then the limits in the right-hand
side of \eqref{kan.eq} reduce to iterated limits of standard
functors of Remark~\ref{st.rem}, so that $\alpha(Q)_*$ exists and
provides an equivalence of categories
\begin{equation}\label{a.Q.1}
\Fun^\natural_{\Dd^o}(\Dd^oQ,\delta^{o*}\C) \cong
\Fun^\natural_{\Delta_a^o}(\Delta_a^oA(Q),\beta^*\C)
\end{equation}
inverse to $\alpha(Q)^*$. Now as in Remark~\ref{qui.rem}, to compute
$P(Q) = \beta_!A(Q)$, one can combine \eqref{be.t} and
Example~\ref{di.exa}, and this provides an identification
\begin{equation}\label{del.P.Q}
\Delta^oP(Q) \cong \Ar^\pm(\Delta)^o_\natural \times_{\Delta_a^o}
\Delta_a^oA(Q),
\end{equation}
while \eqref{s.R} induces a functor
\begin{equation}\label{t.Q}
t(Q):\Delta^o P(Q) \to \Delta_a^o A(Q)
\end{equation}
right-adjoint to $\beta(Q)$. Then by \eqref{adj.I.eq}, the relaive
Kan extension $\beta_!^{\Delta^o}$ provides an equivalence
\begin{equation}\label{a.Q.2}
\Fun^\natural_{\Delta_a^o}(\Delta_a^oA(Q),\beta^*\C) \cong
\Fun^\natural_{\Delta^o}(\Delta^oP(Q),\C) = \Fun^2(\Delta^oP(Q),\C),
\end{equation}
and to finish the proof, it remains to combine \eqref{a.Q.2} and
\eqref{a.Q.1}.
\endproof

Even in the simple case $I=P(Q)$, a lax $2$-functor from a category
$I$ to a $2$-category $\C$ contains much more data than a
$2$-functor (for $I=\ppt$, this is considered below in
Subsection~\ref{mon.subs}). We do not attempt to prove any
classification results similar to Lemma~\ref{qui.le}. However, we do
need one general construction, namely, a $2$-categorical version of
the cylinder construction of Example~\ref{cyl.exa}. Assume given
categories $I_0$, $I_1$ and a functor $\gamma:I_0 \to I_1$, and
consider the cylinder $I = \Cyl(\gamma)$ with its cofibration $\pi:I
\to [1]$. Then we have a functor $\Delta^o(\gamma):\Delta^oI_0 \to
\Delta^oI_1$, and we can also consider the cylinder
$\Cyl(\Delta^o(\gamma))$. This is by definition a cofibration over
$[1] \times \Delta^o$, \eqref{xi.eq} induces a functor
$\Cyl(\Delta(\gamma)) \to I \times \Delta^o$ cocartesian over $[1]
\times \Delta^o_+$, and by \eqref{adj.2.eq}, this corresponds to a
functor
\begin{equation}\label{alpha.eq}
\alpha:\Cyl(\Delta^o(\gamma)) \to \Id^+_{**}(I \times \Delta^o) =
\Delta^o\langle I \rangle
\end{equation}
cocartesian over $[1] \times \Delta^o$. Now assume given a
$2$-category $\C$ and two lax $2$-functors $\phi_l:\Delta^oI_l \to
\C$, $l = 0,1$ equipped with an map $g:\phi_0 \to \phi_1 \circ
\gamma$. Then the triple $\langle \phi_0,\phi_1,g \rangle$ defines a
single functor $\phi:\Cyl(\Delta^o(\gamma)) \to \C$, and we have the
following $2$-categorical version of the cylinder construction.

\begin{lemma}\label{2.cyl.le}
The right Kan extension $\alpha_*\phi:\Delta^{o}\langle I \rangle
\to \C$ with respect to the functor \eqref{alpha.eq} exists, and its
restriction $\Cyl^2(\gamma,g):\Delta^o I \to \C$ to $\Delta^o I
\subset \Delta^o \langle I \rangle$ is a lax $2$-functor.
\end{lemma}

\proof{} Extend $\Delta^o(\gamma)$ to a functor
$\Delta^o(\gamma)^>:\Delta^{o>}I_0 \to \Delta^{o>}I_1$ between
augmented simplicial replacements, and consider the embedding
$\eps:\Cyl(\Delta^o(\gamma)) \to \Cyl(\Delta^o(\gamma)^>)$. Then
$\eps$ is left-closed in the sense of Example~\ref{cl.exa}, so we
have the canonical extensions $\alpha^>:\Cyl(\Delta^o(\gamma)^>) \to
\Delta^o\langle I \rangle^>$, $\phi^>:\Cyl(\Delta^o(\gamma)^>) \to
\C^>$ of the functors $\alpha$ and $\phi$, and $\alpha^> \cong
\eps_*\alpha$, $\phi^> \cong \eps_*\phi$, so that $(\alpha_*\phi)^>
\cong \alpha^>_*\phi^>$, where the right-hand side exists if and
only if so does the left. Thus it suffices to prove that
$\alpha^>_*\phi^>$ exists and restricts to a lax $2$-functor. We
have the embeddings $\sigma_l:I_l \to I$, $l=0,1$, $\sigma_1$ has a
left-adjoint $\tau:I \to I_1$, and the components
$\alpha^>_l:\Delta^{o>}I_l \to \Delta^o\langle I \rangle^>$, $l=0,1$
of the functor $\alpha^>$ are given by $\alpha^>_l =
\Delta^o(\sigma_l)^>$. The functor $\alpha^>_1$ has a left-adjoint
$\alpha^>_{1\dg}$ sending an object $\langle [n],i_\idot \rangle \in
\Delta^o\langle I \rangle^>$ to $\langle [n],\tau \circ
i_\idot\rangle$. Moreover, $\alpha^>_0$ also has a left-adjoint
$\alpha^>_{0\dg} = \mu$ of \eqref{mu.adj.eq}. Then for any object
$\langle [n],i_\idot \rangle$ in $\Delta^o\langle I \rangle^>$, we
have a commutative square
\begin{equation}\label{cyl.fr.sq}
\begin{CD}
\langle [n],i_\idot \rangle @>>> \alpha_0^>(\langle [n_0],i^0_\idot
\rangle)\\ 
@VVV @VVV\\
\alpha_1^>(\langle [n],\tau \circ i_\idot \rangle) @>>>
\alpha_1^>(\langle [n_0],\gamma \circ i^0_\idot \rangle)
\end{CD}
\end{equation}
in the category $\Delta^o\langle I \rangle^>$, where the arrows are
the adjunction maps. This gives a functor $\V \to
\Cyl(\Delta^o(\gamma)^>)^o / \langle [n],i_\idot \rangle$, and these
functors form a framing for $\alpha^o$ in the sense of
Lemma~\ref{kan.le}. If we now compute $\alpha_*$ using this framing,
then the relevant limits over $\V$ reduce to colimits of
half-standard functors in the sense of Remark~\ref{st.rem}, thus
exist, so that $\alpha^>_*\phi^>$ exists. Moreover, by
\eqref{cyl.fr.sq}, it fits into a cartesian square
\begin{equation}\label{cyl.2.sq}
\begin{CD}
\alpha^>_*\phi_> @>>> \phi_0^> \circ \alpha^>_{0\dg}\\
@VVV @VV{g}V\\
\phi_1^> \circ \alpha^>_{1\dg} @>>> \phi_1^> \circ
\Delta^o(\gamma)^> \circ \alpha^>_{0\dg},
\end{CD}
\end{equation}
and since $\alpha^>_{l\dg}$, $l=0,1$ and $\Delta^o(\gamma)^>$ send
anchor maps to anchor maps, $\Cyl^2(\gamma,g)$ is indeed a lax
$2$-functor.
\endproof

\begin{remark}\label{cyl.2.rem}
For an alternative construction of $\Cyl^2(\gamma,g)$ that does not
use the category $\Delta^o\langle I \rangle$, decompose $\gamma$ as
\begin{equation}\label{cyl.deco}
\begin{CD}
I_0 @>{\sigma}>> I @>{\zeta}>> I_1,
\end{CD}
\end{equation}
where $\sigma:I_0 \to I$ is the embedding and $\zeta$ is
left-adjoint to the embedding $I_1 \to I$, and consider the functor
$\nu:\Cyl(\Delta^o(\sigma)) \to \Cyl(\Delta^o(\gamma))$ over $[1]$
given by $\nu_0=\Id$ and $\nu_1=\Delta^o(\zeta)$. Then we have a
commutative diagram
\begin{equation}\label{cyl.2.dia}
\begin{CD}
\Cyl(\Delta^o(\sigma)) @>{\alpha'}>> \Delta^oI\\
@V{\nu}VV @VV{\tau}V\\
\Cyl(\Delta^o(\gamma)) @>{\alpha}>> \Delta^o\langle I \rangle,
\end{CD}
\end{equation}
where $\tau$ is the embedding, and $\alpha'$ is left-adjoint to
the embedding $\Delta^oI = \Cyl(\Delta^o(\sigma))_1 \subset
\Cyl(\Delta^o(\sigma))$. Moreover, the framing \eqref{cyl.fr.sq} for
the functor $\alpha$ lifts to a framing for the functor $\alpha'$,
and therefore the base change map
$$
\Cyl^2(\gamma,g) = \tau^*\alpha_*\phi \to \alpha'_*\nu^*\phi
$$
is an isomorphism, so that its target can be used as a definition of
$\Cyl^2(\gamma,g)$.
\end{remark}

\subsection{Monoidal structures.}\label{mon.subs}

Let us now turn to non-symmetric monoidal structures.  In terms of
Definition~\ref{2cat.def}, these correspond to $2$-categories with a
single object.

\begin{defn}\label{mon.def}
A {\em unital monoidal structure} on a category $\C$ is given by a
$2$-category $B\C$ such that $B\C_{[0]}=\ppt$ is the point category,
and $B\C_{[1]}$ is equipped with an equivalence $B\C_{[1]} \cong
\C$. A {\em lax monoidal structure} on a functor $\gamma:\C \to \C'$
between two categories equipped with unital monoidal structures
$B\C$, $B\C'$ is given by a lax $2$-functor $B\gamma:B\C \to B\C'$
equipped with an isomorphism $B\gamma_{[1]} \cong \gamma$. A lax
monoidal structure $B\gamma$ is {\em monoidal} if it is a
$2$-functor.
\end{defn}

Explicitly, the composition $- \circ -$ in $B\C$ defines the tensor
product $- \otimes -$ in $\C$, and the identity object $\id_{\ppt}
\in \C = B\C(\ppt,\ppt)$ is the unit object $1$ for the tensor
product. A lax monoidal structure on a functor $\gamma$ is given by
the maps \eqref{fi.fu} for the functor $B\gamma$; the essential ones
are the maps
\begin{equation}\label{lax.eq}
\gamma(M) \otimes \gamma(N) \to \gamma(M \otimes N), \qquad 1 \to
\gamma(1)
\end{equation}
corresponding to $m:[1] \to [2]$ and the tautological projection
$[1] \to [0]$. For any monoidal structure $B\C$ on a category $\C$,
the opposite $2$-category $B\C^\iota$ also defines a monoidal
structure on $\C$ (the product is the same but written in the
opposite direction). For consistency, we denote $\C$ with this
monoidal structure by $\C^\iota$.

\begin{exa}\label{end.exa}
For any $2$-category $\C$, and for any object $c \in \C_{[0]}$, the
category $\C(c,c)$ carries a natural monoidal structure $B\C(c,c)
\cong \eps(c)^*\C$, where $\eps(c):\ppt^2 \to \C$ is the embedding
onto $c$.
\end{exa}

\begin{exa}\label{sym.exa}
For any symmetric unital monoidal structure $\Bi\C$ on a category
$\C$ in the sense of Definition~\ref{sym.1.def}, the pullback $B\C =
\Sigma^*\Bi\C$ with respect to the functor \eqref{Si} is a unital
monoidal structure on $\C$ in the sense of Definition~\ref{mon.def},
and we have a canonical identification $B\C \cong B\C^\iota$.
\end{exa}

\begin{exa}\label{point.exa}
For any monoidal category $\C$ and bounded category $I$, the functor
category $\Fun(I,\C)$ carries a natural {\em pointwise} monoidal
structure given by $B\Fun(I,\C) \cong \Fun(I,B\C/\Delta^o)$. If
$\C$ is symmetric, then this is the same structure as in
Example~\ref{point.sym.exa}.
\end{exa}

Definition~\ref{2cat.def} and Definition~\ref{mon.def} are pretty
standard; however, it is more common to use fibrations over $\Delta$
rather than cofibrations over $\Delta^o$. The two notions are
equivalent -- every fibration has its transpose cofibration and vice
versa -- but it is the cofibrations that give the correct notion of
a lax $2$-functor. In particular, lax monoidal functors $\ppt \to
\C$ are the same thing as unital associative algebra objects in
$\C$, and it would be coalgebras were we to use fibrations (the maps
\eqref{lax.eq} would go in the other direction). Our definitions are
also compatible with Definition~\ref{sym.1.def} and
Definition~\ref{sym.2.def}: a unital symmetric monoidal structure
gives a unital monoidal structure via pullback $\Sigma^*$ with
respect to the functor \eqref{Si}, and the same goes for monoidal
and lax monoidal functors.

\begin{exa}\label{del.ten.exa}
The category $\Delta^<$ is a unital monoidal category with respect
to the concatenation product. The empty ordinal is the unit object,
and $[1] \in \Delta^<$ is naturally an algebra object in
$\Delta$. The corresponding cofibration $B\Delta^< \to \Delta^o$ is
obtained by taking $B\Delta^< = \Ar^\pm(\Delta)^o$, with the
projection to $\Delta^o$ opposite to the fibration $s$ of
Example~\ref{ar.exa}, and the Segal condition is \eqref{pm.eq}. In
particular, the fiber $B\Delta_{[1]}$ is by definition the category
$\Delta_{\pm}^o$, and this is canonically identified with $\Delta^<$
by \eqref{p-m-g.eq}. We also have the identity section $\Delta \to
\Ar^\pm(\Delta)$ of the projection $s$ sending $[n]$ to $\id:[n] \to
   [n]$, and the opposite functor $\eta:\Delta^o \to B\Delta$ is a
   lax monoidal functor $\ppt \to \Delta$ corresponding to the
   algebra object $1 \in \Delta^<$.
\end{exa}

\begin{exa}\label{B1.exa}
Let $\Ar^p(\Delta)^o \subset \Ar^\pm(\Delta)^o$ be the full
subcategory spanned by surjective arrows. Then the
projection $s^o:\Ar^p(\Delta)^o \to \Delta^o$ is again a
cofibration. Its fiber $\Ar^p(\Delta)^o_{[1]}$ is natural identified
with the category $[1]$, and $\Ar^p(\Delta)^o \cong B[1]$ defines a
monoidal structure on $[1]$ such that $0 \otimes 0 = 0$ and $0
\otimes 1 = 1 \otimes 1 = 1 \otimes 0 = 1$.
\end{exa}

The lax monoidal functor $\eta$ of Example~\ref{del.ten.exa} is
universal in the following sense. Denote by $s,t:B\Delta =
\Ar^\pm(\Delta)^o \to \Delta^o$ the functors sending a bispecial arrow
$[n] \to [m]$ in $\Delta$ to its source $[n]$ resp.\ its target
$[m]$, and for any $2$-category $\C$, let $P(\C)$ be the product
\begin{equation}\label{path.eq}
P(\C) = \C \times^t_{\Delta^o} B\Delta.
\end{equation}
Alternatively, $P(\C) \subset \Ar(\C)$ is the full subcategory
spanned by cocartesian lifting $c \to c'$ of bispecial arrows in
$\Delta^o$, and we have functors
\begin{equation}\label{st.P}
s,t:P(\C) \to \C
\end{equation}
sending $c \to c'$ to $c'$ resp.\ $c$, and their common section
\begin{equation}\label{eta.P}
\eta:\C \to P(\C)
\end{equation}
sending $c$ to $\id:c \to c$ that is left-adjoint to $t$ and
right-adjoint to $s$. In terms of \eqref{path.eq}, the functor
\eqref{eta.P} is the product of $\id:\C \to \C$ and the lax
$2$-functor $\eta$ of Example~\ref{del.ten.exa}. Now, the projection
$s$ of \eqref{st.P} is a cofibration that turns $P(\C)$ into a
$2$-category, called the {\em path $2$-category} of the $2$-category
$\C$. Then \eqref{eta.P} is a lax $2$-functor, and any lax
$2$-functor $\gamma:\C \to \C'$ uniquely factors as
\begin{equation}\label{P.ga}
\gamma \cong P(\gamma) \circ \eta,
\end{equation}
for a unique $2$-functor $P(\gamma):P(\C) \to \C'$. Explicitly, we
have $P(\gamma) = \eta_!^{\Delta^o}\gamma$, where the relative Kan
extension exists by Example~\ref{adj.I.exa} and is given by
\eqref{adj.I.eq}.

\begin{remark}
Informally, the path $2$-category $P(\C)$ has the same objects as
$\C$, and morphisms in $P(\C)$ are free paths generated by morphisms
in $\C$; this motivates the terminology. Note that if $I$ is a small
category, and $P(I)$ is the path category of the unverlying quiver,
then $P(\Delta^o I) \cong \Delta^oP(I)$.
\end{remark}

For any unital monoidal category $\C$, the opposite category $\C^o$
is also unital monoidal; the corresponding cofibration $B\C^o \to
\Delta^o$ is given by $B\C^o = (B\C)^o_\perp$. In particular, this
applies to $\Delta^{<o} \cong \Delta_{\pm}$. In terms of
$\Delta_{\pm}$, the product is given by the {\em reduced
  concatenation} $[m] * [n] = [m] \copr_{[0]} [n]$, where the
coproduct is taken with respect to the embeddings $t:[0] \to [m]$,
$s:[0] \to [n]$, as in \eqref{seg.sq}. We have the tautological map
$[m] \circ [n] \to [m] * [n]$ from the non-reduced to the reduced
concatenation, and this turns the forgetful functor
$\rho_\flat:\Delta_{\pm} \to \Delta \subset \Delta^<$ into a lax
monoidal functor. In terms of Definition~\ref{mon.def},
$B\Delta_{\pm}$ is the full subcategory $\Tw^\pm(\Delta) \subset
\Tw(\Delta)$ spanned by bispecial arrows, and the functor
\begin{equation}\label{bphi}
B\rho_\flat:B\Delta_{\pm} = \Tw^\pm(\Delta) \to B\Delta^< =
\Ar^\pm(\Delta)^o
\end{equation}
is an embedding both on objects and on morphisms, and sends an arrow
$f:[n] \to [m]$ to an injective arrow $[n] \to [m+n]$. More
precisely, let $\Ar^c(\Delta) \subset \Ar^\pm(\Delta)$ be the
subcategory in $\Ar^\pm(\Delta)$ spanned by injective bispecial
arrows $f:[n] \to [m]$, and maps between them such that the
corresponding square \eqref{ar.dia} is not only commutative but also
cartesian. Then $B\rho_\flat$ induced an equivalence between
$\Tw^\pm(\Delta) \cong B\Delta_{\pm}$ and $\Ar^c(\Delta)^o$.

\subsection{Wreath products.}\label{wr.subs}

As mentioned above, lax $2$-functors from $\Delta^o = \ppt^2$ to
$B\C$ for a monoidal category $\C$ correspond to algebra objects in
$\C$ (this can be proved but we prefer to use it as a
definition). It turns out that the $2$-category
$\Ar^p(\Delta)^o=B[1]$ of Example~\ref{B1.exa} has a similar
universal property: lax $2$-functors $B[1] \to B\C$ classify
morphisms of algebra objects in $\C$.

Namely, for any $2$-category $\C$, a pair of lax $2$-functors
$\gamma_0,\gamma_1:\ppt^2 \to \C$ equipped with a map $\gamma_0 \to
\gamma_1$ give rise to a functor $\gamma:\Delta^o \times [1] \to \C$
over $\Delta^o$ that is cocartesian over anchor maps. The
cofibration $\pi:B[1] \to \Delta^o$ has a left and a right-adjoint
$l,r:\Delta^o \to B[1]$ sending $[n]$ to the arrow $[n] \to [0]$
resp.\ to the arrow $[n] \to [n]$. Both $l$ and $r$ are fully
faithful, and together with the map $l \to r$ adjoint to the
isomorphism $\Id \cong \pi \circ r$ define a functor
\begin{equation}\label{w.eq}
w:\Delta^o \times [1] \to B[1]
\end{equation}
over $\Delta^o$. It is cocartesian over anchor maps, and enjoys the
following universal property.

\begin{lemma}\label{wr.le}
Assume given a $2$-category $\C$, and a functor $\gamma:\Delta^o
\times [1] \to \C$ over $\Delta^o$ that is cocartesian over anchor
maps. Then the right Kan extension $w_*\gamma:B[1] \to \C$ with
respect to \eqref{w.eq} exists and defines a lax $2$-functor, and
the adjunction map $w^*w_*\gamma \to \gamma$ is an
isomorphism. Conversely, for any lax $2$-functor $\gamma':B[1] \to
\C$, the adjunction map $\gamma' \to w_*w^*\gamma'$ is an
isomorphism.
\end{lemma}

\proof{} To compute $w_*$, let us choose a convenient framing of the
opposite functor $w^o:\Delta \times [1]^o \to \Ar^p(\Delta)$. For
any object in $\Ar^p(\Delta)$ represented by a surjective arrow
$g:[n] \to [m]$, consider the corresponding diagram \eqref{m.mb},
and let $v(g):V([m]) \to \Ar^p(\Delta)$ be the functor sending $v
\in V([m])$ to the arrow $[n_v] \to [0]$ (that is, to
$l^o([n_v])$). Then $v(g)$ has a natural augmentation
$v(g)^>:V([m])^> \to \Ar^p(\Delta)$ sending $o$ to $g$, with the
maps $l^o([n_v]) \to g$ induced by the map $e_g$ in
\eqref{m.mb}. Extend $v(g)^>$ to a functor $j(g)^>:[1] \times v(g)^>
\to \Ar^p(\Delta)$ equal to $v(g)^>$ on $1 \times V(g)^>$ and to
$r^o \circ \pi^o \circ v(g)^>$ on $0 \times v(g)^>$, with the
adjunction map $r^o \circ \pi^o \circ v(g)^> \to v(g)^>$, and note
that $[1] \times v(g)^>$ has the largest element $1 \times o$, so
that $[1] \times v(g)^> = J(g)^>$ for the partially ordered set
$J(g) = ([1] \times v(g)^>) \setminus \{1 \times o\}$, and $j(g)^>$
is an augmentation of a functor $j(g):J(g) \to \Ar^p(\Delta)$. This
functor $j(g)$ canonically factors through $w^o$, and the induced
functor $J(g) \to \Ar^p(\Delta)/^{w^o}f$ is a left-admissible full
embedding.

The collection $J(g)$ is our framing, and then for any
$\gamma:\Delta^o \times [1] \to \C$ with components
$\gamma_0,\gamma_1:\Delta^o \to \C$, the expected object
$w_*\gamma(g) \cong \lim_{J(g)}\gamma$ fits into a cartesian square
\begin{equation}\label{w.sq}
\begin{CD}
w_*\gamma(g) @>>> \gamma_1([n])\\
@VVV @VV{e_{g!}^o}V\\
\prod_v \gamma_0([n_v]) @>>> \prod_v \gamma_1([n_v])
\end{CD}
\end{equation}
in the category $\C$. This can be reinterpreted as opposite to an
iterated colimit of half-standard functors $\V^o \to \C$ of
Remark~\ref{st.rem}, as in Lemma~\ref{2.cyl.le}, so that if
$\gamma_0$ and $\gamma_1$ are lax $2$-functors, the limit exists,
and then so does the Kan extension $w_*\gamma$. The equivalence
$w^*w_*\gamma \cong \gamma$ is then obvious from the cartesian
square \eqref{w.sq}, and since any lax $2$-functor
$\gamma':B[1]=\Ar^p(\Delta)^o \to \C$ preserves the limits of
standard functors, we also have $\gamma' \cong w_*w^*\gamma'$.
\endproof

\begin{remark}
In the situation of Lemma~\ref{wr.le}, one can also separate
$\gamma$ into lax $2$-functors $\gamma_0,\gamma_1:\Delta^o \to \C$
equipped with a map $g:\gamma_0 \to \gamma_1$, and consider the
$2$-cylinder $\Cyl^2(\id,g):\Delta^o[1] \to \C$ provided by
Lemma~\ref{2.cyl.le}. The difference between $\Cyl^2(\id,g)$ and the
right Kan extension $w_*\gamma$ of Lemma~\ref{wr.le} is that
$\Delta^o[1]$ does not correspond to a monoidal category: it has two
objects. Since $\Cyl^2(\id,g)$ and $w_*\gamma$ enjoy essentially the
same universal property, we have $e^*w_*\gamma \cong \Cyl^2(\id,g)$,
where $e:\Delta^o[1] \to B[1]$ is the $2$-cylinder of the map $l \to
r$ used to define $w$, but of course $e$ is not an equivalence.
\end{remark}

For a useful application of the same combinatorics to $2$-categories
rather then lax $2$-functors, assume given a $2$-functor
$\gamma:\C_0 \to \C_1$ between some $2$-categories $\C_0$, $\C_1$,
and note that the cylinder $\Cyl(\gamma)$ is then a cofibration over
$\Delta^o \times [1]$ that restricts to $\C_l$ over $\Delta^o \times
l$, $l=0,1$.

\begin{defn}\label{wr.def}
The {\em wreath product} $\C_0 \wr^\gamma \C_1$ of the
$2$-categories $\C_0$, $\C_1$ with respect to the $2$-functor
$\gamma$ is the cofibration
$$
\C_0 \wr^\gamma \C_1 = w_*\Cyl(\gamma)
$$
over $B[1]=\Ar^p(\Delta)^o$, where $w$ is the functor \eqref{w.eq}.
\end{defn}

Explicitly, the wreath product $\C_0 \wr^\gamma \C_1$ can be
computed by the same framing that gives \eqref{w.sq}. This shows
that the fiber $(\C_0 \wr^\gamma \C_1)_g$ over some surjective
$g:[n] \to [m]$ is the category of triples $\langle c_0,c_1,\alpha
\rangle$ of a $2$-functor $c_1:[n] \to C_1$, a $2$-functor
$c_0:[n]_g \to \C_0$, and an isomorphism $\alpha:e_g^*c_1 \cong
c_0$. In particular, $\C_0 \wr^\gamma \C_1$ is a $2$-category, and
we have $\C_0 \wr^\gamma \C_1 \cong \C_0 \wr^{\wgamma}
\gamma^*\C_1$, where $\wgamma$ is the dense component of the
decomposition \eqref{d.f.facto}. We will simplify notation by
writing $\C_0 \wr^\gamma \C_1 = \C_0 \wr \C_1$ when $\gamma$ is
clear from the context.

\begin{exa}\label{wr.exa}
If $\C_1=\ppt^2$ is the point $2$-category, then the only
non-trivial part of a triple $\langle c_0,c_1,\alpha \rangle$
describing an object in $\C_0 \wr \ppt^2$ is $c_0$. Then the
projection $t:\Ar^p(\Delta) \to \Delta$ of \eqref{ar.tw} induces a
projection $\C_0 \wr \ppt^2 \to \Delta^o$ whose fibers are given by
\begin{equation}\label{wr.eq}
(\C_0 \wr \ppt^2)_{[m]} \cong \C_0^{m+1}, \qquad [m] \in \Delta^o.
\end{equation}
This motivates our terminology. If moreover $\C_0 = \Delta^oe(S)$ is
the simplicial replacement of a category $e(S)$, as in
Example~\ref{eS.exa}, then we have
\begin{equation}\label{e.wr}
\Delta^oe(S) \wr \ppt^2 \cong \Delta^oe(S) \times^2 B[1].
\end{equation}
Slightly more generally, for any $2$-category $\C$ and functor
$S:\C_{[0]} \to \Sets$, as in \eqref{C.S.eq}, we have the
$2$-functor $\pi:\C[S] \to \C$, and \eqref{e.wr} induces an
identification
\begin{equation}\label{wr.pi}
\C[S] \wr^\pi \C \cong \C \times^2 B[1].
\end{equation}
On the other hand, if $\C_1$ is arbitrary but $\C_0$ is discrete, we
have
\begin{equation}\label{discr.wr}
\C_0 \wr^\gamma \C_1 \cong \gamma^*\C_1 \times_{\Delta^o}^{t^o} B[1],
\end{equation}
where $t^o:B[1] = \Ar^p(\Delta)^o \to \Delta^o$ is induced by
\eqref{ar.tw}, and $\gamma^*\C_1$ is as in \eqref{2.ga}.
\end{exa}

The wreath product construction is obviously functorial with respect
to $2$-functors: if we are given another $2$-functor $\gamma':\C_0'
\to \C_1'$, and $2$-functors $\phi_0:\C_0 \to \C_0'$, $\phi_1:\C_1
\to \C_1'$ equipped with an isomorphism $\alpha:\gamma' \circ \phi_0
\cong \phi_1 \circ \gamma$, then we have a functor
$\Cyl(\phi):\Cyl(\gamma) \to \Cyl(\gamma')$ cocartesian over
$\Delta^o \times [1]$, and it induces a $2$-functor
\begin{equation}\label{phi.wr}
\phi_0 \wr^\alpha \phi_1 = w_*\Cyl(\phi):\C_0 \wr^\gamma \C_1 \to
\C_0 \wr^{\gamma'} \C_1,
\end{equation}
cocartesian over $B[1]$, where we will again drop $\alpha$ from
notation when it is clear from the context. To extend this to lax
$2$-functors, consider the universal situation: take the path
$2$-categories $P(\C_0)$, $P(\C_1)$, with the functor $P(\C_0) \to
P(\C_1)$ induced by $\gamma$. Then $s$ of \eqref{st.P} induces a
$2$-functor
\begin{equation}\label{s.wr}
s \wr s:P(\C_0) \wr P(\C_1) \to \C_0 \wr \C_1,
\end{equation}
again cocartesian over $B[1]$.

\begin{lemma}\label{lax.le}
The functor \eqref{s.wr} admits a fully faithful right-adjoint lax
$2$-functor $\eta \wr \eta:\C_0 \wr \C_1 \to P(\C_0) \wr P(\C_1)$
over $B[1]$.
\end{lemma}

\proof{} By Lemma~\ref{loc.le}~\thetag{ii}, it suffices to check that
for any $f \in B[1]$, the fiber $(s \wr s)_f:(P(\C_0) \wr P(\C_1))_f
\to (\C_0 \wr \C_1)_f$ of the cocartesian functor \eqref{s.wr}
admits a right-adjoint. But by \eqref{w.sq}, its source resp.\ its
target is a product of categories of the form $P(\C_0)_{[n]}$,
$P(\C_1)_{[m]}$ resp.\ $(\C_0)_{[n]}$, $(\C_1)_{[m]}$ for various
$[n],[m] \in \Delta^o$, and $(s \wr s)_f$ is the
product of the fibers $s_{[n]}$ of the functors \eqref{st.P}. These
have fully faithful right-adjoints induced by \eqref{eta.P}.
\endproof

Now if we have $2$-functors $\gamma:\C_0 \to \C_1$, $\gamma':\C_0'
\to \C_1'$, and lax $2$-functors $\phi_0:\C_0 \to \C_0'$,
$\phi_1:\C_1 \to \C_1'$ equipped with an isomorphism $\gamma' \circ
\phi_0 \cong \phi_1 \circ \gamma$, we can define $\phi_0 \wr \phi_1$
by
\begin{equation}\label{phi.lax.wr}
\phi_0 \wr \phi_1 = (P(\phi_0) \wr P(\phi_1)) \circ (\eta \wr
\eta):\C_0 \wr \C_1 \to \C_0 \wr \C_1,
\end{equation}
where $P(\phi_0)$, $P(\phi_1)$ are as in \eqref{P.ga}, and $\eta
\wr \eta$ is provided by Lemma~\ref{lax.le}. This is a functor
over $B[1]$. If $\phi_0$, $\phi_1$ are actual $2$-functors, then we
have $(s \wr s) \circ (\eta \wr \eta) \cong \Id$ since $\eta \wr
\eta$ is fully faithful, and \eqref{phi.lax.wr} agrees with
\eqref{phi.wr}. In addition to that, if we have a lax $2$-functor
$\gamma:\C_0 \to \C_1$, we can define the wreath product $\C_0
\wr^\gamma \C_1$ by the cartesian square
\begin{equation}\label{wr.sq}
\begin{CD}
\C_0 \wr^{\gamma} \C_1 @>>> P(\C_0) \wr^{P(\gamma)} \C_1\\
@VVV @VV{\id \wr \tau}V\\
\C_0 \wr \ppt^2 @>{\eta \wr \id}>> P(\C_0) \wr \ppt^2,
\end{CD}
\end{equation}
where $\tau:\C_1 \to \ppt^2$ is the tautological projection. This
again agrees with Definition~\ref{wr.def} when $\gamma$ is a
$2$-functor. In all cases, we have a natural projection
\begin{equation}\label{wr.C}
\begin{CD}
\C_0 \wr \C_1 @>{\gamma \wr \id}>> \C_1 \wr^{\id} \C_1 \cong \C_1
\times^2 B[1] @>>> \C_1,
\end{CD}
\end{equation}
where the identification in the middle is \eqref{wr.pi}. It is a lax
$2$-functor, and a $2$-functor if so is $\gamma$.

\subsection{Modules.}\label{C.mod.subs}

Now assume given a category $\C$ equipped with a unital monoidal
structure $B\C$.

\begin{defn}\label{C.mod.def}
A {\em module} over $\C$ is a category $\M$ equipped with a functor
$\mu:\M \to B\C$ such that the composition $\M \to B\C \to \Delta^o$
is a cofibration, with the induced cofibration $\rho^{o*}\M \to
\Delta_+^o$, and the functor
\begin{equation}\label{rho.M}
\zeta([0]) \times \rho^{o*}\mu:\rho^{o}*\M \to M \times
\rho^{o*}B\C, \qquad M = \M_{[0]}
\end{equation}
is an equivalence of categories. A {\em morphism} between two
$\C$-modules $\M$, $\M'$ is a functor $\alpha:\M' \to \M$ over $B\C$
cocartesian over all special maps in $\Delta^o$.
\end{defn}

Explicitly, \eqref{rho.M} provides identifications $\M_{[n]} \cong M
\times B\C_{[n]} \cong M \times \C^n$, $[n] \in \Delta$, and
prescribes the transition functors $f^o_!:\M_{[n]} \to \M_{[n']}$
for all special maps $f:[n'] \to [n]$ in $\Delta$. The essential
part of the structure is the functor
\begin{equation}\label{rho.m}
m = t^o_!:M \times \C \to M
\end{equation}
corresponding to the antispecial map $t:[0] \to [1]$. This functor
defines an action of $\C$ on $M$, and the rest of the structure
encodes the usual associativity and unitality constraints for this
action. A morphism $\alpha$ is then given by a functor
$\alpha_{[0]}:M' = \M'_{[0]} \to \M$ and a map
\begin{equation}\label{fi.fu.mod}
\alpha_{[0]} \circ m' \to m \circ (\alpha_{[0]} \times \id).
\end{equation}
In particular, $\alpha$ can be non-trivial even if $\alpha_{[0]} =
\id$.

\begin{exa}
For any unital monoidal category $\C$, $\kappa^{o*}\C$ with the
projection $a^o_!:\kappa^{o*}\C \to \C$ is a $\C$-module; the
corresponding action \eqref{rho.m} is the action of $\C$ on itself
by left products. For any category $\E$, the product $\E \times B\C$
is a (trivial) $\C$-module; the corresponding action \eqref{rho.m}
is the projection onto the first factor. For any $\C$-module $\M$
with the action map \eqref{rho.m}, and unital monoidal functor
$\gamma:\C' \to \C$ from a unital monoidal category $\C$', the
pullback $\gamma^*\M \to B\C'$ is a $\C'$-module, with the action
map $m \circ (\id \times \gamma)$.
\end{exa}

For any $\C$-module $\M$, the equivalence \eqref{rho.M} immediately
implies that $\rho^{o*}\mu$ is both a fibration and a cofibration,
but $\mu$ itself is neither: in general, it is only a
precofibration. To understand the structure of the category $\M$
better, it is useful to consider the category $\Delta_\idot$ of
\eqref{nu.dot}, with the functor $\nu_\dg:\Delta \to \Delta_\idot$
and its two adjoints $\nu_\idot,\nu_\perp:\Delta_\idot \to
\Delta$. If we denote $B^\hdot\C = \nu^{o*}_\idot B\C$, with the
induced fibration $\nu_\idot:B^\hdot\C \to B\C$, then $\nu_\dg$ and
$\nu_\perp$ induce functors $\nu_\dg:B\C \to B^\hdot\C$,
$\nu_\perp:B^\hdot\C \to B\C$, and we can consider the
precofibration
\begin{equation}\label{M.dot}
\mu_\idot:\M^\hdot = \M \times^{\nu_\perp}_{B\C} B^\hdot\C \to B^\hdot\C.
\end{equation}
Then the composition $\nu_\idot \circ \mu_\idot:\M^\hdot \to B\C$ is
actually a fibration. The fiber of this fibration over some $\langle
[n],c_\idot \rangle \in B\C$, $c_\idot \in \C^n$ is the category
$M(\langle [n],c_\idot \rangle)$ cofibered over $[n]^o$, with all
fibers identified with $M$, and transition functors $m(- \times
c_i):M \to M$, $i = 1,\dots,n$. For any map $f:[n'] \to [n]$ in
$\Delta$, we have $f^*M(\langle [n],c_\idot \rangle) \cong
\M(\langle [n'],f^*c_\idot \rangle$, so that $M(\langle [n],c
\rangle)$ is covariantly functorial with respect to $[n] \in
\Delta$, and \eqref{cyl.eq} immediately shows that it is also
contravariantly functorial with respect to $c_\idot$.

Assume given a $\C$-module $\M$ with bounded $M = \M_{[0]}$, and
some other category $\E$. Then the action functor \eqref{rho.m}
induces a functor
\begin{equation}\label{rho.m.E}
m:\Fun(M,\E) \times \C \to \Fun(M,\E), \qquad m(F \times c)(c') =
F(c \times c'),
\end{equation}
and it turns out that this extends to a module over the opposite
monoidal category $\C^\iota$. To construct it, consider the
precofibration \eqref{M.dot}, take the corresponding precofibration
$\mu(\E):\Fun(\M^\hdot/B^\hdot\C,\E) \to B^\hdot\C$, and note that
the composition $\nu_\idot \circ
\mu(\E):\Fun(\M^\hdot/B^\hdot\C,\E) \to B\C$ is a cofibration, with
fibers $\Fun(M(\langle [n],c \rangle)/[n],\E)$. Therefore we can
consider the transpose fibration $\Fun(\M/B^\hdot\C,E)^\perp \to
B\C$, and $\mu(\E)$ induces a functor
\begin{equation}\label{m.dot.E}
\iota^*\mu(\E)^\perp:\iota^*\Fun(\M/B^\hdot\C,\E)^\perp \to
\iota^*(B^\hdot\C)_\perp \cong B^\hdot\C^\iota,
\end{equation}
where $(B^\hdot\C)_\perp$ is taken with respect to the fibration
$\nu_\idot:B^\hdot \to \C$. Then
\begin{equation}\label{ofun.eq}
\oFun(\M,\E) = \nu_\dg^*\iota^*\Fun(\M/B^\hdot\C,\E)^\perp,
\end{equation}
with the induced functor $\oFun(\M,\E) \to \nu_\dg^*B^\hdot\C^\iota
\cong B\C^\iota$, is a $\C^\iota$-module, with $\oFun(\M,\E)_{[0]}
\cong \Fun(M,\E)$ and the action \eqref{rho.m.E}, while
\eqref{m.dot.E} is identified with the corresponding precofibration
$\oFun(\M,E)^\hdot \to B^\hdot\C^\iota$ of \eqref{M.dot}.

As an application of this construction, let $\Gamma_+$ be the
category of pointed finite sets, with the unital monoidal structure
given by smash product, let it act on itself, and let $\E$ be a
half-additive category in the sense of
Subsection~\ref{adv.stab.subs}. Then the functor \eqref{ga.m}
induced a fully faithful embedding $m^\dg:\E \to \Fun(\Gamma_+,\E)$,
and for any $[n] \in \Delta$, we have a fully faithful embedding
$m^\dg \times \id:\E \times \Gamma_+^n \to \oFun(\Gamma_+,\E)_{[n]}
\cong \Fun(\Gamma_+,\E) \times \Gamma_+^n$. If we let $\E^\otimes
\subset \oFun(\Gamma_+,\E)$ be the full subcategory spanned by the
essential images of these embedding, with the induced projection
$\E^\otimes \to B\Gamma_+$, then by Lemma~\ref{fib.le}~\thetag{i},
$\E^\otimes$ becomes a module over $\Gamma_+ \cong \Gamma_+^\iota$
in the sense of Definition~\ref{C.mod.def}. The corresponding action
functor \eqref{rho.m} is the functor \eqref{ga.m}, and the embedding
$\E^\otimes \to \oFun(\Gamma_+,\E)$ is a morphism of
$\Gamma_+$-modules. Moreover, assume that $\E$ {\em has kernels}, in
the sense that for any map $f:e' \to e$, there exists $\Ker(f) = o
\times_e e'$, where $o \in \E$ is the initial terminal object. Then
essentially as in Lemma~\ref{poly.st.le}, $m^\dg$ admits a
right-adjoint $\Fun(\Gamma_+,\E) \to \E$ sending $E:\Gamma_+ \to \E$
to the kernel of the map $E(\ppt_+) \to E(o)$, and then the
right-adjoint
\begin{equation}\label{I.ot}
\oFun(\Gamma_+,\E) \to \E^\otimes
\end{equation}
to the embedding $\E^\otimes \to \oFun(\Gamma_+,\E)$ provided by
Lemma~\ref{fib.le}~\thetag{i} is also a morphism of
$\Gamma_+$-modules.

Note that the category $\Fun(\Gamma_+,\E)$ is also half-additive, so
that effectively, it has two structures of a $\Gamma_+$-module:
$\oFun(\Gamma_+,\E)$ on one hand, and $\bFun(\Gamma_+,\E) =
\Fun(\Gamma_+,\E)^\otimes$ of \eqref{I.ot} on the other hand
(informally, $\Gamma_+$ can act on $\Fun(\Gamma_+,\E)$ either via
$\Gamma_+$ for via $\E$). To relate the two, let $\Gamma_+$ act on
$\Gamma_+ \times \Gamma_+$ via the left factor. Then the product
functor $m:\Gamma_+ \times \Gamma_+ \to \Gamma_+$ gives rise to a
morphism $\mu$ of $\Gamma_+$-modules, and we have the morphism
\begin{equation}\label{b.o.fun}
\begin{CD}
\oFun(\Gamma_+,\E) @>{m^*}>> \oFun(\Gamma_+ \times \Gamma_+,\E)
 @>>> \bFun(\Gamma_+,\E),
\end{CD}
\end{equation}
where the first arrow is induced by $m$, and the second one is the
morphism \eqref{I.ot} with the identification $\oFun(\Gamma_+ \times
\Gamma_+,\E) \cong \oFun(\Gamma_+,\Fun(\Gamma_+,\E))$.

More generally, if in addition to $\E$ we also have a bounded
category $I$, then the functor category $\Fun(I,\Gamma_+)$
equipped with the pointwise product of Example~\ref{point.exa} acts
on $\Gamma_+ \times I$ by
$$
\Gamma_+ \times I \times \Fun(I,\Gamma_+) \to \Gamma_+ \times I,
\qquad S_+ \times i \times F \mapsto S_+ \wedge F(i) \times i,
$$
and the same construction equips the category $\Fun(I \times \Gamma_+,\E)$
with two structures of a module over $\Fun(I,\Gamma_+)$, while
\eqref{b.o.fun} induced a morphism
\begin{equation}\label{b.o.I.fun}
\oFun(I \times \Gamma_+,\E) \to \bFun(I \times \Gamma_+,\E)
\end{equation}
between the corresponding $\Fun(I,\Gamma_+)$-modules.

\section{Adjunction.}

\subsection{Internal adjunction.}\label{int.subs}

For any $2$-category $\C$ and two objects $c,c' \in \C_{[0]}$, an
{\em adjoint pair} of maps between $c$ and $c'$ is a quadruple
$\langle f,f^\vee,a,a^\vee \rangle$, $f \in \C(c,c')$, $f^\vee \in
\C(c',c)$, $a:\id_c \to f^\vee \circ f$, $a^\vee:f \circ f^\vee \to
\id_{c'}$, subject to the usual relations
\begin{equation}\label{adj.rel}
(a^\vee \circ \id_f) \circ (\id_f \circ a) = \id_f, \quad (\id_{f^\vee}
\circ a^\vee) \circ (a \circ \id_{f^\vee}) = \id_{f^\vee}.
\end{equation}
For any $c,c' \in \C_{[0]}$, adjoint pairs of maps between $c$ and
$c'$ and isomorphisms between them form a category that we denote
$\aAdj(\C)(c,c')$, and the forgetful functor from $\aAdj(\C)(c,c')$
to $\C(c,c')_{\Iso}$ sending $\langle f,f^\vee,a,a^\vee \rangle$ to
$f$ is fully faithful (for the $2$-category of small categories,
this is the standard uniqueness of adjoints, and the same proof
works for a general $2$-category $\C$). We say that $f \in \C(c,c')$
is {\em reflexive} if it extends to an adjoint pair, and we note
that reflexivity is closed under compositions (again, the proof for
small categories works in the general case, or see below in
Subsection~\ref{ite.subs}). Thus we actually have the {\em
  adjunction $2$-category} $\aAdj(\C)$ with the same objects as $\C$
and $\aAdj(-,-)$ as categories of morphisms, and we have a fully
faithful $2$-functor
\begin{equation}\label{adj.c}
  \aAdj(\C) \to \C_\natural
\end{equation}
that is an identity over $[0]$. If $\aAdj(\C)$ is bounded, we denote
its truncation by $\rAdj(\C) = \tau(\aAdj(\C))$. A reflexive morphism
$f \in \C(c,c')$ is an {\em equivalence} if there exists an adjoint
pair $\langle f,f^\vee,a,a^\vee \rangle$ with invertible $a$ and
$a^\vee$.

To understand adjunction in the general $2$-categorical context, it
is useful to look at the universal situation. Let $\adj$ be the
$2$-category with two objects $0$, $1$, and categories of morphisms
\begin{equation}\label{2adj.01}
\adj(0,0) = \Delta^<, \ \adj(0,1)=\Delta_-, \
\adj(1,0)=\Delta_+, \ \adj(1,1) = \Delta_{\pm},
\end{equation}
with compositions
$$
\Delta^< \times \Delta_- \to \Delta_-, \qquad \Delta_+ \times
\Delta^< \to \Delta_+, \qquad \Delta^< \times \Delta^< \to \Delta^<
$$
given by the concatenation product $- \circ -$, and compositions
$$
\Delta_{\pm} \times \Delta_+ \to \Delta_+, \qquad \Delta_- \times
\Delta_{\pm} \to \Delta_-, \qquad \Delta_{\pm} \times \Delta_{\pm} \to
\Delta_{\pm}
$$ given by the reduced concatenation product $- * -$. Then we have
two morphisms $f:0 \to 1$, $f^\vee:1 \to 0$ in $\adj$ corresponding
to the initial objects in $\Delta_+$, $\Delta_-$, and $f^\vee \circ
f \cong [1] \in \Delta^<$, $f \circ f^\vee \cong [1] \in \Delta^{<o}
\cong \Delta_{\pm}$, so that the map $[0] \to [1]$ produces maps
$a:\id_0 \to f^\vee \circ f$, $a^\vee:f \circ f^\vee \to \id_1$. It
is elementary to check that these maps satisfy \eqref{adj.rel}, so
that $\langle f,f^\vee,a,a^\vee \rangle$ is an adjoint pair, and
moreover, it is universal with with property --- for any adjoint
pair $\langle f_1,f_1^\vee,a_1,a_1^\vee \rangle$ of maps between
objects $c$, $c'$ in a $2$-category $\C$, there exists a $2$-functor
$\gamma:\adj \to \C$ sending $0$ to $c$, $1$ to $c'$, $\langle
f,f^\vee,a,a^\vee \rangle$ to $\langle f_1,f_1^\vee,a_1,a_1^\vee
\rangle$, and this $\gamma$ is unique up to a unique
isomorphism. Altogether, morphisms in a $2$-category $\C$ correspond
to $2$-functors from $[1]$ to $\C$, adjoint pairs of morphisms
correspond to $2$-functors $\adj \to \C$, and equivalences
correspond to $2$-functors from $e(\{0,1\})$ to $\C$, where as in
Subsection~\ref{set.sss}, $e(\{0,1\})$ is the category with two
objects $0$, $1$ and exactly one map between any two objects. For
consistency, denote $\nat=\Delta^o[1]$, $\eq=\Delta^oe(\{0,1\})$;
then we have $2$-functors
\begin{equation}\label{1.adj.eq}
\begin{CD}
\nat @>{\delta}>> \adj @>{\nu}>> \eq,
\end{CD}
\end{equation}
where $\nu:\adj \to E\adj_{[0]} \cong \Delta^oe(\{0,1\}) = \eq$ is
as in Example~\ref{E.exa}, and $\delta$ corresponds to the morphism
$f:0 \to 1$ in $\adj$. A morphism is reflexive if and only if it
factors through $\delta$, and an adjoint pair is an equivalence if
and only if it factors through $\nu$. In both cases, the
factorizations are automatically unique.

One can describe the whole cofibration $\adj \to \Delta^o$ of
Definition~\ref{2cat.def} rather explicitly using the packaging of
the monoidal structures on $\Delta^<$ and $\Delta_{\pm}$ given in
Subsection~\ref{mon.subs}. Namely, note that the category $\eq^o =
\Delta e(\{0,1\})$ can be alternatively described as the subcategory
in $\Ar(\Delta^<)$ spanned by injective arrows $a:[n]_1 \to [n]$
with non-empty $[n]$, and with morphisms given by commutative
squares \eqref{ar.dia} that are cartesian (indeed, the complement
$[n]_0 = [n] \setminus [n]_1$ is functorial with respect to such
morphisms, and then an arrow $a:[n]_1 \to [n]$ corresponds to
$\langle [n],e_\idot \rangle \in \Delta e(\{0,1\})$, where
$e_\idot:[n] \to e(\{0,1\})$ sends $[n]_1$ to $1$ and $[n]_0$ to
$0$). Then the category $\adj^o$ opposite to $\adj$ is the category
of cartesian squares
\begin{equation}\label{2adj.dia}
\begin{CD}
[n]_1 @= [m]_1\\
@V{a}VV @VV{b \circ a}V\\
[n] @>{b}>> [m]
\end{CD}
\end{equation}
in $\Delta^<$ such that $[n]$ and $[m]$ are non-empty, $b$ is
bispecial, and $a$ and $b \circ a$ are injective. Morphisms are
given by commutative diagrams
\begin{equation}\label{2.2adj.dia}
\begin{CD}
[n]_1 @>{a}>> [n] @>{b}>> [m]\\
@V{f_1}VV @V{f}VV @VV{g}V\\
[n']_1 @>{a}>> [n'] @>{b}>> [m']
\end{CD}
\end{equation}
such that the leftmost square and the outer rectangle are
cartesian. The $2$-functor $\nu:\adj \to \eq$ of \eqref{1.adj.eq}
sends a diagram \eqref{2adj.dia} to its leftmost vertical arrow
$a:[n]_1 \to [n]$, and then forgetting both $a$ and $[n]_1$ gives
the structural cofibration $\adj \to \Delta^o$.

The correspondence works as follows. The $2$-functor $[n] \to \adj$
represented by a diagram \eqref{2adj.dia} sends $l \in [n]$ to $1$
if $l \in [n]_1 \subset [n]$ and to $0$ if $l \in [n]_0 = [n]
\setminus [n]_1$ (since the leftmost square in \eqref{2.2adj.dia} is
required to be cartesian, both $[n]_1$ and $[n_0]$ are functorial
with respect to the maps \eqref{2.2adj.dia}). Diagrams with empty
$[n]_1$ then describe the full sub-$2$-category $B\Delta^< \subset
\adj$ spanned by $0$, and this is exactly the description given in
Example~\ref{del.ten.exa}. At the opposite extreme, the full
sub-$2$-category $B\Delta_{\pm} \subset \adj$ spanned by $1$ is
described by the diagrams with invertible map $a$, and this is the
description in terms of the embedding \eqref{bphi}. General diagrams
interpolate between the two, and the conditions we impose on
\eqref{2adj.dia} and \eqref{2.2adj.dia} insure that the categories
of morphisms in $\adj$ are those of \eqref{2adj.01}.

\begin{remark}\label{adj.rem}
Say that a map of sets $f:S \to S'$ is {\em strict} on a subset $S_0
\subset S$ if $S_0 = f^{-1}(f(S_0)) \subset S$. Then equivalently,
objects in $\adj$ are pairs of a bispecial arrow $b:[n] \to [m]$ in
$\Delta$ and a subset $[n]_1 \subset [n]$ such that $b$ is strict
and injective on $[n]_1$.
\end{remark}

\subsection{External adjunction.}\label{ext.adj.subs}

We can now apply the universal interpretation of adjunction given in
\eqref{1.adj.eq} to maps between $2$-functors. Namely, for any
$2$-categories $\C$, $\C'$, their $2$-product $\C \times^2 \C'$ is a
$2$-category, an object $c \in \C_{[0]}$ defines a $2$-functor
$\eps(c):\ppt^2 \to \C$, and we denote by $\C' \times c \subset \C'
\times^2 \C$ the essential image of the embedding $\id \times
\eps(c)$. A functor
\begin{equation}\label{alp.eq}
\gamma:\C \times^2 \nat \to \C'
\end{equation}
is called a {\em natural transformation} between its restrictions
$\gamma_0,\gamma_1:\C \to \C'$ to $\C \times 0,\C \times 1 \subset
\C \times_{\Delta^o} \nat$. An {\em adjoint pair} of natural
transformations between $\gamma_0$ and $\gamma_1$ is a $2$-functor
\begin{equation}\label{alp1.eq}
\gamma:\C \times^2 \adj \to \C'
\end{equation}
that restricts to $\gamma_0$ resp.\ $\gamma_1$ on $\C \times 0$
resp.\ $\C \times 1$, and an {\em equivalence} is a $2$-functor
\begin{equation}\label{equi.eq}
  \C \times^2 \eq \to \C'.
\end{equation}
A natural transformation is reflexive if it extends to an adjoint
pair, and an equivalence if it extends to an equivalence.

\begin{exa}\label{equi.exa}
Assume given a $2$-category $\C$ and a functor $S:\C_{[0]} \to
\Sets$, and let $\C[S]$ be the corresponding $2$-category of
\eqref{C.S.eq}, with its $2$-functor $\pi:\C[S] \to \C$. Moreover,
assume that $S(c)$ is non-empty for any $c \in \C_{[0]}$, so that
$\pi$ admits a section $\sigma:\C = \C[\ppt] \to \C[S]$. Then $p
\circ \sigma \cong \Id$, and $\sigma \circ \pi:\C[S] \to \C[S]$ is
equivalent to $\Id$ in the sense of \eqref{equi.eq}. To construct an
equivalence \eqref{equi.eq} between $\sigma \circ \pi$ and $\id$,
note that $\C[S] \times^2 \eq \cong \C[S \times \{0,1\}]$, and
consider the map $S \times \{0,1\} \to S$ equal to $\id$ on $S
\times \{0\}$ and to the compositon $S \to \ppt \to S$ on $S \times
\{1\}$.
\end{exa}

As we see from the explicit description of the $2$-category $\adj$
in terms of \eqref{2adj.dia} and \eqref{2.2adj.dia}, we actually
have a full embedding
\begin{equation}\label{adj.P}
\adj \to P(\eq) \cong \Ar^\pm(\eq)
\end{equation}
into the path $2$-category of $\eq = \Delta^oe(\{0,1\})$. The
projection $\nu:\adj \to \eq$ is induced by the cofibration
$s:P(\eq) \to \eq$ of \eqref{st.P}, thus has a right-adjoint
\begin{equation}\label{eta.eq}
\eta:\eq \to \adj \subset P(\eq)
\end{equation}
induced by \eqref{eta.P}. However, $\eta$ in turn has a
right-adjoint $t:\adj \subset P(\eq) \to \eq \to \Delta^o$ induced
by \eqref{st.P}. For any $2$-category $\C$, we denote $\C\{\adj\} =
\C \times^t_{\Delta^o} \adj \subset P(\C(\eq))$, and we note that
the functors \eqref{st.P} and \eqref{eta.P} induce functors
\begin{equation}\label{st.C}
s,t:\C\{\adj\} \to \C\{\eq\}, \qquad \eta:\C\{\eq\} \to \C\{\adj\},
\end{equation}
where we denote $\C\{\eq\} = \C \times^2 \eq = \C[\{0,1\}]$
for consistency. We further note that $s$ is a cofibration whose
composition with the cofibration $\eq \to \Delta^o$ turns
$\C\{\adj\}$ into a $2$-category. We then define a {\em coadjoint
  pair} of functors from $\C$ to some $2$-category $\C'$ as a
$2$-functor
\begin{equation}\label{coadj.eq}
\gamma:\C\{\adj\} \to \C'.
\end{equation}
If $\C$ is discrete --- that is, $\C \cong \C_{[0]} \times \Delta^o$
--- then $\C\{\adj\} \cong \C_{[0]} \times \adj \cong \C \times^2
\adj$, and a coadjoint pair \eqref{coadj.eq} is the same thing as an
adjoint pair \eqref{alp1.eq}. This is useful since coadjoint pairs
are much easier to construct.

Namely, let $\iota_l:\Delta^{o>} \to \eq^>$, $l=0,1$ be the
embedding induced by the embedding $\ppt \to e(\{0,1\})$ onto $l$,
and note that the embeddings $\iota_0$, $\iota_1$ admit left-adjoint
functors
\begin{equation}\label{s.01.eq}
s_0,s_1:\eq^> \to \Delta^{o>}
\end{equation}
sending an injective arrow $a:[n]_1 \to [n]$ to $[n]_0$
resp.\ $[n]_1$ (we have to pass to the augmented categories since
$[n]_0$ or $[n]_1$ might be empty). Generalizing Lemma~\ref{mu.le},
say that a map $f$ in $\eq$ is {\em $l$-special}, $l=0,1$ if
$s_l(f)$ is invertible. For any $2$-category $\C$ equipped with a
$2$-functor $\C \to \eq$, say that a map $f$ in $\C$ is {\em
  $l$-special} if it is a cocartesian lifting on an $l$-special map
in $\eq$, and say that a lax $2$-functor $\C \to \C'$ to some $\C'$
is {\em $l$-special} if it sends $l$-special maps in $\C$ to maps in
$\C'$ cocartesian over $\Delta^o$. Note that for any $\C/\eq$, the
embeddings $\iota_l^*\C^> \to \C^>$, $l=0,1$ admit left-adjoint
functors
\begin{equation}\label{s.01.C.eq}
s_l:\C^> \to \iota_l^*\C^>, \qquad l = 0,1,
\end{equation}
and $s_l$ inverts $l$-special maps. In these terms, for any
$2$-category $\C$, $\C\{\adj\} \subset P(\C\{\eq\})$ is the full
subcategory spanned by $1$-special maps with respect to the
projection $\C\{\eq\} \to \eq$.

\begin{remark}\label{0.anch.rem}
More generally, for any $2$-category $\C$ equipped with a
$2$-functor $\pi:\C \to \eq$, say that a map $f$ in $\C$ is {\em
  $l$-anchor}, $l=0,1$, if it is cocartesian over $\eq$, and
$s_l(\pi(f))$ is an anchor map. Then any $l$-special lax $2$-functor
$\gamma:\C \to \C'$ also sends $l$-anchor maps to maps cocartesian
over $\Delta^{o>}$, and so does $s_l$.
\end{remark}

\begin{lemma}\label{adj.le}
For any $2$-category $\C$, the functor $\eta$ of \eqref{st.C} is a
$0$-special lax $2$-functor with respect to the projection
$\C\{\eq\} \to \eq$. Moreover, any $0$-special lax $2$-functor
$\gamma:\C\{\eq\} \to \C'$ to some $2$-category $\C$ factors as
\begin{equation}\label{01.adj}
\begin{CD}
\C\{\eq\} @>{\eta}>> \C\{\adj\} @>{\gamma'}>> \C',
\end{CD}
\end{equation}
where $\gamma'$ is a $2$-functor, and the factorization is unique up
to a unique isomorphism.
\end{lemma}

\proof{} As in Lemma~\ref{mu.le}, $0$-special and $1$-special maps
in $\eq^o = \Delta e(\{0,1\})$ form a factorization system ---
namely, any map $f:\langle [m],f^*e_\idot \rangle \to \langle
[n],e_\idot \rangle$ in $\Delta e(\{0,1\})$ can be uniquely factored
as
\begin{equation}\label{01.facto}
\begin{CD}
\langle [m], f^*e_\idot \rangle @>{f_1}>> \langle [m]_f,f_0^*e_\idot
\rangle @>{f_0}>> \langle \langle [n],e_\idot \rangle
\end{CD}
\end{equation}
with $0$-special $f^o_0$ and $1$-special $f^o_1$ (indeed, by virtue
of \eqref{de.m}, it suffices to construct \eqref{01.facto} when
$[n]=[0]$, and then either $f_0=\id$, $f_1=f$, or the other way
around, depending on $e_\idot(0) \in \{0,1\}$). Moreover, we also
have a factorization system given by bispecial and anchor maps, and
if $f$ is bispecial, then so are its components $f_0$, $f_1$. Then
we have $\adj \cong \Ar^{1\pm}(\eq)$, where $1\pm$ is the class of
bispecial $1$-special maps, and the embedding \eqref{adj.P} admits a
left-adjoint functor
\begin{equation}\label{l.eq}
l:P(\eq) \to \adj \subset P(\eq)
\end{equation}
sending an arrow $f$ to the $1$-special component $f_1$ of its
decomposition \eqref{01.facto}. For any $2$-category $\C$,
\eqref{l.eq} induces a functor $l:P(\C\{\eq\}) \to \C\{\adj\}$
left-adjoint to the embedding $\C\{\adj\} \subset P(\C\{\eq\})$.

Now, any lax $2$-functor $\gamma:\C\{\eq\} \to \C'$ has the
canonical decomposition \eqref{P.ga}, and $\gamma$ is $0$-special if
and only if $P(\gamma)$ inverts maps in $P(0)=(s \times t)^*(\Iso
\times 0)$, where $0$ is the class of $0$-special maps in
$\C\{\eq\}$. If $\gamma=\eta$ is the embedding \eqref{01.adj}, then
$P(\gamma)=l$ is the functor \eqref{l.eq} that does invert maps in
$P(0)$. Moreover, \eqref{l.eq} is a localization, with the
adjunction map $c \to l(c)$ in $P(0)$ for any $c \in P(\C\{\eq\})$,
so that for any $\gamma$, $P(\gamma)$ inverts maps in $P(0)$ if and
only if it factors through $l$, and the factorization is unique.
\endproof

\subsection{Combinatorics of adjunction.}\label{combi.subs}

By Lemma~\ref{adj.le}, the study of coadjoint pairs reduces in large
part to the combinatorics of the $2$-category $\eq$. Let us prove
several results in this direction. First, consider the functor
$s_0:\eq^> \to \Delta^{o>}$ of \eqref{s.01.C.eq} for the point
$2$-category $\C = \ppt^2$. Recall that objects in $\eq^>$ can also
be described by injective maps $a:[n_1] \to [n]$ in $\Delta^<$, and
let $\eq_1 \subset \eq \subset \eq^>$ be the full subcategory
spanned by bispecial $a$ (or equivalently, by $\langle [n],e_\idot
\rangle \in \Delta^oe(\{0,1\}) = \eq$ such that
$e_\idot(0)=e_\idot(n)=1$).

\begin{lemma}\label{s0.le}
The functor $s_0^>:\eq^> \to \Delta^{o>}$ is a cofibration whose
every fiber has an initial object, and $\eq_1 \subset \eq^>$ is a
subcofibration that contains all these initial objects.
\end{lemma}

\proof{} Recall that $0$-special and $1$-special maps in $\eq^>$
form a factorization system, in either order (in one of the orders
this is \eqref{01.facto}). Then by Example~\ref{ar.exa}, we have a
cofibration $t:\Ar^0(\eq^<) \to \eq^>$, where $\Ar^0(\eq^>) \subset
\Ar(\eq^>)$ is spanned by $0$-special arrows, and as in
Lemma~\ref{mu.le}, we observe that $\iota_0^*\Ar^0(\eq^>) \cong
\eq^>$, with the equivalence sending $x \in \eq^>$ to the adjunction
arrow $\iota_0(s_0(x)) \to x$. This equivalence identifies $s_0$
with the cofibration $\iota_0^*(t)$, so that $s_0$ is a
cofibration. Its fiber $\eq^>_{[m]}$ for some $[m] \in \Delta^{o>}$
is then equivalent to the category of injective maps $[m] \to [n]$,
with $[m]$ corresponding to $[n_0]$. The functor \eqref{bphi}
identifies the opposite to this category with the left comma-fiber
$\Delta^</[m+1]$, and by \eqref{de.m}, this has a terminal object
$[0]^{m+1}$. Explicitly, the corresponding initial object in
$\eq^>_{[m]}$ is $\langle [2m+2],e_\idot \rangle$, where $e_\idot$
is given by
\begin{equation}\label{e.2m}
e_l = l + 1 \mod 2, \qquad l \in [2m+2] = \{0,\dots,2m+2\},
\end{equation}
so that it manifestly lies in $\eq_1$ for any $m$. Finally, by
construction, a map $f:x \to y$ in $\eq^>$ is cocartesian with
respect to $s_0$ if and only if it is $1$-special, and then if its
source $x$ lies in $\eq_1$, so does its target, so that $\eq_1
\subset \eq^>$ is a subcofibration by Lemma~\ref{fib.le}~\thetag{i}.
\endproof

\begin{remark}\label{adj.path.rem}
According to the universal property of the $2$-category $\adj$, the
dense subcategory $\adj_\natural \subset \adj$ defined by maps
cocartesian over $\Delta^o$ should be the category freely generated
by the morphisms $f$, $f^\vee$ but without the adjuntion maps $a$,
$a^\vee$ --- in other words, $\adj \cong \Delta^o [2]_\Lambda \cong
\Delta^o P([2]_\lambda)$ is the simplicial replacement of the path
category of the wheel quiver $[2]_\lambda$ with two vertices $0$,
$1$ and two arrows $f:0 \to 1$, $f^\vee:1 \to 0$. To see this
explicitly, one can use \eqref{e.2m}. Namely, say that an object
$\langle [n],e_\idot \rangle \in \eq$ is {\em alternating} if
$e_{l+1} = e_l + 1 \mod 2$, $0 \leq l < n$. Then alternating objects
correspond to path in the quiver $[2]_\lambda$, an object $\langle
[n],e_\idot \rangle$ of the form \eqref{e.2m} is alternating, and
for any anchor map $a:[m] \to [n]$, so is the induced object
$\langle [m],a^*e_\idot \rangle$. Now, let $a0$ be the class of maps
in $\eq$ of the form $a \circ f$ with $0$-special $f$ and anchor
$a$. Then we have a factorization system $\langle a0,1\pm \rangle$
on $\eq$, and $\adj \cong \Ar^{1\pm}(\eq)$, so that by
Example~\ref{ar.1.exa}, we have
\begin{equation}\label{adj.n}
\adj_\natural \cong \adj \cap P(\eq)_{t^*a0} \subset P(\eq),
\end{equation}
while \eqref{l.eq} restricts to a functor $l:P(\eq)_{t^*a0} \to
\adj_\natural$ left-adjoint to the full embedding $\adj_\natural \to
P(\eq)_{t^*a0}$. However, Lemma~\ref{s0.le} immediately shows that
$l$ itself has a left-adjoint $p:\adj_\natural \to
P(\eq)_{t^*a0}$. It is automatically fully faithful, and by
\eqref{e.2m}, its essential image is spanned by bispecial arrows in
$\eq$ with alternating target. This is exactly $\Delta^o
P([2]_\lambda)$, in the form \eqref{del.P.Q}, and \eqref{t.Q}
provides a functor
\begin{equation}\label{adj.t}
t:\adj_\natural \to \Delta^o_aA([2]_\lambda),
\end{equation}
where $\Delta^o_aA([2]_\lambda) \subset \eq_a$ is the full
subcategory spanned by  alternating objects.
\end{remark}

Next, consider the wreath product $\eq \wr \ppt^2$ of
Definition~\ref{wr.def}. We have the projection $t^o:\eq \wr \ppt^2
\to \Delta^o$ of Example~\ref{wr.exa}, with fibers \eqref{wr.eq}. By
abuse of notation, let $\eq_1 \wr \ppt^2 \subset \eq \wr \ppt^2$ be
the full subcategory spanned by $\eq_1^{m+1} \subset \eq^{m+1} \cong
(\eq \wr \ppt^2)_{[m]}$, $[m] \in \Delta^o$. We also have $\eq \wr
\ppt^2 \cong \eq \times B[1]$ by \eqref{e.wr}, and if we further
abuse notation by denoting $\eq^> \wr \ppt^2 = \eq^> \times^{s^o}
\Ar(\Delta^<)^o$, then we have full embeddings $\eq_1 \wr \ppt^2
\subset \eq \wr \ppt^2 \subset \eq^> \wr \ppt^2$. The projections
\eqref{ar.tw} and \eqref{ar.tw.co} induce a cofibration $s^o:\eq^>
\wr \ppt^2 \to \eq^>$ and a fibration $t^o:\eq^> \wr \ppt^2 \to
\Delta^{o>}$ with fibers $(\eq^> \wr \ppt^2)_{[m]} \cong
\eq^{>(m+1)}$, $[m] \in \Delta^o$. We also have an embedding
$\iota_0 \times \id:\Ar(\Delta^<)^o \to \eq^> \wr \ppt^2$ cartesian
over $\Delta^o$, and the fiber over each $[m] \in \Delta^{o>}$, it
has a left-adjoint
\begin{equation}\label{psi.wr.m}
\psi'_m \cong s_0^{m+1}:(\eq^> \wr \ppt^2)_{[m]} \cong (\eq^>)^{m+1}
\to \Ar(\Delta^<)^o_{[m]} \cong (\Delta^{o>})^{m+1},
\end{equation}
where the last identification is \eqref{de.m}. By
Lemma~\ref{loc.le}~\thetag{ii}, $\iota \times \id$ then has a
left-adjoint $\psi':\eq^> \times \ppt^2 \to \Ar(\Delta^<)^o$ over
$\Delta^{o>}$ that fits into a commutative square
\begin{equation}\label{psi.sq}
\begin{CD}
\eq^> \wr \ppt^2 @>{s^o}>> \eq^>\\
@V{\psi'}VV @VV{s_0}V\\
\Ar(\Delta^<)^o @>{s^o}>> \Delta^{o>}.
\end{CD}
\end{equation}
Restricting $\psi'$ to $\eq_1 \wr \ppt^2$, we obtain functors
\begin{equation}\label{chi.eq}
\begin{CD}
\eq_1 \wr \ppt^2 @>{\psi}>> \Ar(\Delta^<)^o @>{s^o}>> \Delta^{o>},
\end{CD}
\end{equation}
and we denote their composition by $\chi$.

\begin{lemma}\label{chi.le}
The composition $\chi:\eq_1 \wr \ppt^2 \to \Delta^{o>}$ of the
functors \eqref{chi.eq} is a cofibration, the functor $\psi$ is
cocartesian over $\Delta^{o>}$, and each of its fibers $\psi_{[n]}$,
$[n] \in \Delta^{o>}$ has a fully faithful left-adjoint
$\psi^\dg_{[n]}$.
\end{lemma}

\proof{} First observe that all the statements hold if we replace
$\eq_1 \wr \ppt^2$ with $\eq^> \wr \ppt^2$, and $\psi$, $\chi$ with
$\psi'$, $\chi' = s^o \circ \psi'$. Indeed, $s_0$ in \eqref{psi.sq}
is a cofibration by Lemma~\ref{s0.le}, therefore so is $\chi' \cong
s_0 \circ s^o$. A map $f$ in $\eq^> \wr \ppt^2$ cocartesian with
respect to $\chi$ must also be cocartesian with respect to $s^o$,
this happens if and only if $t^o(f)$ is invertible, and since
$\psi'$ commutes with $t^o$, $t^o(\psi'(f))$ is then invertible, so
that $\psi'(f)$ is cocartesian with respect to $s_o$. Dually, $s^o$
inverts maps cartesian with respect to $t^o$, so that for any $[n]
\in \Delta^{o>}$, both $t^o:(\eq^> \wr \ppt^2)_{[n]} \to
\Delta^{o>}$ and $t^o:\Ar(\Delta^<)^o_{[n]} \cong \Delta^{o>} / [n]
\to \Delta^{o>}$ are fibrations, and since the latter is discrete,
$\psi'_{[n]}$ is also a fibration. But each of its fibers has an
initial object by Lemma~\ref{s0.le}, so it has a fully faithful
left-adjoint by Lemma~\ref{loc.le}~\thetag{ii}.

Now note that since any map $f$ in $\eq^> \wr \ppt^2$ cocartesian
over $\Delta^{o>}$ is inverted by $t^o$, it must lie entirely in one
of the fibers \eqref{psi.wr.m}, and then by Lemma~\ref{s0.le}, if
its source is in $\eq_1 \wr \ppt^2 \subset \eq^> \wr \ppt^2$, then
so is its target. Therefore $\eq_1 \wr \ppt^2 \subset \eq^> \wr
\ppt^2$ is a subcofibration, and $\psi$ is cocartesian over
$\Delta^{o>}$. To finish the proof, it remains to observe that by
\eqref{e.2m}, all the initial objects in the fibers of the
fibrations $\psi_{[n]}$ also lie in $\eq_1 \wr \ppt^2$, and apply
Example~\ref{adm.sub.exa}.
\endproof

Note that by Lemma~\ref{loc.le}~\thetag{ii},\thetag{iv},
Lemma~\ref{chi.le} immediately implies that $\psi$ is a localization
in the sense of Definition~\ref{loc.def}. Since $s^o$ in
\eqref{chi.eq} has a fully faithful right-adjoint $\eta^o$ induced
by \eqref{eta.eq}, it is also a localization by
Example~\ref{adm.exa}, and then so is $\chi = s^o \cong \psi$. Thus
by Example~\ref{kan.loc.exa}, for any functor $E:\eq_1 \wr \ppt^2
\to \E$ to some category $\E$ that inverts maps in $\chi^*\Iso$,
$\chi_!E$ and $\chi_*E$ exist, and we have $\chi_!E \cong \chi_*E$
and $E \cong \chi^*\chi_!E \cong \chi^*\chi_!E$. We also have
$\chi_!E \cong s^o_!\psi_!E \cong \eta^{o*}\psi_!E$, or the same
with $\psi_*E$ instead of $\psi_!E$, but neither $\psi_!E$ not
$\psi_*E$ cannot be expressed as a pullback since $\psi$ does not
have an adjoint. For any map $f:[m] \to [n]$ in $\Delta^{o>}$, we
have the map $\psi^\dg_{[n]} \circ f_!  \to f_! \circ
\psi^\dg_{[m]}$ adjoint to the map \eqref{fi.fu} for $\psi$, but it
goes in the wrong direction, so that $\psi_{[n]}^\dg$ do not form a
functor over $\Delta^{o>}$. On objects, we still have canonical
isomorphisms
\begin{equation}\label{psi.E}
\psi_!E([n]) \cong E(\psi_{[n]}^\dg(\eta^o([n]))), \qquad [n] \in
\Delta^{o>}.
\end{equation}
Explicitly, we can think of objects in $\eq_1 \wr \ppt^2$ as triples
$\langle [n],[n_1],p \rangle$ of an object $[n] \in \Delta^<$, a
subset $[n_1] \subset [n]$, and a surjective map $p:[n] \to [m]$ to
some $[m] \in \Delta^<$. Then for any $[n] \in \Delta^{o>}$,
$\eta^o([n])$ is the identity arrow $\id:[n] \to [n]$, and by
\eqref{e.2m}, $\psi_{[n]}^\dg$ sends it to the triple
\begin{equation}\label{2n.eq}
\beta([n]) = \langle [n] \times [2], [n] \times [1],p \rangle,
\end{equation}
where both $[n] \times [2]$ and $[n] \times [1]$ are equipped with
the lexicographical order, the embedding $[n] \times [1] \subset [n]
\times [2]$ is $\id \times a$, where $a:[1] \to [2]$ is the unique
bispecial embeddding, and $p$ is the projection $p:[n] \times [2]
\to [n]$. However, if one tries to write down explicitly
$\psi_!(E)(f)$ for some map $f:[n] \to [m]$, one ends up with a
zigzag of length $3$.

\subsection{Twisting by adjunction.}

Let us now use the combinatorics of Subsection~\ref{combi.subs} to
prove a useful general result on coadjoint pairs. First, assume
given a $2$-category $\C$, let $\gamma:\C_{[0]}\{\eq\} \to \C_{[0]}
\to \C$ be the composition of the projection $\C_{[0]}\{\eq\} =
\C_{[0]} \times \eq \to \C_{[0]}$ and the tautological embedding
$\C_{[0]} \to \C$, and denote $W(\C) = \C_{[0]}\{\eq\} \wr^\gamma
\C$. We then have $\gamma^*\C \cong \C\{\eq\}$, and \eqref{e.wr}
together with \eqref{discr.wr} provide an identification
\begin{equation}\label{W.C}
W(\C) \cong \C \times^{t^o} B[1] \phantom{\times}^{s^o}\!\!\times \eq.
\end{equation}
In particular, we have a cofibration $W(\C) \to \eq \times^{s^o}
B[1] = \eq \wr \ppt^2$, and we can define a full subcategory
$W(\C)_1 \subset W(\C)$ by
\begin{equation}\label{W.C.1}
W(\C)_1 = W(\C) \times_{\eq \wr\ppt^2} \eq_1 \wr \ppt^2 \cong \C
\times^{t^o} \eq_1 \wr \ppt^2.
\end{equation}
The functors \eqref{chi.eq} then induce functors
\begin{equation}\label{chi.C.eq}
\begin{CD}
W_1(\C) @>{\id \times \psi}>> \Ar^\ddag(\C^>) @>{t}>> \C^>,
\end{CD}
\end{equation}
where $\Ar^\ddag(\C^>) \cong \C^> \times^{t^o} \Ar(\Delta^>)^o$ is
the subcategory in $\Ar(\C^>)$ spanned by arrows cartesian with
respect to the fibration \eqref{ar.tw}, and $t$ is the cofibration
\eqref{ar.tw.co}. Denote the composition of the functors
\eqref{chi.C.eq} by
\begin{equation}\label{chi.C}
\chi:W(\C)_1 \to \C^>.
\end{equation}
Note that $W(\C)$ is a $2$-category equipped with a projection to
$\eq$, so it makes sense to speak of $0$-special maps in $W(\C)$.
By abuse of terminology, say that a map in $W(\C)_1$ is $0$-special
if it is $0$-special as a map in $W(\C)$.

\begin{corr}\label{chi.corr}
The functor \eqref{chi.C} is a localization in the sense of
Definition~\ref{loc.def}, and the class of $0$-special maps is dense
in $\chi^*\Iso$.
\end{corr}

\proof{} Lemma~\ref{chi.le} implies that \eqref{chi.C} is a
cofibration and $\psi \times \id$ in \eqref{chi.C.eq} is cocartesian
over $\C^>$, with the same fibers as $\psi$ in \eqref{chi.eq}. Then
$\psi$ is a localization by
Lemma~\ref{loc.le}~\thetag{ii},\thetag{iv}, and $t$ is a
localization since it has a right-adjoint $\eta:\C^> \to
\Ar^\ddag(\C^>)$ induced by \eqref{eta.eq}. To see that $0$-special
maps are dense, note that the adjunction map $\Id \to \eta \circ t$
is cocartesian over $\Ar(\Delta^<)^o$, so that the subclass $v
\subset t^*\Iso$ of maps that are cocartesian over $\Ar(\Delta^<)^o$
is dense in $t^*\Iso$. Therefore $\psi^*v$ is dense in $\chi^*\Iso$,
and this is exactly the class of $0$-special maps.
\endproof

Now assume given a coadjoint pair of functors $\C'\{\adj\} \to \C$
between $2$-categories $\C'$, $\C$, and let
\begin{equation}\label{adj.ga}
\gamma:\C'\{\eq\} \to \C
\end{equation}
be the corresponding $0$-special lax $2$-functor of
Lemma~\ref{adj.le}. We then have the $2$-categories $\gamma_0^*\C$,
$\gamma_1^*\C$, $\gamma^*\C$ of \eqref{2.ga}, with identifications
$\iota_l^*(\gamma^*\C) \cong \gamma_l^*\C$, $l=0,1$. We also have a
$2$-functor $\gamma_1:\C' \to \C$ and a lax $2$-functor
$\gamma_0:\C' \to \C$, with the decompositions \eqref{d.f.facto} for
$\gamma$, $\gamma_1$ and $\gamma_0$.

\begin{prop}\label{tw.prop}
The component $\wgamma$ of the decomposition \eqref{d.f.facto} of
the lax $2$-functor \eqref{adj.ga} factors as
\begin{equation}\label{dm.facto}
\begin{CD}
\C'\{\eq\} @>{\wgamma_1 \times \id}>> \gamma_1^*\C\{\eq\}
@>{\gamma^\dm}>> \gamma^*\C,
\end{CD}
\end{equation}
where $\gamma^\dm$ is a $0$-special lax $2$-functor over
$\eq/\Delta^o$ equipped with an isomorphism $\iota_1^*(\gamma^\dm)
\cong \id$. Moreover, such a factorization is unique up to a unique
isomorphism.
\end{prop}

\proof{} Denote $W(\C',\gamma) = \C_{[0]}'\{\eq\} \wr^{\gamma_{[0]}}
\C$, where the wreath product with respect to a lax $2$-functor is
defined by \eqref{wr.sq}, and $\gamma_{[0]}$ is the composition of
$\gamma$ and the embedding $\C'_{[0]} \subset \C'$.  Let
$W(\C',\gamma)_1 = W(\C',\gamma) \times_{\eq \wr\ppt^2} \eq_1 \wr
\ppt^2$ as in \eqref{W.C.1}. Then we actually have $W(\C',\gamma)_1
\cong W(\gamma_1^*\C)$, and on the other hand, since $\gamma$ is
$0$-special, the lax $2$-functor $W(\C',\gamma) \to \gamma^*\C$ of
\eqref{wr.C} is also $0$-special. By Remark~\ref{0.anch.rem}, it
then induces a functor
\begin{equation}\label{W.C.ga}
W(\gamma_1^*\C)_1 \cong W(\C',\gamma)_1 \subset W(\C',\gamma) \to
\gamma^*\C
\end{equation}
that sends $0$-special maps to $0$-special maps, and $0$-anchor maps
to $0$-anchor maps. Therefore if we compose \eqref{W.C.ga} with
$s_0:\gamma^*\C \to (\gamma_0^*\C)^>$ to obtain a functor
$\Phi(\gamma):W(\gamma_1^*\C)_1 \to (\gamma_0^*\C)^>$, then
$\Phi(\gamma)$ inverts $0$-special maps, and therefore by
Corollary~\ref{chi.corr}, uniquely
factors through a functor
\begin{equation}\label{Theta.eq}
\Theta(\gamma) = \chi_!\Phi(\gamma):\gamma_1^*\C \to
(\gamma_0^*\C)^>.
\end{equation}
If we let $\pi_l:(\gamma_l^*\C)^> \to \Delta^{o>}$, $l=0,1$ be the
structural cofibrations, then $\pi_0 \circ \Phi(\gamma) \cong \pi_1
\circ \chi$, so again by uniqueness, $\Theta(\gamma)$ is a functors
over $\Delta^{o>}$ (in particular, it factors through $\gamma_0^*\C
\subset (\gamma_0^*\C)^>$). Moreover, if we let $a$ be the class of
cocartesian liftings of anchor maps in $\gamma_1^*\C$, then maps in
$\chi^*a$ cocartesian over $\eq$ are exactly $0$-anchor maps, and
$\Phi(\gamma)$ sends those to maps cocartesian over
$\Delta^o$. Therefore $\Theta(\gamma)$ is a lax $2$-functor from
$\gamma_1^*\C$ to $\gamma_2^*\C$. By virtue of uniqueness, if
$\gamma$ factors through the projection $\C'\{\eq\} \to \C'$, so
that $\gamma_0 \cong \gamma_1$, then $\Theta(\gamma) \cong
\Id$. Moreover, $\Phi(\gamma)$ is obviously functorial with respect
to $\C'$, and then so is $\Theta(\gamma)$: for any $2$-functor
$\delta:\C'' \to \C'$, we have a natural isomorphism
\begin{equation}\label{the.de}
\Theta(\gamma \circ (\delta \times \id)) \cong
\delta^*\Theta(\gamma).
\end{equation}
Moreover, the functor $W(\wgamma_1):W(\C')_1 \to W(\gamma_1^*\C)_1$
induced by $\wgamma_1:\C' \to \gamma_1^*\C$ fits into a commutative
diagram
$$
\begin{CD}
W(\C')_1 @>{W(\wgamma_a)}>> W(\gamma_1^*\C)_1\\
@VVV @VV{\Phi(\gamma)}V\\
\C'\{\eq\} @>{s_0}>> (\gamma_0^*\C)^>,
\end{CD}
$$
where the left vertical arrow is again induced by \eqref{wr.C}, and
again by Corollary~\ref{chi.corr}, this induces an isomorphism
\begin{equation}\label{0.1.the}
\gamma_0 \cong \Theta(\gamma) \circ \wgamma_1
\end{equation}
of lax $2$-functors from $\C'$ to $\gamma_0^*\C$.

Now consider the map $\{0,1\} \times \{0,1\} \to \{0,1\}$ sending $l
\times l'$ to $\max(l,l')$, $l,l' \in \{0,1\}$, and let $m:\eq
\times^2 \eq \to \eq$ be the corresponding $2$-functor. Then
$\gamma_m=\gamma \circ (\id \times m)$ defines a coadjoint pair of
functors between $\C'\{\eq\} = \C'[\{0,1\}]$ and $\C$, so that
\eqref{Theta.eq} provides a lax $2$-functor
\begin{equation}\label{theta.01}
\Theta(\gamma_m):\gamma_1^*\C\{\eq\} \to \gamma^*\C,
\end{equation}
where we identify $\gamma_{m1}^*\C \cong \gamma_1^*\C\{\eq\}$ and
$\gamma_{m0}^*\C \cong \gamma^*\C$. Moreover, say that a map $f$ in
$W(\gamma_{m1}^*\C)_1$ is {\em $00$-special} if it is cocartesian
over $\eq$, and $\chi(f)$ is $0$-special in $\gamma_{m1}^*\C =
\gamma_1^*\C\{\eq\}$; then $\Phi(\gamma_m)$ sends $00$-special maps
to $0$-special maps, and therefore the lax $2$-functor
\eqref{theta.01} is $0$-special. Thus if we take $\gamma^\dm =
\Theta(\gamma \circ (\id \times m))$, with the isomorphism
$\gamma^\dm \circ (\wgamma_1 \times \id) \cong \wgamma$ provided by
\eqref{0.1.the}, we obtain the decomposition
\eqref{dm.facto}. Moreover, applying \eqref{the.de} to the embedding
$\iota_1:\C' \to \C'\{\eq\}$, we obtain an isomorphism
\begin{equation}\label{io.dm}
\iota_1^*(\gamma^\dm) \cong \Theta(\gamma_m \circ
(\iota_1 \times \id)),
\end{equation}
and since $(\id \times m) \circ (\iota_1 \times \id):\C'\{\eq\} \to
\C'\{\eq\}$ factors through the tautological projection $\C'\{\eq\}
\to \C' \cong \iota_1(\C') \subset \C\{\eq\}$, the target of the
isomorphism \eqref{io.dm} is the identity functor.

This proves existence. For uniqueness, assume given some other
decomposition \eqref{dm.facto} of the functor $\wgamma$, and apply
the construction above to the adjoint pair $\gamma^\dm$. Then
$(\gamma^\dm)_1 = i_1^*\gamma^\dm$ is identified with $\id$, so that
we have $\gamma^\dm \cong \Theta(\gamma^\dm \circ (\id \times m))
\circ (\id \times \id) = \Theta(\gamma^\dm \circ (\id \times m))$,
and then \eqref{the.de} for the $2$-functor $\delta = \wgamma_1
\times \id$ provides an identification
$$
\delta^*\gamma^\dm \cong \delta^*\Theta(\gamma^\dm \circ (\id \times
m)) \cong \Theta(\gamma^\dm \circ (\id \times m) \circ (\delta
\times \id)) = \Theta(\gamma_m).
$$
Since $\delta$ is dense, $\delta^* = \id$, so that this reads as
$\gamma^\dm \cong \Theta(\gamma_m)$.
\endproof

We note that in particular, the functor $\gamma^\dm$ in the
decomposition \eqref{dm.facto} induces a lax $2$-functor
\begin{equation}\label{P.al}
\Theta(\gamma) = \iota_0^*(\gamma^\dm):\gamma_1^*\C \to \gamma_0^*\C,
\end{equation}
equipped with an isomorphism $\Theta(\gamma) \circ \wgamma_1 \cong
\wgamma_0$, and this is the essential ingredient of the whole thing
(we actually construct it first, in \eqref{Theta.eq}). We call
$\Theta(\gamma)$ the {\em twisting functor} associated to the
coadjoint pair. To compute its components $\Theta(\gamma)_{[n]}$
more explicitly, one can use \eqref{psi.E} and \eqref{2n.eq}. This
shows that explicitly, the coadjoint pair $\gamma$ defines an
adjoint pair of maps $h:\gamma_0(c) \to \gamma_1(c)$,
$h^\vee:\gamma_1(c) \to \gamma_0(c)$ for any $c \in \C'_{[0]}$, and
for any $[m] \in \Delta$, we have
\begin{equation}\label{adj.fu}
\Theta(\gamma)_{[m]}(c) = \prod_{l=1}^m h^\vee(b^o_{l!}c) \circ
a_{l!}^o(c) \circ h(b^o_{(l-1)!}c), \qquad c \in
(\gamma_1^*\C)_{[m]},
\end{equation}
where the product is the product \eqref{c.m.eq}, and $b_l:[0] \to
[m]$ is the embedding onto $l \in [m]$. The maps \eqref{fi.fu} are
induced by the adjunction maps between $h$ and $h^\vee$, but writing
them down explicitly is hard and probably pointless.

\subsection{Iterated adjunction.}\label{ite.subs}

As an application of Proposition~\ref{tw.prop}, let us give a
somewhat more invariant description of the adjunction $2$-category
\eqref{adj.c} of a $2$-category $\C$.

As in Example~\ref{tw.yo.exa}, let $\Nat = \Tw(\Delta)$, with its
natural cofibration $s \times t:\Nat \to \Delta^o \times \Delta$,
and let $\Eq = (V^o \times V)^*\Tw(\Sets)$, where $V:\Delta \to
\Sets$ is the forgetful functor. Both $\Nat$ and $\Eq$ are
$\Delta$-kernels in the sense of Definition~\ref{ker.def}, and $V$
induces a morphism $\Nat \to \Eq$. For any $[m] \in \Delta$, the
fiber $\Eq_{[n]}$ of the cofibration $t:\Eq \to \Delta$ is given by
$\Eq_{[n]} = \Delta^oe(V([n]))$, and $\Nat_{[n]} \subset \Eq_{[n]}$
is $\Delta^o[n] = (\Delta/[n])^o$, with the embedding induced by the
tautological functor $\id:[n] \to e(V([n]))$. In particular,
$\Eq_{[1]} = \Delta^oe(\{0,1\}) = \eq$ and $\Nat_{[1]} = \nat$, with
the embedding $\nat \to \eq$ given by the composition of the
functors \eqref{1.adj.eq}.

What we want to do is to construct a $\Delta$-kernel $\Adj$ that
fits in between $\Nat$ and $\Eq$ and completes a kernel version of
\eqref{1.adj.eq}. In order to do this, we generalize the description
of the $2$-category $\adj$ in terms of the diagrams \eqref{2adj.dia}
and \eqref{2.2adj.dia} given in Subsection~\ref{int.subs}.

To do this, consider the embedding $\rho:\Delta_+ \to \Delta$ with
its left-adjoint $\lambda:\Delta \to \Delta_+$ and the composition
$\kappa:\Delta \to \Delta$, let $b$ be the class of maps in
$\Delta_+$ that are bispecial, and define a cofibration $C \to
\Delta$ by
\begin{equation}\label{C.eq}
C = \lambda^*\Id^b_{**}\rho^*\Eq,
\end{equation}
where $\Eq$ is also considered as a cofibration over $\Delta$. As in
Definition~\ref{DI.def}, the fibers $\C_{[n]}$ of the cofibration
\eqref{C.eq} can be described as in Example~\ref{id.exa}. By
\eqref{p-m.eq}, we have $[n] \setminus \Delta_+ \cong
(\Delta_-/[n])^o$, and since $\Delta_t \subset \Delta_-$, we then
have an embedding $v_n:\Delta_t/[n] \to \Delta_-/[n]$. This
embedding is fully faithful and left-admissible, with the
left-adjoint functor $v_n^\dg$ sending a special arrow $f:[m] \to
[n]$ to the anchor component $a$ of its decomposition
\eqref{l.r.facto} for the anchor/bispecial fatorization system, and
the adjunction map $\id \to v_n \circ v_n^\dg$ is bispecial. Since
$\Delta_t \cong \N$, we can identify $[n] \cong (\Delta_t/[n])^o$,
and we then have
\begin{equation}\label{C.n}
C_{[n]} =
\Sec^b((\Delta_-/\lambda([n]))^o,p_{\lambda([n])}^*\rho^*\Eq) \cong
\Sec(\lambda([n]),w_n^*\Eq),
\end{equation}
where we denote
\begin{equation}\label{w.n}
w_n = \rho \circ p_{[n]} \circ v^o_n:[n] \to \Delta.
\end{equation}
Moreover, for any special map $f:[n] \to [n']$, we have $v_n^{\dg o}
\circ f^* \cong f^\dg \circ v_{n'}^{\dg o}$, where $f^*$ is the
functor \eqref{comma.tr}, and by adjunction, this induces a map
$$
a(f):w_n \circ f^\dg = \rho \circ p_{[n]} \circ v^o_n \circ f^\dg
\to \rho \circ p_{[n]} \circ f^* \circ v_{n'}^o \cong \rho \circ
p_{[n']} \circ v_{n'}^o = w_{n'}.
$$
Then in terms of \eqref{C.n}, the transition functor $f_!:C_{[n]}
\to C_{[n']}$ for a map $f:[n] \to [n']$ is given by $f_! =
a(\lambda(f))_! \circ \theta(f)^*$, where $\theta$ is the functor
\eqref{theta.eq}, so that $\theta(f) = \lambda(f)^\dg$. Explicitly,
$w_n:[n] \to \Delta$ corresponds to the diagram
\begin{equation}\label{w.n.dia}
\begin{CD}
[n] @>{t_\dg}>> [\nm] @>{t_\dg}>> \dots @>{t_\dg}>> [1] @>{t_\dg}>>
[0]
\end{CD}
\end{equation}
in $\Delta$, with the maps $t_\dg:[l] \to [l-1]$ adjoint to the
embeddings $t:[l-1] \to [l]$, and the fiber $\C_{[n]}$ is opposite
to the category of pairs $\langle \gamma,\pi \rangle$ of a functor
$\gamma:[\np] = \lambda([n]) \to \Delta$ and a map $\pi:V(\gamma)
\to V(w_{n+1})$. For any map $f:[n] \to [n']$, the transition
functor $f_!:C_{[n]} \to C_{[n']}$ sends $\langle \gamma,\pi
\rangle$ to $\theta(f)^*\gamma$ equipped with the composition map
\begin{equation}\label{V.w}
\begin{CD}
V(\theta(f)^*\gamma) @>{V(\theta(f)^*\pi)}>> V(\theta(f)^*w_{n+1})
@>{V(a(\lambda(f)))}>> V(w_{n'+1}).
\end{CD}
\end{equation}
In particular, since $\theta(f)$ is bispecial, both $\gamma(n+1)$
and $\gamma(0)$ with the projection $\pi(0):V(\gamma(0)) \to
V(\kappa([n]))$ are functorial with respect to $\langle \gamma,\pi
\rangle$. Therefore sending $\langle \gamma,\pi \rangle$ to
$\gamma(0) \times \gamma(n+1)$ defines a functor $\chi:C \to
\kappa^*\Eq \times \Delta^o$ cocartesian over $\Delta^o$. We can
then define a category $\Adj'$ by the cartesian square
\begin{equation}\label{Adj.0.eq}
\begin{CD}
\Adj' @>>> C\\
@V{\nu'}VV @VV{\chi}V\\
\Eq @>>> \kappa^*\Eq \times \Delta^o,
\end{CD}
\end{equation}
where the bottom arrow is the product of the full embedding $a_!:\Eq
\to \kappa^*\Eq$ induced by the adjunction map $a:\id \to \kappa$,
and the full embedding $\ppt \to \Delta^o$ onto $[0] \in
\Delta^o$. Then $\nu'$ is obviously a cofibration, so that $\Adj'$
is a $\Delta$-kernel and $\nu'$ is a morphism of $\Delta$-kernels,
while the full embedding $\Adj' \subset C$ is cocartesian over
$\Delta$. For any $[n] \in \Delta$, the subcategory $\Adj'_{[n]}
\subset C_{[n]}$ is spanned by pairs $\langle \gamma,\pi \rangle$
such that $\gamma(n+1) = [0]$, and $\pi$ factors through
$V(w_{n+1}') \subset V(w_{n+1})$, where $w'_n \subset w_n$ is
obtained by replacing $[n]$ in \eqref{w.n.dia} with $t([\nm])
\subset [n]$. In particular, $\Adj'_{[0]}$ is $\Delta^o$, and $\adj'
= \Adj'_{[1]}$ is the arrow category $\Ar(\eq)$.

Now recall that by \eqref{2adj.dia}, $\adj$ is naturally identified
with the full subcategory in $\Ar(\eq)$ spanned by bispecial
$1$-special maps. We can now denote by $\eps:\ppt \to \Delta$ the
embedding onto $[1] \in \Delta$, so that we have $\adj' \cong
\eps^*\Adj'$, and define a category $\Adj$ by the cartesian square
\begin{equation}\label{Adj.eq}
\begin{CD}
\Adj @>{\alpha}>> \Adj'\\
@VVV @VVV\\
\eps_*\adj @>>> \eps_*\adj'
\end{CD}
\end{equation}
where the rightmost vertical arrow is the functor
\eqref{adj.2.1}. Then by definition, $\Adj \subset \Adj'$ is a full
subcategory and a subcofibration over $\Delta$, and its fiber
$\Adj_{[n]} \subset \Adj'_{[n]}$ is spanned by pairs $\langle
\gamma,\pi\rangle$ such that $f_!\langle \gamma,\pi\rangle \in \adj$
for any $f:[n] \to [1]$. Since $\Adj_{[0]} \cong \Adj'_{[0]}$, the
condition is empty if $f$ factors through $[0]$, so that $\Adj_{[0]}
\cong \Delta^o$ and $\Adj_{[1]} \cong \adj$.  The functor $\nu'$ of
\eqref{Adj.0.eq} induces a functor
\begin{equation}\label{Nu.eq}
\nu:\Adj \to \Eq,
\end{equation}
cocartesian over $\Delta$, whose component $\nu_{[1]}$ is identified
with $\nu$ of \eqref{1.adj.eq}.

\begin{lemma}\label{Adj.corr}
The functor \eqref{Nu.eq} is a cofibration, and $\alpha:\Adj \to
\Adj'$ of \eqref{Adj.eq} admits a left-adjoint functor $\beta:\Adj'
\to \Adj$ cocartesian over $\Eq$.
\end{lemma}

\proof{} For any $[n] \in \Delta$, $\langle \gamma,\pi \rangle \in
\Adj'_{[n]}$, $l \in [\np]$, $m \in w'_{n+1}(l)$, denote
$\gamma(l)_m = \gamma(l) \setminus \pi(l)^{-1}(\{0,\dots,m-1\})
\subset \gamma(l)$. This defines a decreasing filtration
$\gamma(l)_\idot$ on $\gamma(l)$, and by the definition of the
functor $w'_{n+1}$, the map $g_l:\gamma(0) \to \gamma(l)$ fits into
a cartesian square
\begin{equation}\label{ga.l.sq}
\begin{CD}
\gamma(0)_l @>{g_l'}>> \gamma(l)_1\\
@VVV @VVV\\
\gamma(0) @>>> \gamma(l),
\end{CD}
\end{equation}
where the vertical arrows are the embeddings. If as in
Subsection~\ref{int.subs}, we interpret $\eq^o$ as the category of
injective arrows $[n]_1 \to [n]$ in $\Delta$, then \eqref{ga.l.sq}
represents a object $\Ar(\eq) \cong \adj'$, and by \eqref{V.w}, this
is exactly the object $f_!\langle \gamma,\pi \rangle$ for the unique
map $f:[n] \to [1]$ such that $l = \theta(f)(1) \in [\np]$. Then
$\langle \gamma,\pi \rangle$ lies in $\Adj_{[n]}$ iff for any $l \in
[\np]$, $g_l:\gamma(0) \to \gamma(l)$ is bispecial and $g'_l$ is
bijective, so that \eqref{ga.l.sq} is a square of the form
\eqref{2adj.dia}.

Now, as in the proof of Lemma~\ref{adj.le}, for any $\langle
\gamma,\pi \rangle \in \Adj_{[n]}$ and $l \in [\np]$, $1 \leq l \leq
n$, we have a unique decomposition
\begin{equation}\label{b.a.l.facto}
\begin{CD}
\gamma(0) @>{b_l}>> \gamma'(l) @>{a_l}>> \gamma(l)
\end{CD}
\end{equation}
of the map $g_l$ such that $b_l$ is bispecial and strict and
injective on $\gamma(0)_l$ in the sense of Remark~\ref{adj.rem}, and
$a_l$ is a composition of an anchor map and a map that is strict and
injective on $\gamma'(l) \setminus b_l(\gamma(0)_l)$. But then $b_l$
is also strict and injective on $\gamma(0)_{l+1} \subset
\gamma(0)_l$, so that by the uniqueness of \eqref{b.a.l.facto}, the
map $\gamma'(l) \to \gamma(l) \to \gamma(l+1)$ factors through
$a_{l+1}:\gamma'(l+1) \to \gamma(l+1)$. Then all the maps
$b_l:\gamma(0) \to \gamma'(l)$ fit together into a functor
$\gamma':[\np] \to \Delta$ equipped with a map $a:\gamma' \to
\gamma$, and $\langle \gamma',\pi \circ V(a) \rangle$ lies in
$\Adj_{[n]} \subset \Adj'_{[n]}$. Therefore sending $\langle
\gamma,\pi \rangle$ to $\langle \gamma',\pi \circ V(a) \rangle$
defines a functor $\beta_{[n]}:\Adj'_{[n]} \to \Adj_{[n]}$
left-adjoint to $\alpha_{[n]}:\Adj_{[n]} \to \Adj'_{[n]}$. By
construction, $\gamma(0) = \gamma'(0)$, so that $\beta_{[n]}$ is a
functor over $\Eq_{[n]}$. Then by Lemma~\ref{fib.le}~\thetag{iii},
$\Adj_{[n]} \to \Eq_{[n]}$ is a cofibration, and $\beta_{[n]}$ is
cocartesian. Finally, again by construction, the functors
$\beta_{[n]}$ commute with the transition functors $f_!$ for all
maps $f$ in $\Delta$, so that to finish the proof, it remains to
apply Lemma~\ref{fib.le}~\thetag{ii}.
\endproof

By Lemma~\ref{Adj.corr}, $\Adj$ is a $\Delta$-kernel, and
\eqref{Nu.eq} is a morphism of $\Delta$-kernels. In particular,
$\Adj_{[n]} \to \Delta^o$ is a cofibration for any $[n] \in
\Delta$. It turns out that these cofibrations can be also
constructed inductively, starting with $\Adj_{[1]} = \adj$. Namely,
assume given a cocartesian square \eqref{seg.sq} with $l=1$, and
denote $p = s^\dg:[n] \to [1]$, $q = t_\dg:[n] \to [\nm]$
(explicitly, $p = \id$ on $s([1]) \subset [n]$ and sends the rest to
$1 \in [1]$, and $q=\id$ on $t([\nm]) \subset [n]$ and sends the
rest to $0 \in [\nm]$). Moreover, note that we have a map
$\tau:t^*w'_{n+1} \to w_n$ equal to $\id$ at $l \in [n]$, $l \geq
1$, and to $q$ at $l = 0$, and sending $\langle \gamma,\pi \rangle$
to $\langle t^*\gamma,\tau \circ t^*\pi \rangle$ defines a functor
$t^*:\Adj_{[n]} \to \Adj_{[n-1]}$.

\begin{lemma}\label{Adj.le}
For any $n \geq 2$, we have a cartesian square
\begin{equation}\label{Adj.sq}
\begin{CD}
\Adj_{[n]} @>{t^* \times p_!}>> \Adj_{[n-1]}\{\adj\}\\
@V{\nu_{[n]}}VV @VV{s \circ (\nu_{[n-1]} \times \id)}V\\
\Eq_{[n]} @>{q_! \times p_!}>> \Eq_{[n-1]} \times^2 \eq,
\end{CD}
\end{equation}
where $s$ is the projection \eqref{st.C}.
\end{lemma}

\proof{} The cocartesian square \eqref{seg.sq} induces a cocartesian
square
$$
\begin{CD}
[0]^> @>{s^>}>> [l]^>\\
@V{t^>}VV @VV{t^>}V\\
[\nl]^> @>{s^>}>> [n]^>,
\end{CD}
$$
and a functor $\gamma:[\np] = [n]^> \to I$ to any category $I$ is
uniquely defined by its restrictions $\gamma' = t^{>*}\gamma:[l]^>
\to I$, $\gamma'' = s^{>*}\gamma:[\nl]^> \to I$, and an isomorphism
$t^{>*}\gamma'' \cong s^{>*}\gamma'$. In particular, we can take $I
= \Delta$. Moreover, $s^>$ and $t^>$ are antispecial, thus have
right-adjoint functors $s^>_\dg$, $t_\dg^>$, and giving two maps
$\pi':V(\gamma') \to V(w'_{l+1})$, $\pi'':V(\gamma'') \to
V(w'_{n-l+1})$ is equivalent to giving a single map
\begin{equation}\label{pi.prod}
\pi = \pi' \times \pi'':V(\gamma) \to (t^>_\dg)^*V(w'_{l+1}) \times
(s^>_\dg)^*V(w'_{n-l+1}).
\end{equation}
For any set $S$, denote by $S^>:[0]^> \to \Sets$ the functor sending
$0 \in [0]^>$ to $S$ and $o \in [0]^>$ to $\ppt$. Then we have
$s^{>*}V(w'_{l+1}) \cong V([l])^>$, $t^{>*}V(w_{n-l+1}) \cong
V([\nl])^>$, and if we let $a = s \circ t = t \circ s:[0] \to [n]$
be the embedding onto $l \in [n]$, and then take $l=1$, we have a
cartesian square
$$
\begin{CD}
V(w'_{n+1}) @>>> (t^>_\dg)^*V(w'_{2}) \times (s^>_\dg)^*V(w'_n)\\
@VVV @VVV\\
(a^>_\dg)^*V([n])^> @>>> (a^>_\dg)^*V([1])^> \times
(a^>_\dg)^*V([\nm])^>,
\end{CD}
$$
where the bottom arrow is induced by $p \times q:V([n]) \to V([1])
\times V([\nm])$. Combined with \eqref{pi.prod}, this gives the
cartesian square \eqref{Adj.sq} with $\Adj'$, $\adj'$ instead of
$\Adj$, $\adj$, and then the explicit description of $\Adj \subset
\Adj'$ in terms of diagrams \eqref{ga.l.sq} shows that it induces
\eqref{Adj.sq} on the nose.
\endproof

Lemma~\ref{Adj.le} immediately implies by induction that for any
$[n]$, the cofibration $\Adj_{[n]} \to \Delta^o$ is a $2$-category,
and $\nu_{[n]}:\Adj_{[n]} \to \Eq_{[n]}$ is a $2$-functor. The
functor \eqref{Nu.eq} has an obvious fully-faithful right-adjoint
\begin{equation}\label{Eta.eq}
\eta:\Eq \to \Adj.
\end{equation}
Explicitly, an object in $\Eq_{[n]} \subset \Eq$ is a pair $\langle
[m],\pi \rangle$, $[m] \in \Delta$, $\pi:V([m]) \to V([n])$, and
$\eta(\langle [m],\pi \rangle)$ is obtained by taking the constant
functor $\gamma:[n] \to \Delta$ with value $[m]$, extending it to a
functor $\gamma^>:[\np] = [n]^> \to \Delta$ sending $o \in [n]^>$ to
$[0]$, and then extending $\pi$ to a map $\pi:V(\gamma^>) \to
V(w_{n+1})$. While $\eta$ is not cocartesian over the whole
$\Delta$, it is obviously cocartesian over all antispecial
maps. Over $[1] \in \Delta$, $\eta_{[1]}:\eq \to \adj$ is the lax
$2$-functor \eqref{eta.eq}, and for any $[n]$, we have the base
change isomorphism
\begin{equation}\label{eta.n}
e \circ \eta_{[n]} \cong (\eta_{[n-1]} \times \id) \circ \eta \circ e,
\end{equation}
induced by \eqref{Adj.sq}, so that by induction, $\eta_{[n]}$ is
a lax $2$-functor for any $n$.

In general, the projection $\nu$ does not have a
left-adjoint. However, if we define a $\Delta$-kernel $\Adj^0$ by
the cartesian square
$$
\begin{CD}
\Adj^0 @>>> \Adj\\
@V{\nu^0}VV @VV{\nu}V\\
\Nat @>>> \Eq,
\end{CD}
$$
then $\nu^0$ does have a left-adjoint $\delta^0:\Nat \to
\Adj^0$. Namely, $\Adj^0 \subset \Adj$ is spanned by pairs $\langle
\gamma,\pi \rangle$ such that $\pi = V(\pi')$ for a map $\pi':\gamma
\to w'_{n+1}$, so that for any $[n] \in \Delta$, we have an object
$\langle w'_{n+1},V(\id) \rangle \in \Adj^0_{[n]}$. Then objects in
$\Nat$ are arrows $f:[m] \to [n]$ in $\Delta$, and $\delta^0$ sends
$f$ to $f_!\langle w'_{n+1},V(\id) \rangle$. By construction,
$\delta^0$ is cocartesian over $\Delta$, and it is easy to check
using Lemma~\ref{Adj.corr} that it is also cocartesian over
$\Delta^o$. Therefore composing $\delta^0$ with the embedding
$\Adj^0 \to \Adj$, we obtain a functor $\delta:\Nat \to \Adj$ that
fits into a diagram
\begin{equation}\label{1.Adj.eq}
\begin{CD}
\Nat @>{\delta}>> \Adj @>{\nu}>> \Eq
\end{CD}
\end{equation}
of morphisms of $\Delta$-kernels. Over $[1] \in \Delta$, it reduces
to \eqref{1.adj.eq}. Now for any $2$-category $\C$, let
\begin{equation}\label{aadj.c}
\aAdj(\C) = (\Adj \otimes_\Delta \C)_\natural \subset \Adj
\otimes_\Delta \C
\end{equation}
be the cofibration over $\Delta^o$ with fibers
$\Fun^2(\Adj_{[m]},\C)$ and transition functors induced by those of
the cofibration $\Adj \to \Delta$. By Example~\ref{tw.yo.exa}, the
morphism $\delta$ of \eqref{1.Adj.eq} induces a functor $\aAdj(\C)
\to \C_\natural$.

\begin{prop}\label{adj.prop}
The cofibration $\aAdj(\C) \to \Delta^o$ is a $2$-category, and the
functor $\aAdj(\C) \to \C_\natural$ induced by the morphism $\delta$
of \eqref{1.Adj.eq} identifies it with the $2$-category
\eqref{adj.c}.
\end{prop}

\proof{} The second claim is by definition true over $[0]$ and $[1]
\in \Delta$, so it follows from the first. Also by definition,
$\aAdj(\C)_{[0]} = \C_{[0]}$ is discrete, so what we need to check
is the Segal condition. By induction on $n$, it suffices to consider
the squares \eqref{seg.sq} with $l=1$. Consider the $2$-category
$\Adj_{[n]}$, with the embeddings
$$
s:\adj=\Adj_{[1]} \to \Adj_{[n]}, \qquad t:\Adj_{[n-1]} \to
\Adj_{[n]},
$$
and the embedding $e:\Adj_{[n]} \to \Adj_{[n-1]}\{\adj\}$ of
\eqref{Adj.sq}. Denote by $S = V([\nm])$ the set of objects of the
$2$-category $\Adj_{[n-1]}$, with the embedding $S =
\Adj_{[n-1],[0]} \to \Adj_{[n-1]}$, and let $p_0:S \times \adj \to
\Adj_{[n]}$ be the adjoint pair given by $s$ on $0 \times \adj$ and
by the constant projection onto $i$ on $i \times \adj$, $i \in S
\setminus \{0\}$. Then by Proposition~\ref{tw.prop} and
Lemma~\ref{adj.le}, $p_0$ extends uniquely to a $2$-functor
$p:\Adj_{[n-1]}\{\adj\} \to \Adj_{[n]}$ such that $p \circ e \cong
\id$. Now assume given a $2$-functor $F:\Adj_{[n-1]} \to \C$ and an
adjoint pair $\gamma:\adj \to \C$ with $\gamma(1)=F(0)$, and let $E$
be the groupoid of $2$-functors $F':\Adj_{[n]} \to \C$ equipped with
isomorphisms $s^*F' \cong \gamma$ and $t^*F' \cong F$. Extend
$\gamma$ to an adjoint pair $\gamma':S \times \adj \to \C$ given by
$\gamma$ on $0 \times \adj$ and by the constant projection onto
$F(i) \in \C_{[0]}$ on $i \times \adj$, $i \in S \setminus
\{0\}$. Let $l:\Adj_{[n-1]} \to \Adj_{[n-1]}\{\adj\}$ be induced by
the embedding $\ppt \to \adj$ onto $1 \in \adj_{[0]}$, let $r:S
\times \adj = S\{\adj\} \to \Adj_{[n-1]}\{\adj\}$ be induced by the
embedding $S \to \Adj_{[n-1]}$, and let $P$ be the category of
$2$-functors $F'':\Adj_{[n-1]}\{\adj\} \to \C$ equipped with
isomorphisms $l^*F'' \cong F$, $r^*F'' \cong \gamma'$. Then the
$2$-functors $e$ and $p$ induce functors $p^*:E \to P$, $e^*:P \to
E$ such that $e^* \circ p^* \cong \Id$, so that $E$ is a retract of
$P$. But by Lemma~\ref{adj.le} and Proposition~\ref{tw.prop}, we
have $P \cong \ppt$.
\endproof

\section{Enriched categories.}

\subsection{Categories and modules.}

Fix a category $\C$ equipped with a unital monoidal structure $B\C$
in the sense of Definition~\ref{mon.def}. Recall that for any set
$S$, we have the small category $e(S)$ of Subsection~\ref{set.sss}
(objects are elements $s \in S$, and there is exactly one morphism
between any two objects).

\begin{defn}\label{enr.def}
A {\em a small $\C$-enriched category} with a set of objects $S$ is
a lax $2$-functor $A$ from $e(S)$ to $B\C$. A {\em functor} between
small $\C$-enriched categories $\langle S,A \rangle$ and $\langle
S',A' \rangle$ is a pair $\langle f,g \rangle$ of a map $f:S \to S'$
and a morphism $g:A \to \Delta^o(e(f))^*A'$.
\end{defn}

\begin{remark}
Explicitly, defining $A$ amounts to giving an object $A(s,s')$ in
$B\C_{[1]} \cong \C$ for any $s,s' \in S$ and equipping them with
the composition and unity maps; thus Definition~\ref{enr.def} is a
repackaging of the usual notion of a category enriched over $\C$.
\end{remark}

\begin{exa}\label{pt.exa}
  If $S = \ppt$ is a one-element set, then a $\C$-enriched category
  $\langle S,A \rangle$ is the same thing as a unital associative
  algebra object in $\C$. In particular, the unit object $1 \in \C$
  is tautologically an algebra, and it defines the point
  $\C$-enriched category $\ppt_\C$. For any $\langle S,A \rangle$
  and object $s \in S$, we have the tautological functor
  $i_s:\ppt_\C \to \langle S,A \rangle$ given by the embedding $\ppt
  \to S$ onto $s$ and the unit map $1 \to A(s,s)$.
\end{exa}

\begin{defn}\label{A.opp.def}
  The {\em opposite category} $\langle S,A \rangle^o$ to a small
  $\C$-enriched category $\langle S,A \rangle$ is the
  $\C^\iota$-enriched category $\langle S,A \rangle^o = \langle
  S,A^\iota \rangle$, where $A^\iota = \iota \circ A$, and
  $\C^\iota$ is $\C$ with the opposite monoidal structure $B\C^\iota
  = \iota^*B\C$.
\end{defn}

Assume given a small $\C$-enriched category $\langle S,A
\rangle$. To define modules over $A$, recall that we have the
functor $\kappa:\Delta \to \Delta$ equipped with the adjunction map
$a:\id \to \kappa$ of \eqref{ka.adj}. The transition functor $a^o_!$
for the opposite map then provides a functor $a^o_!:\kappa^{o*}B\C
\to \C$ cocartesian over $\Delta^o$.

\begin{defn}\label{A.mod.def}
A {\em right module} $M$ over a small $\C$-enriched category
$\langle S,A \rangle$ is a functor $M:\Delta^o e(S) \to
\kappa^{o*}B\C$ over $\Delta$ cocartesian over all left-anchor maps
$s^o$ and equipped with an isomorphism $a^o_! \circ M \cong A$.
\end{defn}

\begin{remark}
More generally, one can replace $\kappa^{o*}B\C$ with a $\C$-module
$\M$ in the sense of Definition~\ref{C.mod.def}, and develop the
theory of $\langle S,A \rangle$-modules with values in $\M$. We will
not need this.
\end{remark}

Recall that $\Delta^o e(S)_{[0]}$ is the set $S$ itself, so that a
right $A$-module $M$ defines an object $M(s) \in B\C_{[1]} \cong \C$
for any $s \in S$; the rest of the structure equips these objects
with the right actions of $A(s,s') \in \C$. The correspondence $M
\mapsto M(s)$ is functorial in $M$, in that we have a functor
\begin{equation}\label{eps.s}
A\amod \to \C, \qquad M \mapsto M(s),
\end{equation}
where $A\amod$ is the category of right $A$-modules, and a similar
functor for left $A$-modules.

\begin{lemma}\label{yo.le}
For any small $\C$-enriched category $\langle S,A \rangle$, the
functor \eqref{eps.s} admits a left-adjoint functor $\C \to A\amod$,
$V \mapsto V_s$.
\end{lemma}

\proof{} For any small category $I$ and object $i \in I$, let
$\Delta_+^oI_i = (\Delta^oI)_i = \Delta_+^oN(I)_i$, where
$(\Delta^oI)_i$ resp.\ $N(I)_i$ are as in \eqref{C.prod}
resp.\ \eqref{ka.de.0}, and recall that if $i \in I$ is initial,
then $\lambda^{o*}\Delta_+^oI_i \cong \Delta^oI$. As in
Lemma~\ref{aug.le} and in \eqref{rho.la.C}, this gives an embedding
$\lambda(i):\Delta^oI \to \Delta_+^oI$ left-adjoint to the embedding
$\rho(i):\Delta_+^oI_i \subset \rho^{o*}\Delta^oI \to \Delta^oI$,
and the opposite functor $\lambda^o(i)$ is then right-adjoint to the
opposite functor $\rho^o(i)$.

Now recall that every $s \in e(S)$ is initial, and denote by
$A\amod_s$ the category of functors $M:\Delta^o_+e(S)_s \to
\rho^{o*}\kappa^{o*}B\C$ over $\Delta^o_+$, cocartesian over all
anchor maps and equipped an isomorphism $\rho^{o*}(a^o_! \circ M)
\cong \rho^{o*}A$. Then the adjunction map $\rho \circ \kappa = \rho
\circ \lambda \circ \rho \to \rho$ provides an identification
$\rho^{o*}\kappa^{o*}B\C \cong \C \times \rho^{o*}B\C$, with the
projection onto the second factor given by $a^o_!$, so that $M \in
A\amod_s$ is completely defined by its first component
\begin{equation}\label{M.s.eq}
M_s:\Delta^o_+e(S)_s \to \C
\end{equation}
that has to be cocartesian over all anchor maps. Since
$\Delta^o_+e(S)_s$ has a terminal object $o$, and the map $e \to o$
is an anchor map for any $e \in \Delta^o_+e(S)_s$, evaluating at $o$
gives an equivalence $A\amod_s \cong \C$. In terms of this
equivalence, the functor \eqref{eps.s} is given by the pullback
functor $\rho^o(s)^*$. Since $\rho^o(s)$ has a right-adjoint
$\lambda^o(s)$, the left Kan extension $\rho^o(s)^{\Delta^o}_!$
exists and is given by \eqref{adj.I.eq} that reads as
\begin{equation}\label{rho.del}
\rho^o(s)^{\Delta^o}_! \cong a^o_!\lambda^o(s)^*.
\end{equation}
Since $a$ is an anchor map and $A$ is a lax $2$-functor, this
implies that
\begin{equation}\label{A.rho}
\rho^o(s)^{\Delta^o}_!\rho^o(s)^*A \cong A,
\end{equation}
and then for any $V \in \C$ that corresponds to some $V_s^+$ under
the equivalence $\C \cong A\amod_s$, the functor
\begin{equation}\label{v.s}
V_s = \rho^o(s)_!^{\Delta^o}V_s^+ = a^o_!\lambda^o(s)^*V_s^+
\end{equation}
comes equipped with an isomorphism
$$
a^o_! \circ V_s \cong \rho^o(s)_!^{\Delta^o}(a^o_! \circ V_s^+) \cong
\rho^o(s)_!^{\Delta^o}\rho^o(s)^*A \cong A,
$$
thus defines an $\langle S,A \rangle$-module. Sending $V$ to $V_s
\in A\amod$ then gives our left adjoint.
\endproof

For an alternative description of right $A$-modules, consider the
nerve $N(e(S))$ of the small category $e(S)$, and note that the left
Kan extension $\kappa^o_!N(e(S))$ exists and is naturally identified
with the nerve $N(e(S)^<)$. The category $e(S)^<$ is in turn
naturally identified with $e(S_+/[1])$, where $S_+ = S \copr \{o\}$
is equipped with the map $S_+ \to [1]$ sending $o$ to $0 \in [1]$
and $S \subset S_+$ to $1 \in [1]$, so we have a projection
$\tau:e(S)^< = e(S_+/[1]) \to [1]$. This projection induces a
projection $\Delta^oe(S_+/[1]) \to \Delta^o[1] = \nat \subset \eq$,
so it makes sense to speak about $1$-special maps in
$\Delta^oe(S_+/[1])$ in the sense of
Subsection~\ref{ext.adj.subs}. Then by adjunction, right $A$-modules
$M$ correspond bijectively to $1$-special lax $2$-functors
\begin{equation}\label{M.adj}
M^\dg:\Delta^oe(S_+/[1]) \to B\C
\end{equation}
equipped with an isomorphism $M^\dg|_{\Delta^oe(S)} \cong A$. In
these terms, an element $s \in S$ with the embedding map $i_s:s \to
S$ defines a section $i_s^<:[1] = [0]^< \to e(S)^< \cong e(S_+/[1])$
of the projection $\tau$, and the functor \eqref{eps.s} is given by
the evaluation at the corresponding object $\langle [1],i_s^<
\rangle \in \Delta^o(e(S_+/[1]))$.

\begin{defn}\label{repr.def}
For any enriched category $\langle S,A \rangle$ over $\C$, a right
$A$-module $M$ is {\em representable} if $M \cong 1_s$ for some $s
\in S$, where $1 \in \C$ is the unit object, and {\em
  polyrepresentable} if $M \cong V_s$, $s \in S$, $V \in \C$, where
$V_s$ and $1_s$ are as in \eqref{v.s}. Moreover, $M$ is {\em
  ind-representable} resp.\ {\em ind-polyrepresentabe} if it is a
filtered colimit of representable resp.\ polyrepresentable modules.
\end{defn}

\begin{exa}\label{k.exa}
Assume given a commutative ring $k$, and let $\C=k\amod$ be the
category of $k$-modules, with its usual tensor structure. Then a
$\C$-enriched small category $A$ is the same thing as a $k$-linear
small category in the sense of Subsection~\ref{hom.subs}, and right
$A$-modules are $k$-linear functors $A \to k\amod$. Representable
modules correspond to representable functors. If $A$ is a
$k$-algebra considered as a $k$-linear category with one object,
then right modules are the right modules in the usual sense, and a
right module is polyrepresentable resp.\ ind-polyrepresentable iff
it is free (that is, of the form $V \otimes_k A$ for some $V \in
k\amod$) resp.\ flat. If we consider the category $P(A)$ of finitely
generated projective left modules over a flat $k$-algebra $A$, then
$P(A)$-modules are the same thing as $A$-modules, and a right
$P(A)$-module is representable resp.\ ind-representable iff the
corresponding $A$-module $M$ is finitely generated projective
resp.\ flat.
\end{exa}

We will also need a covariant version of the Yoneda Lemma similar to
\eqref{yo.cov} and \eqref{yo.M}. Namely, for any small $\C$-enriched
category $\langle S,A \rangle$, define a {\em left module} $N$ over
$\langle S,A \rangle$ as a right module over the opposite category
$\langle S,A \rangle^\iota$ of
Definition~\ref{A.opp.def}. Equivalently, $N$ is given by a functor
$N:\Delta^o e(S) \to \kappa_\iota^{o*}B\C$ over $\Delta^o$
cocartesian over all right-anchor maps and equipped with an
isomorphism $a_{\iota !}^o \circ N \cong A$. Now, the functorial
square \eqref{ka.sq} induces an equivalence $\kappa_\flat^{o*}B\C
\cong \kappa^{o*}B\C \times_{B\C} \kappa_\iota^{*o}B\C$, and for any
right $\langle S,A \rangle$-module $M$, we can consider the functor
\begin{equation}\label{A.box}
M \boxtimes_A N = \zeta(o) \circ (M \times_A N):\Delta^oe(S) \to
\kappa_\flat^{o*}B\C_{o} \cong B\C_{\kappa^o_\flat(o)} \cong \C,
\end{equation}
where we extend $\kappa_\flat$ to $\Delta^<$, and $\zeta(o)$ is the
functor \eqref{beta.eq} for the cofibration $\kappa_\flat^{o*}B\C
\to \Delta^{o>}$. We then define the tensor product $M \otimes_A N
\in \C$ as $\colim_{\Delta^oe(S)}M \boxtimes_A N$, if the colimit
exists.

\begin{lemma}\label{yo.2.le}
For any $\C$-entiched category $\langle S,A \rangle$, left $\langle
S,A \rangle$-module $N$, and ind-polyrepresentable right $\langle
S,A \rangle$-module $M$, the tensor product $M \otimes_A N$ exists,
and if $M \cong V_s$ is polyrepresentable, then $M \otimes_A N \cong
V \otimes N(s)$.
\end{lemma}

\proof{} Since colimits over $\Delta^oe(S)$ commute with filtered
colimits, it suffices to consider the case when $M = V_s$ is
polyrepresentable. By \eqref{v.s}, we have $M \cong
\rho^o(s)_!^{\Delta^o}V_s^+$, where the Kan extension functor
$\rho^o(s)_!^{\Delta^o}$ is given by \eqref{rho.del}. Since the map
$a$ is right-anchor, and both $A$ and $N$ are cocartesian over
right-anchor maps, we have $A \cong
\rho^o(s)_!^{\Delta^o}\rho^o(s)^*A$, $N \cong
\rho^o(s)_!^{\Delta^o}\rho^o(s)^*N$, and
$$
M \boxtimes_A N \cong \rho^o(s)_!^{\Delta^o}(\zeta(o) \circ (V^+_s
\times_{\rho^o(s)^*A} \rho(s)^*N)),
$$
so that
$$
M \otimes_A N \cong \colim_{\Delta^o_+e(S)_s}(\zeta(o) \circ (V^+_s
\times_{\rho^o(s)^*A} \rho(s)^*N)),
$$
Since $\Delta^o_+e(S)_s$ has a terminal object $\langle [0],s
\rangle$, the colimit exists and coincides with $N(s)$.
\endproof

\subsection{Morita $2$-categories: the construction.}

Non-empty small $\C$-en\-riched categories $\langle S,A \rangle$ form
a category that we denote by $\rCat(\C)$. In the situation of
Example~\ref{k.exa}, more is true: small $k$-linear categories and
$k$-linear functors between them form a $2$-category, and one can
define a larger ``Morita $2$-category'' by allowing general
bimodules instead of functors. To construct these $2$-categories
formally and in full generality, it is convenient to encode the
fibers of the corresponding cofibrations over $\Delta^o$ by using
the Grothendieck construction (that is, iterating the cylinder
construction of Example~\ref{cyl.exa}). This uses augmented sets of
Subsection~\ref{set.sss}.

\begin{defn}\label{J.enr.def}
For any monoidal category $\C$ and any pre-ordered set $J$, a {\em
  $J$-augmented $\C$-enriched small category} is the pair $\langle
S/J,A \rangle$ of a proper $J$-augmented set $S$ and a lax
$2$-functor $A:\Delta^oe(S/J) \to B\C$.
\end{defn}

\begin{exa}\label{pb.A.exa}
For any order-preserving map $f:J' \to J$ and small $\C$-enriched
category $\langle S/J,A \rangle$, we have the induced $J'$-augmented
set $f^*S = S \times_J^f J'$, with the natural map $f_S:f^*S \to S$,
and the induced small $\C$-enriched category $f^*\langle S/J,A
\rangle = \langle f^*S/J',\Delta^o(e(f_S))^*A \rangle$.
\end{exa}

We will be mostly interested in $[n]$-augmented enriched categories
for ordinals $[n] \in \Delta$. To consider them all at once, denote
by $\Sets' \subset \Sets$ the full subcategory of non-empty sets,
and consider the cofibrations $E\Sets' \subset E\Sets \to \Delta^o$
of Example~\ref{E.exa}. For any $[n] \in \Delta$, the fiber
$E\Sets_{[n]}$ is naturally identified with the category $\Sets/[n]$
of $[n]$-augmented sets, and $E\Sets'_{[n]} \subset E\Sets_{[n]}$ is
spanned by proper $[n]$-augmented sets. We also have the cofibration
$\phi:E\Sets_+ \to E\Sets'$ of \eqref{C.S.sq}, and for any proper
$[n]$-augmented set $S/[n]$, the right comma-fiber $S
\setminus^{\phi} E\Sets_+$ is naturally identified with the
simplicial replacement $\Delta^oe(S/[n])$ of the corresponding
category $e(S/[n])$. We can then consider the category
$e^*\phi_{**}^a\pi^*B\C$, where $\pi:E\Sets_+ \to \Delta^o$ is the
structural cofibration, $a$ is the class of cocartesian liftings of
anchor maps, and $e:E\Sets'_{\Id} \to E\Sets'$ is induced by the
embedding $\Sets'_{\Id} \to \Sets'$. By definition,
$e^*\phi_{**}^a\pi^*B\C$ is cofibered over $E\Sets'_{\Id}$, hence
also over $\Delta^o$, so we can replace it with its tightening, and
let $\wCat(\C)$ be the reduction of the resulting cofibered category
in the sense of Remark~\ref{red.rem}.

For any $[n] \in \Delta$, the objects in the fiber $\wCat(\C)_{[n]}$
are $[n]$-augmented $\C$-enriched small categories in the sense of
Definition~\ref{J.enr.def}. The transition functor $f_!^o$
corresponding to a map $f:[n'] \to [n]$ in $\Delta$ sends an
$[n]$-augmented $\C$-enriched category $\langle S/[n],A \rangle$ to
$f^*\langle S/[n],A \rangle = \langle f^*S/[n'],\Delta^o(f_S)^*A
\rangle$, where $f^*S = S \times^f_{[n]} [n']$, with the natural map
$f_S:e(f^*S/[n']) = f^*e(S/[n]) \to e(S/[n])$. The evaluation
functor \eqref{adj.2.2} induces a functor
\begin{equation}\label{ev.perp}
\ev:e^*E\Sets_+ \times_{E\Sets'_{\Id}} \wCat(\C) \to B\C
\end{equation}
over $\Delta^o$, and it is cocartesian over all anchor maps.

\begin{defn}\label{2.enr.def}
For any unital monoidal category $\C$, a {\em $\C$-enrichment} of a
$2$-category $\C'$ is the pair $\langle S,A \rangle$ of a functor
$S:\C'_{[0]} \to \Sets'$ and a lax $2$-functor $A:\C'[S] \to B\C$,
where $\C'[S]$ is the $2$-category \eqref{C.S.eq}.
\end{defn}

\begin{exa}
A $\C$-enrichment of the trivial $2$-category $\ppt^2$ is the same
thing as a small $\C$-enriched category $\langle S,A \rangle$ in the
sense of Definition~\ref{enr.def}.
\end{exa}

\begin{lemma}\label{2.enr.le}
For any unital monoidal category $\C$ and $2$-category $\C'$
equip\-ped with a $\C$-enrichment $\langle S,A \rangle$ in the sense
of Definition~\ref{2.enr.def}, there exists a unique pair of a
functor
\begin{equation}\label{Y.eq}
\Y(S,A):\C' \to \wCat(\C),
\end{equation}
cocartesian over $E\Sets_{\Id}$, and an isomorphism $\Y(S,A)^*\ev
\cong A$, where $\ev$ is the functor \eqref{ev.perp}.
\end{lemma}

\proof{} Immediately follows from \eqref{adj.2.eq}. \endproof

Now note that for any partially ordered set $J$ and $J$-augmented
$\C$-enriched category $\langle S/J,A \rangle$ in the sense of
Definition~\ref{J.enr.def}, Definition~\ref{A.mod.def} can be
repeated verbatim with $e(S)$ replaced by $e(S/J)$, so that in
particular, we have the notion of a right $\langle S/J,A
\rangle$-module $M$. For any order-preserving map $f:J' \to J$ and
right $\langle S/J,A \rangle$-module $M$, $f^*M = \Delta^o(f_S)^*M$
is a right module over $f^*\langle S/J',A \rangle$. Moreover, if $J$
has the smallest element $o \in J$, then any $s \in S_o \subset S$
in the fiber over $o \in J$ is still an initial object in $e(S/J)$,
so that Lemma~\ref{yo.le} works with the same proof, and we have the
objects $V_s$, $V \in \C$, $s \in S_o \subset S$. We can then repeat
Definition~\ref{repr.def} verbatim (but limited to $s \in S_o
\subset S$). Now, for any $j \in J'$, we have the full embedding
$\eps_j:j \setminus J \to J$ of the right comma-fiber $j \setminus
J$, and for any $j' \geq j$, we have an embedding $\eps_{j,j'}:j'
\setminus J \to j \setminus J$. Then any $j \in J$ is an initial
object in $j \setminus J$, and we can make the following.

\begin{defn}\label{A.cyl.def}
For any partially ordered set $J$, an $J$-augmented small
$\C$-en\-rich\-ed category $\langle S/J,A \rangle$ is a iterated
{\em cylinder} resp.\ {\em polycylinder} resp.\ {\em ind-cylinder}
resp.\ {\em ind-polycylinder} if for any $j \leq j' \in J$ and
representable right $\eps_j^*\langle S/J,A \rangle$-module $M$, the
right $\eps_{j'}^*\langle S/J,A \rangle$-module $\eps_{j,j'}^*M$ is
representable resp.\ polyrepresentable resp.\ ind-representable
resp.\ ind-polyrepresentable.
\end{defn}

\begin{remark}\label{poly.rem}
Explicitly, Definition~\ref{A.cyl.def} means that if $\langle S/J,A
\rangle$ is an iterated polycylinder, then for any $j \leq j' \in J$
and $s \in S_j$, we have $\eps_{j,j'}^*1_s \cong V_{s'}$ for some $V
\in \C$ and $s' \in S_{j'}$. But then for any $V' \in \C$, we have
$\eps_{j,j'}^*V'_s \cong (V' \otimes V)_s$, so that a
polyrepresentable module also restricts to a polyrepresentable
module. Moreover, $\eps_{j,j'}^*$ commutes with filtered colimits,
so that the same holds for ind-cylinders and ind-representable
modules, and then ind-polycylinders and ind-polyrepresetable
modules.
\end{remark}

\begin{exa}\label{S.cyl.exa}
  If $\C=\Sets$ is the category of sets, with the cartesian product,
  then a $J$-augmented small $\C$-enriched category is the same
  thing as a small category $A$ equipped with a functor $A \to J$
  with non-empty fibers, and it is an iterated cylinder if and only
  if the functor is a cofibration. For $J=[n]$, this is the same
  notion as in Subsection~\ref{cyl.subs}.
\end{exa}

\begin{exa}\label{M.cyl.exa}
For any small $\C$-enriched category $\langle S,A \rangle$ and right
$\langle S,A \rangle$-module $M$, the corresponding lax $2$-functor
\eqref{M.adj} defines a small $[1]$-aug\-ment\-ed $\C$-enriched
category. It is a cylinder resp.\ a polycylinder if $M$ is
representable resp.\ polyrepresentable. We also have a small
$[1]$-augmented $\C$-enriched category defined by
$\iota^*M^\dg:\Delta^oe(S_+/[1])^o \to B\C$. Since all modules over
the unit enriched category $\ppt_\C$ of Example~\ref{pt.exa} are
polyrepresentable, $\iota^*M^\dg$ is always a polycylinder; however,
it is cylinder only if $M(s) \cong 1 \in \C$ for any $s \in S$.
\end{exa}

\begin{prop}\label{enr.prop}
Let $\Cat(\C) \subset \wCat(\C)$ resp.\ $\Mor(\C) \subset \wCat(\C)$
be the full subcategories spanned by cylinders
resp.\ polycylinders. Then the induced projections
$\Cat(\C),\Mor(\C) \to \Delta^o$ are $2$-categories in the sense of
Definition~\ref{2cat.def}. Moreover, assume that $\C$ has all
filtered colimits, and the tensor product preserves them. Then the
same holds for the full subcategories $\iCat(\C)$ and $\iMor(\C)$ in
$\wCat(\C)$ spanned by ind-cylinders and ind-polycylinders.
\end{prop}

The proof of Proposition~\ref{enr.prop} is somewhat technical, so we
postpone it until Subsection~\ref{mor.pf.subs}, and first make two
remarks that do not depend on the specifics of the proof. Firstly, a
lax monoidal functor $\gamma:\C' \to \C$ between unital monoidal
categories $\C$, $\C'$ obviously induces a functor $\wCat(\C') \to
\wCat(\C)$ cocartesian over $\Delta^o$. Moreover, $\gamma$ sends
cylinders to cylinders, and ind-cylinders to ind-cylinders if it is
continuous, so that we have $2$-functors
\begin{equation}\label{ga.ca}
\gamma:\Cat(\C') \to \Cat(\C), \qquad \gamma:\iCat(\C') \to
\iCat(\C).
\end{equation}
If $\gamma$ is not just lax monoidal but monoidal, then it also
sends polycylinders to polycylinders, and ind-polycylinders to
ind-polycylinders if it is continuous, so that we have $2$-functors
\begin{equation}\label{ga.mo}
\gamma:\Mor(\C') \to \Mor(\C), \qquad \gamma:\iMor(\C') \to
\iMor(\C).
\end{equation}
Secondly, for a fixed $\C$ and any $[n] \in \Delta$,
$\C$-enrichments of $e([n]/[n]) \cong [n]$ are lax $2$-functors from
$[n]$ to $B\C$. In particular, a $2$-functor from $[n]$ to $B\C$
gives an enrichment, and these correspond to objects in
$B\C_{[n]}$. This gives a functor $B\C \to \wCat(\C)$ cocartesian
over $\Delta^o$. Any such $[n]$-augmented $\C$-enriched category is
tautologically an iterated polycylinder, so the embedding factors
through $\Mor(\C) \subset \wCat(\C)$, and we get a fully faithful
$2$-functor
\begin{equation}\label{bc.mo}
\ppt_\C:B\C \to \Mor(\C) \subset \iMor(\C).
\end{equation}
Informally, it identifies $B\C$ with the full $2$-subcategory in
$\Mor(\C)$ spanned by the unit enriched category $\ppt_\C$ of
Example~\ref{pt.exa}.

\subsection{Morita $2$-categories: the proof.}\label{mor.pf.subs}

We now prove Proposition~\ref{enr.prop}. For any $[n] \in \Delta$,
$[n]$-augmented $\C$-enriched category $\langle S/[n],A \rangle$,
and $l \in [n]$, we have the $\C$-enriched category $\langle S_l,A_l
\rangle = a_l^*\langle S/[n],A \rangle$, where $a_l:[0] \to [n]$ is
the embedding onto $l$. Assume given such a $\langle S/[n],A
\rangle$, and consider the embedding $t:[\nm] \to [n]$.

\begin{lemma}\label{poly.le}
The following conditions are equivalent.
\begin{enumerate}
\item $\langle S/[n],A \rangle$ is an iterated cylinder,
  resp.\ polycylinder, resp.\ ind-poly\-cy\-linder,
  resp.\ ind-polycylinder. 
\item $t^*\langle S/[n],A \rangle$ is an iterated cylinder,
  resp.\ polycylinder, resp.\ ind-poly\-cy\-linder,
  resp.\ ind-polycylinder, and for any representable $\langle
  S/[n],A \rangle$-module $1_s$, $s \in S_0$, the module
  $M_s=t^*1_s$ is representable, resp.\ poly\-representable,
  resp.\ ind-representable, resp.\ ind-polyrepresentable.
\end{enumerate}
\end{lemma}

\proof{} For any order-preserving map $f:J' \to J$, the pullback
$f^*$ of Example~\ref{pb.A.exa} obviously preserves the cylinder
conditions of Definition~\ref{A.cyl.def}, so
\thetag{i}$\Rightarrow$\thetag{ii} is obvious. In the other
direction, we have to check a condition for any $j \leq j' \in
[n]$. If $j \geq 1$, it follows from the first part of \thetag{ii},
if $j=0$ and $j'=1$, it is the second part on the nose, and if $j =
0$ but $j' \geq 2$, combine both parts and Remark~\ref{poly.rem}.
\endproof

To see explicitly the modules $M_s$ of
Lemma~\ref{poly.le}~\thetag{ii}, it is convenient to use the
cofibration \eqref{mu.adj.eq} associated to the left-closed
embedding $e(S_0) \subset e(S/[n])$. Its fibers are given by
\eqref{mu.fib.eq}, and in particular, for any $s \in S_0$, the fiber
over $\langle [0],s \rangle$ is naturally identified with
$\Delta^{o>}e(t^*S/[\nm])$. This gives an embedding
$j_s:\Delta^oe(t^*S/[\nm]) \to \Delta^oe(S/[n])$ that fits into a
commutative square
\begin{equation}\label{j.s.sq}
\begin{CD}
\Delta^oe(t^*e(S/[\nm])) @>{j_s}>> \Delta^oe(S/[n])\\
@VVV @VVV\\
\Delta^o @>{\kappa^o}>> \Delta^o,
\end{CD}
\end{equation}
where the vertical arrows are structural cofibrations. We then have
\begin{equation}\label{l.r.M}
M_s \cong j_s^*A,
\end{equation}
where $j_s^*A:\Delta^oe(t^*S/[\nm]) \to \kappa^{o*}B\C$ is induced
by the lax $2$-functor $A$.

\medskip

Now consider a square \eqref{seg.sq} and the category $\Delta^o[n]
\cong (\Delta/[n])^o$, and let $\Delta^o[n]^l \subset \Delta^o[n]$
be the full subcategory spanned by maps $f:[m] \to [n]$ whose image
contains $l \in [n]$. Moreover, \eqref{seg.sq} induces left-closed
full embeddings $\Delta^o[l] \to \Delta^o[n]$, $\Delta^o[\nl] \to
\Delta^o[n]$; let $\Delta^o[n]^{\leq l},\Delta^o[n]^{\geq l} \subset
\Delta^o[n]$ be their essential images, and let
\begin{equation}\label{de.l}
\Delta^o[n]^l_0 = \Delta^o[n]^{\leq l} \cup \Delta^o[n]^{\geq l},
\qquad \Delta^o[n]^l_1 = \Delta^o[n]^l_0 \cup \Delta^o[n]^l.
\end{equation}
For any $[n]$-augmented set $S/[n]$, the simplicial replacement
$\Delta^oe(S/[n])$ comes equipped with a projection
$\pi:\Delta^oe(S/[n]) \to \Delta^o[n]$, and we can let
$\Delta^oe(S/[n])^l = \Delta^o[n]^l \times_{\Delta^o[n]}
\Delta^oe(S/[n])$, and similarly for $\Delta^oe(S/[n])^{\leq l}$,
$\Delta^oe(S/[n])^{\geq l}$, $\Delta^oe(S/[n])^l_0$, and
$\Delta^oe(S/[n])^l_1$. If we denote $S^{\leq l} = s^*S$, $S^{\geq
  l} = t^*S$, then we have $\Delta^oe(S/[n])^{\leq l} \cong
\Delta^oe(S^{\leq l}/[l])$, $\Delta^oe(S/[n])^{\geq l} \cong
\Delta^oe(S^{\geq l}/[\nl])$, and
\begin{equation}\label{cap.l}
\Delta^oe(S_l) \cong \Delta^oe(S/[n])^{\leq
  l} \cap \Delta^oe(S/[n])^{\geq l} \subset \Delta^oe(S/[n])^l_0.
\end{equation}
We also have full embeddings
\begin{equation}\label{a.b.eq}
\begin{CD}
\Delta^o e(S/[n])^l_0 @>{\alpha}>> \Delta^oe(S/[n])^l_1 @>{\beta}>>
\Delta^oe(S/[n])
\end{CD}
\end{equation}
over $\Delta^o$. The projection $\pi:\Delta^oe(S/[n])^l_0 \to
\Delta^o$ is a cofibration but not a $2$-category. The projection
$\pi:\Delta^oe(S/[n])^l_1 \to \Delta^o$ is not even a cofibration,
but it is a cofibration over all anchor maps, and the embeddings
\eqref{a.b.eq} are cocartesian over anchor maps. By abuse of
terminology, we will use the term {\em lax $2$-functor} to also
signify a functor from $\Delta^oe(S/[n])^l_0$ or
$\Delta^oe(S/[n])^l_0$ over $\Delta^o$ that is cocartesian over
anchor maps. Note that by Example~\ref{MV.exa}, giving such a lax
$2$-functor $A_0:\Delta^oe(S/[n])^l_0 \to B\C$ is equivalent to
giving lax $2$-functors $A^{\leq l}_0:\Delta^oe(S/[n])^{\leq l} \to
B\C$, $A^{\geq l}_0:\Delta^oe(S/[n])^{\geq l} \to B\C$ and an
isomorphism between their restrictions to the subcategory
\eqref{cap.l}.

\begin{lemma}\label{enr.le}
\begin{enumerate}
\item In the notation and assumptions above, for any lax $2$-functor
  $A_0:\Delta^oe(S/[n])^l_0 \to B\C$, the right Kan extension
  $\alpha_*A_0$ exists and is a lax $2$-functor, while for any lax
  $2$-functor $A_1:\Delta^oe(S/[n])^l_1 \to B\C$, the adjunction map
  $A_1 \to \alpha_*\alpha^*A_1$ is an isomorphism.
\item Assume that $\langle s^*S/[l],A^{\leq l}_0 \rangle$ is an
  iterated ind-polycylinder, and let $A_1=\alpha_*A_0$. Then the
  left Kan extension $\beta^{\Delta^o}_!A_1$ exists and is a lax
  $2$-functor.
\item Finally, for any iterated ind-polyclynider $\langle S/[n],A
  \rangle$, the adjunction map $\beta^{\Delta^o}_!\beta^*A \to A$ is
  an isomorphism
\end{enumerate}
\end{lemma}

\proof{} For \thetag{i}, note that for any $[m] \in \Delta$ and
anchor map $f:[m] \to [n]$, the embedding $\Delta^{o>}e(f^*S/[m])
\to \Delta^{o>}e(S/[n])$ admits a left-adjoint functor $\mu(f)$
provided by Lemma~\ref{mu.le}. In particular, this applies to the
maps $s$, $t$, $a = s \circ t = t \circ s$ in \eqref{seg.sq}, and if
we let $\mu^{\leq l} = \mu(s)$, $\mu^{\geq l}=\mu(t)$, $\mu^l=\mu(a)$,
then any object $x \in \Delta^{o>}e(S/[n])$ fits into a commutative
square
\begin{equation}\label{X1.sq}
\begin{CD}
x @>>> \mu''(x)\\
@VVV @VVV\\
\mu'(x) @>>> \mu^l(x).
\end{CD}
\end{equation}
If $x$ lies in $\Delta^o e(S/[n])^l_1 \setminus \Delta^o
e(S/[n])_0$, then \eqref{X1.sq} is actually a square in
$\Delta^oe(S/[n])^l_1 \subset \Delta^{o>}e(S/[n])$, and all its maps
are anchor maps. Moreover, the oppposite square defines a standard
functor $c:\V \to \Delta e(S/[n])^l_0$ in the sense of
Remark~\ref{st.rem} equipped with an exact $x$-augmentation. The
corresponding embedding $c:\V \to \Delta e(S/[n])^l_0/x$ of
\eqref{aug.lft} is left-ad\-mis\-sible, so we obtain a framing for
the functor $\alpha^o$ in the sense of Lemma~\ref{kan.le}, and if we
use this framing to compute $A_1 = \alpha_*A_0$, we see that it
exists by Remark~\ref{st.rem}. Moreover, the square \eqref{X1.sq}
provides an isomorphism
\begin{equation}\label{A1.eq}
A_1(x) \cong A^{\leq l}_0(\mu^{\leq l}(x)) \times_{A^l_0(\mu^l(x))}
A^{\geq l}_0(\mu^{\geq l}(x)),
\end{equation}
where $A^l:\Delta^oe(S_l) \to B\C$ is the common restriction of
$A^{\leq l}_0$ and $A^{\geq l}_0$. The square \eqref{X1.sq} is
functorial in $x$, and since $\mu^l$, $\mu^{\leq l}$, $\mu^{\geq l}$
send anchor maps to anchor maps, $A_1$ is cocartesian over anchor
maps. Again by Remark~\ref{st.rem}, it is also unique with this
property.

For \thetag{ii}, let $L$ resp.\ $R$ be the classes of those maps in
$\Delta^o[n] \cong (\Delta/[n])^o$ that are bijective over $l \in
[n]$ resp.\ $[n] \setminus \{l\} \subset [n]$, and note that by
Example~\ref{facto.exa}, $\langle L,R \rangle$ is a factorization
system on $\Delta^o[n]$. Moreover, for any $x \in \Delta^oe(S/[n])
\setminus \Delta^oe(S/[n])^l$, the functor $\mu^l$ induces an
equivalence $\Delta^oe(S_l) \cong \Delta^oe(S/[n])^l /_{\pi^*R} x$,
and $\Delta^oe(S/[n])^l /_{\pi^*R} x \subset \Delta^oe(S/[n])^l / x$
is left-admissible by Example~\ref{comma.exa}, so that we obtain a
collection of functors
\begin{equation}\label{sl.eps}
\eps(x):\Delta^{o}e(S_l) \to \Delta^oe(S/[n]), \qquad x \in
\Delta^oe(S/[n]) \setminus \Delta^oe(S/[n])^l
\end{equation}
equipped with admissible augmentations $\eps(x)_>$ sending $o$ to
$x$. By Remark~\ref{fr.rem}, this gives a framing for the embedding
$\Delta^oe(S/[n])^l \subset \Delta^oe(S/[n])$. Moreover, if $x
\not\in \Delta^oe(S/[n])^l_1$, then $\Delta^oe(S/[n])^l_1 / x =
\Delta^oe(S/[n])^l /x$, so that we also obtain a framing for the
embedding $\beta$. Then by Lemma~\ref{kan.le}~\thetag{i} and
\eqref{kan.I.eq}, we have
\begin{equation}\label{be.co}
\beta^{\Delta^o}_!A_1(\langle [m],e_\idot \rangle) =
\colim_{\Delta^oe(S_l)}A_1^\eps(\langle [m],e_\idot \rangle),
\end{equation}
for any $\langle [m],e_\idot \rangle \in \Delta^oe(S/[n]) \setminus
\Delta^oe(S/[n])^l_1$, where we denote
\begin{equation}\label{A.e.m}
A_1^\eps(\langle [n],e_\idot \rangle = \zeta([m]) \circ \eps(\langle
[m],e_\idot \rangle)^*A_1:\Delta^oe(S_l) \to B\C_{[m]},
\end{equation}
and we have to prove that the colimits in \eqref{be.co} exist.

Indeed, by \eqref{c.m.eq}, we have $B\C_{[m]} \cong \C^m$, and it
suffices to prove that colimits in \eqref{be.co} exist after
projecting to each factor $\C$. The projections are given by the
transition functors with respect to anchor maps $a:[1] \to [m]$, and
since $A_1$ is a lax $2$-functor, we have $a^o_! \circ
A_1^\eps(\langle [m],e_\idot \rangle \cong A_1^\eps(\langle
[1],a^*e_\idot \rangle)$, so may assume right away that
$[m]=[1]$. Moreover, the map $f=\pi \circ e_\idot:[m] \to [n]$
factors through $[n'] = f([m]) \cup \{l\} \subset [n]$, and since
iterated ind-polycylinders are preserved by pullbacks, we may replace
$[n]$ with $[n']$ without changing anything, so that we may further
assume that $[n]=[2]$. If $l=0$ or $l=2$, $A_1^\eps(\langle
[1],e_\idot \rangle)$ is constant, and since the simplicial set
$N(e(S_l))$ is contractible, there is nothing to prove. If $l=1$,
then by \eqref{l.r.M} and \eqref{A1.eq}, we have
$$
A_1^*\eps(\langle [1],e_\idot \rangle) \cong M_s \boxtimes_{A^l} N,
$$
where $N$ is a left $\langle S_l,A^l \rangle$-module, $s = e_0$, and
$M_s$ is the right $\langle S_l,A_l \rangle$-module of
Lemma~\ref{poly.le}~\thetag{ii} for $\langle S^{\leq 1},A^{\leq 1}
\rangle$. Then the colimit exists by Lemma~\ref{yo.2.le}.

For \thetag{iii}, by Lemma~\ref{kan.le}~\thetag{ii}, we have to
check that for any $\langle [m],e_\idot \rangle$ in the complement
$\Delta^oe(S/[n]) \setminus \Delta^oe(S/[n])^l_1$, the augmented
functor
\begin{equation}\label{A.e}
A^\eps(\langle [m],e_\idot \rangle) = \zeta([m]) \circ \eps(\langle
[m],e_\idot \rangle)_>^*A:\Delta^{o>}e(S_l) \to B\C_{[m]}
\end{equation}
is exact. It suffices to prove that it is a filtered colimit of
functors contractible in the sense of Definition~\ref{aug.def}. We
again can project to components of the decomposition $B\C_{[m]}
\cong \C^m$, and all projections but one are constant, thus extend
to any contraction of $N(e(S_l))$. Thus as in the proof of
\thetag{ii}, we may assume that $[m]=[1]$, $[n]=[2]$ and $l=1$. Then
if we again let $s=e_0$, the embedding $\eps(\langle
[1],e_\idot\rangle)$ factors through the embedding
$j_s:\Delta^oe(S^{\geq 1}/[1]) \to \Delta^oe(S/[2])$ of
\eqref{l.r.M}, and by Lemma~\ref{poly.le}~\thetag{ii}, the
$\Delta^oe(S^{\geq 1}/[1])$-module $M_s$ is
ind-polyrepresentable. If $M \cong V_{s'}$, $V \in \C$, $s' \in S_1$
is polyrepresentable, then by \eqref{v.s} and \eqref{l.r.M},
$\zeta([1]) \circ j_s^*A$ extends to the category
$\Delta^o_+e(S^{\geq 1}/[1])_{s'}$, and then $A^\eps(\langle
[1],e_\idot \rangle)$ extends to $\Delta^o_+e(S_1)_{s'}$, so that it
is a contractible functor. In the general case, it is a filtered
colimit of such.
\endproof

\begin{corr}\label{A.cyl.corr}
Assume given an iterated ind-polycylinder $\langle S/[n],A \rangle$,
and let $\Delta_p^oe(S/[n]) \subset \Delta_+^oe(S/[n])$ be the
opposite categories to the categories of \eqref{del.p.pm.I}, with
the embedding functor $\beta:\Delta_p^oe(S/[n]) \to
\Delta_+^oe(S/[n])$. Then the adjunction map
$\beta^{\Delta^o}_!\beta^*A \to A$ is an isomorphism.
\end{corr}

\proof{} Use the same induction as in the first claim of
Lemma~\ref{cyl.le}, and replace its second claim with
Lemma~\ref{enr.le}.
\endproof

\proof[Proof of Proposition~\ref{enr.prop}.] Since all the cylinder
conditions of Definition~\ref{A.cyl.def} are stable by pullbacks,
all the subcategories in question are subcofibrations by
Lemma~\ref{fib.le}~\thetag{i}. Their fibers over $[0] \in \Delta^o$
are discrete, so it suffices to check the Segal condition. For any
square \eqref{seg.sq}, any triple of augmented sets $S^{\leq
  l}/[l]$, $S^{\geq l}/[\nl]$ and a set $S_l$ equipped with
isomorphisms $S_l \cong S^{\leq l}_l \cong S^{\geq l}_0$ define an
$[n]$-augmented set $S = S^{\leq l} \copr_{S_l} S^{\geq l}$. Then
lax $2$-functors $A^{\leq l}:\Delta^oe(S^{\leq l}/[l])$, $A^{\geq
  l}:\Delta^oe(S^{\geq l}/[l])$ together with an isomorphism between
their restrictions to the subcategory \eqref{cap.l} define a lax
$2$-functor $A_0:\Delta^oe(S/[n])^l_0 \to B\C$ that uniquely extends
to a lax $2$-functor $A_1:\Delta^oe(S/[n])^l_1 \to B\C$ by
Lemma~\ref{enr.le}~\thetag{i}. We have to shows that if $\langle
S^{\leq l}/[l],A^{\leq l} \rangle$ and $\langle S^{\geq
  l}/[l],A^{\geq l} \rangle$ are iterated cylinders
resp.\ polycylinders resp.\ ind-cylinders resp.\ ind-polycylinders,
then $A_1$ uniquely extends to a lax $2$-functor $A:\Delta^oe(S/[n])
\to B\C$ such that $\langle S/[n],A \rangle$ is an iterated cylinder
resp.\ polycylinder resp.\ ind-cylinder resp.\ ind-polycylinder.

To do this, we note that firstly, uniqueness follows from
Lemma~\ref{enr.le}, and secondly, by induction, it suffices to
consider the case $l=1$. Then by Lemma~\ref{poly.le}, it suffices to
assume that $\langle S^{\leq 1}/[1],A^{\leq 1} \rangle$ is a
cylinder resp.\ polycylinder resp.\ ind-cylinder
resp.\ ind-polycylnder, take the lax $2$-functor
$A=\beta_!^{\Delta^o}A_1$ provided by
Lemma~\ref{enr.le}~\thetag{ii}, and prove that for any $s \in S_0$,
the $\langle S^{\geq l}/[\nm],A^{\geq 1}\rangle$-module $M_s$ of
\eqref{l.r.M} is representable resp.\ polyrepresentable
resp.\ ind-representable resp.\ ind-polyrepresentable.

Indeed, to simplify notation, let $S' = S^{\geq 1}$, $A' = A^{\geq
  1}$, so that $S'_0 \cong S_1$, and consider the subcategory
$$
\Delta^oe(S'/[\nm])^0 = \Delta^oe(S'/[\nm]) \cap
\Delta^oe(S/[n])^l_1 \subset \Delta^oe(S'/[\nm]),
$$
with the embedding functor $\beta:\Delta^oe(S'/[\nm])^0 \to
\Delta^oe(S'/[\nm])$. Then \eqref{sl.eps} also gives a framing for
$\beta$, and the functor $j_s$ of \eqref{j.s.sq} intertwines $\beta$
with the embedding $\beta$ of \eqref{a.b.eq} and identifies the
framings, so that we have $j_s^* \circ \beta^{\Delta^o}_!  \cong
\beta^{\Delta^o}_! \circ j_s^*$. Therefore by \eqref{l.r.M} and
\eqref{A1.eq}, the $\langle S'/[\nm],A' \rangle$-module $M_s$ is
given by
\begin{equation}\label{m.s.a}
M_s = \beta^{\Delta^o}_!(\mu^*(M^{\leq 1}_s) \times_{\mu^*A^1} A'),
\end{equation}
where the projection $\mu:\Delta^{o>}e(S'/[\nm]) \to
\Delta^{o>}e(S_1)$ is left-adjoint to the embedding
$\Delta^{o>}e(S_1) \cong \Delta^{o>}e(S'_0) \to
\Delta^{o>}e(S'/[\nm])$, and $M^{\leq 1}_s$ is the $\langle S_1,A^1
\rangle$-module of Lemma~\ref{poly.le} for the enriched category
$\langle S^{\leq 1}/[1],A^{\leq 1}\rangle$. If this category is a
polycylinder, then $M^{\leq 1}_s \cong V_{s'} \cong
\rho^o(s')_!^{\Delta^o}V^+_s$ for some $V \in \C$, $s' \in S_1 \cong
S'_0$, and then by the same argument as in Lemma~\ref{yo.2.le},
\eqref{m.s.a} yields
$$
M_s \cong \rho^o(s')_!^{\Delta^o}V^+_{s'} \cong V_{s'} \in A^{\geq
  1}\amod,
$$
where $s'$ is now considered as an initial object in
$e(S'/[\nm])$. Therefore $M_s$ is polyrepresentable as required. If
$\langle S^{\leq 1}/[1],A^{\leq 1} \rangle$ is a cylinder, then $V
\cong 1$, so $M_s$ is representable, and in the ind-cases, note that
left Kan extensions commute with filtered colimits.
\endproof

\subsection{Functors and cylinders.}\label{mor.adj.subs}

Let us now describe the relationship between the $2$-categories of
Proposition~\ref{enr.prop} and the category $\rCat(\C)$ of small
$\C$-enriched categories and functors between them.

Assume given $\C$-enriched categories $\langle S_0,A_0 \rangle$,
$\langle S_1,A_1 \rangle$ and a functor $\langle f,g \rangle:\langle
S_0,A_0 \rangle \to \langle S_1,A_1 \rangle$, and let $S_{01} = S_0
\copr S_1$, with the projection $S_{01} \to [1]$ sending $S_l$ to
$l$, $l=0,1$. Then the cylinder $\Cyl(e(f))$ of the functor
$e(f):e(S_0) \to e(S_1)$ is naturally identified with
$e(S_{01}/[1])$, and the $2$-cylinder $\Cyl^2(e(f),g)$ of
Lemma~\ref{2.cyl.le} defines a $[1]$-augmented $\C$-enriched
category
\begin{equation}\label{cyl.f.g}
\Cyl(f,g) = \langle S_{01}/[1],\Cyl^2(e(f),g)\rangle.
\end{equation}
By \eqref{cyl.fr.sq} and \eqref{l.r.M}, for any $s \in S_0$, the
$\langle S_1,A_1 \rangle$-module $M_s$ of Lemma~\ref{poly.le}
corresponding to \eqref{cyl.f.g} is isomorphic to the representable
module $1_{f(s)}$, so that $\Cyl(f,g)$ is a cylinder in the sense of
Definition~\ref{A.cyl.def}. This suggests that there exists a
$2$-functor
\begin{equation}\label{cat.cat}
\Delta^o\rCat(\C) \to \Cat(\C)
\end{equation}
sending $\langle S,A \rangle$ to itself and $\langle f,g \rangle$ to
$\Cyl(f,g)$.

To construct such a $2$-functor, let $S:\rCat(\C) \to \Sets$ be the
forgetful functor sending $\langle S,A \rangle$ to $S$, let
$\rCat^\hdot(\C) \to \rCat(\C)$ be the cofibration with fibers
$e(S)$, and let $\rCat^\hdot_\Delta(\C) \to \rCat(\C)$ be the
cofibration with fibers $\Delta^oe(S)$ (equivalently,
$\rCat^\hdot_\Delta(\C) \to \rCat(C) \times \Delta^o$ is the
discrete cofibration corresponding to $\eps_*S:\rCat(\C) \times
\Delta^o \to \Sets$, where $\eps:\rCat(\C) \to \rCat(\C) \times
\Delta^o$ is the embedding onto $\rCat(\C) \times [0]$). Then the
functors $A$ for various $\langle S,A \rangle \in \rCat(\C)$
together define a single functor $A:\rCat^\hdot_\Delta(\C) \to B\C$
over $\Delta^o$ that is cocartesian over anchor maps. On the other
hand, the functors \eqref{xi.eq} together define a functor
$\xi:\rCat^\hdot_\Delta(\C) \to \rCat^\hdot(\C)$ cocartesian over
$\rCat(\C)$, and then as in Lemma~\ref{2.cyl.le}, \eqref{adj.2.eq}
converts it to a functor
\begin{equation}\label{alha.rcat}
\alpha:\rCat^\hdot_\Delta(\C) \to \Id^+_{**}(\rCat^\hdot \times
\Delta^o) = \Delta^o\langle \rCat^\hdot \rangle
\end{equation}
over $\Delta^o$.

\begin{lemma}
The right Kan extension $\alpha_*A$ with respect to
\eqref{alha.rcat} exists, its restriction to $\Delta^o\rCat[S] =
\Delta^o\rCat^\hdot \subset \Delta^o\langle \rCat^\hdot \rangle$
defines a $\C$-enrichement of the $2$-category $\Delta^o\rCat$, and
the corresponding functor \eqref{Y.eq} factors through a $2$-functor
\eqref{cat.cat}.
\end{lemma}

\proof{} Same as Lemma~\ref{2.cyl.le}. \endproof

Since an iterated cylinder is tautologically an iterated
polycylinder, we have $\Cat(\C) \subset \Mor(\C) \subset \iMor(\C)$,
so that \eqref{cat.cat} also defines a $2$-functor from $\rCat(\C)$
to $\Mor(\C)$ and $\iMor(\C)$. One difference between $\Cat(\C)$ and
the Morita $2$-categories is that in the latter, many morphisms
become reflexive. Namely, recall that for any $\C$-entiched category
$\langle S,A \rangle$, we have the opposite $\C^\iota$-enriched
category $\langle S,A \rangle^o = \langle S,A^\iota \rangle$, and
note that a functor $\langle f,g \rangle:\langle S_0,A_0 \rangle \to
\langle S_1,A_1 \rangle$ defines an opposite functor $\langle
f,g^\iota \rangle$ between the opposite categories. Then
\eqref{cyl.f.g} gives an $[1]$-augmented $\C^\iota$-enriched
category $\Cyl(f,g^\iota)$, and if we let $S_{10}/[1]$ be $S_{01}$
whose augmentation map $S_{01} \to [1]$ is composed with the
order-reversing isomorphism $[1]^o \cong [1]$, then
\begin{equation}\label{cyl.f.g.i}
\Cyl^\iota(f,g) = \Cyl(f,g^\iota)^\iota = \langle
S_{10}/[1],\Cyl^2(e(f),g^\iota) \rangle
\end{equation}
is a $[1]$-augmented $\C$-enriched category that we call the {\em
  dual cylinder} of the functor $\langle f,g \rangle$.

\begin{defn}\label{A.refl.def}
A functor $\langle f,g \rangle:\langle S_0,A_0 \rangle \to \langle
S_1,A_1 \rangle$ between small $\C$-enriched categories is {\em
  reflexive} resp.\ {\em ind-reflexive} if its dual cylinder
\eqref{cyl.f.g.i} is a polycylinder resp.\ an ind-polycylinder.
\end{defn}

\begin{exa}\label{i.s.exa}
By Example~\ref{M.cyl.exa}, the functors $i_s$ of
Example~\ref{pt.exa} are always reflexive in the sense of
Definition~\ref{A.refl.def}.
\end{exa}

\begin{exa}
If $A=k\amod$, as in Example~\ref{k.exa}, and $S_0=S_1=\ppt$, so
that $A_0$ and $A_1$ are $k$-algebras and $g:A_0 \to A_1$ is an
algebra map, then $g$ is reflexive resp.\ ind-reflexive if $A_1$ is
free resp.\ flat as a module over $A_0$.
\end{exa}

\begin{prop}\label{refl.prop}
Assume given a functor $\langle f,g \rangle:\langle S_0,A_0 \rangle
\to \langle S_1,A_1 \rangle$ between small $\C$-enriched categories
that is reflexive resp.\ ind-reflexive in the sense of
Definition~\ref{A.refl.def}. Then the morphism $\Cyl(f,g)$ is
reflexive in $\Mor(\C)$ resp.\ $\iMor(\C)$.
\end{prop}

Before proving Proposition~\ref{refl.prop} in general, let us
consider the situation of Example~\ref{i.s.exa}: we have a small
$\C$-enriched category $\langle S,A \rangle$ and an object $s \in
S$, and we want to extend the functor $i_s:\ppt_\C \to \langle S,A
\rangle$ to an adjoint pair
\begin{equation}\label{ga.s.S}
\gamma_s:\adj \to \Mor(\C).
\end{equation}
We start by recalling the description of the $2$-category $\eq$ in
terms of injective maps $a:[n]_1 \to [n]$. For any such map $a$, the
pullback functor $a^*:\Sets/[n] \to \Sets/[n]_1$ has a fully
faithful right-adjoint $a_*:\Sets/[n]_1 \to \Sets/[n]$ sending $S
\in \Sets/[n]_1$ to $a_*S \in \Sets/[n]$ with $(a_*S)_i = S_i$ for
$i \in [n]_1$ and $(a_*S)_i = \ppt$ for $i \in [n]_0$. For any $S
\in \Sets/[n]$, we then have the $[n]$-augmented set $S_a = a_*a^*S$
and the adjunction map $p:S \to S_a$. Moreover, let
$\Delta_0^oe(S_a/[n]) \subset \Delta^oe(S_a/[n])$ be the full
subcategory spanned by objects $\langle [m],s \rangle$ such that the
composition $[m] \to S_a \to [n]$ is $0$-special. Then by virtue of
\eqref{01.facto}, the embedding $\delta:\Delta_0^oe(S_a/[n]) \to
\Delta^oe(S_a/[n])$ admits a right-adjoint functor
$\delta^\dg:\Delta^oe(S_a/[n]) \to \Delta^o_0e(S_a/[n])$. Therefore
for any lax $2$-functor $A:\Delta^oe(S_a/[n]) \to B\C$, the Kan
extension
\begin{equation}\label{A.a}
A_a = \delta^{\Delta^o}_!\delta^*A:\Delta^oe(S_a/[n]) \to B\C
\end{equation}
exists and is given by \eqref{adj.I.eq}, so that it is a lax
$2$-functor. In particular, if we have an $[n]$-augmented
$\C$-enriched category $\langle S/[n],A \rangle$ and a section
$i:S_a \to S$ of the adjunction map $p:S \to S_a$, then $A_i =
\Delta^oe(i)^*A$ is a lax $2$-functor, and we have the
$[n]$-augmented $\C$-enriched category $\langle S_a/[n],A_{i,a}
\rangle$.

\begin{lemma}\label{S.f.le}
In the assumptions above, if $\langle S/[n],A \rangle$ is an
iterated polycylinder in the sense of Definition~\ref{A.cyl.def},
then so is $\langle S_a/[n],A_{i,a} \rangle$.
\end{lemma}

\proof{} If $n=1$, the claim immediately follows from
Example~\ref{M.cyl.exa}. If $n \geq 2$, choose some $l \in [n]$, $0
< l < n$, consider the square \eqref{seg.sq}, let
$$
[l]_1 = a([n]_1) \cap s([l]) \subset [n], \qquad 
[\nl]_1 = a([n]_1) \cap t([\nl]) \subset [n],
$$
with the embedding maps $a':[l]_1 \to [l]$, $a'':[\nl]_1 \to [\nl]$,
and let $S' = s^*S$, $S'' = t^*S$, with the induced maps $i':S'_{a'}
\to S'$, $i'':S''_{a''} \to S''$. Then if we let $A' = s^*A$, $A'' =
t^*A$, we have
$$
s^*\langle S_a/[n],A_{i,a} \rangle \cong \langle
S'_{a'}/[1],A'_{i',a'} \rangle, \quad t^*\langle S_a/[n],A_{i,a}
\rangle \cong \langle S''_{a''}/[\nm],A''_{i'',a''} \rangle
$$
by \eqref{adj.I.eq}, and both these augmented $\C$-enriched
categories are iterated polycylinders by induction. As in the proof
of Proposition~\ref{enr.prop}, they define an iterated polycylinder
$\langle S_a,A' \rangle$, and we need to check that $A_{i,a} \cong
A'$. Moreover, if we consider the factorization \eqref{a.b.eq} for
the category $\Delta^oe(S_a/[n])$, then $A_{i,a} \cong A'$ on
$\Delta^oe(S_a/[n])^l_0$, and hence also on $\Delta^oe(S_a/[n])^l_1$
by virtue of Lemma~\ref{enr.le}~\thetag{i}. Thus by
Lemma~\ref{enr.le}~\thetag{iii}, it suffices to check that the
adjunction map $\beta^{\Delta^o}_!\beta^*A_{i,a} \to A_{i,a}$ is an
isomorphism. To do this, use the framing \eqref{sl.eps} for $\beta$,
and note that if $l \in [n] \setminus [n]_1$, then $S_{a,l} \cong
S_l$, and $\Delta^oe(i)$ identifies the framings, so that
$\delta^{\Delta^o}_!\delta^*\Delta^oe(i)^*$ commutes with
$\beta^{\Delta^o}_!$, and we are then done by
Lemma~\ref{enr.le}~\thetag{iii}. On the other hand, if $l \in
[n]_1$, then $S_{a,l} = \ppt$ is the single point, and for any
$\langle [m],e_\idot \rangle \in \Delta^oe(S_a/[n]) \setminus
\Delta^oe(S_a/[n])^l_1$, the augmented functor \eqref{A.e} is
constant, hence exact.
\endproof

We can now construct the adjoint pair \eqref{ga.s.S} for any small
$\C$-enriched category $\langle S,A \rangle$ and element $s \in
S$. To do this, let $S:\adj_{[0]} = \eq_{[0]} \to \Sets$ be the
constant functor with value $S$, let $S_+:\eq_{[0]} = \{0,1\} \to
\Sets$ be the functor sending $1$ to $S$ and $0$ to the point, and
let $i^+_s:S_+ \to S$ be the map equal to $\id$ at $1$ and to the
embedding $i_s:\ppt \to S$ onto $s$ at $0$, with the corresponding
$2$-functor $\iota_s:\eq[S_+] \to \eq[S]$. Moreover, let
$\eta:\eq[S_+] \to \adj[S_+]$ be the lax $2$-functor induced by
\eqref{eta.eq}. Since as we saw in the proof of Lemma~\ref{adj.le},
$\eta$ admits a right-adjoint, the Kan extension
$\eta^{\Delta^o}_!E$ exists for any functor $E:\eq[S_+] \to B\C$
over $\Delta^o$, and is given by \eqref{adj.I.eq}. In particular, we
can take the constant $2$-functor $\eq \to \Mor(\C)$ with value
$\langle S,A \rangle$, with the corresponding constant
$\C$-enrichment $E:\eq[S] \to B\C$ of the $2$-category $\eq$, and
consider the functor
\begin{equation}\label{e.pl}
E_+ = \eta^{\Delta^o}_!\iota_s^*E:\adj[S_+] \to B\C.
\end{equation}
Then \eqref{adj.I.eq} immediately shows that this is a lax
$2$-functor, so that $\langle S_+,E_+ \rangle$ is a $\C$-enrichment
of the $2$-category $\adj$. Moreover, for any object $c \in \eq
\subset \adj$ represented by an injective map $a:[n]_1 \to [n]$, the
right comma-fiber $c \setminus \adj$ of the cofibration $\adj[S_+]
\to \adj$ is identified with $\Delta^oe((S \times [n])_a/[n])$, and
the projections \eqref{comma.pr} for the cofibrations $\eq[S_+] \to
\eq$, $\adj[S_+] \to \adj$ fit into a cartesian diagram
\begin{equation}\label{e.dia}
\begin{CD}
\Delta^o_0e((S \times [n])_a) @>{\delta}>> \Delta^oe((S \times
      [n])_a)\\
@V{p_c}VV @VV{p_c}V\\
\eq[S_+] @>{\eta}>> \adj[S_+].
\end{CD}
\end{equation}
Then \eqref{adj.I.eq} shows that the base change map
$\delta_!^{\Delta^o} \circ p_c^* \to p_c^* \circ \eta_!^{\Delta^o}$
is an isomorphism. Therefore the functor $\Y(S_+,E_+)$ of
\eqref{Y.eq} sends the object $c$ in the full subcategory $\eq
\subset \adj$ to the $\C$-enriched $[n]$-augmented category $\langle
(S \times [n])_a/[n], (t(c)^*\iota_s^*E)_a \rangle$. By
Lemma~\ref{S.f.le}, this is an iterated polycylinder, so that
$\Y(S_+,E_+)$ sends $\eq \subset \adj$ into $\Mor(\C) \subset
\wCat(\C)$. Since $\Y(S_+,E_+)$ is cocartesian over $\Delta^o$, it
then sends the whole $\adj$ into $\Mor(\C) \subset \wCat(\C)$, and
gives the adjoint pair \eqref{ga.s.S}.

\subsection{Reflexivity.}

In order to prove Proposition~\ref{refl.prop} in full generality, we
need two technical results. The first is an enriched version of
Lemma~\ref{qui.le}. Assume given a quiver $Q$, with the path
category $P(Q)$. Recall that the simplicial replacement
$\Delta^oP(Q)$ is given by \eqref{del.P.Q}, and let
$$
\Dd^o_pP(Q) = \Ar^p(\Delta)^o_\natural \times_{\Delta_a^o} \Dd^oQ
\subset \Delta^o_pP(Q) = \Ar^p(\Delta)^o_\natural
\times_{\Delta_a^o} \Delta_a^oA(Q) \subset \Delta^oP(Q).
$$
Explicitly, objects in $\Delta^oP(Q)$ are pairs $\langle f,x
\rangle$ of a bispecial arrow $f:[n] \to [m]$ in $\Delta$ and a
quiver map $x:[m]_\delta \to Q$, $\Delta^o_pP(Q) \subset
\Delta^oP(Q)$ is spanned by pairs with surjective $f$, and
$\Dd^o_pP(Q) \subset \Delta^o_pP(Q)$ is spanned by pairs with $m
\leq 1$. Then for any $\C$-enrichment $\langle S,A \rangle$ of
$\Delta^oP(Q)$, the full embeddings $\Dd^o_pP(Q) \subset
\Delta^o_pP(Q) \subset \Delta^oP(Q)$ give rise to full embeddings
\begin{equation}\label{P.S.ab}
\begin{CD}
\Dd^o_pP(Q)[S] @>{\alpha}>> \Delta^o_pP(Q)[S] @>{\beta}>> \Delta^oP(Q)[S]
\end{CD}
\end{equation}
over $\Delta^o$, and $A$ induces functors $A_p =
\beta^*A:\Delta^o_pP(Q)[S] \to B\C$ and $A_d =
\alpha^*A_p:\Dd^o_pP(Q)[S] \to B\C$, again over $\Delta^o$.

\begin{lemma}\label{A.qui.le}
In the assumptions above, assume further that the functor $\Y(S,A)$
of \eqref{Y.eq} takes values in $\iMor(\C) \subset \wCat(\C)$. Then
the adjunction maps $A_p \to \alpha_*A_d$ and $A \to
\beta^{\Delta^o}_!A_p$ are isomorphisms (and both Kan extensions
exist).
\end{lemma}

\proof{} The projection \eqref{t.Q} restricts to a functor
$t:\Delta^oP(Q)[S] \to \Delta^o_aA(Q)$, and by \eqref{p.i.sq}, this
functor is a cofibration. Then a framing for the functor $\alpha^o$
is obtained by taking cocartesian liftings of the right comma-fibers
of the functor $\alpha(Q):\Dd^oQ \to \Delta^o_aA(Q)$ of
\eqref{P.a.Q}, and then as in the proof of Lemma~\ref{qui.le},
computing $\alpha_*$ by means of this framing amounts to taking
iterated limits of standard squares in $B\C$, so that $\alpha_*A_d$
exists and $A_p \cong \alpha_*A_d$. As for the Kan extension
$\beta_!^{\Delta^o}$, note that $\beta$ is a functor over
$\Delta^o_aA(Q)$, and a framing for $\beta$ is given by left
comma-fibers of its fibers $\beta_x$, $x \in \Delta^o_aA(X)$. But
every such fiber is an embedding $\beta$ of Corollary~\ref{A.cyl.corr}.
\endproof

The second result is a corollary of Lemma~\ref{cdm.le}. Consider the
adjunction $2$-category $\adj$, let $\adj^p \subset \adj$ be the
full subcategory spanned by diagrams \eqref{2adj.dia} with
surjective map $b$, and let $\adj^p_\natural = \adj^p \cap
\adj_\natural$. We then have $2$-categories $\adj^p_\natural \subset
\adj^p,\adj_\natural \subset \adj$ with the same objects $0$, $1$,
and for any $S:\{0,1\} \to \Sets$, we have a commutative diagram
\begin{equation}\label{adj.p.dia}
\begin{CD}
\adj^p_\natural[S] @>{\delta}>> \adj_\natural[S]\\
@V{\gamma}VV @VV{\gamma}V\\
\adj^p[S] @>{\delta}>> \adj[S],
\end{CD}
\end{equation}
where all the arrows are fully faithful lax $2$-functors

\begin{lemma}\label{A.adj.le}
For any $\C$-enrichment $\langle S,A \rangle$ of the $2$-category
$\adj^p$, the base change map $\delta^{\Delta^o}_!\gamma^*A \to
\gamma^*\delta^{\Delta^o}_!A$ induced by \eqref{adj.p.dia} is an
isomorphism, and its target exists iff so does its source.
\end{lemma}

\proof{} Let $P^p(\eq) \subset P(\eq)$ be the full subcategory
spanned by surjective bispecial arrows, and consider the diagram
\begin{equation}\label{P.p.dia}
\begin{CD}
P^p(\eq)[S]_{t^*a0} @>{\delta}>> P(\eq)[S]_{t^*a0}\\
@V{\gamma}VV @VV{\gamma}V\\
P^p(\eq)[S] @>{\delta}>> P(\eq)[S],
\end{CD}
\end{equation}
where $a0$ is as in Remark~\ref{adj.path.rem}. Then \eqref{adj.P}
with its left-adjoint functor \eqref{l.eq} induce left-admissible
full embeddings between all the other corresponding vertices of the
diagrams \eqref{adj.p.dia} and \eqref{P.p.dia}, and if we let $A' =
l^*A$, then $A \cong l^{\Delta^o}_!A'$. Therefore it suffices to
prove the claim for the base change map associated to
\eqref{P.p.dia}. To do this, it suffices to prove that for any $c
\in P(\eq)[S]$, the embedding
\begin{equation}\label{adj.p}
P^p(\eq)[S] /^\delta_{t^*a0} c \to P^p(\eq)[S] /^\delta c
\end{equation}
induced by $\delta:P^p(\eq) \to P(\eq)$ is cofinal. By
Lemma~\ref{cof.le}, it suffices to prove this with $a0$ replaced by
$\Iso$. But then, the projection $P(\eq)[S] \to P(\eq) \to P(\ppt^2)
= \Ar^\pm(\Delta)^o$ is a discrete cofibration, and it identifies
the right comma-fibers of the embedding \eqref{adj.p} with those of
\eqref{c.d.m}, so we are done by Lemma~\ref{cdm.le}.
\endproof

\proof[Proof of Proposition~\ref{refl.prop}.] By
Remark~\ref{adj.path.rem}, the $2$-category $\adj_\natural \cong
\Delta^o P([2]_\lambda)$ is the simplicial replacement of the path
category of the quiver $[2]_\lambda$, so that Lemma~\ref{qui.rem}
provides a canonical $2$-functor $\adj_\natural \to \iMor(\C)$ ---
or to $\Mor(\C)$ in the reflexive case --- that sends $0$ resp.\ $1$
to $\langle S_0,A_0 \rangle$ resp.\ $\langle S_1,A_1 \rangle$, and
the two edges of the quiver to $\Cyl(f,g)$ and
$\Cyl^\iota(f,g)$. Let $\langle S,A \rangle$ be the corresponding
$\C$-enrichment of the $2$-category $\adj_\natural$. Then it
suffices to prove that $A:\adj_\natural[S] \to B\C$ extends to
$\adj[S] \supset \adj_\natural[S]$. Moreover, by
Lemma~\ref{A.qui.le}, we have $A \cong \beta^{\Delta^o}_!A_p$, and
the embedding $\beta$ of \eqref{P.S.ab} factors as
$$
\begin{CD}
\Delta^o_pP([2]_\lambda) @>{\eps}>> \adj^p_\natural[S] @>{\delta}>>
\adj_\natural[S] \cong \Delta^oP([2]_\lambda)[S],
\end{CD}
$$
where $\delta$ is as in \eqref{adj.p.dia}. Therefore by
Lemma~\ref{A.adj.le}, it suffices to find a functor $A_w:\adj^p[S]
\to B\C$ over $\Delta^o$ and an isomorphism $\gamma^*A_w \cong
\eps^{\Delta^o}_!A_p$ --- then $\delta^{\Delta^o}_!A_w$ exists and
restricts to $A$ on $\adj_\natural[S] \subset \adj[S]$.

To construct $A_w$, as in Remark~\ref{cyl.2.rem}, decompose the map
$f:S_0 \to S_1$ as
\begin{equation}\label{s.z}
\begin{CD}
S_0 @>{\sigma}>> S_{01} @>{\zeta}>> S_1,
\end{CD}
\end{equation}
where $\sigma:S_0 \to S_0 \copr S_1 = S_{01}$ is the embedding, and
the map $\zeta$ is equal to $f$ on $S_0 \subset S_{01}$ and to $\id$
on $S_1 \subset S_{01}$. Denote $\langle S_{01},A_{01} \rangle =
\zeta^*\langle S_1,A_1 \rangle$. Then we have
$\Delta^oe(\sigma)^*A_{01} \cong \Delta^oe(f)^*A_1$, and the functor
$\langle f,g \rangle$ factors through a functor $\langle \sigma,g
\rangle:\langle S_0,A_0 \rangle \to \langle S_{01},A_{01}
\rangle$. Since $\sigma$ is injective, $\Delta^oe(\sigma)$ is fully
faithful, and the cylinder $\Cyl(\Delta^oe(\sigma))$ is a full
subcategory in $\Delta^oe(S) \times [1]$. Then the functor $w$ of
\eqref{w.eq} restricts to a functor
\begin{equation}\label{w.A.eq}
w:\Cyl(\Delta^oe(\sigma)) \to \adj^p[S] \subset B[1][S],
\end{equation}
while the functor $\langle \sigma,g \rangle$ defines a functor
$A_g:\Cyl(\Delta^oe(\sigma)) \to B\C$ over $\Delta^o$ given by $A_0$
resp.\ $A_{01}$ over $0 \times \Delta^o$ resp.\ $1 \times
\Delta^o$. The same argument as in Lemma~\ref{wr.le} then shows that
the right Kan extension
\begin{equation}\label{A.refl}
A_w=w_*A_g:\adj^p[S] \to B\C
\end{equation}
exists and is a lax $2$-functor.

To construct the isomorphism, consider the restriction
$\alpha^*\eps^*\gamma^*A_w$ of the functor $A_w$ to
$\Dd^o_pP([2]_\lambda)[S]$. Then to obtain an isomorphism
$\alpha^*\eps^*\gamma^*A_w \cong A_d$, it suffices to observe that
$$
\begin{aligned}
\Delta^oe(S_0) \copr \Delta^oe(S_1) &\cong \Delta^oe(S_{01}/[1])
\cap \Delta^oe(S_{10}/[1]),\\
\Dd^o_pP([2]_\lambda)[S] &\cong \Delta^oe(S_{01}/[1])
\cup \Delta^oe(S_{10}/[1]),
\end{aligned}
$$
apply Example~\ref{MV.exa}, and interpret the cylinder and the dual
cylinder in terms of \eqref{cyl.2.dia}. Moreover, we then get the
adjoint map $\eps^*\gamma^*A_w \to A_p \cong \alpha_*A_d$ between
lax $2$-functors, and since lax $2$-functors preserve limits of
standard squares, the same argument as in Lemma~\ref{A.qui.le} shows
that this is also an isomorphism. Thus by Lemma~\ref{A.qui.le}, it
suffices to prove that the adjunction map
\begin{equation}\label{e.A}
\eps^{\Delta^o}_!\eps^*\gamma^*A_w \to \gamma^*A_w
\end{equation}
is an isomorphism. This is again a map of lax $2$-functors, so it
further suffices to prove that it is an isomorphism on the fiber
$\adj^p[S]_{[1]}$ over $[1] \in \Delta^o$.

Now note that the whole fiber $\adj^p[S]_{[1]}$ lies in the
essential image of the functor $w$ of \eqref{w.A.eq} (this amounts
to observing that a surjective bispecial map $[n] \to [m]$ in
$\Delta$ with $[n]=[1]$ is either a bijection, or the projection
onto $[0]$). Therefore it suffices to prove that \eqref{e.A} becomes
an isomorphism after applying $w^*$. Moreover, $w$ sends
$\Delta^oe(S_0) = \Cyl(\Delta^oe(\sigma))_0$ into
$\Delta^o_pP([2]_\lambda)[S] \subset \adj^p_\natural[S]$ where
\eqref{e.A} is tautologically an isomorphism, so it suffices to
further restrict to $\Cyl(\Delta^oe(\sigma))_1 = \Delta^oe(S_{01}) =
\eq[S]$. The corresponding component $w_1=\eta:\eq[S] \to \adj^p[S]$
of the functor $w$ is then the embedding \eqref{eta.eq}. It is fully
faithful, and we have
\begin{equation}\label{A.A.A}
\eta^*A_w \cong A_{01} \cong \nu^{\Delta^o}_!A_w,
\end{equation}
where $\nu:\adj^p[S] \to \eq[S]$ is left-adjoint to $\eta$ induced
by \eqref{1.adj.eq}. Moreover, $\eq[S] \cap \adj^p_\natural[S]
\subset \eq[S]$ is the dense subcategory $\eq[S]_{a0}$, the functors
$\eta$ and $\nu$ induce an adjoint pair of functors between
$\eq[S]_{a0}$ and $\adj^p_\natural[S]$, and if we let
$\gamma:\eq[S]_{a0} \to \eq[S]$ be the embedding, then \eqref{A.A.A}
also holds for the restrictions $\gamma^*A_w$ and
$\gamma^*A_{01}$. Finally, if we denote $\eq_p[S]_{a0} = \eq[S]_{a0}
\cap \Delta^o_pP([2]_\lambda)[S]$, with the embedding
$\eps:\eq_p[S]_{a0} \to \eq[S]_{a0}$, then $\eta$ and $\nu$ further
induce an adjoint pair of functors between $\eq_p[S]_{a0}$ and
$\Delta^o_pP([2]_\lambda)[S]$, \eqref{A.A.A} still holds, and
$\eta^*$ sends the map \eqref{e.A} to the adjunction map
\begin{equation}\label{e.A.A}
\eps^{\Delta^o}_!\eps^*\gamma^*A_{01} \to \gamma^*A_{01}.
\end{equation}
But $A_{01}$ only depends on $f:S_0 \to S_1$ and $A_1$, so that from
now on, we may assume that $g:A_0 \to \Delta^oe(f)^*A_1$ is an
isomorphism. Moreover, if we let $S_0' = S_1' = S_1$ and
$A_0'=A_1'=A_1$, then the map $\zeta$ of \eqref{s.z} factors through
a map $h = f \copr \id:S_{01} \to S'_{01} = S'_0 \copr S'_1$, and
the functor $h:\eq[S]_{a0} \to \eq[S']_{a0}$ identifies the left
comma-fibers of the embeddings $\eps:\eq_p[S]_{a0} \to \eq[S]_{a0}$
and $\eps:\eq_p[S']_{a0} \to \eq[S']_{a0}$, so that $h^* \circ
\eps^{\Delta^o}_! \cong \eps^{\Delta^o}_! \circ h^*$. Therefore the
map \eqref{e.A.A} is the pullback of the same map for the
functor $\langle \id,\id \rangle:\langle S_1,A_1 \rangle \to \langle
S_1,A_1 \rangle$, and we may assume right away that $\langle f,g
\rangle = \langle \id,\id \rangle$. But then $A_w$ trivially extends
to the whole $\adj[S]$, and \eqref{e.A}, hence also \eqref{e.A.A} is
an isomorphism by Lemma~\ref{A.qui.le}.
\endproof

\subsection{Admissible bases.}\label{ev.subs}

It turns out that the adjoint pairs \eqref{ga.s.S} allow one to
characterize the Morita $2$-category by a universal property
refining Lemma~\ref{2.enr.le}. To do this, it is convenient to
change notation: assume given a $2$-category $\C$, and assume that
$\C$ is {\em pointed}, in the sense that we have a distinguished
object $o \in \C_{[0]}$. Let $\C_o = \C(o,o)$ be the endomorphisms
category of the object $o$, with its natural unital monoidal
structure of Example~\ref{end.exa}.

\begin{defn}\label{base.def}
A {\em base} of a pointed $2$-category $\langle \C,o \rangle$ is a
collection of non-empty sets $S(c)$, $c \in \C_{[0]}$ of reflexive
morphisms $i_s \in \C(o,c)$, $s \in S(c)$, such that $s(o)$ consists
of the unit object $1 \in \C_o = \C(o,o)$. A base is {\em
  admissible} if for any $c,c' \in \C_{[0]}$, $f \in \C(c,c')$ and
$s \in S(c)$ we have $f \circ i_s \cong i_{s'} \circ V$ for some $s'
\in S(c')$, $V \in \C_o$, and {\em strongly admissible} if one can
arrange that $V \cong 1$. If $\C$ has filtered colimits, then a base
is {\em ind-admissible} resp.\ {\em strongly ind-admissible} if for
any $c,c' \in \C_{[0]}$, $f \in \C(c,c')$ and $s \in S(c)$, $f \circ
i_s$ is a filtered colimit of objects of the form $i_{s'} \circ V$
resp.\ $i_{s'}$.
\end{defn}

For any pointed $2$-category $\langle \C,o \rangle$ with a base $S$,
the adjoint pairs \eqref{ga.s.S} taken together define an adjoint
pair
\begin{equation}\label{ga.S}
\gamma:\C[S]_{[0]} \times \adj \to \C.
\end{equation}
Since $\C[S]_{[0]}$ is discrete, $\gamma$ is also a coadjoint pair,
so that by Proposition~\ref{tw.prop}, it gives rise to the twisting
lax $2$-functor $\Theta(\gamma):\gamma_1^*\C \to \gamma_0^*\C$ of
\eqref{P.al}. Note that we have $\gamma_1^*\C \cong \C[S]$, while
$\gamma_0^*\C \cong \C[S]_{[0]} \times B\C_o$ is the disjoint union
of several copies of the category $B\C_o$. If we let
$\tau:\C[S]_{[0]} \to \ppt$ be the tautological projection, then we
have the lax $2$-functor
\begin{equation}\label{the.2}
\Theta = (\tau \circ \id) \circ \Theta(\gamma):\C[S] \to B\C_o,
\end{equation}
and the pair $\langle S,\Theta \rangle$ is a $\C_o$-enrichment of
the $2$-category $\C$ in the sense of Definition~\ref{2.enr.def}.

\begin{lemma}\label{base.le}
Assume that the base $S$ of a pointed $2$-category $\C$ is
admissible resp.\ strongly admissible. Then the functor
$\Y=\Y(S,\Theta)$ of \eqref{Y.eq} corresponding to the enrichment
\eqref{the.2} factors through $\Mor(\C_o)$ resp.\ $\Cat(\C_o)
\subset \wCat(\C_o)$. Moreover, assume that $\C$ has filtered
colimits. Then if $S$ is ind-admissible resp.\ strongly
ind-admissible, $\Y$ factors through $\iMor(\C_o)$
resp.\ $\iCat(\C_o)$.
\end{lemma}

\proof{} Note that the description of right $A$-modules in terms of
$1$-special lax $2$-functors \eqref{M.adj} also works in in the
augmented setting: for any $[n] \in \Delta$ and $\C_o$-enriched
$[n]$-augmented category $\langle S/[n],A \rangle$, we have
$e(S/[n])^< \cong e(S_+/[\np])$, with the augmentation map $S_+ \to
[\np] = [n]^<$ sending $S$ to $[n] = t([n]) \subset [\np]$ and $o
\in S_+$ to $0$, and then right $A$-modules $M$ correspond
bijectively to lax $2$-functors $M^\dg:\Delta^oe(S_+/[\np]) \to
B\C_o$, equipped with an isomorphism $M^\dg|_{\Delta^oe(S/[n])}
\cong A$ and $1$-special with respect to the projection $S_+ \to
\{0,1\}$ sending $o$ to $0$ and $S$ to $1$. Now take some object $c
\in \C_{[n]} \subset \C$ with $b_{0!}(c) = c_0 \in \C_{[0]}$, an
element $s \in S(c)_0 \subset S(c)$ with the corresponding reflexive
morphism $i_s \in \C(o,c_0)$, and an object $V \in \C_o$. Denote
$V_s = i_s \circ V \in \C(o,c_0) \subset \C_{[1]}$, and consider the
product
$$
V_s \times c \in \C_{[n+1]} \cong \C_{[1]} \times_{\C_{[0]}}
\C_{[n]}.
$$
Then since $S(o)$ consists of the single element $1 \in \C_o$, we
have an identification $S(V_s \times c) \cong S(c)_+$, and $\Y(V_s
\times c):\Delta^oe(S(c)_+/[\np]) \to B\C_o$ is $1$-special by
virtue of \eqref{adj.fu}, thus corresponds to a right module over
$\Y(c)$. Moreover, again by \eqref{adj.fu}, the functor
\eqref{eps.s} sends this module to $i_s^\vee \circ V_s \in \C_o$,
the adjunction map $V \to i_s^\vee \circ i_s \circ V \cong i_s^\vee
\circ V_s$ induces a map
$$
\V^\dg_s \to \Y(V_s \times c),
$$
where $V^\dg_s$ corresponds to the polyrepresentable right
$\Y(c)$-module $V_s$, and \eqref{adj.fu} immediately implies that
this map is an isomorphism. This finishes the proof: since $\Y$ is
cocartesian over $\Delta^o$, we have $f_!\Y(V_s \times c) \cong
\Y(f_!(V_s \times c))$, and then the four cylinder conditions of
Definition~\ref{A.cyl.def} directly translate to the corresponding
admissibility conditions of Definition~\ref{base.def}.
\endproof

The $2$-category $\iMor(\C_o)$ itself is pointed by the point
enriched category $o=\ppt_{\C_o}$ of Example~\ref{pt.exa}, with
$\iMor(\C_o)_o \cong \C_o$ by \eqref{bc.mo}, and it is equipped with
the standard base $S_{triv}$ given by all the reflexive maps $i_s$
of \eqref{ga.s.S}. This base is ind-admissible by definition, and
its restriction to $\Mor(\C_o) \subset \iMor(\C_o)$ is
admissible. Let us now show that the base $S_{triv}$ enjoys a
universal property analogous to Lemma~\ref{2.enr.le}.

\begin{prop}\label{base.prop}
For any pointed $2$-category $\langle \C,o \rangle$ equipped with an
admissible base $S$, there exists a $2$-functor
\begin{equation}\label{2.yo}
\Y:\C \to \Mor(\C_o)
\end{equation}
equipped with an isomorphism $\phi:S \cong \Y^*S_{triv}$ and
isomorphisms $i_{\phi(s)} \cong \Y(i_s)$ for any $s \in S(c)$, $c
\in \C_{[0]}$. Moreover, such a functor is unique up to a unique
isomorphisms, and if $\C$ has filtered colimits, then the same
statement holds for an ind-admissible base $S$ and the $2$-category
$\iMor(\C_o)$.
\end{prop}

\begin{remark}
The $2$-functor \eqref{2.yo} is a sort of a $2$-Yoneda embedding for
the pointed $2$-category $\C$: we realize objects $c \in \C_{[0]}$
by $\C_o$-enriched categories formed by objects $i_s \in
\C(o,c)$. This motivates our notation.
\end{remark}

\proof{} By virtue of Lemma~\ref{base.le}, it suffices to consider
the universal situation: we take the $\C_o$-enrichment
$\Theta_{triv}$ of the $2$-category $\Mor(\C_o)$ corresponding to
the base $S_{triv}$, and we need to construct an isomorphism
$\Y(S_{triv},\Theta_{triv}) \cong \Id$ (the same argument will also
work for $\iMor(\C_o)$ in the ind-admissible case). By
Lemma~\ref{2.enr.le}, it further suffices to construct an
isomorphism
\begin{equation}\label{theta.ev}
\Theta_{triv} \cong \ev,
\end{equation}
where $\ev:\Mor(\C_o)[S_{triv}] \to B\C_o$ is induced by the
evaluation functor \eqref{ev.perp}. To do this, we construct
explicitly the decomposition \eqref{dm.facto} for the disjoint union
\begin{equation}\label{ga.mor}
\gamma:\Mor(\C_o)[S_{triv}]_{[0]} \times \adj \to \Mor(\C_o)
\end{equation}
of the adjoint-coadjoint pairs \eqref{ga.s.S}. We again use
Lemma~\ref{S.f.le}. To simplify notation, denote $\C =
\Mor(\C_o)[S_{triv}] \cong \gamma_1^*\Mor(\C_o)$, with the natural
projection $\pi:\C\{\eq\} \to \C = \Mor(\C_o)[S_{triv}] \to
\Mor(\C_o)$, and denote also $S = \pi^*S_{triv}:\C\{\eq\}_{[0]} =
\C_{[0]} \times \eq_{[0]} \to \Sets$. Moreover, let $S_+$ be the
functor from $\C\{\eq\}_{[0]}$ to sets equal to $S$ over $0 \in
\{0,1\} = \eq_{[0]}$, and sending $\C_{[0]} \times 1$ to the
point. Then let $i^+:S_+ \to S$ be equal to $\id$ over $\C_{[0]}
\times 0$ and to the embedding $i_s$ on $\langle S/[n],A,s \rangle
\times 1 \in \C_{[0]} \times \eq_{[0]}$, with the corresponding
induced $2$-functor $\iota:\C\{\eq\}[S_+] \to
\C\{\eq\}[S]$. Furthermore, denote by $\eta:\C\{\eq\}[S_+] \to
\C\{\adj\}[S_+]$ the lax $2$-functor induced by \eqref{st.C}, and
let $E:\C\{\eq\}[S] \to B\C_o$ be the $\C_o$-enrichment of the
$2$-category $\C\{\eq\}[S]$ corresponding to the projection
$\pi$. Finally, let
$$
E_+ = \eta^{\Delta^o}_!\iota^*E:\C\{\adj\}[S_+] \to B\C_o,
$$
where as in \eqref{e.pl}, the Kan extension exists by
\eqref{adj.I.eq}, and defines a lax $2$-functor, so that $\langle
S_+,E_+ \rangle$ is a $\C_o$-enrichment of the $2$-category $\C
\times \adj$. Moreover, for any object in $\eq$ represented by an
injective arrow $a:[n]_1 \to [n]$, and any object $\langle S/[n],A,s
\rangle \in \C_{[n]}$, we have a counterpart of the cartesian
diagram \eqref{e.dia}, and we conclude by base change that
$\Y(S_+,E_+)$ sends the product $\langle S/[n],A,s \rangle \times a
\in \C \times_{\Delta^o} \eq$ to the $\C_0$-enriched $[n]$-augmented
category $\langle S_a/[n],(i_{a,s}^*A)_a \rangle$, where
$i_{a,s}:S_a \to S$ is equal to $\id$ over $[n]_1$, and sends the
single element in $(S_a)_l$, $l \in [n]_0$ to the value $s(l) \in
S_l$ of the section $s:[n] \to S$. Then by Lemma~\ref{S.f.le}, this
is an iterated polycylinder, so that $\Y(S_+,E_+):\C\{\adj\} \to
\wCat(\C_o)$ factors through $\Mor(\C_o) \subset \wCat(\C_o)$ and
then restricts to a $0$-special lax $2$-functor
\begin{equation}\label{ga.dm.mor}
\gamma^\dm:\C\{\eq\} \to \Mor(\C_o)
\end{equation}
that fits into a decomposition \eqref{dm.facto} for the adjoint pair
\eqref{ga.mor}.

To finish the proof in the admissible case, it remains to recall
that by Proposition~\ref{tw.prop}, the decomposition
\eqref{dm.facto} is unique. Therefore $\Theta_{triv}$ in
\eqref{theta.ev} is naturally identified with the component
$\gamma^\dm_0:\C = \Mor(\C_o)[S_{triv}] \to B\C_o \subset
\Mor(\C_o)$ of the lax $2$-functor \eqref{ga.dm.mor}. The latter
sends $\langle S/[n],A,s \rangle$ to $A(s)$ and is canonically
identified with $\ev$.

In the ind-admissible case, the proof is the same, with $\Mor(\C_o)$
replaced by $\iMor(\C_o)$, and Lemma~\ref{S.f.le} replaced with its
obvious generalization for iterated ind-polycylinders.
\endproof

\subsection{Dualizability.}

Now as in Proposition~\ref{enr.prop}, let $\C$ be a unital monoidal
category, and note that in general, the $2$-categories of
Proposition~\ref{enr.prop} form a diagram
\begin{equation}\label{T.eq}
\begin{CD}
  \Cat(\C) @>{T}>> \Mor(\C)\\
  @VVV @VVV\\
  \iCat(\C) @>{T}>> \iMor(\C),
\end{CD}
\end{equation}
where all the $2$-functors are faithful (the bottom row only exists
if $\C$ has filtered colimits preserved by the tensor product). As
it happens, under an additional assumption on $\C$, one can also
construct a $2$-functor going in the other direction, namely,
\begin{equation}\label{P.eq}
  P:\iMor(\C) \to \iCat(\C).
\end{equation}
Here is the assumption.

\begin{defn}\label{dua.def}
  An object $V$ in a unital monoidal category $\C$ is {\em
    (left-)dua\-lizable} if it is reflexive as a morphism in the
  $2$-category $B\C$. The category $\C$ is {\em well-generated} if
  it has filtered colimits, the full subcategory $D(\C) \subset \C$
  spanned by dualizable objects is essentially small, and every
  object in $\C$ is a filtered colimit of dualizable objects.
\end{defn}

\begin{exa}\label{k.dua.exa}
For any commutative ring $k$, the category $k\amod^{fl}$ of flat
$k$-modules with its usual tensor structure is well-generated, and
the dualizable objects are finitely generated projective
$k$-modules.
\end{exa}

Assume that the category $\C$ is well-generated in the sense of
Definition~\ref{dua.def}. Then informally, we want our $2$-functor
\eqref{P.eq} to send a small $\C$-enriched category $\langle S,A
\rangle$ to the category of right $A$-modules of the form $V_s$, $s
\in S$, $V \in D(\C)$, with an appropriate enrichement. To achieve
this formally, we use the $2$-functor $\Y$ of
Proposition~\ref{base.prop}, but for a different admissible base.

Choose a set $D$ of representatives of the isomorphism classes of
dualizable objects, with the functor $\phi:D \to D(\C) \subset \C$
(where $D$ is treated as a discrete category). For any $\C$-enriched
category $\langle S,A \rangle$ and $s \in S$, let $i_{d \times s} =
i_s \circ \phi(d) \in \Mor(\C)(\ppt_\C,\langle S,A \rangle)$. This
is a composition of reflexive morphisms, thus reflexive. Now
consider the $2$-category $\iMor(\C)[\{0,1\}] = \iMor(\C) \times^2
\eq$, with the distinguished object $o = \ppt_\C \times 0$, let $f
\in \eq(0,1)$ be unique isomorphism, and define a base $S_+$ for the
pointed $2$-category $\iMor(\C)[\{0,1\}]$ as follows:
\begin{itemize}
\item for any $\C$-enriched category $\langle S,A \rangle$, the set
  $S_+(\langle S,A \rangle \times 0)$ consists of the maps $i_s$, $s
  \in S$, and the set $S_+(\langle S,A \rangle \times 1)$ consists
  of maps $i_{d \times s} \circ (\id \times f)$, $s \in S$, $d \in
  D$.
\end{itemize}
This is obviously an ind-admissible base in the sense of
Definition~\ref{base.def}, so that Proposition~\ref{base.prop}
provides a $2$-functor
\begin{equation}\label{P.pl}
P_+:\iMor(\C)[\{0,1\}] \to \iMor(\C).
\end{equation}
Moreover, since $\C$ is well-generated, the base $S_+$ is strongly
ind-admissible on $\iMor(\C) \times 1 \subset \iMor(\C)[\{0,1\}]$,
so that by Lemma~\ref{base.le}, the restriction of the $2$-functor
\eqref{P.pl} to $\iMor(\C) \times 1$ factors through $\iCat(\C)
\subset \iMor(\C)$ and induces a $2$-functor \eqref{P.eq}.

We tautologically have $P \circ T \cong \Id$, and one can probably
show that $P$ is fully faithful and right-adjoint to $T$ of
\eqref{T.eq} in the appropriate $2$-categorical sense, but we will
not need it; the following trivial observation will be sufficient
for our purposes.

\begin{lemma}\label{dua.le}
  For any well-generated unital monoidal category $\C$, the
  composition $T \circ P:\iMor(\C) \to \iMor(\C)$ of the
  $2$-functors \eqref{P.eq} and \eqref{T.eq} is equivalent to the
  identity in the sense of \eqref{equi.eq}.
\end{lemma}

\proof{} An equivalence is given by the $2$-functor \eqref{P.pl}.
\endproof

\section{Trace theories.}\label{tr.sec}

\subsection{Cyclic nerves.}\label{tr.subs}

Consider now the cyclic category $\Lambda$ formed by the categories
$[n]_\Lambda$, $n \geq 1$, and horizontal functors $f:[n]_\Lambda
\to [m]_\Lambda$. Denote by $r:\Lambda^\hdot \to \Lambda$ be the
fibration with fibers $\Delta^o[n]_\Lambda$ and transition functors
$\Delta^o(f_\dg):\Delta^o[m]_\Lambda \to \Delta^o[n]_\Lambda$, where
$f_\dg$ is as in \eqref{la.du}. The structural cofibrations
$\Lambda^\hdot_{[n]} = \Delta^o[n]_\Lambda \to \Delta^o$, $[n] \in
\Lambda$ then define a functor $l:\Lambda^\hdot \to \Delta^o$
cartesian over $\Lambda$.

\begin{defn}\label{cycl.def}
The {\em cyclic nerve} $\Lambda\C$ of a $2$-category $\C/\Delta^o$
is the cofibration $\Lambda\C = r_*l^*\C$.  The {\em cyclic nerve}
$\C^\hash$ of a unital monoidal category $\C$ is the cyclic nerve
$\C^\hash = \Lambda B\C$ of the corresponding $2$-category $B\C$.
\end{defn}

By definition, the cofibration $\Lambda\C \to \Lambda$ has fibers
$\Lambda\C_{[n]} = \Fun^2([n]_\Lambda,\C)$, with transition functors
$f_! = f_\dg^*:\Fun^2([n]_\Lambda,\C) \to \Fun^2([m]_\Lambda,\C)$.
To see these fibers explicitly, one can use the cocartesian squares
\eqref{de.la}. In effect, we have the functor $j:\Delta \cong ([1]
\setminus \Lambda)^o \to \Lambda^o \cong \Lambda$, and the pullback
functors $\omega^*_n$ for the functors $\omega_n$, $n \geq 1$ of
\eqref{de.la} provide a fully faithful embedding
\begin{equation}\label{om.eq}
\omega^*:j^*\Lambda\C \to \theta^*\C \subset \rho_\flat^{o*}\C \cong
P(\C)_{[1]},
\end{equation}
where $\theta$ is the functor \eqref{theta.eq} and $P(\C)$ is the
path $2$-category of \eqref{path.eq}. In particular, for any $c \in
\C_{[0]}$, we have the natural embedding
\begin{equation}\label{i.c}
i_c:\C(c,c) \to \Lambda\C_{[1]} \subset j^*\Lambda\C \subset \Lambda\C
\end{equation}
sending $f \in \C(c,c)$ to the corresponding path of length $1$ from
$c$ to itself. If $\C$ has one object, then \eqref{om.eq} is an
equivalence (but even in this good case, we lose the full cyclic
structure on $\Lambda\C$). If $\C = \Delta^oI$ is the simplicial
replacement of a category $I$, we will simplify notation by writing
$\Lambda I = \Lambda\Delta^oI$.

\begin{exa}\label{del.la.exa}
Consider the category $\Delta^<$ with the unital monoidal structure
$B\Delta^< = \Ar^\pm(\Delta)^o$ of Example~\ref{del.ten.exa}. Then
$\Delta^{<\hash} = \Lambda B\Delta^<$ is the full subcategory
$\Ar^<(\Lambda) \subset \Ar(\Lambda^<)$ spanned by arrows whose
target is in $\Lambda \subset \Lambda^<$.
\end{exa}

Cyclic nerves are obviously functorial with respect to $2$-functors
between the underlying $2$-categories, in that a $2$-functor
$\gamma:\C \to \C'$ induces a functor
\begin{equation}\label{la.ga.0}
\Lambda\gamma:\Lambda\C \to \Lambda\C'
\end{equation}
cocartesian over $\Lambda$. To extend this to lax $2$-functors, one
can use path $2$-ca\-te\-gories in the same way as in
Subsection~\ref{wr.subs}. Namely, for any $2$-category $\C$, the
$2$-functor $P(\C) \to \C$ of \eqref{st.P} has a fully faithful
right-adjoint $\eta:\C \to P(\C)$ of \eqref{eta.P}, and by
Lemma~\ref{loc.le}~\thetag{ii}, $\Lambda s$ of \eqref{la.ga.0} then
has a fully faithful right-adjoint $\Lambda\eta$ over
$\Lambda$. Then for any lax $2$-functor $\gamma:\C \to \C'$, we
consider its decomposition \eqref{P.ga}, and we let
\begin{equation}\label{la.ga}
\Lambda\gamma = \Lambda P(\gamma) \circ \Lambda\eta:\Lambda\C \to
\Lambda\C'.
\end{equation}
This is a functor over $\Lambda$, not necessarily cocartesian. If
$\gamma$ is a $2$-functor, then $P(\gamma) \cong \gamma \circ s$,
and since $\Lambda\eta$ is fully faithful, we have $\Lambda s \circ
\Lambda \eta \cong \Id$, so that \eqref{la.ga.0} and \eqref{la.ga}
are consistent. By construction, we also have a functorial
isomorphism
\begin{equation}\label{e.la}
\omega^* \circ \Lambda\gamma \cong \theta^*\gamma \circ \omega^*,
\end{equation}
where $\omega^*$ is the embedding \eqref{om.eq}. This also works in
families, in the following sense. For any category $I$, say that a
cofibration $\C \to \Delta^o \times I$ is a {\em family of
  $2$-categories} over $I$ if $\C_i = \C|_{\Delta^o \times i}$ is a
$2$-category for any $i \in I$. For any two such families $\C$,
$\C'$, define a {\em $2$-functor} resp.\ {\em lax $2$-functor}
$\gamma:\C \to \C'$ as a functor over $\Delta^o \times I$ whose
restriction $\gamma_i:\C_i \to \C'_i$ is a $2$-functor resp.\ a lax
$2$-functors for any $i \in I$. Then for any family $\C$, path
$2$-categories $P(\C_i)$ form a family $P(\C) \to \Delta^o \times
I$, with the $2$-functor $s:P(\C) \to \C$ and a lax $2$-functor
$\eta:\C \to P(\C)$, and just as in the absolute case, a lax
$2$-functor $\gamma:\C \to \C'$ canonically factors as $\gamma =
P(\gamma) \circ \eta$ for a $2$-functor $\gamma':P(\C) \to \C'$. For
any $I$, we can consider the diagram
$$
\begin{CD}
\Lambda \times I @<{r \times \id}<< \Lambda^\hdot \times I @>{l
  \times \id}>> \Delta^o \times I,
\end{CD}
$$
and define the {\em relative cyclic nerve} $\Lambda(\C/I)$ of a
family of $2$-categories $\C$ over $I$ by $\Lambda(\C/I) = (r \times
\id)_*(l \times \id)^*\C$. With this definition, a $2$-functor
$\gamma:\C \to \C'$ induces a functor $\Lambda(\gamma):\Lambda(\C/I)
\to \Lambda(\C'/I)$ over $\Lambda \times I$ that restricts to
$\Lambda(\C_i)$ on any $\Lambda \times i$, $i \in I$, and then
\eqref{la.ga} extends it to lax $2$-functors. Moreover, if $\gamma$
is cocartesian over a map $f$ in $I$, then so is $\Lambda(\gamma)$.

\begin{defn}\label{tr.def}
A {\em trace theory} on a $2$-category $\C$ with values in a
category $\E$ is a functor $E:\Lambda\C \to \E$ cocartesian over
$\Lambda$. A {\em trace functor} from a unital monoidal category $\C$
to some $\E$ is a trace theory on $B\C$ with values in $\E$.
\end{defn}

Definition~\ref{tr.def} is a generalization of \cite[Definition
  2.8]{trace}. In the most basic example, an algebra object $A$ in a
unital monoidal category $\C$ corresponds to a lax $2$-functor $\ppt
\to B\C$, and then \eqref{la.ga} provides a canonical section
$A_\hash:\Lambda \to \C_\hash$ of the cofibration $\C^\hash \to
\Lambda$. If we also given a trace functor $E_\hash:\C^\hash \to \E$
to some category $\E$, we can define a cyclic object $E_\hash
A_\hash \in \E$; this was essentially the main construction of
\cite{trace}. In general, a trace theory $E$ on a $2$-category $\C$
provides a collection of functors $i_c^*E:\C(c,c) \to \E$, where $c
\in \C_{[0]}$ is an object in $\C$, and $i_c$ is the embedding
\eqref{i.c}. For a monoidal category $\C$, we have only one object
$o \in B\C$, and $i_o^*E$ is a simply a functor from $\C =
\C^\hash_{[1]}$ to $\E$; we call it the {\em underlying functor} of
the trace functor $E$.

\begin{remark}\label{la.rem}
Definition~\ref{cycl.def} suggests that $\Lambda\C$ really should be
thought of as a cofibration over $\Lambda^o$ rather than $\Lambda$,
and we force it to be a cofibration over $\Lambda$ by applying
\eqref{la.du}. We do it for consistency with earlier definitions of
the objects $F_\hash A_\hash$ (including the original object
$A_\hash$ of \cite{connes}).
\end{remark}

If a $2$-category $\C$ is bounded, then bounded trace theories on
$\C$ with values in some $\E$ form a well-defined category denoted
$\Tr(\C,E)$; in terms of \eqref{na.na}, we have $\Tr(\C,E) =
\Fun_\natural(\Lambda\C/\Lambda,\E)$. We also have the fibration
transpose to $\Lambda\C \to \Lambda$ that we denote by
$\Lambda^\perp\C \to \Lambda^o$. We denote $\Tr^\perp(\C,E) =
\Fun^\natural(\Lambda^\perp\C/\Lambda^o,\E)$, and we note that
\eqref{na.na} provides a natural equivalence $\Tr(\C,\E) \cong
\Tr^\perp(\C,E)$ sending $E \in \Tr(\C,\E)$ to
\begin{equation}\label{tw.na}
E^\perp \cong l!q^*r^*E \in \Tr^\perp(\C,\E),
\end{equation}
where $l$, $r$ and $q$ are as in \eqref{rtl.dia} for the cofibration
$\Lambda\C \to \Lambda$.

If $\gamma:\C' \to \C$ is a $2$-functor between $2$-categories, then
the functor $\Lambda\gamma$ of \eqref{la.ga} is cocartesian over
$\Lambda$, the transpose functor
$\Lambda^\perp\gamma:\Lambda^\perp\C' \to \Lambda^\perp\C$ is
cartesian over $\Lambda$, and for any trace theory $E$ on $\C$, the
pullback $\Lambda\gamma^*E$ is a trace theory on $\C'$, and we have
a canonical isomorphism $(\Lambda\gamma^*E)^\perp \cong
\Lambda^\perp\gamma^*E$. If the $2$-categories $\C$ and $\C'$ are
bounded, we obtain a pullback functor $\Lambda\gamma^*:\Tr(\C,\E)
\to \Tr(\C',\E)$ for any target category $\E$.

\begin{lemma}\label{S.le}
  Assume given a bounded $2$-category $\C$ and a functor $S:\C_{[0]}
  \to \Sets$ with values in non-empty sets, as in
  Example~\ref{equi.exa}, and let $\pi:\C[S] \to \C$ be the
  corresponding $2$-functor \eqref{C.S.eq}. Then for any target
  category $\E$, the pullback functor $\Lambda\pi^*:\Tr(\C,\E) \to
  \Tr(\C[S],\E)$ is an equivalence of categories.
\end{lemma}

\proof{} We need to construct the inverse equivalence. Informally,
the idea is to consider the left Kan extension $\Lambda\pi_!E$ of a
trace theory $E \in \Tr(\C[S],\E)$ and prove that, while it is not
necessarily cocartesian, it nevertheless becomes cartesian after
applying the functor $l_!q^*r^*$ of \eqref{tw.na}. Unfortunately,
$\Lambda\pi_!E$ need not even exist, so we repackage the same
argument slightly differently. Namely, we define a category
$\Tw(\Lambda\C/\Lambda,S)$ by the cartesian product
\begin{equation}\label{tw.S.dia}
\begin{CD}
\Tw(\Lambda\C,S) @>{\pi'}>> \Tw^\perp(\Lambda\C/\Lambda)\\
@V{R}VV @VV{r \circ q}V\\
\Lambda\C[S] @>{\Lambda\pi}>> \Lambda\C,
\end{CD}
\end{equation}
where $r$ and $q$ are as in \eqref{tw.na}, we let $\phi = l \circ
\pi':\Tw(\Lambda\C,S) \to \Lambda^\perp\C$, and we observe
that it suffices to prove that for any $E \in \Tr(\C[S],\E)$, the
left Kan extension $\phi_!R^*E$ exists, is cartesian over
$\Lambda^o$, and the adjunction map $R^*E \to \phi^*\phi_!R^*E$ is
an isomorphism. Indeed, then the base change map \eqref{bc.eq} and
\eqref{tw.na} induce an isomorphism $\phi_!R^*\Lambda\pi^*E \cong
E^\perp$ for any $E$ in $\Tr(\C,\E)$, and an isomorphism
$\Lambda^\perp\pi^*\phi_!R^*E \cong E^\perp$ for any $E$ in
$\Tr(\C[S],\E)$.

Since $\phi$ is a cofibration, Kan extensions $\phi_!$ can be
computed by the framing \eqref{fr.co}, that is, by \eqref{kan.eq}
with the comma-fibers replaced by the usual fibers. Moreover, since
$E$ is a trace theory, $R^*E$ is locally constant after restriction
to each of these fibers. To describe the fiber over some $c \in
\Lambda^\perp\C_{[n]}$, choose a map $f:[n] \to [1]$, and let
$\eps(f):[\nm] \to [n]_\Lambda$ be the corresponding embedding
\eqref{la.de}. Then Example~\ref{delta.n.exa} and \eqref{la.de.n}
provide an identification
\begin{equation}\label{la.fib}
\Tw(\Lambda\C,S)_{\langle [n],c \rangle} \cong \Delta^o
e(E_S(\eps(f)^*c))/[\nm]),
\end{equation}
where $E_S = E(S) \circ \nu:\C \to E\Sets$ is as in
\eqref{C.S.sq}. But since for any $[m] \in \Delta$ and
$[m]$-augmented set $S$, any object $s \in e(S_0) \subset e(S/[m])$
is initial, a locally constant functor $F:\Delta^oe(S/[m]) \to \E$
is constant by Lemma~\ref{aug.le}, and $\colim_{\Delta^oe(S/[m])}F
\cong F(s)$. Moreover, the same is true for $e(S_0) \subset
e(S/[m])$, so that in particular, the natural map
$\colim_{\Delta^oe(S_0)}F \to \colim_{\Delta^oe(S/[m])}F$ is an
isomorphism. Applying this to $E_S(\eps(f)^*c)/[\nm]$ and the
restriction of the functor $R^*E$ to the fiber \eqref{la.fib}, we
conclude that $\phi_!R^*E$ exists and is cartesian along all maps of
the form $[1] \to [n]$, while the adjunction map $R^*E \to
\phi^*\phi_!R^*E$ is an isomorphism. Since any map $[n] \to [m]$ can
be composed with a map $[1] \to [m]$, $\phi_!R^*E$ inverts all
cartesian maps.
\endproof

Lemma~\ref{S.le} immediately implies that an equivalence
\eqref{equi.eq} between some $2$-functors $\gamma_0$ and $\gamma_1$
induces an isomorphism $\Lambda\gamma_0^* \cong \Lambda\gamma_1^*$
between the corresponding pullback functors (just take $S$ to be the
constant functor with value $\{0,1\}$).

\subsection{Functoriality by adjunction.}\label{fun.adj.subs}

Explicitly, for any $2$-category $\C$, the fiber $\Lambda\C_{[1]}$
of the cyclic nerve $\Lambda\C \to \Lambda$ is the category of pairs
$\langle c,f\rangle$ of an object $c \in \C_{[0]}$ and an
endomorphism $f \in \C(c,c)$, with the embedding \eqref{i.c} sending
$f$ to $\langle c,f \rangle$, while the fiber $\Lambda\C_{[2]}$
consists of quadruples $\langle c,c',f,f' \rangle$, $c,c' \in
\C_{[0]}$, $f \in \C(c,c')$, $f' \in \C(c',c)$. For any trace theory
$E$, we then have isomorphisms
$$
E(\langle c,c',f,f' \rangle) \cong E(\langle c,f' \circ f \rangle),
\qquad E(\langle c,c',f,f' \rangle) \cong E(\langle c',f \circ f'
\rangle)
$$
provided by the two maps $[2] \to [1]$ in $\Lambda$, and these
provide an isomorphism
\begin{equation}\label{tau.ff}
\tau_{f,f'}:E(\langle c,f' \circ f \rangle) \cong E(\langle c',f
\circ f' \rangle),
\end{equation}
a sort of a categorified trace property for $E$. This explains our
terminology. Note that if $f$ and $f'$ form an adjoint pair $\langle
f,f',a,a' \rangle$, then \eqref{tau.ff} allows to define a natural
map
\begin{equation}\label{adj.dia}
\begin{CD}
E(\langle c,\id_c \rangle) @>{E(a)}>> E(\langle c,f' \circ f
\rangle) \cong E(\langle c',f \circ f' \rangle) @>{E(a')}>>
E(\langle c',\id_{c'} \rangle.
\end{CD}
\end{equation}
The map is obviously invariant under automorphisms of $f$, so that
if one checks that it is compatible with compositions, and $\C$ is
bounded, then a trace theory $E$ defines a functor
\begin{equation}\label{radj.E}
\rAdj(E):\rAdj(\C) \to \E
\end{equation}
sending $c$ to $E(\langle c,\id \rangle)$, and an adjoint pair
$\langle f,f',a,a' \rangle$ to the map \eqref{adj.dia}. If $a$ and
$a'$ are invertible, \eqref{adj.dia} is an isomorphism, so that the
functor $\rAdj(E)$ inverts equivalences.

In principle, it is not hard to check that \eqref{adj.dia} is
compatible with compositions by a direct computation. However, for
homotopical applications, we will give a more invariant argument
based on the description of the $2$-category $\aAdj(\C)$ given in
Proposition~\ref{adj.prop}. We will need the following technical
result.

\begin{lemma}\label{01.le}
Assume given a $2$-category $\C$ equipped with a $2$-functor $\C \to
\eq$, and a $0$-special functor $E:\C \to \C$ to some category
$\E$. Then the map
$$
\colim_{i_1^*\C}i_1^*E \to \colim_{\C}E
$$
induced by the embedding $i_1:i_1^*\C \to \C$ is an isomorphism, and
its source exists iff so does its target.
\end{lemma}

\proof{} Let $\eq_1 \subset \eq$ be the full subcategory spanned by
injective maps $a:[n]_1 \to [n]$ with non-empty $[n]_1$, denote
$\C_1 = \C \times_{\eq} \eq_1$, and note that the embedding $i_1^*\C
\to \C$ factors as
$$
\begin{CD}
i_1^*\C @>{\alpha}>> \C_1 @>{\beta}>> \C.
\end{CD}
$$
Then $\alpha$ is fully faithful and has a left-adjoint $s_1:\C_1 \to
\C$ induced by \eqref{s.01.C.eq}, so that $i_1^*\C \subset \C_1$ is
left-admissible, and then by \eqref{adm.eq}, it suffices to prove
that the left Kan extension $\beta_!\beta^*E$ exists, and the
adjunction map $\beta_!\beta^*E \to E$ is an isomorphism. But the
full embedding $\beta$ has a framing given by the full subcategories
$\beta(c) \subset \C_1/c$, $c \in \C$ spanned by $0$-special maps,
and to finish the proof, it remains to observe that for any $c \in
\C_{[n]}$ that is not in $\C_1$, $\beta(c)$ is equivalent to the
category $(\Delta/[\np])^o \cong \Delta^o[\np]$, so that $\beta(c)^>
\cong (\Delta^{o>})^{n+1}$. Then since $E$ is $0$-special,
$E|_{\beta(c)^>}$ is locally constant, thus constant, and then exact
by Lemma~\ref{aug.le}.
\endproof

Consider the relative cyclic nerve $\Lambda(\Eq/\Delta)$. By
definition, it comes equipped with a cofibration
$\pi:\Lambda(\Eq/\Delta) \to \Lambda$, so we can define a category
$\Tw(\Lambda,\Eq)$ by the cartesian square
\begin{equation}\label{tw.Eq.sq}
\begin{CD}
\Tw(\Lambda,\Eq) @>{\pi'}>> \Tw(\Lambda)\\
@V{R}VV @VV{t}V\\
\Lambda(\Eq/\Delta) @>{\pi}>> \Lambda,
\end{CD}
\end{equation}
where $t:\Tw(\Lambda) \to \Lambda$ is as in \eqref{ar.tw.co}. Then
$s$ of \eqref{ar.tw} induces a cofibration $\phi = s \circ
\pi':\Tw(\Lambda,\Eq) \to \Lambda^o$, and $R$ composed with the
cofibration $\Lambda(\Eq/\Delta) \to \Delta$ gives rise to a
cofibration $\Tw(\Lambda,\Eq) \to \Delta$. For any $[n] \in \Delta$,
its fiber $\Tw(\Lambda,\Eq)_{[n]} \cong \Tw(\Lambda,V([n])$ is the
category \eqref{tw.S.dia} for $\C = \ppt^2$ and $S = V([n])$. Both
$R$, $\phi$ are cocartesian over $\Delta$, and restrict to the
eponymous functors on the fibers $\Tw(\Lambda,\Eq)_{[n]} \cong
\Tw(\Lambda,V([n]))$.

Now assume given a $2$-category $\C$, and consider the adjunction
$2$-category $\aAdj(\C)$ of \eqref{aadj.c}, with the transpose
fibration $\aAdj(\C)^\perp \to \Delta$. Note that since $\aAdj(\C)
\to \Delta^o$ is semidiscrete, we actually have $\aAdj(\C)^\perp
\cong \aAdj(\C)^o$. Then the evaluation functor \eqref{adj.2.2} and
the embedding \eqref{Eta.eq} induce a functor
$$
\begin{CD}
\Eq \times_\Delta \aAdj(\C)^o \cong \Eq \times_\Delta
\aAdj(\C)^\perp @>{\eta \times \id}>> \Adj \times_\Delta
\aAdj(\C)^\perp @>>> \C
\end{CD}
$$
that we denote by $\ev$, and then \eqref{la.ga} provides a functor
\begin{equation}\label{la.ev}
\Lambda\ev:\Lambda(\Eq \times_\Delta
\aAdj(\C)^\perp/\aAdj(\C)^o) \cong \Lambda(\Eq/\Delta)
\times_\Delta \aAdj(\C)^o \to \Lambda\C.
\end{equation}
On the other hand, \eqref{tw.Eq.sq} gives rise to functors
\begin{equation}\label{st.la}
\begin{aligned}
\phi \times \id:&\Tw(\Lambda,\Eq) \times_\Delta \aAdj(\C)^o \to
\Lambda^o \times \aAdj(\C)^o\\ 
R \times \id:&\Tw(\Lambda,\Eq) \times_\Delta \aAdj(\C)^o \to
\Lambda(\Eq/\Delta) \times_\Delta \aAdj(\C)^o.
\end{aligned}
\end{equation}
Say that a map $f$ in $\aAdj(\C)^o$ is {\em antispecial} if so is
its image in $\Delta$, and say that a functor $\aAdj(\C)^o \to \E$
to some category $\E$ is {\em antispecial} if it inverts all
antispecial maps.

\begin{lemma}\label{aa.le}
Let $\Lambda\ev$ be as in \eqref{la.ev}, and let $R$, $\phi$ be as
in \eqref{st.la}. Then for any trace theory $E:\Lambda\C \to \E$ on
$\C$ with values in some category $\E$, the left Kan extension
$\aa(E) = (\phi \times \id)_!(R \times \id)^*\Lambda\ev^*E:\Lambda
\times \aAdj(\C)^o \to \E$ exists. Moreover, for any $c \in
\aAdj(\C)^o$, the restriction $\aa(E)_c = \aa(E)|_{\Lambda \times
  c}$ is locally constant, and for any $[n] \in \Lambda$, the
restriction $\aa(E)_{[n]} = E|_{[n] \times \aAdj(\C)^o}$ is
antispecial.
\end{lemma}

\proof{} For any $[n] \in \Delta$ and $c \in \aAdj(\C)^o_{[n]}
\subset \aAdj(\C)^o$, the evaluation functor $\ev$ induces a lax
$2$-functor $\ev_c:\Eq_{[n]} = \Delta^oe(V([n])) \to \C$. For any
$[l] \in \Lambda$, denote by $\F(n,l) = \Tw(\Lambda,V([n]))_{[l]}$
the fiber of the cofibration $\phi:\Tw(\Lambda,V([n])) \to
\Lambda^o$, and let $E_c = R^*\Lambda(\ev_c)^*E|_{\F(n,l)}$. As in
the proof of Lemma~\ref{S.le}, the fiber $\F(n,l)$ is explicitly
given by \eqref{la.fib} that reads as
$$
\F(n,l) \cong \Delta^oe(V([n]) \times [l-1]/[l-1]),
$$
and we have
$$
\aa(E)([l] \times c) = \colim_{\F(n,l)}E_c.
$$
If $[n]=0$, then the colimit exists by Lemma~\ref{S.le}, and we in
fact have $\aa(E)_c \cong \Lambda\eps(c)^*E^\perp$, where
$\eps(c):\ppt^2 \to \C$ is the embedding onto $c \in \C_{[0]} =
\aAdj(\C)_{[0]}$. In particular, it is locally constant. Moreover,
the functor \eqref{Eta.eq}, hence also $\ev$ and $\Lambda\ev$ are
cocartesian over antispecial maps. Thus it suffices to prove that
for any $[n]$ and $c$, and any antispecial map $f:[m] \to
[n]$ in $\Delta$, the map
$$
\colim_{\F(m,l)}\F(f)^*E_c \to \colim_{\F(n,l)}E_c
$$
induced by the embedding $\F(f):\F(m,l) \to \F(n,l)$ is an
isomorphism. Moreover, it obviously suffices to prove it for
$m=0$. But then, we can consider the map $p:V([n]) \to \{0,1\}$
sending $n$ to $1$ and the rest to $0$, with the induced projection
$p:\Eq_{[n]} \to \eq$, and \eqref{eta.n} immediately implies by
induction that $\ev_c$ is $0$-special with respect to $p$. Then
$\Lambda\ev_c$ is $0$-special with respect to the induced projection
$\F(n,l) \to \eq$, and we are done by Lemma~\ref{01.le}.
\endproof

As a corollary of Lemma~\ref{aa.le}, we see that any trace theory
$E:\Lambda\C \to \E$ on a bounded $2$-category $\C$ gives rise to an
antispecial functor
$$
\aAdj(E)^o:\aAdj(\C)^o \to \Tr(\ppt^2,\E) \cong \E, \quad c \mapsto
\aa(E)_c.
$$
If we take the opposite functor $\aAdj(E):\aAdj(\C) \to \E^o$ and
compose it with $\iota:\iota^*\aAdj(\C) \to \aAdj(\C)$, then the
resulting functor $\aAdj(E) \circ \iota:\iota^*\aAdj(\C) \to \E^o$
is special in the sense of Definition~\ref{sp.def}, thus factors
through $\rAdj(\C)^o \cong \tau(\iota^*\aAdj(\C))$ by
Corollary~\ref{sp.corr}, and this gives rise to the functor
\eqref{radj.E} that we set out to construct. Explicitly, we have
\begin{equation}\label{rAdj.E}
\rAdj(E) = \xi^\perp_!\aAdj(E)^o,
\end{equation}
where for any category $I$, we let $\xi^\perp:\Delta I \to I$ be the
functor sending $\langle [n],i_\idot \rangle$ to $i_\idot(n) \in I$,
and we extend it to $2$-categories by composing with the
truncation functor $\C^\perp \to \Delta\tau(\C)$.

\subsection{Expansion.}

Now assume given a bounded $2$-category $\C$, and assume that it is
pointed in the sense of Subsection~\ref{ev.subs}, with the object $o
\in \C_{[0]}$, the embedding $j_o:\Delta^o \to \C$ and the unital
monoidal category $\C_o=\C(o,o)$, $B\C_o \cong j_o^*\C$. Then $j_o$
has the factorization \eqref{d.f.facto}, $\overline{j}_o:B\C_o \to
\C$ is a fully faithful embedding, and we can define the {\em
  reduction functor}
\begin{equation}\label{red.eq}
\Red = \Lambda\overline{j}_o^*:\Tr(\C,\E) \to \Tr(B\C_o,\E)
\end{equation}
for an arbitrary target category $\E$. It turns out that in many
cases, one can also define a functor going in the other direction,
and in fact reconstruct a trace theory $E \in \Tr(\C,E)$ from its
reduction $\Red(E)$.

Firstly, assume that the target category $\E$ is cocomplete. Then
assume given a base $S$ of the pointed $2$-category $\langle \C,o
\rangle$ in the sense of Definition~\ref{base.def}, with the
corresponding coadjont pair \eqref{ga.S} and the lax $2$-functor
$\Theta$ of \eqref{the.2}, consider the category
$\Tw(\Lambda\C,S)$ of \eqref{tw.S.dia} with its projections
$$
\begin{CD}
\Lambda\C[S] @<{R}<< \Tw(\Lambda\C,S) @>{\phi}>> \Lambda^\perp\C,
\end{CD}
$$
and define the {\em expansion functor} $\Exp:\Fun(\C^\hash_o,\E)
\to \Fun(\Lambda^\perp\C,\E)$ by
\begin{equation}\label{exp.eq}
\Exp(E) =  \phi_!R^*\Lambda\Theta^*E.
\end{equation}
Moreover, let $\Tw(\Lambda\C,S,\eq) = \Tw(\Lambda\C,
S \times \{0,1\})$, with the projections
$$
\begin{CD}
\Lambda\C[S]\{\eq\} \cong \Lambda\C[S \times \{0,1\}] @<{R_+}<<
\Tw(\Lambda\C,S,\eq) @>{\phi^+}>> \Lambda^\perp\C,
\end{CD}
$$
and let $\Theta_+ = \ogamma \circ \gamma^\dm:\C[S]\{\eq\} \to \C,$
where $\ogamma$ and $\gamma^\dm$ are the components of the
decompositions \eqref{d.f.facto} and \eqref{dm.facto} of the
coadjoint pair \eqref{ga.S}. Then for any $E \in \Tr(\C,\E)$, the
embeddings $i_0,i_1:\C[S] \to \C[S]\{\eq\}$ induce maps
$$
\begin{CD}
i_0^*\phi_!^+R_+^*\Lambda\Theta_+^*E @>{a_0}>>
\phi_!^+R_+^*\Lambda\Theta_+^*E @<{a_1}<<
i_1^*\phi_!^+R_+^*\Lambda\Theta_+^*E.
\end{CD}
$$
We have $i_0^*\Theta_+ \cong \Lambda j_o \circ \Theta$, so that
$i_0^*\phi_!^+R_+^*\Lambda\Theta_+^*E \cong \Exp(\Red(E))$, and
$i_1^*\Theta_+$ is the projection $\pi:\C[S] \to \C$, so that
$i_1^*\phi_!^+R_+^*\Lambda\Theta_+^*E \cong E^\perp$ by
Lemma~\ref{S.le}. Moreover, the fibers \eqref{la.fib} of the map
$\phi^+$ carry natural projections to $\eq$ induced by $\C[S]\{\eq\}
\to \eq$, and then Lemma~\ref{01.le} shows that the map $a_1$ is an
isomorphism. Altogether, we obtain a functorial map
\begin{equation}\label{exp.red}
\Exp(\Red(E)) \to E^\perp
\end{equation}
for any $E \in \Tr(\C,\E)$. Note that since by definition, $S(o)$
consists of a single map $\id:o \to \id$, we also have a functorial
isomorphism
\begin{equation}\label{red.exp}
\Red(\Exp(E)) \cong E^\perp
\end{equation}
for any $E \in \Tr(B\C_o,\E)$ induced by \eqref{tw.na} and
\eqref{bc.eq}.

\begin{lemma}\label{exp.le}
Assume that the base $S$ is ind-admissible in the sense of
Definition~\ref{base.def}. Then the expansion $\Exp(E)$ of any trace
theory $E \in \Tr(B\C_o,\E)$ lies in $\Tr^\perp(\C,\E) \subset
\Fun(\Lambda^\perp\C,\E)$.
\end{lemma}

\proof{} Denote by $\F([n],c)$ the fiber of the cofibration $\phi$
over an object $\langle [n],c \rangle \in \Lambda^\perp\C_{[n]}$,
and let $E([n],c)$ be the restriction of the functor
$R^*\Lambda\Theta^*E$ to this fiber. As in the proof of
Lemma~\ref{S.le}, it suffices to prove that $\Exp(E)$ is cartesian
along any map $f:[1] \to [n]$, $[n] \in \Lambda$. This amounts to
checking that the map
$$
\colim_{\F([1],c)}E([1],f^*c) \cong \colim_{\F([1],c)}F^*E([n],c)
\to \colim_{\F([n],c)}E([n],c)
$$
induced by the embedding $F:\F([1],f^*c) \to \F([n],c)$ is an
isomorphism. Fix a map $h:[n] \to [1]$ to obtain an identification
\eqref{la.fib} of the fiber $\F([n],c)$, with the corresponding
identification of the fiber $\F([1],f^*c)$ induced by the
composition $h \circ f:[1] \to [1]$, and let $S = E_S(\eps(h)^*c)
\in \Sets/[\nm]$ and $l = f(0) \in [\nm]$. Then as in the proof of
Lemma~\ref{01.le}, $F$ factors as 
$$
\begin{CD}
\F([1],f^*c) @>{\alpha}>> \F([n],c)_f @>>> \F([n],c),
\end{CD}
$$
where $\F([n],c)_f = \Delta^oe(S/[\nm])^l \subset \F([n],c) =
\Delta^oe(S/[\nm])$ is as in \eqref{sl.eps}, and $\alpha$ is a
left-admissible embedding, so that it suffices to prove that the
adjunction map $\beta^{\Delta^o}_!\beta^*E([n],c) \to E([n],c)$ is
an isomorphism.

If $n=1$, $\F([n],c)_f = \F([n],c)$, and there is nothing to
prove. If not, we can use the framing \eqref{sl.eps} for $\beta$; we
then have to prove that for any object $\langle [\mm],e_\idot
\rangle \in \F([n],c) \setminus \F([n],c)_f$, with the corresponding
augmented functor $\eps = \eps(\langle [\mm],e_\idot \rangle)$ of
\eqref{sl.eps}, the augmented functor $\eps^*E([n],c)$ is exact.

To do this, let $p:S \to [\nm]$ be the augmentation map, and denote
by $q:[m] \to [n]$ the map in $\Lambda$ corresponding to $p \circ
s:[m-1] \to [\nm]$ under \eqref{la.de.n}. Choose a map $g:[1] \to
[m]$, and denote by $\omega_m:[m] \to [m]_\Lambda$
resp.\ $\omega_n:[n] \to [n]_\Lambda$ the functors \eqref{de.la}
corresponding to the maps $g:[1] \to [m]$ resp.\ $q \circ g:[1] \to
[n]$. Then the $2$-functor $\Y(S,\Theta)$ of
Proposition~\ref{base.prop} sends $\omega_n^*c \in \C_{[n]}$ to an
iterated ind-polycylinder $\langle S_\idot/[n],A \rangle$ in
$\iMor(\C_o)$, the augmented functor \eqref{sl.eps} naturally lifts
to an embedding
$$
\eps_\delta:\Delta^{o>}e(S_l) \to
\theta^*\Delta^oe(S_\idot/[n]), \qquad \eps_\delta(o)= \langle
      [m],\omega_m^*s \rangle,
$$
where $\theta$ is the functor \eqref{theta.eq}, and the lax
$2$-functor $A$ restricts to a functor
$$
A_\delta:\theta^*\Delta^oe(S_\idot/[n]) \to \theta^*B\C_o \cong
j^*\C^\hash_o \to \C^\hash_o,
$$
where we identify $\theta^*B\C_o \cong j^*\C^\hash_o$ by
\eqref{om.eq}. Furthermore, we have $\Theta \cong \Theta_{triv}
\circ \Y(S,\Theta)$, and by virtue of \eqref{e.la} and
\eqref{theta.ev}, we then have an isomorphism $\eps^*E([n],c) \cong
\eps_\delta^*A_\delta^*E$, so it suffices to prove that
$\eps_\delta^*A_\delta^*E$ is exact. Since $E$ is a trace theory, we
have $\eps_\delta^*A_\delta^*\zeta([m])^*E \cong
\eps_\delta^*A_\delta^*E$, where $\zeta([m])$ is the functor
\eqref{beta.eq} for the cofibration $B\C_o$. But since $A$ is an
iterated ind-polycylinder, then as we saw in the proof of
Lemma~\ref{enr.le}~\thetag{iii}, $\zeta([m]) \circ A_\delta \circ
\eps_\delta$ is a filtered colimit of contractible augmented
functors. Since $E$ commutes with filtered colimits, the same then
holds for $\eps_\delta^*A_\delta^*\zeta([m])^*E$, so that it is
exact.
\endproof

\subsection{Reconstruction.}

We can now state and prove our reconstruction theorem. Recall that
for any $2$-category $\C$ and object $c \in \C_{[0]}$, we have the
embedding \eqref{i.c}.

\begin{defn}\label{tr.exa.def}
For any two objects $c,c' \in \C_{[0]}$ in a pointed $2$-category
$\langle \C,o \rangle$ equipped with a base $S$, and any
contractible simplicial set $X$, a functor $g:\Delta^{o>}X \to
\C(c,c')$ is {\em $S$-contractible} if for any $s \in S(c)$, the
composition
$$
\begin{CD}
\Delta^{o>}X @>{g}>> \C(c,c') @>{- \circ i_s}>> \C(o,c')
\end{CD}
$$
is contractible in the sense of Definition~\ref{aug.def}. A trace
theory $E$ on the $2$-category $\C$ is {\em $S$-exact} if for any $c
\in \C_{[0]}$, contractible simplicial set $X$, and $S$-contractible
functor $g:\Delta^{o>}X \to \C(c,c)$, the functor $g^*i_c^*E$ is
exact.
\end{defn}

We note that by definition, the embedding \eqref{i.c} factors
through the fiber $(\Lambda\C)_{[1]} \subset \Lambda\C$ over $[1]
\in \Lambda$, and we have $(\Lambda\C)_{[1]} \cong
(\Lambda^\perp\C)_{[1]}$, so \eqref{i.c} also defines an embedding
\begin{equation}\label{i.c.p}
i_c^\perp:\C(c,c) \to \Lambda^\perp\C.
\end{equation}
Then for any $E \in \Tr(\C,\E)$, we have $i_c^*E \cong
i_c^{\perp*}E^\perp$, so that $E$ is $S$-exact if and only if
$g^*i_c^{\perp*}E^\perp$ is exact for any $S$-contractible
$g:\Delta^{o>}X \to \C(c,c)$.

\begin{theorem}\label{rec.thm}
Assume given a bounded pointed $2$-category $\langle \C,o \rangle$
that admits an ind-admissible base $S$ in the sense of
Definition~\ref{base.def}. Then for any cocomplete target category
$\E$, the reduction functor $\Red$ of \eqref{red.eq} has a fully
faithful left-adjoint functor
\begin{equation}\label{exp.mor.eq}
\Exp:\Tr(B\C_o,\E) \to \Tr(\C,\E)
\end{equation}
whose essential image consists of trace theories that are $S$-exact
in the sense of Definition~\ref{tr.exa.def}.
\end{theorem}

\proof{} By Lemma~\ref{exp.le}, sending $E \in \Tr(B\C,\E)$ to its
expansion $\Exp(E)$ of \eqref{exp.eq} defines a functor
\eqref{exp.mor.eq}, and we have the isomorphism \eqref{red.exp} and
the functorial map \eqref{exp.red}. To prove the theorem, it then
suffices to shows that firstly, for any $E \in \Tr(B\C_o,\E)$, the
trace theory $\Exp(E)$ is $S$-exact, and secondly, for any $S$-exact
trace theory $E \in \Tr(\C,\E)$, the map \eqref{exp.red} is an
isomorphism. For the first claim, take some $c \in \C_{[0]}$, and
define a category $\Tw(\Lambda\C,S)_c$ by the cartesian product
$$
\begin{CD}
\Tw(\Lambda\C,S)_c @>{\phi'}>> \C(c,c)\\
@V{i_c'}VV @VV{i^\perp_c}V\\
\Tw(\Lambda\C,S) @>{\phi}>> \Lambda^\perp\C,
\end{CD}
$$
where $i^\perp_c$ is the embedding \eqref{i.c.p}. Then
\eqref{la.fib} provides an identification $\Tw(\Lambda\C,S)_c
\cong \Delta^oe(S(c)) \times \C(c,c)$, with $\phi'$ given by the
projection to the second factor, and we have a projection
$$
\Theta(c) = \Lambda\Theta \circ R \circ i'_c:\Delta^oe(S(c)) \times \C(c,c)
\cong \Tw(\Lambda\C,S)_c \to \C^\hash_o
$$
such that $i_c^{\perp *}\Exp(E) \cong \phi'_!\Theta(c)^*E$ for any
$E \in \Tr(B\C_o)$ by virtue of the base change isomorphism
\eqref{bc.eq}. By definition, $\Theta(c)$ factors through
$j^*\C^\hash_o \subset \C^\hash_o$, and for any $\langle
      [n],s\rangle \in \Delta^oe(S(c))$, its composition
$$
\Theta(c)|_{\langle [n],s \rangle \times \C(c,c)}:\C(c,c) \to
  (B\C_o)_{[n]} \subset \theta^* B\C_o
$$
with the embedding \eqref{om.eq} is explicitly given by
\eqref{adj.fu}. In particular, its first component $\C(c,c) \to
\C_o$ with respect to the decomposition \eqref{c.m.eq} sends $f \in
\C(c,c)$ to $i_{s(n)}^\vee \circ f \circ i_{s(o)}$, and the other
components do not depend on $f$ at all. Therefore for any functor
$g:\Delta^{o>}X \to \C(c,c)$ that is $S$-contractible in the sense
of Definition~\ref{tr.exa.def}, the pullback $(\id \times
g)^*\Theta(c)^*E$ restricts to a contractible, hence exact functor
$\Delta^{o>}X \to \E$ over any $\langle [n],s \rangle \in
\Delta^oe(S(c))$, and then $g^*\Exp(E) \cong g^*\phi'_!\Theta(c)^*E
\cong \phi'_!(\id \times g)^*\Theta(c)^*E$ is also exact.

For the second claim, take some $S$-exact trace theory $E \in
\Tr(\C,\E)$, and note that since both $E^\perp$ and $\Exp(\Red(\E))$
are cartesian over $\Lambda$, it suffices to prove that
\eqref{exp.red} is an isomorphism over $[1] \in \Lambda$. Take an
object $c \in \C_{[0]}$, and define a category
$\Tw(\Lambda\C,S,\eq)_c$ by the cartesian diagram
$$
\begin{CD}
\Tw(\Lambda\C,S,\eq)_c @>{\phi'}>> \C(c,c)\\
@V{i_c^+}VV @VV{i^\perp_c}V\\
\Tw(\Lambda\C,S,\eq) @>{\phi_+}>> \Lambda^\perp\C.
\end{CD}
$$
Then again, \eqref{la.fib} induces an identification
$$
\Tw(\Lambda\C,S,\eq)_c \cong \Delta^oe(S(c))\{\eq\} \times
\C(c,c),
$$
and we have the projection
$$
\Theta(c)_+ = \Lambda\Theta_+ \circ R_+ \circ i^+_c:\Delta^oe(S(c))
\times \C(c,c) \to \Lambda \C.
$$
Moreover, let $\sigma:\Delta^{o>} \to \eq$ be the functor sending
$[n] \in \Delta^<$ to the injective map $s:[0] \to [\np] =
\kappa([n])$. Then $\sigma^*\Delta^oe(S(c))\{\eq\} \cong S(c) \times
\Delta^{o>}e(S(c))$, and since $S(c)$ is by assumption non-empty, we
can choose an element in $S(c)$ and lift $\sigma$ to an embedding
$\sigma_c:\Delta^{o>}e(S(c)) \to \Delta^oe(S(c))\{\eq\}$. Let
$\theta(c,f) = \Theta(c)_+ \circ (\sigma_c \times
j_f):\Delta^{o>}e(S(c)) \to \Lambda\C$, where $j_f:\ppt \to \C(c,c)$
is the embedding onto some $f \in \C(c,c)$, and denote $F = (\id
\times j_f)^*\Theta(c)_+^*E$. Then by the same argument as in
Lemma~\ref{01.le}, the map
$$
F(o) = \colim_{\Delta^{o>}e(S(c))}\sigma_c^*F \to
\colim_{\Delta^oe(S(c))\{\eq\}}F
$$
induced by $\sigma_c$ is an isomorphism, and then \eqref{exp.red} is
an isomorphism at $f \in \C(c,c) \subset \Lambda^\perp\C$ if and
only if the augmented functor $\sigma_c^*F = \theta(c,f)^*E$ is
exact. Moreover, since $E$ is cocartesian over $\Lambda$, we have
$\theta(c,f)^*E \cong \theta(c,f)^*\zeta(\theta(c,f)(o))^*E$, where
$\zeta(\theta(c,f)(o))$ is the functor \eqref{beta.eq} for the
cofibration $\Lambda\C \to \Lambda$. Since $\theta(c,f)(o)$ lies in
the image of the embedding \eqref{i.c}, we have
$\zeta(\theta(c,f)(o)) \circ \theta(c,f) \cong i_c \circ g$ for some
functor $g:\Delta^{o>}e(S(c)) \to \C(c,c)$, and what we need to
check is the exactness of $g^*i_c^*E$. However, $E$ is by assumption
$S$-exact, and \eqref{adj.fu} immediately shows that for any $s \in
S(c)$, $g \circ i_s$ extends to the contraction $\Delta^o_+e(S(c))_s
\supset \Delta^oe(S(c))$ of Lemma~\ref{yo.le}, so that $g$ is
$S$-contractible.
\endproof

\subsection{Enriched categories.}

As in Theorem~\ref{rec.thm}, assume given a monoidal category $\C_o$
and a trace functor $E \in \Tr(B\C_o,\E)$ with values in some
cocomplete category $\E$.  By virtue of Proposition~\ref{base.prop},
it is actually sufficient to describe the expansion $\Exp(E)$ in the
universal case $\C = \iMor(\C_o)$ --- indeed, it is obvious from
\eqref{exp.eq} that the expansion commutes with the pullback
$\Lambda(\Y)^*$ with respect to the $2$-functor \eqref{2.yo}. One
problem with this is that the $2$-category $\iMor(\C_o)$ is not
bounded (it has too many objects and too few morphisms). Therefore
it is necessary to replace $\iMor(\C_o)$ with a sufficiently large
but bounded full $2$-subcategory $\iMor(\C_o)_b \subset \iMor(\C_o)$
and enlarge it if necessary (this does not change the
expansion). With this convention in mind, let us give a more
explicit description of the trace theory $\Exp(E)$ and the
corresponding functor \eqref{radj.E}.

To simplify notation, let $\C=\C_o$. Consider the embedding
$j^o:\Delta^o \to \Lambda$, $j^o([n]) = [\np]$ and the cofibration
$\C^\hash = \Lambda B\C \to \Lambda$, and note that the functorial
map $\eps=\eps(f)$ of \eqref{la.de} induces a functor
\begin{equation}\label{eps.hash}
\eps^*:j^{o*}\C^\hash \to B\C
\end{equation}
cocartesian over $\Delta^o$.

\begin{defn}\label{A.bimod.def}
A {\em bimodule} over a small $\C$-enriched category $\langle S,A
\rangle$ is a functor $M:\Delta^oe(S) \to j^{o*}\C^\hash$ over
$\Delta^o$, cocartesian over anchor maps and equipped with an
isomorphism $\eps^* \circ M \cong A$.
\end{defn}

For any cocomplete category $\E$, an $\E$-valued trace functor $E
\in \Tr(B\C,\E)$ restricts to a functor $j^{o*}E:j^{o*}\C^\hash \to
\E$, and for any bimodule $M$ over a small $\C$-enriched category
$\langle S,A \rangle$, we can consider the object
\begin{equation}\label{M.A.hash}
E(M/A)_\hash = \pi_!M^*j^{o*}E \in \Fun(\Delta^o,\E),
\end{equation}
where $\pi:\Delta^oe(S) \to \Delta^o$ is the structural cofibration.
Let us then define the {\em $E$-twisted trace} of $M$ by
\begin{equation}\label{E.tr}
\Tr^E_A(M) = \colim_{\Delta^o}E(M/A)_\hash \cong
\colim_{\Delta^oe(S)}M^*j^{o*}E.
\end{equation}
This is obviously functorial with respect to $M$, so we obtain a
functor from the category $A\bimod$ of $\langle S,A
\rangle$-bimodules to $\E$.

Now, by definition, objects $c \in \Lambda\iMor(\C)_{[1]}$ are
represented by $2$-functors from $[1]_\Lambda$ to $\iMor(\C)$, and
by Lemma~\ref{2.enr.le}, such a $2$-functor defines a
$\C$-enrichment $\langle S(c),A(c) \rangle$ of the $2$-category
$\Delta^o[1]_\Lambda$. This consists of a set $S=S(c)$ and a lax
$2$-functor $A(c):\Delta^o[1]_\Lambda[S] \to B\C$. Then by
\eqref{la.ga}, the $2$-functor $A(c)$ defines a functor
$\Lambda(A(c)):\Lambda [1]_\Lambda[S] \to \C^\hash = \Lambda
B\C$. The cyclic nerve $\Lambda [1]_\Lambda$ is rather
large. However, since $\Lambda$ was defined as a subcategory in
$\rCat$, we have the Yoneda embedding $\Y:\Delta^o \cong
(\Lambda/[1])^o \to \Lambda [1]_\Lambda$, and a commutative diagram
\begin{equation}\label{yo.la.dia}
\begin{CD}
\Delta^oe(S(c)) @>>> \Lambda [1]_\Lambda[S]
@>{\Lambda(A(c))}>> \C_o^\hash\\
@VVV @VVV @VVV\\
\Delta^o @>{\Y}>> \Lambda [1]_\Lambda @>>> \Lambda.
\end{CD}
\end{equation}
The composition of the two bottom arrows is the embedding
$j^o:\Delta^o \to \Lambda$, so that the composition of the top two
arrows induces a functor
\begin{equation}\label{A.ga}
M(c):\Delta^oe(S(c)) \to j^{o*}\C_o^\hash.
\end{equation}
If we let $A = \eps^* \circ M(c):\Delta^oe(S(c)) \to B\C$, then
$\langle S,A \rangle$ is a small $\C$-enriched category, and $M(c)$
is an $A$-bimodule.

\begin{lemma}\label{A.ga.le}
For any trace functor $E \in \Tr(B\C,\E)$ and any object $c$ in the
fiber $\Lambda\iMor(\C)_{[1]}$, we have a natural identification
\begin{equation}\label{exp.ga}
\Exp(E)(c) \cong \Tr_A^E(M(c)),
\end{equation}
where $M(c)$ is the functor \eqref{A.ga}, and the right-hand side is
the $E$-twisted trace of \eqref{E.tr}.
\end{lemma}

\proof{} Combine \eqref{exp.eq}, \eqref{la.fib}, \eqref{theta.ev},
and the definition of $\ev$.
\endproof

\begin{remark}
The notation in \eqref{M.A.hash} is chosen for
consistency with \cite{trace} where we worked out to some extent the
particular case of Example~\ref{k.exa}.
\end{remark}

Now let us turn to the functor \eqref{radj.E}. For any small
$\C$-enriched category $\langle S,A \rangle$, \eqref{la.ga} provides
a functor $\Lambda(A):\Lambda\Delta^oe(S) \cong \Lambda[S] \to
\C^\hash$, and for any trace functor $E \in \Tr(B\C,\E)$, we can
define the object
\begin{equation}\label{A.hash}
EA_\hash = \pi_!\Lambda(A)^*E \in \Fun(\Lambda,\E),
\end{equation}
where as in \eqref{M.A.hash}, $\pi:\Lambda[S] \to \Lambda$ is the
structural cofibration. This is obviously functorial with respect to
$A$ and also with respect to $S$, so that we obtain a functor
$\rCat(\C) \to \Fun(\Lambda,\E)$. Composing it with the projection
$\tw_\Lambda:\Fun(\Lambda,\E) \to \Tr(\ppt,\E)$ of Lemma~\ref{av.le}
then gives a functor
\begin{equation}\label{CC.E}
\CC(E):\Refl(\C) \to \Tr(\ppt,\E).
\end{equation}
On the other hand, let $\Refl(\C) \subset \rCat(\C)$ be the dense
subcategory defined by the class of functors reflexive in the sense
of Definition~\ref{A.refl.def}. We then have the functor
\begin{equation}\label{CC.1.E}
\CC'(E):\Refl(\C) \to \Tr(\ppt,\E)
\end{equation}
obtained by composing the natural functor $\Refl(\C) \to
\rAdj(\iMor(\C))$ and the functor $\rAdj(\Exp(E))$ of
\eqref{radj.E}.

\begin{lemma}\label{A.ha.le}
The restriction of the functor $\CC(E)$ of \eqref{CC.E} to the
subcategory $\Refl(\C) \subset \rCat(\C)$ is isomorphic to the
functor $\CC'(E)$ of \eqref{CC.1.E}.
\end{lemma}

\proof{} For any categories $I$ and $\E$, to construct an
isomorphism $E \cong E'$ between any two functors $E,E' \in
\Fun(I,\E)$, it suffices to construct isomorphisms $E(i) \cong
E'(i)$ and all objects $i \in I$ that are compatible with all maps
$f:i \to i'$. Compatibility means that $E(f)=E'(f)$, or more
generally, that for any diagram \eqref{ar.dia}, we have
\begin{equation}\label{ar.compa}
E(g) \circ E'(f_0) = E'(f_1) \circ E(g').
\end{equation}
In our case, objects are small $\C$-enriched categories $\langle S,A
\rangle$, and isomorphisms
\begin{equation}\label{CC.A}
\tw_\Lambda(EA_\hash) = \CC(E)(\langle S,A \rangle) \cong
\CC'(E)(\langle S,A \rangle)
\end{equation}
are provided by the same argument as in Lemma~\ref{A.ga.le}. What we
have to check is \eqref{ar.compa}.

By definition, morphisms $f:\langle S_0,A_0 \rangle \to \langle
S_1,A_1 \rangle$ in $\rAdj(\iMor(\C))$ are represented by adjoint
pairs in the $2$-category $\iMor(\C)$ that correspond to
$\C$-enrichments $\langle S(f),A(f) \rangle$ of the $2$-category
$\adj$. Such an enrichment gives rise to a $\C$-enrichment $\langle
S(f),A(f) \circ \eta \rangle$ of the $2$-category $\eq$, and we can
consider the object $EA(f)_\hash = \pi_!\Lambda(A(f))^*E \in
\Fun(\Lambda,\E)$, where we again let $\pi:\Lambda\eq[S(f)] \to
\Lambda$ be natural discrete cofibration. We then have the diagram
\begin{equation}\label{i.01.la}
\begin{CD}
\tw_\Lambda(EA_{0\hash}) @>{i_0}>>
\tw_\Lambda(EA(f)_\hash) @<{i_1}<<
\tw_\Lambda(EA_{1\hash}),
\end{CD}
\end{equation}
the map $i_1$ is invertible by Lemma~\ref{01.le}, and in terms of
\eqref{CC.A}, the map $\rAdj(\Exp(E))(f)$ is given by $i_1^{-1}
\circ i_0$. Then if we consider $f$ as an object in the arrow
category $\Ar(\Refl(\C))$, both $EA_{0\hash}$ and $EA_{1\hash}$ in
\eqref{i.01.la} are functorial with respect to $f$, and to prove
\eqref{ar.compa}, it suffices to show the same for $EA(f)_\hash$ and
the maps $i_0$, $i_1$. Moreover, the whole diagram \eqref{i.01.la}
depends functorially on the enrichment $\langle S(f),A(f) \rangle$,
so it suffices to check that any commutative diagram
\begin{equation}\label{cat.sq}
\begin{CD}
\langle S_0,A_0 \rangle @>{\langle f,g \rangle}>> \langle S_1,A_1
\rangle\\
@V{\langle h_0,r_0 \rangle}VV @VV{\langle h_1,r_1 \rangle}V\\
\langle S'_0,A'_0 \rangle @>{\langle f',g' \rangle}>> \langle S'_1,A'_1
\rangle
\end{CD}
\end{equation}
in $\rCat(\C)$ with reflexive $\langle f,g \rangle$ and $\langle
f',g' \rangle$ gives rise to maps $h:\adj[S(f)] \to \adj[S(f')]$ and
$r:A(f) \to h^*A(f')$ that restrict to $\langle h_0,r_0 \rangle$
resp.\ $\langle h_1,r_1 \rangle$ over the objects $0$ resp.\ $1$ in
$\adj$.

Now, for any reflexive functor $\langle f,g \rangle:\langle S_0,A_0
\rangle \to \langle S_1,A_2 \rangle$ between two small $\C$-enriched
categories, the corresponding adjoint pair has been constructed in
Proposition~\ref{refl.prop}. Namely, we take $S(f) = S_0 \copr S_1$,
construct the decomposition \eqref{s.z}, and consider the diagram
\begin{equation}\label{w.d.dia}
\begin{CD}
\Cyl(\Delta^oe(\sigma)) @>{w}>> \adj^p[S_f] @>{\delta}>> \adj[S(f)],
\end{CD}
\end{equation}
where $w$ is the functor \eqref{w.A.eq}, and $\delta$ is as in
\eqref{adj.p.dia}. We then construct the functor $A_g$, and take
$A_w = w_*A_g$ and $A(f) = \delta^{\Delta^o}_!A_w$. Both $S(f)$,
\eqref{w.d.dia} and $A_g$ are obviously functorial in $f$, so that a
square \eqref{cat.sq} gives rise to a map $A_g \to h^*A'_g$, and we
have a diagram
$$
\begin{CD}
w_*A_g @>>> w_*h^*A_g @<<< h^*w_*A'_g,
\end{CD}
$$
where the map on the right is the base change map. However, it is
obvious from the explicit description of $w_*$ given in
Lemma~\ref{wr.le} that the base change map is in fact an isomorphism
and can be inverted. Therefore $A_w$ is also functorial in $f$,
and then again by base change, so is $A(f)$.
\endproof

\begin{corr}\label{inva.corr}
For any trace functor $E$, the functor \eqref{CC.E} inverts
equivalences, and identifies reflexive functors that are isomorphic
as morphisms in $\Cat(\C)$.
\end{corr}

\proof{} Clear. \endproof

\begin{remark}
In practice, one is often only interested in the functor
\eqref{CC.E} induced by a trace functor $E$, and can define it
directly without going through all the machinery of
Theorem~\ref{rec.thm}, Lemma~\ref{aa.le} and the rest of the
material in this section. However, Corollary~\ref{inva.corr} then
becomes rather cumbersome to prove.
\end{remark}

\subsection{Additional structures.}\label{more.subs}

Let us now describe some additional structures trace theories can
carry, and show that these are preserved by the expansion functor of
Theorem~\ref{rec.thm}.

\subsection{Expansion in families.}

Assume given a $2$-functor $\gamma:\C \to \C'$ between pointed
bounded $2$-categories, and assume that $\gamma$ is pointed (that
is, $\gamma(o) = o$). Define an {\em ind-admissible base} for
$\gamma$ as a pair of a functor $S':\C'_{[0]} \to \Sets$, with
non-empty values, and an adjoint pair $\iota:\adj \times
\C[\gamma^*S'] \to \C$ such that $\langle \gamma^*S',\iota \rangle$
is an ind-admissible base for $\C$, and $\langle S',\gamma \circ
\iota \rangle$ is an ind-admissible base for $\C'$. Then being
pointed, $\gamma$ restricts to a $2$-functor $\gamma_o:B\C_o \to
B\C_o'$, and if we let $\Red$, $\Red'$ resp.\ $\Exp$, $\Exp'$ be the
functors \eqref{red.eq} resp.\ \eqref{exp.mor.eq} for $\C$, $\C'$,
we have an obvious isomorphism $\Lambda(\gamma_o)^* \circ \Red'
\cong \Red \circ \Lambda(\gamma)^*$ that gives rise to the base
change map
\begin{equation}\label{exp.bc}
\Exp \circ \Lambda(\gamma_o)^* \to \Lambda(\gamma)^* \circ \Exp.
\end{equation}
A moment's reflection shows that \eqref{exp.bc} is also an
isomorphism: indeed, the functor $\Lambda(\gamma)$ identifies the
fibers \eqref{la.fib} for the cofibration $\phi$ in \eqref{exp.eq},
so that \eqref{exp.bc} reduces to the base change isomorphism
\eqref{bc.eq}.

Alternatively, one can consider the cylinder $\Cyl(\gamma)$ of the
functor $\gamma$. Then it is a family of $2$-categories over $[1]$
in the sense of Subsection~\ref{tr.subs}, and giving an
ind-admissible base for the $2$-functor $\gamma$ is equivalent to
giving a functor $S:\Cyl(f)_{[0]} \to \Sets$, with non-empty values
and cocartesian over $[1]$, together with a functor $\iota:\adj
\times \Cyl(f)_{[0]}[S] \to \Cyl(f)$, cocartesian over $\Delta^o
\times I$ and such that $\langle S,\iota \rangle$ restrict to an
admissible base on $\C = \Cyl(f)_0$ and $\C' = \Cyl(f)_1$.

Now more generally, say that a family of $2$-categories $\C \to
\Delta^o \times I$ over a bounded category $I$ is {\em pointed} if
it is equipped with a cocartesian section $o:I \to \C_{[0]}$ of the
discrete cofibration $\C_{[0]} \to I$. An {\em ind-admissible base}
for a pointed family $\langle \C/I,o \rangle$ is a pair of a functor
$S:\C_{[0]} \to \Sets$, witn non-empty values and cocartesian over
$I$, and a functor $\iota:\adj \times \C_{[0]} \to \C$, cocartesian
over $\Delta^o \times I$ and such that $\langle S,\iota \rangle$
restricts to an ind-admissible base for each $2$-category $\C_i$, $i
\in I$. We can then consider the relative cyclic nerve
$\Lambda(\C/I) \to \Lambda \times I$, and say that for any category
$\E$, an {\em $E$-valued trace theory} of $\C/I$ is a cocartesian
functor $\Lambda(\C/I) \to \E$. If $\C$ and $I$ are bounded, these
form a category $\Tr(\C,\E)$. The section $o:I \to \C_{[0]}$ gives
rise to a family of $2$-categories $B\C_o \to \Delta^o \times I$
over $I$ and the cocartesian full embedding $\overline{j}_o:B\C_o
\to \C$, and we have the functor
\begin{equation}\label{red.I}
\Red = \Lambda\overline{j}_o^*:\Tr(\C,\E) \to \Tr(B\C_o,\E),
\end{equation}
a relative version of \eqref{red.eq}. Moreover, we can consider the
subcofibration $\Tr(\C/I,\E) \subset \Fun(\C/I,\E)$ over $I$ spanned
by $\Tr(\C_i,\E) \subset \Fun(\Lambda\C_i,\E)$, and we have
$\Tr(\C,\E) \cong \Sec^\natural(I^o,\Tr(\C/I,\E))$, and similarly
for $\Tr(B\C_o,\E)$. The functor \eqref{red.I} is then induced by a
cocartesian functor
\begin{equation}\label{Red.I}
\Red_I:\Tr(\C/I,\E) \to \Tr(B\C_o/I,\E)
\end{equation}
whose fibers are the functors \eqref{red.eq} for the $2$-categories
$\C_i$. If $\E$ is cocomplete, these have left-adjoint expansion
functors \eqref{exp.eq} of Theorem~\ref{rec.thm}, and crucially, for
any map $f:i \to i'$, the adjunctions maps \eqref{exp.bc} for the
transition functor $\gamma = f_!:\C_i \to \C_{i'}$ are
isomorphisms. Therefore $\Red_I$ admits a cocartesian left-adjoint
\begin{equation}\label{exp.I}
\Exp_I:\Tr(B\C_o,\E) \to \Tr(\C,\E)
\end{equation}
over $I^o$, and taking the global sections, we also obtain a
left-adjoint $\Exp$ to the functor \eqref{red.I}.

\subsubsection{Multiplication.}

Next, we want to discuss the relationship between trace theories and
symmetric monoidal structures of Subsection~\ref{mono.subs}. We
start with a relative version of Definition~\ref{sym.1.def} and
Definition~\ref{sym.2.def}.

\begin{defn}\label{sym.I.def}
A {\em unital symmetric monoidal structure} on a cofibration $\C \to
I$ is given by a cofibration $\Bi\C \to \Gamma_+ \times I$ equipped
with an equivalence $\Bi\C|_{\ppt_+ \times I} \cong \C$ such that
for any $i \in I$, $\Bi\C_{\Gamma_+ \times i}$ is a unital symmetric
monoidal structure on $\C_i$ in the sense of
Definition~\ref{sym.1.def}. A {\em lax monoidal structure} on a
functor $\gamma:\C \to \C'$ over $I$ between two cofibrations
$\C,C'/I$ equipped with unital symmeric monoidal structures
$\Bi\C,\Bi\C'/\Gamma_+ \times I$ is a functor $\Bi\gamma:\Bi\C \to
\Bi\C'$ over $\Gamma_+ \times I$ such that for any $i \in I$,
$\Bi\gamma|_{\Gamma_+ \times i}$ is a lax monoidal structure on
$\gamma_i:\C_i \to \C'_i$ in the sense of Definition~\ref{sym.2.def}.
\end{defn}

\begin{exa}\label{bi.bi.exa}
For any unital symmetric monoidal structure $\Bi\C$ on a category
$\C$, $\Bi\C \to \Gamma_+$ carries a natural symmetric monoidal
structure $\Bi\Bi\C = m^*\Bi\C$, where $m:\Gamma_+ \times \Gamma_+
\to \Gamma_+$ is the smash product functor.
\end{exa}

In the situation of Example~\ref{bi.bi.exa}, $\Bi\C$ induces a
non-symmetric unital monoidal structure $B\C = \Sigma^*\Bi\C$ on
$\C$, and then $\Bi\Bi\C$ induces a unital symmetric monoidal
structure $\Bi B\C = (\id \times \Sigma)^*\Bi\Bi\C$ on
$B\C/\Delta$. Moreover, for any partially ordered set $J$ and
$\C$-enriched $J$-augmented categories $\langle S/J,A \rangle$,
$\langle S'/J,A' \rangle$, we can define the product $A \boxtimes
A'$ as the composition
$$
\begin{CD}
\Delta^oe(S \times_J S') \subset \Delta^oe(S/J) \times_{\Delta^o}
\Delta^oe(S'/J) @>{A \times A'}>> B\C \times_{\Delta^o} B\C \to B\C,
\end{CD}
$$
where the last functor is the product on $B\C/\Delta^o$. Then
the cofibration $\wCat(\C)/\Delta^o$ also carries a natural unital
symmetric monoidal structure, with the unit given by section
$\Delta^o \to \wCat(\C)$ sending $[n] \in \Delta^o$ to $[n] \times
\ppt_\C$, and product given by
$$
\langle S/[n],A \rangle \otimes \langle S'/[n],A' \rangle = \langle
S \times_{[n]} S'/[n],A \boxtimes A' \rangle.
$$
This induces unital symmetric monoidal structures on all the
$2$-categories of Proposition~\ref{enr.prop}.

Now, for any $2$-category $\C$, a unital symmetric monoidal
structure $\Bi\C$ on $\C/\Delta^o$ is a family of $2$-categories
over $\Gamma_+$, and moreover, it is pointed by the unit section
$\Gamma_+ \to \Bi\C_{[0]}$ of the commutative monoid
$\Bi\C_{[0]}$. We then have the unital symmetric monoidal structure
$\Bi\Lambda(\C) = \Lambda(\Bi\C/\Gamma_+)$ on its cyclic nerve
$\Lambda\C/\Lambda$, and for any category $\E$ equipped with a
unital symmetric monoidal structure $\Bi\E$, we can define a {\em
  multiplicative structure} $\Bi E$ on a functor $E:\Lambda\C \to
\E$ as a lax monoidal structure on the product $E \times
\pi:\Lambda\C \to \E \times \Lambda$, where $\pi:\Lambda\C \to
\Lambda$ is the projection. A {\em multiplicative trace theory} is a
trace theory equipped with a multiplicative structure, and if $\C$
is bounded, we denote the category of $\E$-valued multiplicative
trace theories on $\C$ by $\Tri(\C,\E)$. The reduction functor
\eqref{red.I} then induces a reduction functor
\begin{equation}\label{red.mult}
\Red:\Tri(\C,\E) \to \Tri(B\C_o,\E).
\end{equation}
Moreover, \eqref{Bi.FF} makes sense in the relative setting, so that
we have a natural unital symmetric monoidal structure on the
category $\Fun(\Bi^s \Lambda(\C)/\Gamma,\E)$, and the full
subcategory $\Tr(\Bi^s \C/\Gamma,\E) \subset
\Fun(\Bi^s\Lambda(\C)/\Gamma,\E)$ is obviously a monoidal
subcategory. Therefore \eqref{mo.ff} induces identifications
\begin{equation}\label{mo.tr}
\begin{aligned}
\Tri(\C,\E) &\cong \Sec^\natural_\infty(\Gamma,\Tr(\Bi^s
\C/\Gamma,\E)),\\
\Tri(B\C_o,\E) &\cong \Sec^\natural_\infty(\Gamma,\Tr(\Bi^s
B\C_o/\Gamma,\E)),
\end{aligned}
\end{equation}
and the functor \eqref{red.mult} is induced by a lax monoidal
structure $\Bi\Red_\Gamma$ on the functor $\Red_\Gamma$ of
\eqref{Red.I}.

\begin{lemma}\label{exp.mult.le}
Assume that the unital symmetric monoidal category $\E$ is
cocomplete, and $e \otimes -:\E \to \E$ preserves colimits for any
$e \in \E$. Then for any bounded symmetric monoidal $2$-category
$\C$, the expansion functor $\Exp_\Gamma$ of \eqref{exp.I} admits a
monoidal structure $\Bi\Exp_\Gamma$ left-adjoint over $\Gamma_+$ to
the lax monoidal structure $\Bi\Red_\Gamma$.
\end{lemma}

\proof{} We need to check that for any map $f$ in $\Gamma_+$, the
base change maps adjoint to the maps \eqref{fi.fu} for the functor
$\Bi\Red_\Gamma$ are isomorphism. By induction, it suffices to
consider the map $m:\{0,1\}_+ \to \ppt_+$. Moreover, it suffices to
consider the universal case $\C = \iMor(\C_0)$, and since we are
dealing with trace theories, it suffices to prove that the maps are
isomorphism after evaluation at any object $c \in
\Lambda\C_{[0]}$. Then by Lemma~\ref{A.ga.le}, this amound to
checking that the for any two sets $S_0$, $S_1$, the diagonal
embedding $\Delta^o e(S_0 \times S_1) \to \Delta^oe(S_0) \times
\Delta^oe(S_1)$ is cofinal, and this immediately follows from
Lemma~\ref{shuf.le}.
\endproof

By virtue of \eqref{mo.tr}, Lemma~\ref{exp.mult.le} immediately
implies that for any trace functor $E$ on $\C_o$ equipped with a
monoidal resp.\ lax monoidal structure, the expansion $\Exp(E)$
carries a natural monoidal resp.\ lax monoidal
structure. Explicitly, if $\C=\iMor(\C_o)$, then $\Exp(E)$ is given
by \eqref{exp.ga} and \eqref{E.tr}, and for any two $\C_o$-enriched
categories $\langle S,A \rangle$, $\langle S',A' \rangle$
equipped with bimodules $M$, $M'$, we have natural maps
\begin{equation}\label{mult.exp.dia}
\begin{CD}
\Tr^E_A(M) \otimes \Tr^E_{A'}(M') \cong \colim_{\Delta^o
  \times \Delta^o}E(M/A)_\hash \boxtimes E(M'/A')\hash\\
@AAA\\
\colim_{\Delta^o} E(M/A)_\hash \otimes E(M'/A')\hash\\
@VVV\\
\colim_{\Delta^oe} E(M \otimes M'/A \otimes A')_\hash \cong
\Tr^E_{A \otimes A'}(M \otimes M')
\end{CD}
\end{equation}
where the bottom map is induced by the monoidal structure on the
trace functor $E$, and the top map is invertible by
Lemma~\ref{shuf.le}. This is the multiplication map for the monoidal
structure on $\Exp(E)$.

\begin{remark}
If one is only interested in a non-symmetric monoidal structure on
the categories of Proposition~\ref{enr.prop}, then it suffices to
have a non-symmetric monoidal structure on $B\C/\Delta^o$, and for
this, the monoidal structure on $\C$ does not have to symmetric: it
suffices to ask for it to be braided. The expansion then still sends
monoidal resp.\ lax monoidal functors to monoidal resp.\ lax
mono\-idal ones. We do not go into this for lack of interesting
examples.
\end{remark}

\subsubsection{Extra functoriality.}\label{extra.subs}

Assume now given two bounded unital mono\-idal categories $\C$,
$\C'$. Then a monoidal functor $\gamma:\C' \to \C$ induces
$2$-functors \eqref{ga.mo}, and by \eqref{exp.bc}, the pullbacks
$\Lambda\gamma^*$ commute with expansion. If $\gamma$ is only lax
monoidal, then the $2$-functors \eqref{ga.mo} do not exists, but we
still have $2$-functors \eqref{ga.ca}. If we further assume that
$\C$ is well-generated in the sense of Definition~\ref{dua.def},
then we have the $2$-functor $P \circ \gamma \circ T:\iMor(\C') \to
\iMor(\C)$. We still have the base change isomorphism
\begin{equation}\label{PT1.eq}
\Lambda \gamma^*\Lambda T^*\Exp(E) \cong \Lambda
T^*\Exp(\Lambda\gamma^*E)
\end{equation}
but neither $\Lambda\gamma^*E$ nor $\Exp(\Lambda\gamma^*E)$ are trace
theories. However, assume in addition that we have given a trace
theory $E' \in \Tr(\C',\E)$ and a morphism $\alpha:\Lambda\gamma^* E
\to E'$. We then have the induced morphism
\begin{equation}\label{PT2.eq}
\Lambda P^*\Lambda T^* \Exp(\Lambda\gamma^*E) \to \Lambda
P^*\Lambda T^*\Exp(E'),
\end{equation}
and by Lemma~\ref{dua.le}, $P \circ T$ is equivalent to the
identity, so that by Lemma~\ref{S.le}, the target of this morphism
is naturally identified with $\Exp(E')$. Combining \eqref{PT1.eq}
and \eqref{PT2.eq}, we obtain a functorial morphism
\begin{equation}\label{PT.eq}
\Lambda P^*\Lambda\gamma^*\Lambda T^*\Exp(E) \to \Exp(E')
\end{equation}
of $\E$-valued trace theories on $\iMor(\C')$. If $\C$, $\C'$ and
$\E$ are unital symmetric monoidal, $\gamma$ is lax monoidal, and
$E$, $E'$ and $\alpha$ are multiplicative, then \eqref{PT.eq} is a
multiplicative map.

A typical example of such a situation occurs in the following
case. Consider the functor $V:\Lambda \to \Gamma \subset \Gamma_+$
sending $[n] \in \Lambda$ to the set $V([n]_\lambda)$ of its
vertices. Then for any bounded category $\C$ equipped with a unital
symmetric monoidal structure $\Bi\C$, we have $\C^\hash = \Lambda
B\C \cong V^*\Bi\C$. Since $V$ factors through $\Gamma \subset
\Gamma_+$, and $\ppt \in \Gamma$ is the terminal object,
\eqref{beta.eq} for the cofibration $\Bi\C|_{\Gamma_+}$ induced then
a natural functor
\begin{equation}\label{i.triv}
\Id_{triv}:\C^\hash \to \C
\end{equation}
canonically identified with $\Id$ on $\C \cong \C^\hash_{[1]}$. Thus
the identity functor $\C \to \C$ trivially extends to a trace
functor $\Id_{triv} \in \Tr(B\C,\C)$, and more generally, any
functor $E:\C \to \E$ to some category $\E$ lifts to a trace functor
$E_{triv} = E \circ \Id_{triv} \in \Tr(B\C,\E)$. Moreover,
$\Id_{triv}$ is multiplicative, so that if $\E$ is symmetric
monoidal and $E$ is multiplicative, then $E_{triv}$ is also
multiplicative.

Then if we are given another bounded unital symmetric monoidal
category $\C'$ and a lax monoidal functor $\gamma:\C' \to \C$, we
can also consider the composition $\gamma^*E:\C' \to \E$, and then
$\gamma$ induces a natural map $\alpha:\gamma^*E_{triv} \to
(\gamma^*E)_{triv}$. Plugging it into the construction of the map
\eqref{PT.eq}, we obtain a functorial map
\begin{equation}\label{PT.triv}
\Lambda P^*\Lambda\gamma^*\Lambda T^*\Exp(E_{triv}) \to
\Exp((\gamma^*E)_{triv})
\end{equation}
of $\E$-valued trace theories on $\iMor(\C')$. Again, if $\E$ is
unital symmetric monoidal and $E$ is lax monoidal, then
\eqref{PT.triv} is a multiplicative map.

\subsubsection{Edgewise subdivision.}\label{edge.subs}

The last additional structure on trace theories that we will need is
induced by the edgewise subdivision functors \eqref{i.pi}. Namely,
observe that for any $2$-category $\C$, these functors fit into a
commutative diagram
\begin{equation}\label{i.pi.dia}
\begin{CD}
\Lambda\C @<{\pi_l(\C)}<< \Lambda_l\C @>{i_l(\C)}>> \Lambda\C\\
@VVV @VVV @VVV\\
\Lambda @<{\pi_l}<< \Lambda_l @>>> \Lambda,
\end{CD}
\end{equation}
where the square on the left is cartesian, and the functor $i_l(\C)$
over some $v:[nl] \to [n]$ is given by the pullback
$v^*:\Fun^2([n]_\Lambda,\C) \to \Fun^2([n]_\Lambda,\C)$. Fix a
target category $\E$.

\begin{defn}\label{F.def}
An {\em $\langle F,l \rangle$-structure} on an $\E$-valued trace
theory $E$ on $\C$ is a map
$$
F:\pi_l(\C)^*E \to i_l(\C)^*E.
$$
\end{defn}

If $\E$ has finite limits, then since $\pi_l$ is a bifibration wth
finite fibers, so is $\pi_l(\C)$, and there exists the right Kan
extension $\pi_l(\C)_*i_l(\C)^*E$. Moreover, by \eqref{bc.eq}, this
is again a trace theory on $\C$, and by adjunction, giving an
$\langle F,l \rangle$-structure on $E$ is equivalent to giving a map
\begin{equation}\label{F.dg}
F^\dg:E \to \pi_l(\C)_*i_l(\C)^*E.
\end{equation}

\begin{defn}\label{edge.def}
The trace theory $\pi_l(\C)_*i_l(\C)^*E$ is called the {\em $l$-th
  edgewise subdivision} of the trace theory $E$.
\end{defn}

If $\C/\Delta^o$ and $\E$ are symmetric monoidal, then
$\Lambda_l\C/\Lambda_l$ is also symmetric monoidal by pullback, so
it makes sense to say that an $\langle F,l \rangle$-structure $F$ on a
multiplicative trace theory $E$ is multiplicative. If $\E$ has
finite limits, then the $l$-th edgewise subdivision
$\pi_l(\C)_*i_l(\C)^*E$ of a multiplicative trace theory $E$ is
multiplicative, and for any multiplicative $\langle F,l
\rangle$-structure $F$ on $E$, the corresponding map \eqref{F.dg} is
multiplicative as well.

Finally, for any functor $S:\C_{[0]} \to \Sets$, the diagram
\eqref{i.pi.dia} obviously induces an analogous diagram for the
category $\Tw(\Lambda\C,S)$, so that for any monoidal
category $\C$, an $\langle F,l \rangle$-structure $F$ on a trace
functor $E$ on $\C$ induces an $\langle F,l \rangle$-structure
$\Exp(F)$ on its expansion $\Exp(E)$. If $F$ is multiplicative, then
so is $\Exp(F)$. We note that if $\C$ is symmetric monoidal, then
giving an $\langle F,l \rangle$-structure on $E_{triv} \in
\Tr(B\C,\E)$ for some functor $E:\C \to \E$ is equivalent to giving
a map
\begin{equation}\label{F.triv}
E \to \delta_l^*m_l^*E,
\end{equation}
where $m_l:\C^l \to \C$ is the $l$-fold product, and $\delta_l:\C
\to \C^l$ is the diagonal embedding. Explicitly, this amounts to
giving a map $E(c) \to E(c^{\otimes l})$ for any $c \in \C$,
functorially with respect to $c$.

\begin{remark}
Since limits do not commute with colimits, edgewise subdivision of
Definition~\ref{edge.def} in general does {\em not} commute with the
expansion functor \eqref{exp.eq}.
\end{remark}

\section{Homotopy trace theories.}\label{ho.tr.sec}

\subsection{Definitions.}

Recall that for any bounded category $\C$, we have the category
$\Ho(\C)$ of functors from $\C$ to $\Top_+$ localized with respect
to pointwise weak equivalences, as in Subsection~\ref{ho.subs}, and
for any commutative ring $k$, we also have the category $\Ho(\C,k)$
of functors from $\C$ to $\Delta^ok\amod$, again localized with
respect to pointwise weak equivalences, as in
Subsection~\ref{add.subs}. For any class $v$ of maps in $\C$, we
have the full subcategory $\Ho^v(\C) \subset \Ho(\C)$ spanned by
objects $E$ such that $\ho(E):\C \to \Ho$ inverts maps in $v$, and
similarly for $\Ho(-,k)$. For any cofibration $\C \to I$ of bounded
categories, say that an object $E$ in $\Ho(\C)$ or $\Ho(C,k)$ is
{\em homotopy cocartesian over $I$} if so is $\ho(E):\C \to \Ho$,
and similarly, for any fibration $\C \to I$, say that $E$ is {\em
  homotopy cartesian over $I$} if so is $\ho(E)$. Then
$\Ho^\natural(\C) \subset \Ho(\C)$ resp.\ $\Ho^\ddag(\C) \subset
\Ho(\C)$ is the full subcategory spanned by homotopy cocartesian
resp.\ homotopy cartesian functors, and similarly for
$\Ho(-,k)$. Assume given a cofibration $\C \to I$ of bounded
categories, with the transpose fibration $\C^\perp \to I^o$, and
consider the corresponding diagram \eqref{rtl.dia}.

\begin{lemma}\label{perp.le}
In the assumptions above, the functor
\begin{equation}\label{rtl.ho}
l_! \circ q^* \circ r^*:\Ho(\C) \to \Ho(\C^\perp)
\end{equation}
induces an equivalence $\Ho^\natural(\C) \cong \Ho^\ddag(\C^\perp)$,
and similarly for $\Ho(-,k)$.
\end{lemma}

\proof{} Denote $\phi = r \circ q:\Tw^\perp(\C/I) \to \C$. Note that
for any $i \in I$ and $c \in \C^\perp_i \subset \C^\perp$, the fiber
$\Tw^\perp(\C/I)_c$ of the cofibration $l:\Tw^\perp(\C/I) \to
\C^\perp$ is naturally identified with the right comma-category $i
\setminus I$. In particular, it is homotopy contractible, with the
initial object $i$, so that if we denote $\dm = l^*(\ddag)$, the
functors $l_!$ and $l^*$ induce an equivalence
$\Ho^{\dm}(\Tw^\perp(\C/I)) \cong \Ho^\ddag(\C^\perp)$. Moreover,
$\phi^*(\natural) \subset \dm$, so that $\phi^*$ induces a functor
\begin{equation}\label{ho.perp}
\Ho^\natural(\C) \to \Ho^{\dm}(\Tw^\perp(\C/I)) \cong
\Ho^\ddag(\C^\perp), \qquad E \mapsto E^\perp = l_!\phi^*E.
\end{equation}
We have $\ho(E^\perp)(i) \cong \ho(E)(i)$ for any $i \in I$, and
$\ho(E^\perp) \cong \ho(E)^\perp$. To go in the other direction, let
$\C' = \C^o_\perp \to I$ be the cofibration transpose to $\C^o \to
I^o$, consider the diagram \eqref{rtl.dia} for the cofibration $\C'
\to I$, and pass to the opposite functors to obtain the diagram
\begin{equation}\label{rtl.2.dia}
\begin{CD}
  \C^\perp @<{\phi^o}<< \Tw^\perp(\C'/I)^o @>{l^o}>> \C.
\end{CD}
\end{equation}
Then again, the fibers of the fibration $l^o:\Tw^\perp(\C'/I) \to
\C$ are homotopy contractible, and $\phi^{o*}$ induces a functor
$$
\Ho^\ddag(\C^\perp) \to \Ho^\natural(\C), \qquad E' \mapsto E'_\perp =
l^o_*\phi^{o*}E'
$$
such that $\ho(E_\perp) \cong \ho(E)_\perp$. To finish the proof, it
remains to construct functorial isomorphisms $E \cong
(E^\perp)_\perp$ and $E' \cong (E'_\perp)^\perp$ for any $E \in
\Ho^\natural(\C)$ and $E' \in \Ho^\ddag(\C^\perp)$.

To do this, consider the product $\Tw(\C,\C') = \Tw^\perp(\C/I)
\times_{\C^\perp} \Tw^\perp(C'/I)$. Then objects in $\Tw(\C,\C')$
are composable pairs $c \to c' \to c''$ of arrows in $\C$
cocartesian over $I$, and sending such an arrow to $c \to c''$
defines a functor $\mu:\Tw(\C,\C') \to \Aa(\C)$, where $\Aa(\C)
\subset \Ar(\C)$ is spanned by cocartesian arrows. We also have the
projections $\phi,l:\Tw(\C,\C') \to \C$ sending the pair to $c$
resp.\ $c''$, $l^*$ induces an equivalence $\Ho^\natural(\C) \cong
\Ho^{\dm}(\Tw(\C,\C'))$, where we let $\dm = l^*\natural$, and
$\phi^*$ then induces a functor
$$
\Ho^\natural(\C) \to \Ho^{\dm}(\Tw(\C,\C')) \cong
\Ho^\natural(\C), \qquad E \mapsto (E^\perp)_\perp.
$$
Thus to identify $E \cong (E^\perp)_\perp$, it suffices to construct
a functorial isomorphism $\phi^*E \cong l^*E$ for any $E \in
\Ho^\natural(\C)$. But we have $\phi = s \circ \mu$, $l = t \circ
\mu$, where $s,t:\Aa(\C) \to \C$ send $c \to c'$ to $c$ resp.\ $c'$,
and we obviously have $s^*E \cong t^*E$ for any $E \in
\Ho^\natural(\C)$. This provides the required isomorphism $E \cong
(E^\perp)_\perp$. The argument for $E'$ is dual, and then the
argument for the categories $\Ho(-,k)$ is identically the same.
\endproof

Now assume given a bounded $2$-category $\C$. Then its cyclic nerve
$\Lambda\C$ is also bounded, so we may consider the homotopy
categories $\Ho(\Lambda\C)$ and $\Ho(\Lambda\C,k)$, $k$ a
commutative ring. Recall that we have the functor \eqref{ho.Ho}, and
similar functors for $\Ho(\Lambda\C,k)$.

\begin{defn}\label{ho.tr.def}
A {\em homotopy trace theory} resp.\ a {\em $k$-valued homotopy
  trace theory} on a bounded $2$-category $\C$ is an object $E$ in
$\Ho(\Lambda\C)$ resp.\ $\Ho(\Lambda\C,k)$ homotopy cocartesian over
$\Lambda$. A {\em homotopy trace functor} on a unital monoi\-dal
category $\C$ is a homotopy trace theory on $B\C$.
\end{defn}

Equvalently, $E$ is a homotopy trace theory iff $\ho(E)$ is a trace
theory in the sense of Definition~\ref{tr.def}. We denote by
$\Ho_{tr}(\C) = \Ho^\natural(\Lambda\C) \subset \Ho(\lambda\C)$,
$\Ho_{tr}(\C,k) = \Ho^\natural(\Lambda\C,k) \subset
\Ho(\Lambda\C,k)$ the full subcategories spanned by homotopy trace
theories in the sense of Definition~\ref{ho.tr.def}.

Now, the reader is invited to observe that with this definition, and
with Lemma~\ref{perp.le}, all the material of Section~\ref{tr.sec}
extends to homotopy trace theories, with identical proofs: all one
has to do is to replace colimits with homotopy colimits, Kan
extensions with homotopy Kan extensions, and reinterpret $E^\perp$
for a homotopy cocartesian functor $E$ in terms of
\eqref{ho.perp}. In particular, for any pointed bounded $2$-category
$\langle \C,o \rangle$ equipped with an admissible base $S$, we have
a pair of adjoint functors
\begin{equation}\label{ho.red.exp}
\Red^{ho}:\Ho_{tr}(\C) \to \Ho_{tr}(B\C_o), \qquad
\Exp^{ho}:\Ho_{tr}(B\C_o) \to \Ho_{tr}(\C)
\end{equation}
given by the homotopy counterparts of \eqref{red.eq} and
\eqref{exp.eq}, and similar for $\Ho_{tr}(-,k)$. We also have the
notion of an {\em $S$-homotopy exact} homotopy trace theory obtained
by replacing ``exact'' in Definition~\ref{tr.exa.def} with
``homotopy exact'', and the essential image of the full expansion
functor $\Exp$ of \eqref{ho.red.exp} consists of homotopy exact
homotopy trace theories.

All the additional structures of Subsection~\ref{more.subs} also
have obvious homotopy counterparts. In particular, for monoidal
structures, observe that the fibration $\Tr(\Bi^s\C/\Gamma,\Top_+)
\to \Gamma$ of \eqref{mo.tr} induces a fibration
$\ho(\Bi^s\C/\Gamma) \to \Gamma$ with fibers $\Ho(\C \times^2 \dots
\times^2 \C)$ that is also a cofibration and carries a symmetric
monoidal structure, and define a multiplicative homotopy trace
theory on $\C$ as an object in
$\Sec^\natural_\infty(\Gamma,\ho(\Bi^s\C/\Gamma))$. Note that
Lemma~\ref{exp.mult.le} then holds with the same proof.

Moreover, as in Lemma~\ref{aa.le}, a homotopy trace theory $E$ on a
bounded $2$-category $\C$ gives rise to an object $\aa(E)$ in
$\Ho(\aAdj(\C) \times \Lambda)$ that induces a functor
\begin{equation}\label{ho.radj}
\rAdj(E):\rAdj(\C) \to \Ho_{tr}(\ppt^2) = \Ho^\forall(\Lambda),
\end{equation}
and similarly for $\Ho(-,k)$. However, since the category $\Lambda$
is not homotopy contractible, it is no longer true that
$\Ho_{tr}(\ppt^2)$ is equivalent to $\Ho$.

In effect, the category $\Delta^o$ is homotopy contractible, so that
we at least have $\Ho^\forall(\Delta^o)\cong \Ho$, and the embedding
$j^o:\Delta^o \to \Lambda$ provides a pullback functor
$j^{o*}:\Ho^\forall(\Lambda) \to \Ho^\forall(\Delta^o)$. Moreover,
for any small category $I$, we can define a homotopy version $\tw_I$
of the twist functor \eqref{av.eq} by replacing Kan extensions with
homotopy Kan extensions, and then Lemma~\ref{av.le} holds with the
same proof, so that $\tw_\Lambda$ sends the whole $\Ho(\Lambda)$
into $\Ho^\forall(\Lambda^o)$. If we now let $\Av_\Lambda =
\tw_{\Lambda^o} \circ \tw_\Lambda$, then the same argument as in
Lemma~\ref{perp.le} provides a functorial map $\Id \to \Av_\Lambda$
whose evaluation $E \to \Av_\Lambda(E)$ at some $E \in \Ho(\Lambda)$
is an isomorphism whenever $E$ is locally constant. Therefore our
averaging functor $\Av_\Lambda:\Ho(\Lambda) \to
\Ho^\forall(\Lambda)$ is left-adjoint to the full embedding
$\Ho^\forall(\Lambda) \subset \Ho(\Lambda)$. Explicitly, for any $E
\in \Ho(\Lambda)$, we have
\begin{equation}\label{HH.eq}
j^{o*}\Av_\Lambda(E) \cong \hocolim_{\Delta^o}j^{o*}E \in
\Ho^\forall(\Delta^o) \cong \Ho,
\end{equation}
and the pullback $j^{o*}$ then has a left-adjoint given by
\begin{equation}\label{K.Av}
\Av_\Lambda \circ j^o_!:\Ho^\forall(\Delta^o) \to
\Ho^\forall(\Lambda).
\end{equation}
By adjunction, the composition $\K = \Av_\Lambda \circ j^o_! \circ
j^{o*}$ is a comonad on $\Ho^\forall(\Lambda)$, and the composition
$\K^\dg = j^{o*} \circ \Av_\Lambda \circ j^o_!$ is a monad on
$\Ho^\forall(\Delta^o) \cong \Ho$. To describe the comonad $\K$,
note that since homotopy colimits over $\Delta^o$ preserve finite
products, the functor $\tw_\Lambda$, hence also $\Av_\Lambda$ is
monoidal with respect to pointwise smash-product, and by the
projection formula, for any $E \in \Ho^\forall(\Lambda)$, we have
$j^o_!j^{o*}E \cong E \wedge j^o_!\ppt_+$, where $\ppt_+:\Delta^o
\to \Sets_+ \subset \Top_+$ is the constant functor with value
$\ppt_+$. Therefore if we let $\K = \Av_\Lambda(j^o_!\ppt_+)$, then
$\K$ is a coalgebra object in $\Ho^\forall(\Lambda)$, and we have
$\K(E) = E \wedge \K$. We also note that for any $E \in
\Ho(\Lambda)$, we have functorial isomorphisms
\begin{equation}\label{HH.HC}
\begin{aligned}
\hocolim_\Lambda E \wedge K &\cong \hocolim_\Lambda \Av_\Lambda(E
\wedge K) \cong \hocolim_\Lambda \K(\Av_\Lambda(E)) \cong\\
&\cong \hocolim_{\Delta^o} j^{o*}\Av_\Lambda(E) \cong
\hocolim_{\Delta^o}j^{o*}E,
\end{aligned}
\end{equation}
where the last isomorphism is \eqref{HH.eq}. Dually,
$\K^\dg=j^{o*}\K$ is an algebra object in $\Ho^\forall(\Delta^o)$,
and we have $\K^\dg(E) = E \wedge \K^\dg$ for any $E \in
\Ho^\forall(\Delta^o)$.

Now, up to an equivalence, $j^o:\Delta^o \cong [1] \setminus \Lambda
\to \Lambda$ is a discrete cofibration, and $j^o_!\ppt_+$ is easy to
compute; in particular, $j^{o*}j_o!\ppt_+ \in \Ho(\Delta^o)$ is
canonically identified with $\Sigma_+$, where $\Sigma$ is the
standard simplicial circle $\Sigma:\Delta^o \to \Sets_+$ of
\eqref{Si}. Then $\K^\dg \cong S^1_+$ in $\Ho^\forall(\Delta^o)
\cong \Ho$, with the algebra structure induced by the group
structure on the circle $S^1$, and for any $E \in
\Ho^\forall(\Lambda)$, the pullback $j^{o*}E$ comes equipped with an
action map
\begin{equation}\label{S1.eq}
S^1_+ \wedge j^{o*}E \to j^{o*}E.
\end{equation}
While $j^{o*}E \in \Ho^\forall(\Delta^o) \cong \Ho$ is simply a
homotopy type, the action \eqref{S1.eq} might well be non-trivial
--- for example, if $E=\K$, then \eqref{S1.eq} is the multiplication
in the algebra $\K^\dg$.

The same procedure works for $k$-valued homotopy trace theories, but
in this case, one can say more. Namely, if $k=\Z$, $\K \in
\Ho^\forall(\Lambda,\Z)$ can be represented by an explicit complex
$\K_\idot$ in $\Fun(\Lambda,\Z)$ that fits into a four-term exact
sequence
\begin{equation}\label{4.term}
\begin{CD}
0 @>>> \Z @>>> \K_1 @>>> \K_0 @>>> \Z @>>> 0,
\end{CD}
\end{equation}
where $\Z$ stands for the constant functor with value $\Z$ (see
e.g.\ \cite[Lemma 1.6]{hw}). For any $k$ and $E \in
\Ho^\forall(\Lambda,k)$, we then have $\K(E) \cong E \otimes_{\Z}
\K_\idot$. The algebra $\K^\dg(k) = j^{o*}\K_\idot \otimes_{\Z} k$
is then a DG algebra over $k$ quasiisomorphic to the homological
chain complex $C_\idot(S^1,k)$ of the circle equipped with the
Pontryagin product. This DG algebra is formal --- that is, we have
$\K^\dg(k) \cong k[B]$, the free graded-commutative algebra in one
generator $B$ of homological degree $1$ with trivial differential
--- and \eqref{S1.eq} can be refined to an equivalence of categories
\begin{equation}\label{k.B}
\D^\forall(\Lambda,k) \cong \D(k[B]).
\end{equation}
Explicitly, for any $E \in \Fun(\Lambda,k)$, the pullback
$j^{o*}\Av_\Lambda(E)$ is identified with the homology object
$C_\idot(\Delta^o,j^{o*}E) \in \D(k)$ by \eqref{HH.eq}, this can be
computed by the standard complex $C_\idot(j^{o*}E)$ of the
simplicial $k$-module $j^{o*}E$, and the generator $B$ in
\eqref{k.B} then acts on $C_\idot(j^{o*}E)$ by the Connes-Tsygan
differential (also known as Rinehart differential, see \cite[Chapter
  2]{Lo} and bibliographical comments therein).

\subsection{Stabilization.}

Let us now discuss homotopy trace theories from the point of view of
the stabilization formalism of Section~\ref{stab.sec}.

\subsubsection{Stable trace theories.}

Say that a $2$-category $\C$ is {\em half-additive} if for any $c,c'
\in \C_{[0]}$, the category $\C(c,c')$ is half-additive -- that is,
pointed with finite coproducts -- and for any $f \in \C(c,c')$, all
the composition functors $f \circ -$ and $- \circ f$ preserve finite
coproducts. In particular, $\C(c,c)$ is half-additive for any $c \in
\C_{[0]}$.

\begin{defn}
A homotopy trace theory or a $k$-valued homotopy trace theory $E$ on
a half-additive bounded $2$-category $\C$ is {\em stable} if for any
$c \in \C$, its restriction $i_c^*E$ with respect to \eqref{i.c} is
stable.
\end{defn}

In particular, for any homotopy trace theory $E$ on $\C$, object $c
\in \C_{[0]}$, and its endomorphism $f \in \C(c,c)$, we have the
object $\ho(E)(\langle c,f \rangle) \in \Ho$ induced by
\eqref{ho.Ho}, and if $E$ is stable, we also have the stable
$\Gamma$-space
\begin{equation}\label{ho.tr}
\ho^{st}(E)(\langle c,f \rangle) \in \Ho^{st}(\Gamma_+)
\end{equation}
induced by \eqref{ho.st}. If $f=\id$ is the identity endomorphism,
then we have the embedding $\eps(c):\ppt^2 \to \C$ onto $c$, and if
we let $E(c) = \Lambda\eps(c)^*E$ in $\Ho_{tr}(\ppt^2) =
\Ho^\forall(\Lambda)$, then $\ho(E)(\langle c,\id \rangle) \in \Ho$
corresponds to $j^{o*}E(c)$ under the equivalence
$\Ho^\forall(\Delta^o) \cong \Ho$. To obtain a similar
interpretation of the stable $\Gamma$-space \eqref{ho.tr}, note that
if we equip $\Gamma_+$ with the monoidal structure given by smash
product, then $B\Gamma_+$ is half-additive, and for any
half-additive $2$-category $\C$, we have a $2$-functor
\begin{equation}\label{2.ga.m}
m:B\Gamma_+ \times^2 \C \to \C
\end{equation}
that induces the functors \eqref{ga.m} for the half-additive
categories $\C(c,c')$. We then have a $2$-functor $\eps(c)_+ = m
\circ (\id \times \eps(c)):B\Gamma_+ \to \C$, and we can consider
the object
\begin{equation}\label{E.c.pl}
E(c)_+ = \Lambda\eps(c)_+^*E \in \Ho^{st}_{tr}(\Gamma_+).
\end{equation}
By definition, this is a stable homotopy trace functor on $\Gamma_+$
but it turns out that this is the same thing as a locally constant
family of stable $\Gamma$-spaces over $\Lambda$. Namely, since
$\Gamma_+$ is symmetric monoidal, it carries the tautological
$\Gamma_+$-valued trace functor \eqref{i.triv}, and taking its
product with the projection $\Gamma_+^\hash \to \Lambda$ provides a
functor
\begin{equation}\label{m.hash}
m_\hash:\Gamma_+^\hash \to \Gamma_+ \times \Lambda.
\end{equation}
Let $\Ho_\natural^{st}(\Gamma_+, \Lambda) \subset
\Ho^{st}(\Gamma_+, \Lambda)$ be the full subcategory spanned
by objects cocartesian with respect to the projection $\Gamma_+
\times \Lambda \to \Lambda$. Then \eqref{m.hash} is cocartesian over
$\Lambda$, so that it induces a functor
\begin{equation}\label{m.hash.pb}
m_\hash^*:\Ho^{st}_\natural(\Gamma_+, \Lambda) \to
\Ho^{st}_{tr}(\Gamma_+).
\end{equation}

\begin{lemma}\label{st.pt.le}
The functor \eqref{m.hash.pb} is an equivalence.
\end{lemma}

\proof{} By definition, the fibers of the cofibration
$\Gamma_+^\hash \to \Lambda$ are products $\Gamma_+^n$, $n \geq
1$. Let $\Ho^{st}(\Gamma_+^\hash) \subset \Ho(\Gamma_+^\hash)$ be
the full subcategory spanned by objects whose restriction to each of
these fibers is polystable in the sense of
Lemma~\ref{poly.st.le}. Then \eqref{m.hash} also induces a functor
\begin{equation}\label{m.hash.pb.1}
m_\hash^*:\Ho^{st}(\Gamma_+, \Lambda) \to
\Ho^{st}(\Gamma_+^\hash)
\end{equation}
that has a left-adjoint $\Stab \circ m_{\hash!}$. On each fiber
$\Gamma_+^n \subset \Gamma_+^\hash$, \eqref{m.hash} restricts to the
functor $m_n$ of \eqref{ga.m.n} (with $I = \ppt$), \eqref{m.hash.pb}
restricts to \eqref{m.n.pb}, and by \eqref{bc.eq} and
Remark~\ref{hoco.stab.rem}, $\Stab \circ m_{\hash!}$ restricts to
$\Stab \circ m_{n!}$. Then by Lemma~\ref{poly.st.le}, the adjunction
maps for the adjoint pair $m_\hash^*$, $\Stab \circ m_{\hash!}$ are
isomorphisms, so that \eqref{m.hash.pb.1} is an equivalence. To prove
that \eqref{m.hash.pb} is also an equivalence, it remains to check
that $m_\hash^*E$ is homotopy cocartesian over $\Lambda$ only if so
is $E \in \Ho^{st}(\Gamma_+, \Lambda)$. But this is obvious:
any map $f$ in $\Gamma_+ \times \Lambda$ cocartesian over
$\Lambda$ lifts to a cocartesian map in $\Gamma_+^\hash$.
\endproof

By virtue of Lemma~\ref{st.pt.le}, the object $E(c)_+$ of
\eqref{E.c.pl} defines a structure of a stable $\Gamma$-space on
$E(c) \in \Ho^\forall(\Lambda)$, and this can be used to refine
slightly the map \eqref{S1.eq} for $j^{o*}E(c)$. Namely, let
$\Ho^{st}_\natural(\Gamma_+, \Delta^o) \subset \Ho^{st}(\Gamma_+,
\Delta^o)$ be the full subcategory spanned by objects cocartesian
with respect to $\Gamma_+ \times \Delta^o \to \Delta^o$, and note
that since $\Delta^o$ is homotopy contractible, we have
$\Ho^{st}_\natural(\Gamma_+, \Delta^o) \cong
\Ho^{st}(\Gamma_+)$. Then $j^{o*}E(c)_+$ defines a structure of a
stable $\Gamma$-space on $j^{o*}E(c)$. For any small category $I$,
the construction of the map \eqref{S1.eq} works in exactly the same
way relatively over the fibers of the projection $I \times \Lambda
\to I$, so that in particular, the stable $\Gamma$-space
$j^{o*}E(c)_+$ also comes equipped with such a map. By adjunction,
we then have a based map
\begin{equation}\label{L1.eq}
j^{o*}E(c)_+ \to L(j^{o*}E(c))_+
\end{equation}
to the free loop space $L(j^{o*}E(c)_+)$. Moreover, we have natural
embeddings $\Omega(j^{o*}E(c)_+),j^{o*}E(c)_+ \to L(j^{o*}E(c)_+)$
in $\Ho^{st}(\Gamma_+)$ onto based resp.\ constant loops. But then
we can use the product map \eqref{G.prod} to combine these
embeddings into a weak equivalence $\Omega(j^{o*}E(c)_+) \times
j^{o*}E(c)_+ \to L(j^{o*}E(c)_+)$. The map \eqref{L1.eq} then splits
as $B \times \Id$, where $B$ is a map
\begin{equation}\label{B.eq}
B:j^{o*}E(c)_+ \to \Omega j^{o*}E(c)_+.
\end{equation}
in $\Ho^{st}(\Gamma_+)$. In the $k$-linear case, this is the
Connes-Tsygan differential that appears in \eqref{k.B}. The
construction is obviously functorial, so that for any half-additive
$2$-category $\C$ and stable homotopy trace theory $E$ on $\C$, the
functor
\begin{equation}\label{hh.radj}
j^{o*} \circ \rAdj(E):\rAdj(\C) \to \Ho^{st}(\Gamma_+)
\end{equation}
induced by \eqref{ho.radj} comes equipped with a map
\eqref{B.eq}. Moreover, if a stable homotopy trace theory $E \in
\Ho^{st}_{tr}(\Gamma_+)$ is multiplicative, then $j^{o*}E \in
\Ho^{st}(\Gamma_+)$ is an algebra object with respect to the unital
symmetric monoidal structure \eqref{bi.ho}, and in general, for any
multiplicative stable homotopy trace theory $E$ on a unital
symmetric monoidal half-additive $2$-category $\C$, the same applies
to the value $j^{o*}\rAdj(E)(o)$ of the functor \eqref{hh.radj} on
the unit object $o \in \C_{[0]}$.

\begin{remark}
One can show that the equivalence of Lemma~\ref{st.pt.le} is
multiplicative, and a multiplicative stable homotopy trace theory $E
\in \Ho^{st}_{tr}(\Gamma_+)$ actually defines a locally constant
family of algebra objects in $\Ho^{st}(\Gamma_+)$ parametrized by
$\Lambda$. In particular, the algebra structure is compatible with
the map \eqref{B.eq} in a certain natural way. However, we will not
need this.
\end{remark}

\begin{prop}\label{la.stab.prop}
For any bounded half-additive $2$-category $\C$, the full subcategory
$\Ho^{st}_{tr}(\C) \subset \Ho_{tr}(\C)$ spanned by stable homotopy
trace theories is left-admissible, with the stabilization functor
$$
\Stab_{\C}:\Ho_{tr}(\C) \to \Ho_{tr}^{st}(\C)
$$
adjoint to the embedding $\Ho^{st}_{tr}(\C) \to \Ho_{tr}(\C)$, and
for any object $c \in \C_{[0]}$, we have $i_c^* \circ \Stab_{\C}
\cong \Stab_{\C(c,c)} \circ i_c^*$. The same holds for the category
$\Ho_{tr}(\C,k)$ for any commutative ring $k$.
\end{prop}

Unfortunately, the stabilization procedure used in
Proposition~\ref{stab.prop} is not compatible with monoidal
structures, and cannot be applied immediately in the setting of
Proposition~\ref{la.stab.prop}. This is a well-known problem with a
well-known solution that underlies the theory of symmetric spectra
of \cite{HSS} --- roughly speaking, one needs to replace the colimit
\eqref{stab.bo} with a colimit over a more appropriate indexing
category. For symmetric spectra, this is the category of finite sets
and injective maps. For trace theories, we need something slightly
more complicated, so we first discuss the construction in an
abstract setting and exhibit its essentially $2$-categorical nature.

\subsubsection{Abstract stabilization.}\label{rec.subs}

For any unital monoidal category $J$, the pullback $\kappa^{o*}BJ
\to \Delta^o$ satisfies the Segal condition, so that its reduction
$(\kappa^{o*}BJ)^{red}$ of Remark~\ref{red.rem} is a
$2$-category. Say that $J$ is {\em essentially discrete} if the
cofibration $(\kappa^{o*}BJ)^{red} \to \Delta^o$ is discrete. In
this case, $(\kappa^{o*}BJ)^{red} \cong \Delta^oJ_{red}$ is the
simplicial replacement of a category $J_{red}$ that we call the {\em
  reduction} of the unital monoidal category $J$. The composition of
the embedding $(\kappa^{o*}BJ)^{red} \to \kappa^{o*}BJ$ and the
projection $a^o_!:\kappa^{o*}BJ \to BJ$ gives a $2$-functor
\begin{equation}\label{ee}
e:\Delta^oJ_{red} \cong (\kappa^{o*}BJ)^{red} \to BJ.
\end{equation}
Explicitly, objects in $J_{red}$ are objects in $J$, morphisms from
$j$ to $j'$ are pairs $\langle j'',a \rangle$ of an object $j'' \in
I$ and an isomorphism $a:j \otimes j'' \cong j'$, and the
$2$-functor \eqref{ee} sends such a morphism to $j''$.  In
particular, $J_{red}$ only depends on the isomorphism groupoid
$\overline{J} \subset J$ of the category $J$ with the induced
monoidal structure, so that $\overline{J}_{red} \cong J_{red}$, and
the reduction $J_{red}$ has an initial object $o \in J_{red}$
corresponding to the unit object $1 \in J$.

\begin{exa}\label{de.red.exa}
The unital monoidal category $\Delta^<$ is essentially discrete, and
its reduction $\Delta^<_{red}$ is the partially ordered set $\N$ of
non-negative integers.
\end{exa}

\begin{exa}\label{ga.red.exa}
The category $\Gamma$ of finite sets with cocartesian product is
essentially discrete, and its reduction $\Gamma_{red}$ is the
category of finite sets and injective maps. The monoidal functor
$V:\Delta^< \to \Gamma$ sending an ordinal to the set of its
elements induces a functor $V_{red}:\N \to \Gamma_{red}$ sending
maps in $\N$ to maps $V(s):V([n]) \to V([n'])$ in $\Gamma_{red}$.
\end{exa}

Assume given an essentially discrete unital monoidal category $J$,
and a $J$-module $\M$ in the sense of Definition~\ref{C.mod.def},
with $M = M_{[0]}$ and the action functor \eqref{rho.m}. Then since
$\Delta^oJ_{red} \to \Delta^o$ is discrete, the precofibration
$e^*\M \to \Delta^oJ_{red}$ induced by $\mu:\M \to BJ$ is a
cofibration, and since \eqref{rho.M} is an equivalence, the
cofibration is special. Then by Lemma~\ref{sp.le}, we have $e^*\M
\cong \xi^*\M^{red}$ for a unique cofibration $\M^{red} \to
J_{red}$. Explicitly, all fibers $\M^{red}_j$, $j \in J_{red}$ are
identified with $M=\M_{[0]}$, and the transition functor
corresponding to $\langle j'',a \rangle$ is $m(-,j'')$. In
particular, the fiber $\M^{red}_o$ over the initial object $o \in
J_{red}$ is identified with $M$, and \eqref{s.C} gives a functor
\begin{equation}\label{s.m}
\sigma:M \times J_{red} \to \M^{red}, \quad x \times j \mapsto
\langle j,\id \rangle_!x \cong m(x,j), \ x \in M, j \in J
\end{equation}
over $J_{red}$. Moreover, say that $\M$ is {\em reflexive} if for
any $j \in J$, the functor $m(-,j):M \to M$ has a right-adjoint
$j^*:M \to M$. Then in this case, the cofibration $\M^{red} \to
J_{red}$ is a bifibration, and \eqref{beta.eq} provides a functor
\begin{equation}\label{o.m}
\omega = \zeta(o):\M^{red} \to M
\end{equation}
sending $x \in M = M_j \subset \M^{red}$, $j \in J$ to $j^*x \in
M$. If we have two $I$-modules $\M_0$, $\M_1$, then a morphism
$\alpha:\M_0 \to \M_1$ of $I$-module induces a functor
$\alpha^{red}:\M_0^{red} \to \M_1^{red}$ over $I_{red}$. If $\M_1$
is also reflexive, we can consider the composition
$$
\begin{CD}
M_0 \times J_{red} @>{\sigma}>> \M_0^{red} @>{\alpha^{red}}>>
  \M_1^{red} @>{\omega}>> M_1
\end{CD}
$$
that gives rise to a functor
\begin{equation}\label{sp.m}
M_0 \to \Fun(J_{red},M_1)
\end{equation}
as soon as $J_{red}$ is bounded, so that the right-hand side is
well-defined.

To apply this abstract machinery to stabilization, it is convenient
to use both the category $\Top_+$ of pointed compactly generated
topological spaces and the category $\Delta^o\Sets_+$ of pointed
simplicial sets. Both are unital symmetric monoidal with respect to
the smash product, and we have an adjoint pair
\begin{equation}\label{real.eq}
\begin{CD}
\Top_+ @>{\Sing}>> \Delta^o\Sing_+ @>{\Real}>> \Top_+
\end{CD}
\end{equation}
of the singular complex functor $\Sing$ and the geometric
realization functor $\Real$. By Milnor Theorem, the latter is
symmetric monoidal (see \cite{dr} for a modern proof). We also have
the tautological embedding $\Gamma_+ \to \Sets_+$ that induces an
embedding $\eps:\Delta^o\Gamma_+ \to \Delta^o\Sets_+$.

Now, \eqref{b.o.I.fun} with $\E=\Sets_+$ and $I=\Delta^o$ provides a
morphism
\begin{equation}\label{b.o.de}
\oFun(\Gamma_+,\Delta^o\Sets_+) \to \bFun(\Gamma_+,\Delta^o\Sets_+)
\end{equation}
of modules over the unital monoidal category $\Delta^o\Gamma_+$, and
its target is equivalent to $\Fun(\Gamma_+,\eps^*\Delta^o\Sets_+)$,
where $\Delta^o\Sets_+$ is considered as module over itself. Then the
realization functor \eqref{real.eq} induces a
$\Delta^o\Gamma_+$-module morphism
\begin{equation}\label{re.de}
\Fun(\Gamma_+,\eps^*\Delta^o\Sets_+) \to
\Fun(\Gamma_+,\eps^*\Real^*\Top_+),
\end{equation}
and for any unital monoidal category $J$ equipped with a monoidal
functor $\gamma:J \to \Delta^o\Gamma_+$, the composition of
morphisms \eqref{b.o.de} and \eqref{re.de} induces a $J$-module
morphism
\begin{equation}\label{a.a}
\alpha(\gamma):\gamma^*\oFun(\Gamma_+,\Delta^o\Sets_+) \to
\Fun(\Gamma_+,\gamma^*\eps^*\Real^*\Top_+).
\end{equation}
Moreover, $\Top_+$ is reflexive as a module over itself. Therefore
the target of the morphism \eqref{a.a} is also reflexive, and if $J$
is essentially discrete, \eqref{a.a} gives rise to the corresponding
functor \eqref{sp.m}. Precomposing it with $\Sing$ of
\eqref{real.eq} then gives a functor
\begin{equation}\label{sp.I}
\Sp(J,\gamma):\Fun(\Gamma_+,\Top_+) \to \Fun(J_{red} \times
\Gamma_+,\Top_+).
\end{equation}
By construction, \eqref{sp.I} is functorial with respect to pairs
$\langle J,\gamma \rangle$ --- if we have another essentially
discrete unital monoidal category $J'$ and a monoidal functor
$\phi:J' \to J$, then we have a canonical isomorphism
\begin{equation}\label{sp.phi}
\phi_{red}^*\Sp(J,\gamma) \cong \Sp(J',\phi^*\gamma),
\end{equation}
where $\phi_{red}:J'_{red} \to J_{red}$ is induced by $\phi$. In
particular, we can always replace $J$ with its isomorphism groupoid
$\overline{J}$, and by abuse of notation, we will denote
$\Sp(J,\gamma) = \Sp(\overline{J},\gamma)$ even when $\gamma$ is
only defined on $\overline{J}$. Also by construction, the functor
\eqref{sp.I} respects homotopy equivalences, so we can define an
endofunctor $\Stab(J,\Sigma)$ of the category $\Ho(\Gamma_+)$ by
setting
\begin{equation}\label{stab.I}
\Stab(J,\gamma) = \pi_! \circ \Sp(J,\gamma),
\end{equation}
where $\pi:J_{red} \times \Gamma_+ \to \Gamma_+$ is the projection.
Then \eqref{sp.phi} provides a  map
\begin{equation}\label{stab.phi}
\Stab(J',\phi^*\gamma) \to \Stab(J,\gamma).
\end{equation}
In the trivial case $J' = \ppt$, this is simply a functorial map
$\Id \to \Stab(J,\gamma)$.

If we take $J = \Delta^<$, as in Example~\ref{de.red.exa}, and let
$\gamma:\overline{\Delta^<} \to \Delta^o\Gamma_+$ be the only
monoidal functor that sends $[0]$ to the simplicial circle $\Sigma$
of \eqref{Si}, then $\Stab(\Delta^<,\gamma)$ is precisely the
stabilization functor $\Stab$ of \eqref{stab.bo} constructed in
Proposition~\ref{stab.prop}. However, $\Gamma_+$ is symmetric
monoidal; therefore we can also take $J=\Gamma$, as in
Example~\ref{ga.red.exa}, with the monoidal functor
$\gamma:\overline{\Gamma} \to \Delta^o\Gamma_+$ of
Example~\ref{ga.iso.exa} sending $V([0]) = \ppt \in \Gamma$ to
$\Sigma$.

\begin{lemma}\label{I.stab.le}
The map $\Stab = \Stab(\Delta^<,V^*\gamma) \to \Stab(\Gamma,\gamma)$
of \eqref{stab.phi} induced by the embedding $V:\Delta^< \to \Gamma$
of Example~\ref{ga.red.exa} is an isomorphism.
\end{lemma}

\proof{} By \eqref{sp.phi}, for any $X \in \Ho(\Gamma_+)$,
$V_{red}^*\Sp(\Gamma,\gamma)(X) \in \Ho(\N \times \Gamma_+)$ is
represented by the inductive system $\Omega^nB^nX$ of
\eqref{stab.bo}, so that in particular, it is cocartesian over $\N$
for stable $X$. Since every map in $\Gamma_{red}$ is isomorphic to a
map in the image of $V_{red}$, $\Sp(\Gamma,\gamma)(X)$ then is
cocartesian over $\Gamma_{red}$, and if we let $j:\ppt \to
\Gamma_{red}$ be the embedding onto the initial object $o \in
\Gamma_{red}$, the adjunction map $\Sp(\Gamma,\gamma)(X) \to (j
\times \id)_!X$ is an isomorphism. Therefore for a stable $X$, the
canonical map $X \to \Stab(\Gamma,\gamma)(X)$ is an isomorphism, and
to finish the proof, it remains to show that
$\Stab(\Gamma,\gamma)(X)$ is stable for any $X$. Moreover, again by
\eqref{sp.phi}, the functor $\Ho(\Sp(\Gamma,\gamma)(X)):\Gamma_{red}
\to \Ho(\Gamma_+)$ sends a set $S \in \Gamma_{red}$ of cardinality
$n$ to $\Omega^nB^nX$, and the latter is $n$-stable.  Let
$\Gamma^{\geq n}_{red} \subset \Gamma_{red}$ be the subcategory of
sets of cardinality $\geq n$, and assume for the moment that we know
the following:
\begin{itemize}
\item for any $n \geq 0$ and any $Y \in \Ho(\Gamma_{red})$, the
  natural map
$$
\hocolim_{\Gamma_{red}^{\geq n}}Y|_{\Gamma_{red}^{\geq
      n}} \to \hocolim_{\Gamma_{red}}Y
$$
is an isomorphims.
\end{itemize}
Then we are done: by the base change isomorphism \eqref{bc.eq},
$\Stab(\Gamma,\gamma)(X)$ is $n$-stable for any $n$, thus stable.

The argument for \thetag{$\bullet$} is non-trivial but well-known
(apparently it goes back to \cite[Proposition
  VI.4.6.12]{ill}). Since it is also very short, we reproduce it for
the convenience of the reader. By induction on $n$, it suffices to
show that the embedding $\Gamma_{red}^{\geq n+1} \subset
\Gamma_{red}^{\geq n}$ is homotopy cofinal for any $n$. By
Remark~\ref{cof.rem}, we only have to consider the right comma-fiber
over the single object in $\Gamma_{red}^{\geq n}$ that is not in
$\Gamma_{red}^{\geq n+1}$, and for any $n$, this right comma-fiber
is equivalent to the category $I = \Gamma^{\geq 1}_{red}$ of
non-empty finite sets and injective maps. To prove that it is
homotopy contractible, consider the product $I \times I$, the
diagonal embedding $\delta:I \to I \times I$ and the coproduct
functor $\mu:I \times I \to I$. We then have a natural map $\Id \to
\delta \circ \mu$, and also a map $\Id \to \mu \circ \delta$ (say
given by the embedding $S \to S \copr S$ onto the first copy of
$S$). Then by Lemma~\ref{Qui.le}, $\delta$ is a weak
equivalence. Then so are its one-sided inverses $p_1,p_2:I \times I
\to I$ given by the projection onto the two factors, and then also
the section $i:I \to I \times I$ of the projection $p_1$ given by
the embedding onto $I \times S$ for some fixed $S \in I$. But then
$p_2 \circ i$ factors through the embedding $\ppt \to I$ onto $S$.
\endproof

\subsubsection{Stabilization for trace theories.}

We now observe that the construction of the functor \eqref{sp.m}
works in the relative setting. Namely, as in
Definition~\ref{sym.I.def}, define a unital monoidal structure on a
cofibration $\C \to I$ as a cofibration $B\C \to \Delta^o \times I$
whose restriction to $\Delta^o \times i$, $i \in I$ is a unital
monoidal structure on $\C_i$. Moreover, define a $\C$-module as a
cofibration $\M \to \Delta^o \times I$, with the fiber $M = \M|_{[0]
  \times I}$, equipped with a functor $\mu:\M \to B\C$ over
$\Delta^o \times I$ that is cocartesian over all maps $f \times f'$,
$f$ special, and such that the functor
$$
\zeta([0]) \times (\rho^o \times \id)^*\mu:M \times_I (\rho^o \times
\id)^*B\C
$$
is an equivalence of categories. Then for any $i \in I$, the
restriction $\M|_{\Delta^o \times i}$ is a $\C_i$-module in the
sense of Definition~\ref{C.mod.def}, and we have the action functor
\begin{equation}\label{rho.I.m}
m:M \times_I \C \to M.
\end{equation}
Moreover, if $M$ is bounded, and we have a cocomplete category $\E$,
then the fibration $\Fun(M/I,\E)^\perp \to I$ is also a cofibration,
with transition functors given by left Kan extension $(f_!)_!$ with
respect to transition functors $f_!$ of the cofibration $M \to I$,
and we have the action functor
\begin{equation}\label{rho.I.m.E}
m:\Fun(M/I,\E)^\perp \times_I \C \to \Fun(M/I,\E)^\perp
\end{equation}
whose fiber $m_i$ over some $i \in I$ is the functor \eqref{rho.m.E}
for the $\C_i$-module $M_i$, and the maps \eqref{fi.fu} are the base
change maps \eqref{bc.eq} induced by the corresponding maps for the
functor \eqref{rho.I.m}. Then identically the same construction as
in Subsection~\ref{C.mod.subs} provides a $\C$-module
$\oFun(\M/I,\E)$ with $\oFun(\M,\E)|_{[0] \times I} \cong
\Fun(M/I,\E)^\perp$ and the action functor
\eqref{rho.I.m.E}. Moreover, say that a unital monoidal category
$J/I$ is essentially discrete if so is each of the fibers $J_i$, $i
\in I$; then for such a $J$, we have the cofibration $J_{red} \to
I$, and if $J_{red}$ is bounded, then for any morphism $\alpha:\M_0
\to \M_1$ of $J$-modules with reflexive $\M_1$, we have the relative
version
\begin{equation}\label{sp.I.m}
M_0 \to \Fun(J_{red}/I,M_1)
\end{equation}
of the functor \eqref{sp.m}.

To apply this to stabilization, take an integer $n \geq 1$, let
$m:\Gamma_+^n \to \Gamma_+$ be the iterated smash product functor,
and note that by induction on $n$, for any half-additive category
$\E$ that has kernels, the composition $m^* \circ m^\dg$ of the
embedding $m^\dg:\E \to \Fun(\Gamma_+,\E)$ and the pullback functor
$m^*:\Fun(\Gamma_+,\E) \to \Fun(\Gamma_+^n,\E)$ is fully faithful
and admits a right-adjoint. Therefore as in
Subsection~\ref{C.mod.subs}, we have an embedding $m^*\E^\otimes
\subset \Fun(\Gamma_+^n,\E)$ of $\Gamma_+^n$-modules and its
right-adjoint. We also have the morphisms \eqref{b.o.fun} and
\eqref{b.o.I.fun} with $\Gamma_+$ replaced by $\Gamma_+^n$. More
generally, consider the cyclic nerve $\Gamma^\hash_+ \to \Lambda$ of
the category $\Gamma_+$, with the monoidal structure induced by the
symmetric monoidal structure on $\Gamma_+$, and let
$\Id_{triv}:\Gamma_+^\hash \to \Gamma_+$ be the functor
\eqref{i.triv}. Then as soon as $\E$ is cocomplete, we have a fully
faithful embedding
$$
\Id_{triv}^* \circ m^\dg:\Id_{triv}^*\E \to
\oFun(\Gamma_+^\hash/\Lambda,\E)
$$
of $\Gamma_+^\hash$-modules, and for any bounded $I$ and
half-additive cocomplete $\E$ with kernels, we have a morphism
\begin{equation}\label{b.o.I.la}
\oFun(I \times \Gamma^\hash_+/\Lambda,\E) \to \Id_{triv}^*\bFun(I
\times \Gamma_+,\E)
\end{equation}
of $\Fun(I,\Gamma^\hash_+/\Lambda)$-modules, a relative version of
the morphism \eqref{b.o.I.fun}. We can now take $\E = \Sets_+$ and
$I = \Delta^o$, let $J =\Gamma^\hash$, where $\Gamma$ is the
essentially discrete symmetric monoidal category of
Example~\ref{ga.red.exa}, equip it with the functor
$\gamma^\hash:\overline{\Gamma}^\hash \to
\Fun(\Delta^o,\Gamma^\hash_+/\Lambda)$ induced by
$\gamma:\overline{\Gamma} \to \Delta^o\Gamma_+$ of
Lemma~\ref{I.stab.le}, replace \eqref{sp.m} with \eqref{sp.I.m}, and
repeat the construction of Subsection~\ref{rec.subs} to obtain a
functor
\begin{equation}\label{sp.I.la}
\Sp(\Gamma^\hash,\gamma^\hash):\Fun(\Gamma^\hash_+/\Lambda,\Top_+)^\perp
\to \Fun(\Gamma^\hash_{red} \times_\Lambda
\Gamma_+^\hash/\Lambda,\Top_+)^\perp
\end{equation}
over $\Lambda$, a relative version of the functor \eqref{sp.I}. It
still respects weak equivalences, thus descends to a functor
$\ho(\Gamma^\hash_+/\Lambda) \to \ho(\Gamma_{red}^\hash
\times_\Lambda \Gamma^\hash_+/\Lambda)$, and we can further define a
functor
\begin{equation}\label{stab.la}
\Stab^\Lambda = \pi_! \circ
\Sp(\Gamma^\hash,\gamma^\hash):\ho(\Gamma_+/\Lambda) \to
\ho(\Gamma_+/\Lambda)
\end{equation}
over $\Lambda$, where $\pi:\Gamma_{red}^\hash \times \Gamma^\hash_+
\to \Gamma^\hash_+$ is the projection.

\begin{lemma}\label{stab.la.le}
The functor \eqref{stab.la} is cartesian over any map $f:[n] \to
[1]$.
\end{lemma}

\proof{} For any $[q] \in \Lambda$, we have $(\Gamma^\hash_+)_{[q]}
\cong \Gamma_+^q$, and under these identifications, the transition
functor $f_!$ corresponding to the map $f$ is isomorphic to the
iterated product functor $m:\Gamma_+^n \to \Gamma_+$.  For any $X
\in \ho(\Gamma^\hash_+/\Lambda)_{[1]} \cong \Ho(\Gamma_+)$, we have
$\Stab^\Lambda(X) = \Stab(\Gamma,\gamma)(X)$, where
$\Stab(\Gamma,\gamma)$ is as in Lemma~\ref{I.stab.le}. Moreover, by
construction, the functor \eqref{sp.I.la} is cocartesian, so that
$\Stab^\Lambda(m^*X) \cong m^*\Stab(\Gamma^n,\phi^*\gamma)(X)$,
where $\phi:\Gamma^n \to \Gamma$ is the iterated product functor for
the monoidal structure on $\Gamma$ (that is, the iterated coproduct
of finite sets). What we have to check, then, is that the map
\eqref{stab.phi} induced by $\phi$ is an isomorphism. But by
Lemma~\ref{I.stab.le}, its target is isomorphic to $\Stab$, and
exactly the same argument proves that the same is true for its
source.
\endproof

\proof[Proof of Proposition~\ref{la.stab.prop}.] Consider the
product $B\Gamma_+ \times^2 \C$, with the corresponding $2$-functor
\eqref{2.ga.m} and the $2$-functor $t = \eps(o) \times \id:\C \to
B\Gamma_+ \times^2 \C$, and consider the corresponding diagram
$$
\begin{CD}
\Lambda\C @>{\Lambda t}>> \Gamma^\hash_+ \times_\Lambda \Lambda\C
@>{\Lambda m}>> \Lambda\C,
\end{CD}
$$
of cyclic nerves and functors cocartesian over $\Lambda$. Define an
endofunctor $\Stab$ of the homotopy category $\Ho(\Gamma^\hash_+
\times_\Lambda \Lambda\C)$ as the composition
$$
\begin{CD}
\Ho(\Gamma_+^\hash \times_\Lambda \Lambda\C)
@>{\Sp(\Gamma^\hash,\gamma^\hash)}>> \Ho(\Gamma_{red}^\hash
\times_\Lambda \Gamma_+^\hash \times_\Lambda \Lambda\C) @>{\pi_!}>>
\Ho(\Gamma_+^\hash \times_\Lambda \Lambda\C)
\end{CD}
$$
where $\pi:\Gamma_{red}^\hash \times_\Lambda \Gamma_+^\hash
\times_\Lambda \Lambda\C \to \Gamma_+^\hash \times_\Lambda
\Lambda\C$ is the projection, and $\Sp(\Gamma^\hash,\gamma^\hash)$
is induced by the functor \eqref{stab.la} for the cofibration
$\Gamma_+^\hash \times_\Lambda\Lambda\C \to \Lambda\C$. Let
$$
\Stab_{\C} = \Lambda t^* \circ \Stab \circ \Lambda m^*:\Ho(\Lambda\C)
\to \Ho(\Lambda(\C)).
$$
Then we have a functorial map $\Id \to \Stab_{\C}$, and by
definition, we have $i_c^* \circ \Stab_{\C} \cong \Stab_{\C(c,c)}
\circ i_c^*$ for any $c \in \C_{[0]}$. Thus to prove the claim, it
suffices to check that for any homotopy trace theory $E$,
$\Stab_{\C}(E)$ is a homotopy trace theory. But since every map $[n]
\to [m]$ in $\Lambda$ can be composed with a map $[m] \to [1]$, it
suffices to check that $\Ho(\Stab_{\C}(E))$ is cocartesian over maps
$[n] \to [1]$, and this immediately follows from
Lemma~\ref{stab.la.le}. The argument for $\Ho_{tr}(-,k)$ is
identically the same.
\endproof

\subsection{Homotopy invariance.}\label{inva.subs}

For any bounded half-additive pointed $2$-category $\langle \C,o
\rangle$, the restriction functor \eqref{red.eq} commutes with
stabilization, thus induces a functor $\Red^{st}:\Ho_{tr}^{st}(\C)
\to \Ho_{tr}^{st}(B\C_o)$, and similarly for $\Ho(-,k)$. In the
$k$-linear case, stabilization commutes with homotopy colimits,
hence also with expansion, so that \eqref{exp.eq} provides a fully
faithful left-adjoint functor $\Exp^{st} \cong \Exp^{ho}$ to
$\Red^{st}$. In the absolute case, $\Red^{st}$ still admits a fully
faithful left-adjoint given by
\begin{equation}\label{st.exp.eq}
\Exp^{st} = \Stab_{\C} \circ \Exp^{ho}:\Ho_{tr}^{st}(B\C_o) \to
\Ho^{st}_{tr}(\C),
\end{equation}
while \eqref{phi.stab} provides an isomorphism
\begin{equation}\label{phi.exp}
\phi \circ \Ext^{st} \cong \Exp^{st} \circ \phi.
\end{equation}
However, the absolute $\Stab_{\C}$ no longer commutes with
$\Exp^{ho}$, and it is not true that the essential image of the
functor \eqref{st.exp.eq} consist of stable $S$-homotopy exact trace
theories. To describe this essential image, we need to modify the
notion of homotopy exactness. It turns out that the appropriate
condition in the stable case is actually shorter.

\begin{defn}\label{st.tr.exa.def}
For any bounded pointed half-additive $2$-category $\langle \C,o
\rangle$ equipped with a base $S$, a stable homotopy trace theory $E
\in \Ho_{tr}^{st}(\C)$ is {\em stably $S$-homotopy exact} if for any
object $c \in \C_{[0]}$ and $S$-contractible functor $g:\Delta^{o>}
\to \C(c,c)$, the pullback $g^*i_c^*E \in \Ho^{st}(\Delta^{o>})$ is
homotopy exact.
\end{defn}

The difference with Definition~\ref{tr.exa.def} is that we do not
consider an additional contractible simplicial set $X$. Let us show
that this is enough.

\begin{lemma}\label{st.exa.le}
In the assumptions of Definition~\ref{st.tr.exa.def}, assume further
that $\C$ is large, and the base $S$ is ind-admissible. Then a
stable homotopy trace theory $E \in \Ho^{st}_{tr}(\C)$ is stably
$S$-homotopy exact if and only if the natural map $E \to
\Exp^{st}(\Red^{st}(E))$ is an isomorphism.
\end{lemma}

\proof{} For any small category $I$, with the tautological
projection $\tau:I \to \ppt$, the embedding $e:I \to I^>$, and the
embedding $i:\ppt \to I^>$ onto $o \in I^>$, say that an object $E
\in \Ho^{st}(\Gamma_+, I^>)$ is {\em stably homotopy exact} if
the augmentation map $\Stab((\id \times \tau)_!(\id \times e)^*E)
\to (\id \times i)^*E$ is an isomorphism in
$\Ho^{st}(\Gamma_+)$. Since $\hocolim_{\Delta^o}$ commutes with
finite products, an object $E \in \Ho^{st}(\Gamma_+,
\Delta^o)$ is stably homotopy exact if and only if its restriction
to $\ppt_+ \times \Delta^o$ is homotopy exact. In particular, for
any bounded half-additive category $\E$ and small category $I$, a
functor $\gamma:I^> \to \E$ uniquely extends to a half-additive
functor $\gamma_+:\Gamma_+ \times I^> \to \E$, and if $I =
\Delta^{o>}$, then $\gamma_+^*E$ for some $E \in \Ho^{st}(\E)$ is
stably homotopy exact if and only if $\gamma^*E$ is homotopy exact.

If the bounded half-additive category $\E$ is also large, then it
has filtered colimits, thus all coproducts. Then for any simplicial
set $X$ with the projection $\pi:\Delta^{o>}X = (\Delta^o X)^> \to
\Delta^{o>}$, and any functor $g:\Delta^{o>}X \to \E$, the left Kan
extension $(\id \times \pi)_!g_+ = ((\id \times \pi_!g)_+:\Gamma_+
\times \Delta^{o>} \to \E$ exists. If $X$ is contractible and $g$ is
contractible in the sense of Definition~\ref{aug.def}, then $(\id
\times \pi)_!g$ is also contractible. Moreover, for any $X$ and any
$E \in \Ho^{st}(\E)$, we can consider the natural map
\begin{equation}\label{pi.sta}
\Stab((\id \times \pi)_!g_+^*E) \to ((\id \times \pi)_!g_+)^*E,
\end{equation}
and by the same argument as in the proof of
Proposition~\ref{mult.prop}, \eqref{Sig.st} shows that
\eqref{pi.sta} is an isomorphism in $\Ho^{st}(\Gamma_+,
\Delta^{o>})$. Therefore $(\id \times g_+)^*E$ is stably homotopy
exact iff $((\id \times \pi)_!g_+)^*E$ is homotopy exact.

Applying this in the setting of the Lemma, we conclude that $E \in
\Ho^{st}_{tr}(\C)$ is stably $S$-homotopy exact if and only if for
any contractible $X \in \Delta^o\Sets_+$, any $c \in \C_{[0]}$, and
any functor $g:\Delta^{o>}X \to \C(c)$ that is $S$-contractible in
the sense of Definition~\ref{tr.exa.def}, the object $g_+^*i_c^* \in
\Ho^{st}(\Gamma_+, \Delta^{o>}X)$ is stably homotopy exact. To
finish the proof, it remains to apply the same argument as in
Theorem~\ref{rec.thm}, with $\Exp$ replaced by $\Exp^{st}$ of
\eqref{st.exp.eq}.
\endproof

One can also rewrite Definition~\ref{st.tr.exa.def} in the following
way. For any bounded $2$-category $\C$, the cofibration $\C^\Delta =
\Fun(\Delta^o,\C/\Delta^o) \to \Delta^o$ is a $2$-category --- its
objects are the objects $c \in \C_{[0]}$, and morphism categories
are $\C^\Delta(c,c') = \Fun(\Delta^o,\C(c,c'))$, so that morphisms
in $\C^\Delta$ are simplicial objects in the morphism categories
$\C(c,c')$. We have the evaluation functor $\ev:\Delta^o \times
\C^\Delta \to \C$ and the projection $\pi:\Delta^o \times \C^\Delta
\to \C^\Delta$, these induces the corresponding functors on cyclic
nerves, and for any homotopy trace theory $E \in \Ho_{tr}(\C)$, we
can define its {\em simplicial extension} $E^\Delta$ by
$$
E^\Delta = \pi_!\ev^*E \in \Ho_{tr}(\C^\Delta).
$$
If $\C$ is half-additive, then $\C^\Delta$ is also half-additive,
and if $E$ is stable, then $E^\Delta$ is stable. Moreover, assume
that $\C$ is pointed, and $\langle \C,o \rangle$ is equipped with a
base $S$. Then any functor $g:\Delta^{o>} \to \C(c,c)$, $c \in
\C_{[0]}$ defines objects $g_0 = g|_{\Delta^o} \in
\Fun(\Delta^o,\C(c,c))$, $g(o) \in \C(c,c)$, and an augmentaion map
\begin{equation}\label{aug.C}
a:g_0 \to g(o)
\end{equation}
in $\Fun(\Delta^o,\C(c,c)) \subset \Lambda\C^\Delta_{[1]} \subset
\Lambda\C^\Delta$ whose target is the constant simplicial object
with value $g(o)$. Say that an augmentation map \eqref{aug.C} is
{\em $S$-contractible} if so is the functor $g:\Delta^{o>} \to
\C(c,c)$. Then $E$ is stably $S$-homotopy exact if and only if
$\ho(E^\Delta)$ inverts all $S$-contractible augmentation maps
\eqref{aug.C}.

This interpretation is especially useful if we take the unital
monoidal category $k\amod^{fl}$ of flat modules over a commutative
ring $k$, and let $\C = \aAlg(k) \subset \iMor(k\amod^{fl})$ be the
full subcategory in the corresponding Morita $2$-category spanned by
enriched categories with one object. Then explicitly, objects in
$\aAlg(k)$ are flat $k$-algebras $A$, and morphisms from $A_0$ to
$A_1$ are modules $M$ over $A_0^o \otimes_k A_1$ that are flat over
$A_1$. In particular, $\aAlg(k)(A,A) \cong A\bimod^{fl} \subset
A\bimod$ is the category of $A$-bimodules that are flat as right
$A$-modules (``right-flat''). By virtue of the Dold-Kan equivalence,
morphisms in $\aAlg(k)^\Delta$ are then chain complexes of such
module concentrated in non-negative homological degrees. The
$2$-category $\aAlg(k)$ is pointed, with $o$ being $k$ itself, and it
has a standard ind-admissible base $S$ formed by the unit maps $k
\to A$, $A \in \aAlg(k)_{[0]}$.

\begin{defn}\label{inv.tr.def}
A homotopy trace theory $E \in \Ho_{tr}(\aAlg(k))$ is {\em
  homoto\-py-invariant} if for any $A \in \aAlg(k)_{[0]}$,
$\Ho(i_A^*E^\Delta):C_\idot(A\bimod) \to \Ho$ inverts
quasiisomorphisms.
\end{defn}

\begin{lemma}\label{inv.le}
A stable homotopy trace theory $E \in \Ho_{tr}^{st}(\aAlg(k))$ is
homo\-to\-py-invariant in the sense of Definition~\ref{inv.tr.def}
if and only if it is stably $S$-homotopy exact with respect to the
standard base $S$.
\end{lemma}

\proof{} By the Dold-Kan equivalence, an augmented simplicial object
$\Delta^{o>} \to A\bimod$ for some $A \in \aAlg(k)_{[0]}$ is the same
thing as a triple $\langle M_\idot,M,a \rangle$ of a complex
$M_\idot \in C_{\geq 0}(A\bimod)$, an object $M \in A\bimod$ and a
map $a:M_\idot \to M$, where $M$ is treated as a complex placed in
degree $0$. Such a triple is $S$-contractible if and only if the
cone $M'_\idot$ of the map $a$ is contractible as a complex of right
$A$-modules (``right-contractible''). In this case, $a$ is a
quasiisomorphism, so that a homotopy-invariant stable homotopy trace
theory $E$ on $\aAlg(k)$ is stably $S$-homotopy exact.

To prove the converse, note that every short exact sequence
\begin{equation}\label{short}
\begin{CD}
0 @>>> M @>{f}>> M' @>>> M'' @>>> 0
\end{CD}
\end{equation}
in $A\bimod^{fl}$ gives rise to a cocartesian square
$$
\begin{CD}
M @>>> 0\\
@V{f}VV @VV{f'}V\\
M' @>>> M''
\end{CD}
$$
that we can interpret as a functor $\gamma:[1]^2 = \V^> \to
A\bimod^{fl} \subset \Lambda\aAlg(k)$. If \eqref{short} is split,
then for any $E \in \Ho^{st}_{tr}(\aAlg(k))$, $\gamma^*E$ is stably
homotopy exact in the same sense as in the proof of
Lemma~\ref{st.exa.le}. Moreover, since $\hocolim_{\Delta^o}$
commutes with stabilization, the same is true for $E^\Delta$ and a
sequence \eqref{short} in $C_{\geq 0}(A\bimod^{fl})$, and it
suffices to require that the sequence is termwise split. However, it
is well-known that homotopy cocartesian squares in
$\Ho^{st}(\Gamma_+)$ coincide with homotopy cartesian
ones. Therefore for any termwise split sequence \eqref{short} in
$C_{\geq 0}(A\bimod)$, $\ho(E^\Delta)$ inverts the map $f$ iff it
inverts the map $f'$ (that is, annihilates $M''$).

Now assume given a stable homotopy trace theory $E \in
\Ho^{st}_{tr}(\aAlg(k))$ that is stably $S$-homotopy exact, and note
that for any map $f:M_\idot \to M_\idot'$ in $C_{\geq
  0}(A\bimod^{fl})$, with cone $C_\idot$ and cylinder $N_\idot$, we
have termwise-split short exact sequences
$$
\begin{CD}
0 @>>> M_\idot @>>> N_\idot @>>> C_\idot @>>> 0,\\
0 @>>> M'_\idot @>>> N_\idot @>>> C'_\idot @>>> 0,
\end{CD}
$$
where the cone $C'_\idot$ of the identity map $\id:M_\idot \to
M_\idot$ is contractible. Therefore $\ho(E^\Delta)$ inverts $f$ if
and only if it annihilates $C_\idot$. But if $f$ is a
quasiisomorphism, $C_\idot$ is an acyclic complex of right-flat
$A$-bimodules. Therefore it is a filtered colimit of acyclic
complexes of right-projective $A$-bimodules, and these are
right-contractible, thus annihilated by $\ho(E^\Delta)$.
\endproof

By virtue of Lemma~\ref{inv.le} and Lemma~\ref{st.exa.le}, any
homotopy trace functor on the unital monoidal category $k\amod^{fl}$
extends uniquely and canonically to a homotopy-invariant stable
homotopy trace theory $E$ on the $2$-category $\aAlg(k)$, with its
simplicial extension $E^\Delta \in
\Ho^{st}_{tr}(\aAlg(k)^\Delta)$. Moreover, we can invert
quasiisomorphisms in the morphism categories of $\aAlg(k)^\Delta$ to
obtain a $2$-category $\aAlg(k)^\D$ whose objects are again flat
$k$-algebras, and whose morphism categories are given by
$$
\aAlg(k)^\D(A_0,A_1) = \D^{\leq 0}(A_0^o \otimes_k A_1) \subset
\D(A_0^o \times_k A_1),
$$
where $\D(A_0^o \otimes_k A_1)$ is the derived category of left
$A_0^o \otimes_k A_1$-modules. Then being homotopy-invariant,
$E^\Delta$ descends to a homotopy trace theory $E^\D$ on
$\aAlg(k)^\D$. This is useful to know even if one is only interested
in $k$-algebras and bimodules since $\aAlg(k)^\D \supset \aAlg(k)$
has more adjoint pairs. In effect, any multiplicative, possibly
non-unital map $f:A \to B$ between two $k$-algebras gives rise to a
morphism in $\aAlg(k)$ given by $f(1)B \in A^o \otimes_k
B\bimod^{fl}$, and this morphism becomes reflexive in
$\aAlg(k)^\D$. Thus in particular, the functor \eqref{ho.radj}
associated to $E^\D$ induces a functor
\begin{equation}\label{e.alg}
\Alg(k) \subset \rAdj(\aAlg(k)^\D) \to \Ho_{tr}(\ppt^2),
\end{equation}
where $\Alg(k)$ is the category of flat $k$-algebras and non-unital
maps. Explicitly, any morphism in $\aAlg(k)$ can be represented by a
simplicial pointwise-flat $A_0^o \otimes_k A$-module, and by
Proposition~\ref{refl.prop}, these correspond to morphisms in
$\Cat(k\amod^{fl})$ that are reflexive in the sense of
Definition~\ref{A.refl.def}. Therefore by Lemma~\ref{A.ha.le}, the
functor \eqref{e.alg} is the homotopical version of the functor
\eqref{CC.E}: it sends an algebra $A$ to the homotopical averaging
$\Av_\Lambda(EA_\hash)$ of the object $EA_\hash \in \Ho(\Lambda)$.

Finally, we note that Definition~\ref{inv.tr.def},
Lemma~\ref{inv.le} and the discussion above have an obvious
counterpart for the categories $\Ho(-,k')$, where $k'$ is a
commutative ring (possibly different from $k$). The proofs are
identically the same, and we leave the details to the reader.

\section{Topological Hochschild Homology.}\label{thh.sec}

\subsection{Bifunctor homology.}

Fix a commutative ring $k$. The shortest way to define Hochschild
homology modules $HH_\idot(A/k,M)$ of a flat $k$-algebra $A$ with
coefficients in an $A$-bimodule $M$ is by using the functor $\Tr_A$
of Example~\ref{tr.A.exa}: one considers the object
\begin{equation}\label{ch.M}
CH_\idot(A/k,M) = L^\hdot\Tr_A(M) \in \D(k),
\end{equation}
and defines $HH_\idot(A/k,M)$ as its homology modules. If we have
another commutative ring $k'$ with a map $f:k' \to k$, a flat
$k'$-algebra $A$, and a bimodule $M$ over $A \otimes_{k'} k$, then
we have a taugolocal quasiisomorphism
\begin{equation}\label{CH.bc}
CH_\idot(A/k',f_*M) \cong CH_\idot(A \otimes_{k'} k/k,M),
\end{equation}
where as in Example~\ref{norma.exa}, $f_*M$ is $M$ considered as an
$A$-bimodule by restriction of scalars. If we have two flat
$k$-algebras $A$, $A'$, then $\Tr_{A \otimes_k A'} \cong
\Tr_A \circ \Tr_{A'} \cong \Tr_{A'} \circ \Tr_A$, and for any
bimodules $M$, $M'$ flat over $k$, we have the K\"unneth isomorphism
\begin{equation}\label{kun.hh}
CH_\idot(A,M) \lotimes_k CH_\idot(A',M') \cong CH_\idot(A \otimes_k
A',M \otimes_k M').
\end{equation}
If one uses the standard bar resolution of the bimodule $M$ to
compute the derived functor, one obtains an explicit {\em Hochschild
  complex} $CH_\idot(A/k,M)$ representing \eqref{ch.M} with terms
$CH_i(A/k,M) = A^{\otimes_k i} \otimes_k M$ and a certain
differential $b$. If $M=A$ is the diagonal bimodule, one shortens
$HH_\idot(A/k,A)$ to $HH_\idot(A/k)$ and $CH_\idot(A/k,A)$ to
$CH_\idot(A/k)$. For more details, we refer the reader to \cite{Lo}
or \cite{FTadd}.

More generally, if $A_\idot$ is a flat DG algebra over $k$, with a
DG bimodule $M_\idot$, one defines $CH_\idot(A_\idot/k,M_\idot)$ by
taking the Hochschild complex $CH_\idot(-,-)$ termwise, and then
taking the product-total complex of the resulting bicomplex, and one
denotes by $HH_\idot(A_\idot/k,M_\idot)$ the homology modules of the
complex $CH_\idot(A_\idot,M_\idot)$ (for details, see \cite{kel}).

{\em Topological Hochschild Homology} $THH$ was introduced by
B\"okstedt \cite{bo} as a generalization of Hochschild homology to
the absolute setting: for any associative ring spectrum $A$ and
$A$-bimodule spectrum $M$, one constructs a simplicial spectrum
$CH(A,M)$, a generalization of the complex $CH_\idot(-,-)$, and
defines the spectrum $THH(A,M)$ as the homotopy colimit
$\hocolim_{\Delta^o}CH(A,M)$. To make this work, one needs to
understand a ``ring spectrum'' in an appropriate sense, and there
are various ways to do it. However, in any approach, an associative
algebra $A$ over a commutative ring $k$ gives an example of an
associative ring spectrum, so it makes sense to speak of the
spectrum $THH(A,M)$. Moreover, in such a situation, $THH(A,M)$ is
automatically a module over $k$, so that {\em a posteriori}, it is
not a spectrum but a complex $TCH_\idot(A,M) \in \D^{\leq 0}(k)$ of
$k$-modules. Its homotopy groups $THH_\idot(A,M) =
\pi_\idot(THH(A,M))$ are then homology groups of the complex
$TCH_\idot(A,M)$ (in particular, they are $k$-modules in a natural
way).

It was realized pretty soon after \cite{bo} that in the particular
case of an algebra $A$ over a ring $k$, $THH(A,M)$ actually admits
several alternative definitions, some of them manifestly
$k$-linear. One such uses bifunctor homology of
Subsection~\ref{hom.subs}. Say that a right module $P$ over a flat
$k$-algebra $A$ is {\em finite free} if $P \cong V \otimes_k A$ for
a finitely generated projective $k$-module $V$, and let $P(A)$ be
the category of finite free left $A$-modules. Then $P(A)$ is
$k$-linear and small, and any $A$-bimodule $M$ defines a
$k$-bilinear functor $P(M):P(A)^o \times P(A) \to k\amod$ sending $V
\times V'$ to $\Hom_A(V,V' \otimes_A M)$. The correspondence $M
\mapsto P(M)$ is an equivalence between $A\bimod$ and $\Fun_k(P(A)^o
\times P(A),k)$ that intertwines $\Tr_A$ with $\Tr_{P(A)}$, so that
we have $HH_\idot(P(A)/k,P(M)) \cong HH_\idot(A/k,M)$. The absolute
bifunctor homology is then naturally identified with $THH(A,M)$ --
we have
\begin{equation}\label{thh.fu}
CH_\idot(P(A),P(M)) \cong TCH_\idot(A,M) \in \D^{\leq 0}(k) \cong
\Ho^{st}(\Gamma_+,k),
\end{equation}
and the map \eqref{c.i.m.aug} provides a functorial map
$TCH_\idot(A,M) \to CH_\idot(A,M)$. The homology groups
$HH_\idot(P(A),P(M))$ are thereby identified with the homotopy
groups $THH_\idot(A,M)$. If we have another flat $k$-algebra $A'$
equipped with a $k$-flat bimodule $M'$, then we have
$$
TCH_\idot(A,M) \lotimes_k TCH_\idot(A',M') \cong CH_\idot(P(A)
\times P(A'),P(M) \boxtimes_k P(M')),
$$
and the tensor product functor $- \otimes_k -:P(A) \times P(A') \to
P(A \otimes_k A')$ induces a K\"unneth map
\begin{equation}\label{kun.thh}
TCH_\idot(A,M) \lotimes_k TCH_\idot(A',M') \to TCH_\idot(A \otimes_k
A',M \otimes_k M').
\end{equation}
Unlike \eqref{kun.thh}, the map \eqref{kun.thh} is usually not an
isomorphism.

To make bifunctor homology more amenable to computation, one can use
the stabilization functor \eqref{stab.k} and the equivalence
\eqref{Q.A.eq}. One observes that the equivalence $\Fun_k(P(A),k)
\cong A\amod$ extends to an equivalence $\D(Q_\idot(P(A)/k)) \cong
\D(Q_\idot(A/k))$, and similarly for $P(A)^o$. For bimodules, we
then obtain an equivalence
\begin{equation}\label{Q.P.A}
\Ho^{st}(P(A)^o \times P(A),k)
\cong \D^{\leq 0}(Q_\idot(A/k)^o \otimes_k Q_\idot(A/k)),
\end{equation}
where by abuse of notation, $\Ho^{st}(P(A)^o \times P(A),k) \subset
\Ho(P(A)^o \times P(A),k)$ is the full subcategory spanned by objects
stable with respect to both $P(A)$ and $P(A)^o$. 
Taking stabilization with respect to both $P(A)^o$ and $P(A)$,
in any order, and using \eqref{Q.P.A}, we obtain a functor
$$
\Stab:\Fun(P(A)^o \times P(A),k) \to \D^{\leq 0}(Q_\idot(A/k)^o
\otimes_k Q_\idot(A/k)).
$$
If we compute $\Stab$ by the resulution \eqref{Q.reso}, then for any
$M \in \Fun(P(A)^o \times P(A),k)$, the degree-$0$ term $\Stab(M)_0$
of the resulting complex $\Stab(M)_\idot$ is $M$ itself, and we have
an obvious map $\Tr_{P(A)}(M) \to HH_0(Q_\idot(A),\Stab(M)_\idot)$ that
gives rise to a map
\begin{equation}\label{C.Q}
CH_\idot(P(A),M) = L^\hdot\Tr_{P(A)}M \to
CH_\idot(Q_\idot(A),\Stab(M)_\idot)
\end{equation}
in $\D^{\leq 0}(k)$. For any object $P \in P(A)$ and functor $N \in
\Fun_k(P(A),k)$, the map \eqref{C.Q} for $M = k_P \boxtimes N$ is a
quasiisomorphism by \eqref{yo.M} and \eqref{yo.A}, and since every
$M \in \Fun_k(P(A)^o \times P(A),k)$ admits a resolution by functors
of this form, \eqref{C.Q} is also a quasiisomorphism for such
$M$. In particular, for any $M \in A\bimod$, we have a natural
isomorphism
\begin{equation}\label{C.Q.1}
HH_\idot(P(A),P(M)) \cong HH_\idot(Q_\idot(A/k),M),
\end{equation}
where we identify $\Stab(P(M)) \cong P(M)$ since $P(M)$ is already
additive. The target of this isomorphism is known as {\em Mac Lane
  Homology} and denoted $HM_\idot(A,M)$; it has been introduced by
Mac Lane \cite{Mc} back in 1956. Interpretation of Mac Lane Homology
in terms of bifunctor homology is the great discovery of Jibladze
and Pirashvili \cite{JP}, and the identification between either of
them and $THH$ is \cite{P.et.al} and \cite{PW}. The product
\eqref{kun.thh} corresponds under \eqref{C.Q.1} to a product induced
by \eqref{kun.hh} and the lax monoidal structure on the functor
$Q_\idot$.

\begin{remark}
As noticed in \cite{PW}, as soon one has some sufficiently
well-developed ``brave new algebra'' that allows one to work ``over
the sphere spectrum $\Ss$'', an identification $THH_\idot(A,M) \cong
HM_\idot(A,M)$ becomes obvious: this is simply \eqref{CH.bc} with
$k' = \Ss$.
\end{remark}

\subsection{Mac Lane cohomology.}\label{mcl.subs}

One advantage of working with Mac Lane Homology it that it has a
cohomology theory attached to it: for any flat $k$-algebra $A$ and
$A$-bimodule $M$, one defines the {\em Hochschild Cohomology}
$HH^\hdot(A,M)$ as $\RHom^\hdot_{A^o \times_k A}(A,M)$, with the
obvious generalization to DG algebras, and then {\em Mac Lane
  Cohomology} is given by
$$
HM^\hdot(A,M) = HH^\hdot(Q_\idot(A),M).
$$
Note that if $A$ is commutative, the lax monoidal structure on
$Q_\idot(-)$ and the product map $A \otimes_k A \to A$ turn
$HM_\idot(A)$ into a commutative algebra, and then if $A=k$ is a
field, $HM_\idot(k)$ is actually a commutative cocommutative Hopf
algebra, with the dual algebra given by $HM^\hdot(k)$. In this case,
the bifunctor $P(k):P(k)^o \times P(k) \to k\amod$ is isomorphic to
$T^* \boxtimes_k T$, where $T:P(k) \to k\amod$ sends $V \in P(K)$ to
itself, and $T^*:P(k)^o \to k\amod$ sends it to the dual vector
space $V^*$, and then \eqref{lot.I} provides an isomorphism
$$
HM_\idot(k) \cong \Tor^\hdot_{P(k)}(T^*,T),
$$
where $\Tor^\hdot_I$ denotes the homology modules of the derived
tensor product \eqref{lot.I}. For cohomology, we then have
$$
HM^\hdot(k) \cong \Ext^\hdot_{P(k)}(T,T),
$$
the $\Ext$-groups computed in the abelian category $\Fun(P(k),k)$
(or equivalently, in the derived category $\D(Q_\idot(k,k)) \cong
\Ho^{st}(P(k),k) \subset \D_{\geq 0}(P(k),k)$).

Another advantage is computational. Since the complex computing the
Mac Lane Homology $HM_\idot(A,M)$ is obtained by totalizing a
bicomplex, and similarly for $HM^\hdot(A)$, we have spectral
sequences
\begin{equation}\label{sp.hm}
\begin{aligned}
HH_\idot(HQ_\idot(A/k),M) &\Rightarrow HM_\idot(A,M), \\
HH^\hdot(HQ_\idot(A/k),M) &\Rightarrow HM^\hdot(A,M),
\end{aligned}
\end{equation}
where $HQ_\idot(A/k)$ is the homology algebra of the DG algebra
$Q_\idot(A/k)$ (with zero differential). If $A=M=k$ is a perfect
field of some positive odd characteristic $p = \cchar k \geq 3$,
then \eqref{st.odd} induces an isomorphism
\begin{equation}\label{hq.st}
HQ_\idot(k/k) = (k \otimes k)[\beta,\tau_i,\xi_i],
\end{equation}
and the $E_2$-pages of the sequences \eqref{sp.hm} can be described
very explicitly. Namely, given some formal generators
$x_0,\dots,x_n$ of some homological degrees $\deg x_0,\dots,\deg
x_n$, let $k\{x_0,\dots,x_n\}$ be the free divided power algebra
generated by $x_\idot$ (that is, the free graded-commutative algebra
generated by symbols $x_i^j=x_i^{[p^j]}$, $j \geq 0$, $0 \leq i \leq
n$ of degrees $\deg x_i^j = p^j\deg x_i$, modulo the relations
$(x_i^j)^p=0$). Then the standard Hochschild-Kostant-Rosenberg
argument provides algebra isomorphisms
$$
\begin{aligned}
HH_\idot(k[x_0,\dots,x_i],k) &\cong k\{y_0,\dots,y_i\},\\
HH^\idot(k[x_0,\dots,x_i],k) &\cong k[y'_0,\dots,y'_i],
\end{aligned}
$$
where the free graded-commutative algebra $k[x_0,\dots,x_n]$ acts on
$k$ via the augmentation map (that is, all $x_i$ act by $0$), and
$y_i$ resp.\ $y_i'$ are generators of homological
resp.\ cohomological degree $\deg x_i+1$. Plugging this into
\eqref{hq.st}, and observing that $HH_\idot(k \otimes k,k) \cong
HH^\hdot(k \otimes k,k) \cong k$ since $k$ is perfect, we conclude
that the page $E^2_{\idot,\idot}(k)$ resp. $E_2^{\hdot,\hdot}(k)$ of
the spectral sequence \eqref{sp.hm} for $HM_\idot(k)$
resp.\ $HM^\hdot(k)$ is given by
\begin{equation}\label{e2}
E^2_{\hdot,\hdot}(k) = k\{\sigma,A_i,B_i\}, \quad
E_2^{\hdot,\hdot}(k) = k[\eps,C_i,D_i], \qquad i \geq 1,
\end{equation}
where $\sigma$, $A_i$, $B_i$ resp.\ $\eps$, $C_i$, $D_i$ have
homological resp.\ cohomological bidegrees $(1,1)$,
$(2p^i-2,1)$, $(2p^i-1,1)$. If $p=2$, then a similar computation
produces a spectral sequence
\begin{equation}\label{sp.hm.2}
k\{\sigma\} \Rightarrow HM_\idot(k), \qquad \deg \sigma = 2
\end{equation}
that degenerates for obvious dimension reasons. For \eqref{e2}, this
is not true: again for dimension reasons, the sequences degenerate
up to the $E^{p-1}$ resp. $E_{p-1}$-page, but then there are
differentials $d^{p-1}$ resp.\ $d_{p-1}$ given by
\begin{equation}\label{C.D.i}
d^{p-1}(B_{i-1}^{[p]}) = c_i^{-1}A_i, \qquad d_{p-1}(C_i) = c_iD_{i-1}^p
\end{equation}
for some possibily non-trivial elements $c_i \in k$.

As it happens, all the elements $c_i$ in \eqref{C.D.i} are indeed
non-trivial (and then again for dimension reasons, the spectral
sequences degenerate at $E^p$ resp.\ $E_p$). For Mac Lane Homology,
this was first discovered by Breen \cite{breen}, and for $THH$, this
is in the original paper of B\"okstedt \cite{bo} (the two arguments
are completely different, and the identification between $THH$ and
$HM_\idot$ was in any case only discovered later). In 1994, there
appeared another proof by Franjou, Lannes and Schwartz \cite{slf}
using cohomology rather than homology, and homological rather than
homotopical methods. A good overview of this proof can be found in
\cite{LP}, together with a very nice description of the
state-of-the-art at the time of its writing. Since then, there has
been a lot of progress (such as the introduction of strictly
polynomial functors in \cite{FS}, and a lot of subsequent work in
the area including the groundbreaking paper of Touz\'e
\cite{tou}). With all these modern advances, and with the ground we
have prepared in Subsection~\ref{div.subs}, we can now cut the proof
down to several lines. Let us do it for the convenience of the
reader.

\medskip

The numbers $c_i$ in \eqref{C.D.i} obviously only depend on $p$, so
that we may assume right away that $k=\Ff_p$ is a prime field. For
any integer $n \geq 1$, let $P^n(k)$ be the $k$-linear category
whose objects are finite-dimensional vector spaces, and whose
morphisms are given by
$$
P^n(k)(V,V') = D^n(P(k)(V,V')),
$$
where $D^n$ is lax monoidal divided power functor \eqref{S.D}. A
{\em strictly polynomial functor} of degree $n$ is a $k$-linear
functor $P^n(k) \to k\amod$. These form an abelian category
$\Fun_k(P^n(k),k)$ with its derived category $\D(P^n(k))$. We will
need the following general fact.

\begin{lemma}\label{efi.le}
For any $n \geq 1$, the category $\Fun_k(P^n(k),k)$ of strictly
polynomial functors of degree $n$ has finite homological dimension.
\end{lemma}

\proof{} By an easy and standard argument, for any $N \geq n$,
$\D(P^n(k))$ is equivalent to the derived category of degree-$n$
polynomial representations of the algebraic group $GL(N,k)$
(\cite[Theorem 2.4]{BFS}, or see \cite[Theorem 4.7]{ef} and the
following discussion). The latter is a highest weight category and
therefore has finite homological dimension (\cite[Theorem
  3.7.2]{CPS}, \cite{Do}); alternatively, by \cite[Theorem 1.4]{ef},
it is a full subcategory in the derived category of quasicoherent
sheaves on a Grassmann variety $G(N,M)$ for a sufficiently large $M$.
\endproof

Now, while the category $P^n(k)$ is not half-additive, we still have
a functor $m:\Gamma_+ \times P^n(k) \to P^n(k)$ induced by the
functor \eqref{ga.m} for the category $P(k)$. We say that an object
$E \in \D(P^n(k))$ is {\em stable} if so is $m^*E \in \D(\Gamma_+
\times P^n(k),k)$, and we denote by $\D^{st}(P^n(k)) \subset
\D(P^n(k))$ the full subcategory spanned by stable objects. By the
same argument as in Subsection~\ref{add.exa.subs}, we then have an
equivalence
\begin{equation}\label{D.Q.st}
\D^{st}(P^n(k)) \cong \D(Q_\idot^n(k)),
\end{equation}
where $Q_\idot^n(k)$ is the DG algebra \eqref{Q.n}. In particular,
the category is trivial unless $n$ is a power of $p$. If $n=p^d$ is
a power of $p$, we have at least one non-zero stable object given by
the ``iterated Frobenius'' functor $T^{(d)}:P^n(k) \to k\amod$ that
sends a vector space to itself, and acts on morphisms via the
$d$-fold iteration of the map \eqref{D.fr}. The equivalence
\eqref{D.Q.st} sends $T^{(d)}$ to the trivial DG module $k$. If we
let $HM^\hdot(k;d)$ be the algebra given by
$$
HM^\hdot(k;d) = \Ext^\hdot_{P^n(k)}(T^{(d)},T^{(d)}) \cong
\Ext^\hdot_{Q^{(d)}_\idot(k)}(k,k),
$$
where $Q^{(d)}_\idot$ is as in \eqref{Q.i}, then \eqref{HQ.i}
provides a spectral sequence
\begin{equation}\label{e2.d}
k[\eps,C_i,D_i] \Rightarrow HM^\hdot(k;d),
\end{equation}
where the generators are numbered by $i$ with $1 \leq i \leq d$ and
have the same degrees as in \eqref{e2}. Moreover, these spectral
sequences form an inductive system with respect to $d$, and by
\eqref{Q.lim}, the cohomological sequence of \eqref{e2} is the
colimit of the sequences \eqref{e2.d} over $d$. In particular, they
all degenerate up to $E_{p-1}$, and their differentials $d_{p-1}$
are given by \eqref{C.D.i}, with the same numbers $c_i$.

\begin{lemma}\label{bok.le}
All the elements $c_i \in k$ in \eqref{C.D.i} are non-trivial, the
spectral sequences \eqref{e2} degenerate at $E^p$ resp.\ $E_p$, and
we have $HM_i(k) \cong HM^i(k) \cong k$ if $i$ is even and $0$
otherwise.
\end{lemma}

\proof{} Assume that $c_d=0$ for some $d$, take the smallest such
$d$, and consider the spectral sequences \eqref{e2.d} for
$HM^\hdot(k;d)$ and $HM^\hdot(k;d-1)$. For dimension
reasons, the latter degenerates at $E_p$, and $HM^\hdot(k;d-1)$
is an algebra concentrated in degrees between $0$ and $2p^d-2$. But
since $c_d=0$, this means that again for dimension reasons, the
former also degenerates at $E_p$. Then $D_d$ survives into the
$E_\infty$-term and provides an element $D \in HM^\hdot(k;d)$
such that $D^n \neq 0$ for any $n \geq 1$. This contradicts
Lemma~\ref{efi.le}.
\endproof

\subsection{Trace theories.} Let us now interpret Hochschild
homology and its generalizations as homotopy trace theories in the
sense of Section~\ref{ho.tr.sec}.

Hochschild homology itself is in the fact the prototypical example
of a trace theory. To see it as such, consider the category
$k\amod^{fl}$ of flat $k$-modules, and as in
Subsection~\ref{inva.subs}, let $\aAlg(k) \subset
\iMor(k\amod^{fl})$ be the $2$-category of flat $k$-algebras, with
morphisms given by right-flat bimodules. As in
Subsection~\ref{extra.subs}, since $k\amod^{fl}$ is symmetric
monoidal, the embedding $\I(k):k\amod^{fl} \subset k\amod$ has a
trivial structure of a trace functor that we denote by $\I(k)_{triv}
\in \Tr(Bk\amod^{fl},k\amod)$. Then for any flat $k$-algebra $A$ and
right-flat $A$-bimodule $M$, we have a natural identification
$\Exp(\I(k)_{triv})(\langle A,M \rangle) \cong \Tr_A(M)$. The
homotopy expansion $\Exp^{ho}(\I(k)_{triv})$ in
$\Ho_{tr}(\aAlg(k),k)$ is stable, and if we denote
$$
\HH(k) = \Exp^{ho}(\I(k)_{triv}) \cong \Exp^{st}(\I(k)_{triv}) \in
\Ho_{tr}(\aAlg(k),k),
$$
then we have a functorial identification
$$
\ho(\HH(k))(\langle A,M \rangle) \cong CH_\idot(A/k,M) \in \Ho(k) \cong
\D^{\leq 0}(k),
$$
where $\ho(-)$ is the functor \eqref{ho.Ho}. Moreover, we also have
the simplicial $k$-module $(M/A)_\hash = \I(k)_{triv}(M/A)_\hash \in
\Delta^ok\amod$ of \eqref{M.A.hash}, and its standard complex
$C_\idot(\Delta^o,(M/A)_\hash)$ gives the Hochschild complex
$CH_\idot(A,M)$. The homotopy trace theory $\HH(k)$ is
multiplicative, and by \eqref{mult.exp.dia}, the corresponding
porduct is the usual K\"unneth product \eqref{kun.hh}.

The first non-trivial example of a trace functor is obtained by
using the edgewise subdivision of Subsection~\ref{edge.subs} (it is
described in some detail in \cite{trace}). For any positive integer
$p \geq 1$, shorten notation by setting $\pi_p(k) =
\pi_p(Bk\amod^{fl})$, $i_p(k)=i_p(Bk\amod^{fl})$, and denote
$$
\I^{(p)}(k) = \pi_p(k)_*i_p(k)_{triv}^*\I(k) \in
\Tr(Bk\amod^{fl},k\amod)
$$
the $p$-th edgewise subdivision of the trace functor $\I(k)_{triv}$
in the sense of Definition~\ref{edge.def}. We then consider its
stable homotopy expansion, and define
\begin{equation}\label{hh.p}
\begin{aligned}
&\HH^{(p)}(k) = \Exp^{st}(I^{(p)}(k)) \in \Ho_{tr}^{st}(\aAlg(k),k),\\
&CH^{(p)}_\idot(A/k,M) = \HH^{(p)}(k)(\langle A,M \rangle) \in \Ho(k)
\end{aligned}
\end{equation}
for a flat $k$-algebra $A$ and a right-flat $A$-bimodule $M$. We let
$HH^{(p)}_\idot(A/k,M)$ be the homology modules of the complex
$CH^{(p)}_\idot(A/k,M)$. If $M=A$ is the diagonal bimodule, we let
$CH_\idot^{(p)}(A/k) = CH_\idot^{(p)}(A/k,A)$, and we recall that
$CH_\idot^{(p)}(A/k) \in \D^\forall(\Delta^o,k) \cong \D(k)$ extends
to the category $\Lambda$: if we denote
\begin{equation}\label{cc.p}
CC^{(p)}(A/k) = \rAdj(\Exp^{st}(I^{(p)}(k)))(A) \in
\D^\forall(\Lambda,k),
\end{equation}
then $CH^{(p)}(A/k) \cong j^{o*}CC^{(p)}(A)$. Alternatively, we can
first take the stabilization $\Stab(\I^{(p)}(k)) \in
\Ho^{st}_{tr}(Bk\amod^{fl},k)$, and then identify $\HH^{(p)}(k) =
\Exp^{ho}(\Stab(\I^{(p)}(k)))$. Note that the underlying functor
$k\amod^{fl} \to k\amod$ of the trace functor $I^{(p)}(k)$ is the
$p$-th cyclic power functor of \eqref{C.M.eq}, so that
$\Stab(\I^{(p)}(k))$ can be computed by the same argument as in
Lemma~\ref{tate.le}. Namely, for any bifibration $\pi:I' \to I$ of
bounded categories whose fibers are connected groupoids with finite
automorphisms groups, we have a natural trace map $\tr:\pi_! \to
\pi_*$ between the Kan extension functors $\pi_!,\pi_*:\Fun(I',k)
\to \Fun(I,k)$ (if $I=\ppt$, so that $I'$ has a single object $o$,
$\tr$ is the averaging over the automorphism group $\Aut(o)$, and in
the general case, see e.g.\ \cite[Subsection 1.1]{hw}). Moreover,
see e.g.\ \cite[Subsection 3.1]{dege}, one can define the relative
Tate cohomology functor $\pi^\flat$ that fits into a functorial
exact triangle
\begin{equation}\label{pi.fl}
\begin{CD}
R^\hdot\pi_*(E) @>>> \pi^\flat(E) @>>> L^\hdot\pi_!(E)[1] @>{\tr}>>
\end{CD}
\end{equation}
for any $E \in \Fun(I',k)$, where $[1]$ stands for the homological
shift. We then have
\begin{equation}\label{stab.I.p}
\Stab(\I^{(p)}(k)) \cong \tau_{\geq 0}\pi_p(k)^\flat
i_p(k)^*\I(k)_{triv}.
\end{equation}
In particular, $\Stab(I^{(p)}(k))$ fits into a distinguished
triangle
\begin{equation}\label{I.p.tria}
\begin{CD}
\I^{(p)}(k) @>>> \Stab(I^{(p)}(k)) @>>>
L^\hdot\pi_p(k)_!i_p(k)^*\I(k)[1] @>>>
\end{CD}
\end{equation}
in the triangulated category $\D(\Lambda Bk\amod^{fl},k)$ induced by
\eqref{pi.fl}. The object $CH_\idot^{(p)}(A,M) \in \D(k)$ is given
by
$$
CH^{(p)}_\idot(A,M) \cong C_\idot(\Delta^o,\tau_{\geq
  0}\pi_p(k)^\flat i_p(k)^*(M/A)_\hash),
$$
and the homology of the category $\Delta^o$ can be again computed by
taking the standard complexes of simplicial $k$-modules, and
totalizing the resulting bicomplex. Since $I^{(p)}(k)$ is the
edgewise subdivision of a multiplicative trace functor, it is itself
multiplicative, and then $\HH^{(p)}(k)$ is a multiplicative homotopy
trace theory; the product is obtained by composing the K\"unneth
product \eqref{kun.hh} and the standard lax monoidal structure on
the relative Tate cohomology functor $\pi_p(k)^\flat$.

Now consider the category $\Sets$ of sets, with the cartesian
product, and the category $\Sets_+$ of pointed sets, with the smash
product. The forgetul functor $F:\Sets_+ \to \Sets$ has a
left-adjoint $\I:\Sets \to \Sets_+$ sending $S \in \Sets$ to the
coproduct $S_+ = S \copr \{o\} \in \Sets_+$, and $\I$ is symmetric
monoidal, so that $F$ is lax monoidal by adjunction. Then again,
since $\Sets$ is symmetric monoidal, $\I$ has a trivial trace
functor structure $\I_{triv} \in \Tr(B\Sets,\Sets_+)$, and we can
consider its stable expansion $\Exp^{st}(\I_{triv}) \in
\Ho^{st}_{tr}(\iMor(\Sets))$. However, we also have the forgetful
functor \eqref{phi.eq}. Denote by $\phi_- = F \circ \phi:k\amod^{fl}
\to \Sets$ its composition with the forgetful functor $F$, and note
that both $F$ and $\phi$ are lax monoidal, so that $\phi_-$ is lax
monoidal as well. The monoidal category $k\amod^{fl}$ is
well-generated by Example~\ref{k.dua.exa}, so that we have the
$2$-functors $T$, $P$ of \eqref{T.eq} and \eqref{P.eq}, and the map
\eqref{PT.triv} of Subsection~\ref{extra.subs} whose stabilized
homotopical counterpart reads as
\begin{equation}\label{thh.hh}
\Lambda P^*\Lambda\phi_-^*\Lambda T^*\Exp^{st}(\I_{triv}) \to
\Exp^{st}((\phi_-^*\I)_{triv}).
\end{equation}
Note that we have $\I \circ \phi_- \cong \phi \vee \ppt_+$, where
$\ppt_+$ is the constant functor sending everything to $\ppt_+ \in
\Sets$. Then $\Stab(\ppt_+) \cong \emptyset_+$, the constant
functors with value $\emptyset_+$, and the map $\Stab(\phi_-^*\I)
\cong \Stab(\phi \vee \ppt_+) \cong \Stab(\phi) \times \Stab(\ppt_+)
\to \Stab(\phi)$ induced by the adjunction map $\I \circ \phi_- \to
\phi$ is an isomorphism. Combining this with \eqref{phi.exp}, we
obtain identifications
$$
\Exp^{st}((\phi_-^*\I)_{triv}) \cong \Exp^{st}(\phi) \cong \phi
\circ \Exp^{st}(\I(k)),
$$
so that after restriction to $\aAlg(k) \subset \iMor(k)$, the target
of the map \eqref{thh.hh} becomes naturally identified with $\phi
\circ \HH(k)$. Let us denote its source by
$$
\THH(k) = \Lambda P^*\Lambda\phi_-^*\Lambda
T^*\Exp^{st}(\I_{triv})|_{\aAlg(k)} \in \Ho^{st}_{tr}(\aAlg(k)),
$$
and justify the notation by the following result.

\begin{lemma}\label{thh.le}
For any flat $k$-algebra $A$ and right-flat $k$-bimodule $M$, we
have a functorial identification
$$
\ho^{st}(\THH(k))(\langle A,M \rangle) \cong THH(A,M) \cong
\phi(THH_\idot(A,M)) \in \Ho^{st}(\Gamma_+),
$$
where $\ho^{st}(\THH(k))$ is as in \eqref{ho.tr}.
\end{lemma}

\proof{} By definition, objects in $\iMor(\Sets)$ are small
categories, and morphisms frm $I$ to $I'$ are functors $F:I^o \times
I' \to \Sets$ satisfying the appropriate representability condition
(namely, $F(i \times -):I' \to \Sets$ is ind-polyrepresentable in
the sense of Definition~\ref{repr.def} for any $i \in I$). In
particular, endomorphisms of a small category $I$ are bifunctors
$F:I^o \times I \to \Sets$. But we have the
functor \eqref{xi.fl}, and then as in the proofs of
Lemma~\ref{exp.le} and Lemma~\ref{S.le}, \eqref{la.fib} provides a
natural identification
$$
\ho(\Exp^{ho}(\I_{triv}))(\langle I,F \rangle) \cong
\hocolim_{\Delta^o I} \I \circ \xi_\flat^*(s \times t)^*F
$$
for any small category $I \in \iMor(\Sets)_{[0]}$ and its
endomorphism $F$. By virtue of \eqref{phi.stab} and \eqref{stab.pb},
if we apply this to $I=P(A)$ and $F=P(M)$, we obtain an
identification
$$
\ho^{st}(\THH(k))(\langle A,M \rangle) \cong \phi \circ
C_\idot(\Delta^oP(A),\xi_\flat^*(s \times t)^*P(M)),
$$
and then by \eqref{thh.fu} and \eqref{CH.Tw}, it only remains to
prove that the functor \eqref{xi.fl} is homotopy cofinal. The
argument for this is exactly the same as in Lemma~\ref{sp.le}, with
special maps replaced by bispecial maps.
\endproof

\begin{remark}
Lemma~\ref{thh.le} suggest that the unstable homotopy trace theory
$\Exp^{ho}(\I_{triv})$ is actually more fundamental. If we take as
$F$ the identity endomorphism of a small category $I$, then
$\ho(\Exp^{ho}(\I_{triv}))(\langle I,\id \rangle)$ is simply the
geometric realization of the cyclic nerve $N^{cy}(I)$ of the
category $I$. Considering the cyclic nerve of the category $P(A)$
has been suggested by Goodwillie \cite{G.W} back in 1988; among
other things, he observed that this cyclic nerve has an infinite
loop space structure induced by the direct sum in $P(A)$, and gives
rise to a ``cyclic $K$-theory spectrum'' of the algebra
$A$. Morally, this should be a ``master theory'' that gives all other
interesting theories by formal procedures (e.g.\ stabilization gives
$THH$). Unfortunately, the cyclic nerve is very hard to compute; it
seems that even when $A = \Ff_p$ is a prime field, the homotopy
groups of $|N^{cy}(P(A))|$ are not known.
\end{remark}

\begin{remark}
Our notation seems to imply that $\THH(k)$ somehow depends on $k$
(which would be strange since $THH$ is an absolute theory). This an
illusion: $k$ in $\THH(k)$ simply denotes the domain of definition
of the trace theory, and for any commutative ring map $f:k' \to k$,
we have a natural identification $f^*\THH(k) \cong \THH(k')$.
\end{remark}

With out new notation, \eqref{thh.hh} is simply a map $\THH(k) \to
\phi(\HH(k))$, and under the identication of Lemma~\ref{thh.le}, it
induces a functorial map
\begin{equation}\label{thh.aug}
THH(A,M) \to \phi(CH_\idot(A,M))
\end{equation}
for any flat $k$-algebra $A$ and right-flat $A$-bimodule $M$. In
terms of bifunctor homology, this is the augmentation map
\eqref{c.i.m.aug}. Moreover, $\THH(k)$ is a multiplicative homotopy
trace theory, and again by \eqref{mult.exp.dia}, the product
coincides with the functor homology product \eqref{kun.thh}, while
the augmentation map \eqref{thh.aug} is a multiplicative map.

In addition, since the product in $\Sets$ is cartesian, a natural
map \eqref{F.triv} for the identity functor $\Id$ and any integer
$l$ is given by the diagonal maps $S \to S^l$, $S \in \Sets$, and
this then induces a map \eqref{F.triv} for any functor $E:\Sets \to
\E$ to some category $\E$. In particular, this applies to $\I$, so
that the expansion $\Exp^{ho}(\I_{triv})$ and its stabilization
$\Exp^{st}(\I_{triv})$ acquire a natural $\langle F,l
\rangle$-structure in the sense of Definition~\ref{F.def}. The map
\eqref{thh.hh} then induces a natural map
\begin{equation}\label{thh.p}
\THH(k) \to \phi \circ \HH^{(p)}(k),
\end{equation}
where $\HH^{(p)}(k)$ is the homotopy trace theory \eqref{hh.p}. For
any flat $k$-algebra $A$ and right-flat $A$-bimodule $M$, we then
have functorial maps
\begin{equation}\label{thh.psi}
\begin{aligned}
\psi:THH(A,M) &\to \phi(CH^{(p)}_\idot(A/k,M)),\\
\psi_\idot:THH_\idot(A,M) &\to HH^{(p)}_\idot(A/k,M).
\end{aligned}
\end{equation}
Here $\psi$ is a map in $\Ho^{st}(\Gamma_+)$ (that is, a map of
connective spectra), and $\psi_\idot$ is the induced map on homotopy
groups. Since \eqref{thh.p} is a map of trace theories, $\psi$
commutes with the corresponding maps $B$ of \eqref{B.eq}. Moreover,
\eqref{thh.p} is a map of multiplicative trace theories, so that in
the case $M=A=k$, $\psi$ and $\psi_\idot$ are multiplicative
maps. However, neither of them needs to be $k$-linear (and in
interesting example, they are not).

\subsection{B\"okstedt periodicity.}\label{bok.subs}

Now fix a prime $p$, and from now on, assume that $k$ is a perfect
field of characteristic $p$. Then $THH_\idot(k) \cong HM_\idot(k)$
is the abutment of the spectral sequence \eqref{sp.hm} with the
$E_2$-page \eqref{e2}, or \eqref{sp.hm.2} if $p=2$, and in either
case, the sequence for $THH_i(k)$, $i = 0,1,2$ degenerates for
dimension reasons, so that $THH_0(k) \cong k$, $THH_1(k)=0$, and
$THH_2(k) \cong k$ is the free $k$-module generated by a single
element $\sigma$. Morover, $THH_\idot(k)$ is a commutative
cocomutative Hopf algebra, and in particular, we have a natural
algebra map
\begin{equation}\label{bok.eq}
\alpha:k[\sigma] \to THH_\idot(k),
\end{equation}
where $k[\sigma]$ is the free commutative $k$-algebra generated by
$\sigma$, $\deg \sigma = 2$.

\begin{prop}\label{bok.prop}
The map \eqref{bok.eq} is an isomorphism.
\end{prop}

Our proof of Proposition~\ref{bok.prop} uses the maps
\eqref{thh.psi} --- specifically, we take $A=M=k$, and consider the
corresponding maps
\begin{equation}\label{psi.k}
\psi:THH(k) \to \phi(CH_\idot^{(p)}(k/k)), \quad \psi_\idot:THH_\idot(k)
\to HH_\idot^{(p)}(k/k).
\end{equation}
Since the underlying functor $k\amod^{fl} \to k\amod$ of the trace
functor $\I^{(p)}(k)$ is the cyclic power functor $\C$ of
\eqref{C.M.eq}, for any $M \in k\amod^{fl}$, we have
\begin{equation}\label{ch.p.r}
\begin{aligned}
CH^{(p)}_\idot(k/k,M) &\cong R(M) \cong CH_\idot(k/k,R(M))
\cong\\
&\cong CH_\idot(P(k)/k,P(M)^*R),
\end{aligned}
\end{equation}
where $R \in \Ho(k\amod,k)$ is the stabilization \eqref{R.M.eq} of
the functor $C$. The map \eqref{psi.k} is then obtained by stabilizing
the map
$$
\hocolim_{\Tw(P(k))}\phi(s \times t)^*P(M)) \to
\hocolim_{\Tw(P(k))}\phi((s \times t)^*P(M)^*C)
$$
induced by \eqref{sp.psi.eq}, taking $M=k$, and composing the
resulting map
$$
THH(k) \cong \phi(CH_\idot(P(k),P(k))) \to \phi(CH_\idot(P(k),P(k)^*R))
$$
with the map \eqref{c.i.m.aug} for the category $P(k)$. In
particular, since the map \eqref{sp.psi.eq} is obviously
Frobenius-semilinear, so are the maps $\psi_\idot$ of \eqref{psi.k}:
we have $\psi_i(\lambda a) = \lambda^p\phi_i(a)$ for any $i \geq 0$,
$\lambda \in k$ and $a \in THH_i(k)$.

Now recall that we also have the object $CC^{(p)}_\idot(k/k)$ of
\eqref{cc.p} such that $CH_\idot^{(p)}(k/k) \cong
j^{o*}CC^{(p)}_\idot(k/k) \in \D^\forall(\Delta^o,k) \cong
\D(k)$. By \eqref{stab.I.p}, it is given by $CC^{(p)}_\idot(k/k)
\cong \tau_{\geq 0}\pi_p^\flat k$, where $k$ is the constant functor
with value $k$. Then \eqref{I.p.tria} induces an exact triangle
\begin{equation}\label{cc.tria}
\begin{CD}
k @>>> CC^{(p)}_\idot(k/k) @>>> L^\hdot\pi_{p!}k @>>>
\end{CD}
\end{equation}
in $\D^\natural(\Lambda,k)$ that restricts to a triangle
\begin{equation}\label{p.R}
\begin{CD}
k @>>> R(k) @>{r}>> \overline{R}(k) \cong C_\idot(\Z/p\Z,k) @>>>
\end{CD}
\end{equation}
in $\D^\forall(\Delta^o)(k) \cong \D(k)$ (where $R(k) \cong
CH_\idot^{(p)}(k/k)$ by \eqref{ch.p.r}). All terms in \eqref{p.R}
carry the Connes-Tsygan differential \eqref{B.eq}, and the
differential for $\overline{R}(k)$ is easy to compute: it has been
shown in \cite[Lemma 3.2]{cart} that
\begin{equation}\label{K.u.l}
L^\hdot\pi_{p!}k \cong \bigoplus_{i \geq 0}\K(k)[2i],
\end{equation}
where $\K(k) \cong \K_\idot \otimes k$ is the complex
\eqref{4.term}, and then $H_i(\Z/p\Z,k) \cong k$ in any non-negative
degree $i \geq 0$, and $B=B_i:H_i(\Z/p\Z,k) \to H_{i+1}(\Z/p\Z,k)$
is given by
\begin{equation}\label{B.p}
B_i = \begin{cases} \id, &\quad i = 2j,\\
0, &\quad i = 2j + 1.
\end{cases}
\end{equation}
As an algebra, $HH_\idot^{(p)}(k/k) \cong \tau^{\leq
  0}\vH^\hdot(\Z/p\Z,k)$ is the free graded-com\-muta\-tive algebra
given by
\begin{equation}\label{hh.p.alg}
HH^{(p)}_\idot(k/k) \cong k[\beta,\sigma], \qquad \deg \sigma =2, \deg
\beta = 1
\end{equation}
if $p$ is odd. If $p=2$, then graded-commutative means commutative,
and we have $\eps^2=\sigma$ instead of $\eps^2=0$.

\begin{lemma}\label{psi.2.le}
The component $\psi_2:THH_2(k) \to HH^{(p)}_2(k)$ of the map
\eqref{psi.k} is an isomorphism.
\end{lemma}

\proof{} Since $THH_1(k)=0$, the map \eqref{B.eq} for $THH(k)$
factors through the canonical truncation $\tau_{\geq 1}\Omega THH(k)
\cong \Omega \tau_{\geq 2}THH(k)$, and by \eqref{B.p}, the same is
true for $\overline{R}(k)[1]$. Applying further the canonical
truncation $\tau_{\leq 1}$, we obtain a diagram
$$
\begin{CD}
\phi(k) \cong \tau_{\leq 1}THH(k) @>{B}>> \phi(k[1]) \cong \tau_{\leq
  1}\tau_{\geq 1}\Omega THH(k)\\
@V{\psi'}VV @VV{\psi_2}V\\
\phi(k[1]) \cong \tau_{\leq 1}\phi(\overline{R}(k)[1]) @>{B}>>
\phi(k[1]) \cong \tau_{\leq 1}\tau_{\geq 1}\phi(\overline{R}(k))
\end{CD}
$$
in $\Ho^{st}(\Gamma_+)$, where $\psi' = r \circ \psi$ is the
composition of $\psi$ of \eqref{psi.k} and the projection $r$ in
\eqref{p.R}. If we further compose everything with the embedding
$\tau_{\leq 1}\tau_{\geq 1}\psi(\overline{R}(k)) \to \tau_{\leq
  1}\psi(\overline{R}(k))$, the diagram becomes commutative, and
since the embedding is split, it was commutative to begin with. But
$B$ in the bottom line is an isomorphism, and $\psi'$ is not equal
to $0$ by Lemma~\ref{boks.le}. Therefore $\psi_2:k \to k$ is not
equal to $0$ either; being Frobenius-semilinear, it must be an
isomorphism.
\endproof

\proof[Proof of Proposition~\ref{bok.prop}.] Lemma~\ref{psi.2.le}
and \eqref{hh.p.alg} immediately show that the map \eqref{bok.eq} is
injective: $\phi_{2n}\alpha(\sigma^n) = \phi_2(\sigma)^n \neq 0$ for
any $n$, so $\alpha(\sigma^n) \neq 0$ as well. If $p=2$, then the
degenerate spectral sequence \eqref{sp.hm.2} shows that the source
and the target of \eqref{bok.eq} are $k$-vector spaces of the same
dimension, so this finishes the proof. If $p \geq 3$ is odd, use
Lemma~\ref{bok.le}.
\endproof

\begin{remark}
The original proof of Proposition~\ref{bok.prop} in \cite{bo} used
non-trivial techniques from homotopy theory such as Dyer-Lashof
operations. Very recently, a very short topological proof appeared
in \cite{KN} where the authors deduce the theorem from a more
fundamental fact: $Q_\idot(k)$ considered as an $E_2$-algebra is
free on one generator $\beta$ of homological degree $1$. The
deduction is extremely easy but the fact itself again requires
Dyer-Lashof operations and other homotopy theory techniques
(including some classical but non-trivial computations). It would be
interesting to see if one can find a direct proof in terms of the
description of $Q_\idot(k)$ given in Subsection~\ref{steen.subs}.
\end{remark}

\subsection{Hochschild-Witt Homology.}

Proposition~\ref{bok.prop} shows that if $k$ is a perfect field of
characteristic $p \geq 0$, then the map $\psi_\idot$ of
\eqref{psi.k} is injective. It is certainly not an isomorphism,
since its target is $k$ in all degrees including the odd ones. It
turns out that one can modify the target so that the map becomes an
isomorphism, thus giving a purely algebraic model for
$THH(A,M)$. This uses Hochschild-Witt Homology of \cite{hw} and
polynomial Witt vectors of \cite{witt}. Let us briefly recall the
construction.

Recall that for any integer $m \geq 1$, we have the ring $W_m(k)$ of
$m$-truncated $p$-typical Witt vectors, with the restriction and
Frobenius ring maps $R,F:W_{m+1}(k) \to W_m(k)$. If $m=1$, then
$W_1(k) \cong k$, and the map $F:W_2(k) \to k$ is the composition of
the map $R$ and the Frobenius endomorphism of $k$, and since $k$ is
assumed to be perfect, the ideal $\Ker R \subset W_m(k)$ for any $m$
is generated by $p$. Now for any $m \geq 0$, simplify notation by
writing $i_{p^m}=i_{p^m}(W_{m+1}(k))$,
$\pi_{p^m}=\pi_{p^m}(W_{m+1}(k))$,
$\I_{triv}=\I(W_{m+1}(k))_{triv}$, and define trace functors
$Q_m,Q'_m \in \Tr(BW_{m+1}(k)\amod^{fl},W_{m+1}(k))$ and a map
$R:Q'_m \to Q_m$ by constructing the diagram
\begin{equation}\label{Q.dia}
\begin{CD}
\pi_{p^m!}i_{p^m}^*\I_{triv} @>{p\tr}>>
\pi_{p^m*}i_{p^m}^*\I_{triv} @>>> Q'_m @>>> 0\\
@V{p\id}VV @V{\id}VV @VV{R}V\\
\pi_{p^m!}i_{p^m}^*\I_{triv} @>{\tr}>>
\pi_{p^m*}i_{p^m}^*\I_{triv} @>>> Q_m @>>> 0
\end{CD}
\end{equation}
with exact rows. For $m=0$, we have $Q_0=0$ and $Q'_o =
\I_{triv}$. As it turns out (\cite[Proposition 2.3]{witt}), for any
$m \geq 1$ and flat $W_{m+1}(k)$-module $M$, $Q_m(M)$ and $Q'_m(M)$
only depend on the quotient $M/p$, so that we have $Q_m \cong
q^*W_m$, $Q'_m=q^*W'_m$ for some trace functors $W_m,W'_m \in
\Tr(Bk\amod,W_m(k))$, where $q:W_{m+1}(k)\amod^{fl} \to k\amod$
sends $M$ to $M/p$. One further shows (\cite[Corollary 2.7]{witt})
that $W_{m+1} \cong R_*W_m$ (where as in Example~\ref{norma.exa},
for any ring map $f:k' \to k$, we denote by $f_*:k\amod \to k\amod$
the restriction of scalars). Therefore $W_1(M) \cong M$, and the map
$R$ in \eqref{Q.dia} induces a functorial restriction map $R:W_{m+1}
\to R_*W_m$. All the functors $W_m$ are also symmetric lax monoidal
(\cite[Proposition 3.10]{witt}), and the maps $R$ are lax monoidal,
so that the limit $W=\lim W_\idot$ is a symmetric lax monoidal trace
functor on $k\amod$ with values in $p$-adically complete modules
over $W(k) = \lim W_\idot(k)$. Moreover, the Frobenius maps
$F:W_{m+1}(k) \to W_m(k)$ lift to lax monoidal trace functor maps
$F:W_{m+1} \to F_*W_m$, and we also have trace functor maps
$C:R_*W_m \to W_{m+1}$ that fit into short exact sequences
\begin{equation}\label{W.F.seq}
\begin{CD}
0 @>>> R^m_*W_n @>{C^m}>> W_{m+n} @>{F^n}>>
F^n_*\pi_{p^m!}i_{p^m}^*W_m @>>> 0
\end{CD}
\end{equation}
for any $m,n \geq 0$. We recall the functor $W_2$ already appeared
in the proof of Lemma~\ref{bok.le}, and \eqref{W.F.seq} for $m=n=1$
is \eqref{w2.2.seq}. The Frobenius maps $F$ commute with the
restriction maps $R$, thus define an $\langle F,p\rangle$-structure
in the sense of Definition~\ref{F.def} on the limit trace functor
$W$, and the map
\begin{equation}\label{F.W}
F:W \to \pi_{p*}i_p^*W
\end{equation}
of \eqref{F.dg} is actually an isomorphism. This makes it very easy
to compute the stabilization $W^{st}=\Stab(W)$ of the trace functor
$W$. Namely, for any $k$-vector space $M$, $W(M^{\otimes_k p})$
carries a natural $\Z/p\Z$-action induced by the trace functor
structure on $W$, with the $\Z/p\Z$-equivariant map
\begin{equation}\label{r.M}
W(M^{\otimes_k p}) \to W_1(M^{\otimes_k p}) \cong M^{\otimes_k p}
\end{equation}
induced by the restriction map $W \to W_1$, and we have the
following.

\begin{lemma}\label{W.p.le}
The map $\vH^i(\Z/p\Z,W(M^{\otimes_k p})) \to
\vH^i(\Z/p\Z,M^{\otimes_k p})$ induc\-ed by \eqref{r.M} is bijective
if $i$ is even, and its source vanishes if $i$ is odd.
\end{lemma}

\proof{} This is \cite[Corollary 3.9]{witt}. \endproof

By the same argument as in Lemma~\ref{tate.le}, Lemma~\ref{W.p.le}
together with the isomorphism \eqref{F.W} provide an identification
$$
W^{st} \cong \tau_{\geq 0}\pi_p(k)^\flat i_p(k)^*W,
$$
where $\pi^\flat$ is as in \eqref{stab.I.p}, so that $W^{st}$ fits
into an exact triangle
\begin{equation}\label{W.tria}
\begin{CD}
W @>>> W^{st} @>>> L^\hdot\pi_{p!}i_p^*W[1] @>>>
\end{CD}
\end{equation}
in the category $\Ho_{tr}(Bk\amod,W(k))$. As in \eqref{hh.p}, we can
now take the stable homotopy expansion, and define
\begin{equation}\label{whh.st}
\begin{aligned}
&\WHH^{st} = \Exp^{st}(W) \cong \Exp(W^{st}) \in
  \Ho_{tr}^{st}(\aAlg(k),W(k)),\\
&WCH^{st}_\idot(A/k,M) = \WHH(\langle A,M \rangle) \in \Ho(W(k))
\end{aligned}
\end{equation}
for a flat $k$-algebra $A$ and a right-flat $A$-bimodule $M$. By
\eqref{exp.ga} and \eqref{E.tr}, we have
\begin{equation}\label{whh.st.tr}
WCH^{st}_\idot(A/k,M) \cong C_\idot(\Delta^o,W^{st}(M/A)_\hash) \in \D(W(k))
\end{equation}
for a functorial object $W^{st}(M/A)_\hash \in \Ho(\Delta^o,W(k))$,
and we denote by $WHH^{st}(A/k,M)$ the homology modules of the
object \eqref{whh.st.tr}. If $M=A$ is the diagonal bimodule, then
$W^{st}(A/A)_\hash = j^{o*}W^{st}(A)_\hash$ for a cyclic object
$W(A)_\hash \in \Ho(\Lambda,W(k))$, and this is again functorial
with respect to $A$ in the sense of Lemma~\ref{A.ha.le}.

If we do not stabilize, the same procedure gives a homotopy trace
theory $\WHH = \Exp(W) \in \Ho_{tr}(\aAlg(k),W(k)$, functorial
objects $W(M/A)_\hash \in \Fun(\Delta^o,W(k))$, $W(A)_\hash \in
\Fun(\Lambda^o,W(k))$, and then the {\em Hochschild-Witt Homology}
modules $WHH_\idot(A,M)$ the {\em Hochschild-Witt complex}
$WCH_\idot(A,M)$ introduced and studied in some detail in
\cite{hw}. We also have the functorial stabilization map
\begin{equation}\label{whh.stab}
WHH_\idot(A,M) \to WHH^{st}_\idot(A,M)
\end{equation}
induced by a map of homotopy trace theories $\WHH \to
\WHH^{st}$. By \eqref{W.tria}, the corresponding map $W(M/A)_\hash
\to W^{st}(M/A)_\hash$ fits into an exact triangle
\begin{equation}\label{W.A.tria}
W(M/A)_\hash \longrightarrow W^{st}(M/A)_\hash \longrightarrow
L^\hdot\pi_{p!}i_p^*W(M/A)_\hash[1] \longrightarrow
\end{equation}
in $\D(\Delta^o,W(k))$, and if $M=A$ is the diagonal bimodule, the
triangle comes from a triangle in $\D(\Lambda,W(k))$. Moreover, both
$\WHH$ and $\WHH^{st}$ are multiplicative homotopy trace theories,
so that we have product maps
$$
\begin{aligned}
WHH_\idot(A,M) \otimes_{W(k)} WHH_\idot(A',M') &\to WHH_\idot(A
\otimes_k A',M \otimes_k M')\\
WHH^{st}_\idot(A,M) \otimes_{W(k)} WHH^{st}_\idot(A',M') &\to
WHH^{st}_\idot(A \otimes_k A',M \otimes_k M'),
\end{aligned}
$$
and the map \eqref{whh.stab} is compatible with these maps. In
addition to this, the restriction map $W(M) \to W_1(M) \cong M$
induces a map of homotopy trace theories $r:\WHH \to \HH(k)$ that
factors through \eqref{whh.stab} since $\HH(k)$ is stable, and by
virtue of the isomorphism \eqref{F.W}, its $p$-fold edgewise
subdivision defines a multiplicative map
\begin{equation}\label{whh.hh.p}
\WHH^{st} \to \HH^{(p)}(k),
\end{equation}
where $\HH^{(p)}(k)$ is as in \eqref{hh.p}.

To compare $WHH^{st}(A,M)$ to Topological Hochschild Homology, we
need one more ingredient introduced in \cite[Subsection 3.3]{witt}:
the non-additive {\em Teichm\"uller map} $T:\phi(M) \to
\phi(W_m(M))$. It is induced by the map $N \to Q_m(N)$, $N \in
W_m(k)$ sending $n \in N$ to $n^{\otimes p^m}$, and it in fact
factors through $M = N/p$ (\cite[Lemma 2.2]{witt}). We have $\phi(R)
\circ T = T$, so that the individual Teichm\"uller maps combine
together to a map $T:\phi(M) \to \phi(W(M))$. Moreover, this is a
map of symmetric lax monoidal trace functors, thus gives rise to a
map \eqref{PT.triv} that in this case reads as a map
\begin{equation}\label{T.whh}
T:\THH(k) \to \phi(\WHH^{st})
\end{equation}
of multiplicative stable homotopy trace theories on $k\amod$. Thus
for any flat $k$-algebra $A$ and $k$-flat $A$-bimodule $M$, we have
a map
\begin{equation}\label{thh.whh}
THH_\idot(A,M) \to WHH^{st}(A,M),
\end{equation}
and these maps are functorial in $M$ and compatible with the
products.

\begin{prop}\label{whh.prop}
The map \eqref{T.whh} is an isomorphism, so that the map
  \eqref{thh.whh} is an isomorphism for any $A$ and $M$.
\end{prop}

\proof{} By the universal property of expansion, it suffices to
check that the map
\begin{equation}\label{T.thh}
T_\idot:THH_\idot(k) \to WHH^{st}(k)
\end{equation}
induced by \eqref{T.whh} is an isomorphism. Indeed, \eqref{whh.hh.p}
induces a map
$$
S_\idot:WHH^{st}(k) \to HH^{(p)}(k)
$$
whose composition $S_\idot \circ T_\idot$ with \eqref{T.thh} is the
map $\psi_\idot$ of \eqref{psi.k}. But by Proposition~\ref{bok.prop}
and Lemma~\ref{W.p.le}, both $S_\idot$ and $\psi_\idot = S_\idot
\circ T_\idot$ are injective and have the same image.
\endproof

\subsection{Bott periodicity.}

Proposition~\ref{whh.prop} shows that for algebras over a perfect
field $k$ of positive characteristic $p$, Topological Hochschild
Homology is $W(k)$-linear as a homotopy trace theory: we have
$\THH(k) \cong \phi(\WHH^{st})$, and $\WHH^{st}$ is $W(k)$-linear by
construction. In particular, the map $B$ of \eqref{B.eq} for
$\THH(k)$ is $W(k)$-linear. This allows us to prove a comparison
theorem for {\em Periodic Topological Cyclic Homology} $TP(A)$
introduced in \cite{tp}. In general, this is defined by promoting
$THH(A)$ to a functor with values in $\Ho^{st}_{tr}(\Gamma_+,
\Lambda)$, as in \eqref{ho.radj}, refining \eqref{S1.eq} to
interpret the latter as the category of connective spectra equipped
with an action of the circle $S^1$, and taking the Tate fixed points
with respect to the circle. Tate fixed points normally take one out
of the world of connective spectra, so this is beyond the scope of
this paper. However, there is an alternative construction that works
perfectly well in the homological setting. The cohomology
$H^\hdot(\Lambda,\Z) = \Ext^\hdot_\Lambda(\Z,\Z)$ of the category
$\Lambda$ with integer coefficients is the algebra $\Z[u]$ with a
generator $u$ of degree $2$ known as the ``inverse Bott periodicity
generator''. Explicitly, the generator is represented by the complex
\eqref{4.term}. For any ring $R$ and object $E \in \Fun(\Lambda,R)$,
the homology $H_\idot(\Lambda,R)$ is naturally a module over
$H^\hdot(\Lambda,\Z)$, and one can consider the object
\begin{equation}\label{cp.e}
CP_\idot(E) = R^\hdot\lim_n C_\idot(\Lambda,E)[-2n] \in \D(R),
\end{equation}
where the limit is taken with respect to the inverse Bott
periodicity maps $u$. Its homology objects are denoted
$HP_\idot(E)$.

For the usual Hochschild Homology, $HP_\idot(A/k) =
HP_\idot(A_\hash)$ is the periodic cyclic homology of a flat
$k$-algebra $A$, and for the Hoshchild-Witt Homology, the periodic
theory $WHP_\idot(A) = HP_\idot(W(A)_\hash)$ has been introduced and
studied in \cite{hw}. We can also consider the stabilized theory
\begin{equation}\label{whp.st}
WHP^{st}_\idot(A) = HP_\idot(W^{st}(A)_\hash).
\end{equation}
The isomorphism \eqref{T.whh} of Proposition~\ref{whh.prop} then
also induces an isomorphism of the corresponding functors of
Lemma~\ref{A.ha.le}, and since the trace theory $\WHH^{st}$ is
obviously homotopy invariant in the sense of
Subsection~\ref{inva.subs}, the functors are defined on the whole
category $\Alg(k)$ of $k$-algebras and non-unital maps. Therefore
\eqref{whp.st} is functorially isomorphic to the Topological
Periodic Cyclic Homology groups $TP_\idot(A)$ of \cite{tp}, and
\eqref{whh.stab} induces a functorial map
\begin{equation}\label{whp.tp}
WHP_\idot(A) \to TP_\idot(A) \cong WHP^{st}_\idot(A).
\end{equation}
If $A$ is commutative and finitely generated over $k$, and $X =
\Spec A$ is smooth, then it has been proved in \cite[Theorem
  6.15]{hw} that the Hochschild-Witt Homology groups $WHH_\idot(A)$
are naturally identified with the spaces of de Rham-Witt forms on
$X$, and the differential \eqref{B.eq} corresponds to the
differential in the de Rham-Witt complex. It relatively
straightforward to deduce from this that if $p$ is large enough, or
if one inverts $p$, then $WHP_\idot(A)$ is the $2$-periodized
version of the cristalline cohomology $H^\hdot_{cris}(X)$. For the
$TP$-theory, analogous results were established in \cite{tp}, and
this seems contradictory since \eqref{whh.stab} is not an
isomorphism. However, it becomes one if we pass to the periodic
versions.

\begin{lemma}
For any perfect field $k$ and $k$-algebra $A$, the map
\eqref{whp.tp} is an isomorphism.
\end{lemma}

\proof{} Since $W(M)$ is $p$-adically complete for any $k$-vector
space $M$, the object $W(A)_\hash \in \Fun(\Lambda,W(k))$ is also
$p$-adically complete, and so is its restriction $i_p^*W(A)_\hash
\in \Fun(\Lambda_p,W(k))$. Then by \eqref{W.A.tria}, it suffices to
show that $CP_\idot(L^\hdot\pi_{p!}E)=0$ for any $p$-adically
complete $E \in \Fun(\Lambda_p,E)$. Indeed, we have
$$
C_\idot(\Lambda,L^\hdot\pi_{p!}E) \cong C_\idot(\Lambda_p,E),
$$
and the transition maps in the inverse system \eqref{cp.e} for
$CP_\idot(L^\hdot\pi_{p!}E)$ are given by the action of the element
$\pi_p^*u \in H^2(\Lambda_p,\Z)$. But $\pi_p^*u = pi_p^*u$, so this
element is divisible by $p$, and to finish the proof, it suffices to
check that $C_\idot(\Lambda_p,E)$ is $p$-adically complete. Indeed,
it can be computed by \eqref{C.lot}, and it is well-known (see
e.g.\ \cite[Appendix, A3]{FTadd}) that $k \in \Fun(\Lambda^o,k)$
admits a projective resolution $P_\idot$ whose every term is the sum
of a finite number of representable functors. Therefore as in
Remark~\ref{stab.prod.rem}, $C_\idot(\Lambda^o_p,-)$ commutes with
arbitrary products in $\D^{\leq 0}(\Lambda_p,k)$ and sends
$p$-adically complete objects to $p$-adically complete ones.
\endproof

Another useful comparison result suggested to us by A. Mathew
expresses Topological Hochschild Homology groups $THH_\idot(A)$ in
terms of the so-called {\em co-periodic cyclic homology}
$\bHP_\idot(A/k)$ introduced in \cite{coper}. To define it, recall
that for any commutative ring $k$ and flat $k$-algebra $A$, that
Connes-Tsygan differential $B$ of \eqref{B.eq} can be represented by
an explicit functorial endomorphism of the standard Hochschild
complex $CH_\idot(A/k)$, and then the periodic cyclic homology
$HP_\idot(A/k)$ can be computed by the complex $CP_\idot(A/k) =
CH_\idot(A/k)((u))$ with the differential $b+uB$, where $b$ is the
differential in the Hochschild complex, and $u$ is a generator of
cohomological degree $2$. Then $\bHP_\idot(A)$ is defined as the
homology of the complex $\bCP_\idot(A/k) = CH_\idot(A/k)((u^{-1}))$,
with the same differential. If $k$ contains $\Q$, these groups
vanish taulologically, but in other cases they are non-trivial and
provide an interesting homological invariant. In particular, if $k$
is a perfect field of positive characteristic $p$, then it has been
shown in \cite[Proposition 4.4]{coper} that we have a functorial
isomorphism
\begin{equation}\label{cp.b}
\bCP_\idot(A) \cong C_\idot(\Lambda,\pi_p^\flat i_p^*A_\hash[1])
\end{equation}
in the derived category $\D(k)$. For any object $E \in
\Fun(\Lambda_p,k)$, the relative Tate cohomology complex
$\pi_p^\flat E$ is equipped with an isomorphism $u:\pi_p^\flat E \to
\pi_p^\flat E[-2]$ that corresponds to the generator $u \in
H^2(\Z/p\Z,k)$, and in terms of \eqref{cp.b}, this corresponds to
the $k((u^{-1}))$-module structure on $\bCP_\idot(A/k)$. Moreover,
$\pi_p^\flat E[1]$ is equipped with the increasing {\em conjugate
  filtration} $V_\idot$, $V_i\pi_p^\flat E[1] = \tau^{\leq
  2i}(\pi_p^\flat E[1])$, it induces a generalized conjugate filtration
\begin{equation}\label{V.eq}
V_i\bCP_\idot(A) \cong C_\idot(\Lambda,V_i\pi_p^\flat i_p^*A_\hash[1])
\end{equation}
on $\bCP_\idot(A)$, and by \cite[Proposition 6.4]{coper}, this gives
rise to a convergent {\em conjugate spectral sequence}
$HH_\idot(A/k)((u^{-1})) \Rightarrow \bHP_\idot(A)$.

Now consider the Hochschild-Witt trace theory and the corresponding
object $W(A)_\hash \in \Fun(\Lambda,W(k))$. Then the restriction
maps $R:W(M) \to M$, $M \in k\amod$ induces a functorial
$W(k)$-linear map $R:W(A)_\hash \to R_*A_\hash$, and this gives rise
to a functorial $W(k)$-linear map
\begin{equation}\label{W.V.0}
W^{st}(A)_\hash[1] \cong \tau^{\leq 0}\pi_p^\flat W(A)_\hash[1] \to
V_0\pi_p^\flat i_p^*R_*A_\hash[1].
\end{equation}
Moreover, relative Tate cohomology is a lax monoidal functor, so
that for any $E \in \Fun(\Lambda_p,k)$, we have a natural product
map $\pi_p^\flat k \otimes_k \pi_p^\flat E \to \pi_p^\flat E$. By
\eqref{K.u.l}, we have $\pi_p^\flat k \cong \K(k)[-1]((u))$, so that
in particular, we have a product map $\K(k)[-1] \otimes_k
\pi_p^\flat E \to \pi_p^\flat E$. This allows to refine
\eqref{W.V.0} to a functorial $W(k)$-linear map
\begin{equation}\label{W.V.1}
\K_\idot \otimes W^{st}(A)_\hash \to V_0\pi_p^\flat
i_p^*R_*A_\hash[1],
\end{equation}
where $\K_\idot$ is the complex \eqref{4.term}.

\begin{lemma}\label{hpb.le}
The map \eqref{W.V.1} is an isomorphism in $\D(\Lambda,W(k))$.
\end{lemma}

\proof{} For any $n \geq 1$, let $e_n:\ppt \to \Lambda$ be
the embedding onto $[n] \in \Lambda$. It suffices to check that
\eqref{W.V.1} becomes a quasiisomorphism after applying $e_n^*$ for
all $n$. We have
$$
e_n^*\pi_p^\flat i_p^*A_\hash \cong \vC^\hdot(\Z/p\Z,M^{\otimes p}),
\quad e_n^*\pi_p^\flat i_p^*W(A)_\hash \cong
\vC^\hdot(\Z/p\Z,W(M^{\otimes p})),
$$
where $M=A^{\otimes n}$, and the $\Z/p\Z$-action on $W(M^{\otimes
  p})$ is twisted as in \cite[Corollary 3.9]{witt}. As we saw in
Lemma~\ref{tate.le}, $\vH^i(\Z/p\Z,M^{\otimes p}) \cong M^{(1)}$ for
any integer $i$, while by \cite[Corollary 3.9]{witt},
$\vH^i(\Z/p\Z,W(M^{\otimes p})) \cong M^{(1)}$ if $i$ is even and
$0$ if $i$ is odd. The map \eqref{W.V.1} is an isomorphism in
homological degree $0$, and since it commutes with $u$, it suffices
to check that for any map, the map
\begin{equation}\label{e.pp}
\vH^0(\Z/p\Z,M^{\otimes p}) \to \vH^1(\Z/p\Z,M^{\otimes p})
\end{equation}
given by the action of the generator $\eps \in H^1(\Z/p\Z,k) \cong
k$ is an isomorphism. But for $M=k$ this is true tautologically, and
both the source and the target of \eqref{e.pp} commute with
arbitrary sums.
\endproof

\begin{remark}
The main argument in the proof of Lemma~\ref{hpb.le} is essentially
\cite[Proposition 5.6]{coper} but there is one difference: in
\cite{coper}, we used relative Tate homology instead of relative
Tate cohomology, and the coaction map rather than the action
map. The difference between Tate homology and Tate cohomology is
simply the homological shift $[1]$, and if $p \neq 2$, both
arguments work, but the homological argument requires the equality
$\eps^2=0$ that is not true for $p=2$. The cohomological argument
works for all $p$ (and it can be used to improve \cite{coper}
obviating the need for \cite[Section 4]{dege}).
\end{remark}

\begin{corr}\label{hpb.corr}
There exists a functorial isomorphism
\begin{equation}\label{thh.hpb}
THH(A) \cong V_0\bCP_\idot(A) \in \D(k),
\end{equation}
where $V_0$ is the conjugate filtration \eqref{V.eq}
\end{corr}

\proof{} Apply $C_\idot(\Lambda,-)$ to the isomorphism
\eqref{W.V.1}, and recall that by \eqref{HH.HC}, we have
$C_\idot(\Delta^o,j^{o*}E) \cong C_\idot(\Lambda,\K \otimes E)$ for
any $E \in \D(\Lambda,k)$.
\endproof

We note that $THH(A)$ is a module over $THH(k) \cong k[\sigma]$, and
by Lemma~\ref{bok.le}, possibly after rescaling by a non-zero
constant, the isomorphism \eqref{thh.hpb} intertwines the B\"okstedt
periodicity generator $\sigma$ with the Bott periodicity generator
$u^{-1}$. Therefore we also have
$$
THH(A) \otimes_{k[\sigma]} k(\sigma,\sigma^{-1}) \cong \bCP_\idot(A),
$$
and the conjugate spectral sequence corresponds to the spectral
sequence connecting $THH_\idot(A)$ and $HH_\idot(A)[\sigma]$.

\begin{remark}
Co-periodic cyclic homology is a non-commutative counterpart of the
non-standard derived de Rham cohomology obtained by considering the
de Rham complex of a free simplicial resolution of a commutative
algebra, and then taking its ``wrong'' totalization. The latter
formed a crucial part in the approach to $p$-adic Hodge theory of
Beilinson \cite{Be} and Bhatt \cite{Bh}, and it always puzzled the
author why this completely disappeared in the integral $p$-adic
Hodge theory developed later in \cite{BMS}. Corollary~\ref{hpb.corr}
explains the mystery: \cite{BMS} consistently uses Topological
Hochschild Homology, and this has the same effect as taking the
``wrong'' totalization.
\end{remark}

{\small\noindent
Affiliations:
\begin{enumerate}
\renewcommand{\labelenumi}{\arabic{enumi}.}
\item Steklov Mathematics Institute, Algebraic Geometry Section
  (main affiliation).
\item Laboratory of Algebraic Geometry, National Research University
Higher\\ School of Economics.
\end{enumerate}}

{\small\noindent
{\em E-mail address\/}: {\tt kaledin@mi-ras.ru}
}

\end{document}